\documentclass[oneside]{amsart}
\usepackage{esvect} 
\usepackage{geometry}
\usepackage{amsmath,amssymb,amsfonts, accents}
\usepackage{tikz}
\usepackage{tikz-cd}
\usepackage{hyperref}
\hypersetup{bookmarksdepth=2}
    \newgeometry{vmargin={30mm,40mm}, hmargin={34.5mm,34.5mm}}
\usepackage{xcolor}

\usepackage{cleveref} 
\usepackage[small,nohug,heads=vee]{diagrams}

\usepackage{xstring}
\usepackage{xcolor}
\usepackage[mathscr]{euscript}   
\usepackage[cmtip,all,2cell]{xy}

\usepackage{graphicx}
\usepackage{tipa}
\usepackage{enumitem} 
\usepackage[colorinlistoftodos]{todonotes}

\newcommand{\dgob}{\mathbf}   
\DeclareMathOperator{\Surj}{\mathbf{Sur}}

\usepackage[makeroom]{cancel}

\usepackage{tikzit} 
\usepackage{enumitem}

\tikzstyle{bullet}=[fill={rgb,255: red,6; green,6; blue,6}, draw=black, shape=circle, minimum size=.02cm, scale=0.6]

\tikzstyle{new edge style 0}=[-]

\entrymodifiers={+!!<0pt,\fontdimen22\textfont2>}

\newcommand{\ZZ}{\mathbb{Z}}

\newcommand{\EE}{\mathbb{E}}

\newcommand{\LL}{\mathbb{L}}
\DeclareMathOperator{\Fun}{Fun}
\newcommand{\Com}[0]{\mathbf{Com}}

\DeclareMathOperator{\Ch}{\mathbf{Ch}}

\DeclareMathOperator{\res}{res}
\newcommand{\rt}[0]{\longrightarrow}
\newcommand{\lt}[0]{\longleftarrow}
 \DeclareMathOperator{\weirdleq}{\substack{  _{} \vspace{-3.8pt} \\ \eqslantless }}
 
\DeclareMathOperator{\Pro}{Pro}
\DeclareMathOperator{\fin}{fin}
\DeclareMathOperator{\QC}{QC}
\DeclareMathOperator{\Op}{Op}
\DeclareMathOperator{\pd}{pd}
\DeclareMathOperator{\gen}{gen}

\DeclareMathOperator{\Fin}{\mathscr{F}in}
\DeclareMathOperator{\Alg}{Alg}
\DeclareMathOperator{\Barr}{Bar}
\DeclareMathOperator{\dgBar}{B}
\DeclareMathOperator{\aug}{aug}
\newcommand{\mb}[1]{\mathbf{#1}}
\newcommand{\mm}[1]{\mathrm{#1}}
\newcommand{\mf}[1]{\mathfrak{#1}}
\newcommand{\mc}[1]{\mathcal{#1}}
\newcommand{\ul}[1]{\underline{#1}}
\newcommand{\ol}[1]{\overline{#1}}
\newcommand{\cat}[1]{
\StrLen{#1}[\mystrlen]
\ifnum\mystrlen=1 \mathscr{#1}
\else \mathrm{#1}
\fi}

\usepackage{ upgreek } 
\DeclareMathOperator{\cooliota}{\upiota}
\DeclareMathOperator{\coolU}{\upupsilon}

\DeclareMathOperator{\End}{End}

\DeclareMathOperator{\Hom}{Hom}
\DeclareMathOperator{\Sur}{Sur}
\DeclareMathOperator{\Vect}{Vect}

\DeclareMathOperator{\Lie}{Lie}
\DeclareMathOperator{\Mod}{Mod}

\DeclareMathOperator{\sSeq}{sSeq}
\DeclareMathOperator{\Perf}{Perf}
\DeclareMathOperator{\APerf}{APerf}
\DeclareMathOperator{\Coh}{Coh}
\DeclareMathOperator{\id}{id}

\DeclareMathOperator{\Rings}{Rings}
\DeclareMathOperator{\Tot}{Tot}
\DeclareMathOperator{\simcos}{sc}

\DeclareMathOperator{\Sp}{Sp}

\DeclareMathOperator{\N}{N}

\DeclareMathOperator{\Map}{Map}

\DeclareMathOperator{\KD}{KD}
\DeclareMathOperator{\op}{op}
\DeclareMathOperator{\coOp}{coOp}

\DeclareMathOperator{\colim}{colim}

\usepackage{eqparbox}
\newcommand{\bcirc}{\bar{\circ}}
\newcommand{\Boxed}[2][A]{\fbox{{$\mathstrut#2$}}}
 
\DeclareMathOperator{\conv}{conv}
\DeclareMathOperator{\reg}{reg}
\DeclareMathOperator{\An}{\cat{S}}

\DeclareMathOperator{\StPrSigma}{\mathscr{P}r^{St,\Sigma}}
 \DeclareMathOperator{\SCR}{SCR}

\newcommand{\Del}{\Delta}
\DeclareMathOperator{\PrL}{\cat{P}r^L}
\DeclareMathOperator{\CAlg}{CAlg}
\DeclareMathOperator{\ev}{ev}
\newcommand{\PP}[0]{\mathscr{P}}
\newcommand{\CC}[0]{\mathscr{C}}
\DeclareMathOperator{\lev}{lev} 
\DeclareMathOperator{\Orb}{\mc{O}}

\DeclareMathOperator{\Coalg}{coAlg}
\DeclareMathOperator{\Coop}{coOp}

\DeclareMathOperator{\Cobar}{coBar}

\newcommand{\Cl}[0]{\cat{C}_-}
\newcommand{\Cm}[0]{\cat{C}_{\mf{m}}}
\newcommand{\Cr}[0]{\cat{C}_+}
\newcommand{\Dl}[0]{\cat{D}_-}
\newcommand{\Dm}[0]{\cat{D}_{\mf{m}}}
\newcommand{\Dr}[0]{\cat{D}_+}
\newcommand{\Ml}[0]{\cat{M}_-}
\newcommand{\Mm}[0]{\cat{M}_{\mf{m}}}
\newcommand{\Mr}[0]{\cat{M}_+}
\DeclareMathOperator{\LMod}{LMod}
\DeclareMathOperator{\RMod}{RMod}
\DeclareMathOperator{\LComod}{LComod}
\DeclareMathOperator{\Tw}{Tw}
\DeclareMathOperator{\triv}{triv}
\newcommand{\sLL}[0]{\scalebox{0.6}{$\LL$}}
\newcommand{\ROS}[0]{{R[\Orb_{\Sigma}]}}
\newcommand{\scMod}{\mathbf{Mod}^{\simcos}}
\newcommand{\cMod}{\mathbf{Mod}^{\mm{c}}}
\DeclareMathOperator{\dgOp}{\mathbf{Op}}
\DeclareMathOperator{\dgAlg}{\mathbf{Alg}}
\DeclareMathOperator{\dgCoalg}{\mathbf{coAlg}}
\newcommand{\scsSeq}{\mathbf{sSeq}^{\simcos}}

\DeclareMathOperator{\sd}{sd}
\DeclareMathOperator{\sdBar}{Bar^{sd}}
\DeclareMathOperator{\sdK}{K^{sd}}
\newcommand{\nondeg}[1]{((#1))}
\newcommand{\ope}[1]{\mathbf{#1}}
\newcommand{\partLieR}{\dgob{Lie}^{\pi}_{R,\Delta}}
\newcommand{\sG}{\scalebox{0.7}{$G$}}
\newcommand{\ulR}{\scalebox{0.7}{$\ul{R}$}}
\newcommand{\pmu}{\pm_{(\mb{u}, \alpha)}}
\newcommand{\pmc}{\pm_{||}}
\newcommand{\Caes}[1]{{#1} {}_{||}}

\setlist[itemize]{itemsep=3pt, topsep=3pt, labelindent=1ex, itemindent=0mm, labelsep=1.2ex, leftmargin=*}
\setlist[enumerate]{itemsep=3pt, topsep=3pt, labelindent=1ex, itemindent=0mm, labelsep=1.2ex, leftmargin=*}

\usepackage{amsthm}

    \theoremstyle{plain}

 \theoremstyle{plain}
    \newtheorem{theorem}{Theorem}[section]
    \newtheoremstyle{TheoremNum}
        {}{}              
        {\itshape}                      
        {}                              
        {\bfseries}                     
        {.}                             
        { }                             
        {\thmname{#1}\thmnote{ \bfseries #3}}
    \theoremstyle{TheoremNum}
    \newtheorem{theoremn}{Theorem}

 \theoremstyle{plain}
    
    \newtheoremstyle{LemmaNum}
        {}{}              
        {\itshape}                      
        {}                              
        {\bfseries}                     
        {.}                             
        { }                             
        {\thmname{#1}\thmnote{ \bfseries #3}}
    \theoremstyle{TheoremNum}

 \theoremstyle{plain}
    
    \newtheoremstyle{DefinitionNum}
        {}{}              
        {\itshape}                      
        {}                              
        {\bfseries}                     
        {.}                             
        { }                             
        {\thmname{#1}\thmnote{ \bfseries #3}}
    \theoremstyle{TheoremNum}
    \newtheorem{definitionn}{Definition}

 \theoremstyle{plain}
    
    \newtheoremstyle{CorollaryNum}
        {}{}              
        {\itshape}                      
        {}                              
        {\bfseries}                     
        {.}                             
        { }                             
        {\thmname{#1}\thmnote{ \bfseries #3}}
    \theoremstyle{TheoremNum}

\theoremstyle{plain}

\newtheorem{lemma}[theorem]{Lemma}  
\newtheorem{proposition}[theorem]{Proposition} 
    \newtheoremstyle{PropositionNum}
        {}{}              
        {\itshape}                      
        {}                              
        {\bfseries}                     
        {.}                             
        { }                             
        {\thmname{#1}\thmnote{ \bfseries #3}}
    \theoremstyle{PropositionNum}

\newtheorem*{proposition*}{Proposition}
\newtheorem*{corollary*}{Corollary}

\theoremstyle{definition}
\newtheorem{definition}[theorem]{Definition}
\newtheorem{cons}[theorem]{Construction}
\newtheorem{example}[theorem]{Example}
\newtheorem{notation}[theorem]{Notation}

\newtheorem{corollary}[theorem]{Corollary}
\newtheorem{reminder}[theorem]{Reminder}

\newtheorem{remark}[theorem]{Remark}
 \newtheorem{sgnrule}[theorem]{Sign Rule}

\newtheorem{warning}[theorem]{Warning}

\begin{document}  
 
\title{PD Operads and Explicit Partition Lie Algebras}  

\author[Lukas Brantner]{Lukas Brantner }
\address{Lukas Brantner, Oxford University,   Universit\'{e} Paris--Saclay (CNRS)} 
\email{lukas.brantner@universite-paris-saclay.fr, lukas.brantner@maths.ox.ac.uk}

\author[Ricardo Campos]{Ricardo Campos} 
\address{Ricardo Campos, Universit\'{e}   Toulouse III Paul Sabatier (CNRS)}
\email{ricardo.campos@math.univ-toulouse.fr}

\author[Joost Nuiten]{Joost Nuiten} 
\address{Joost Nuiten, Universit\'{e}   Toulouse III Paul Sabatier}
\email{joost.nuiten@math.univ-toulouse.fr}  
\maketitle

\begin{abstract} \hspace{-5pt}
Infinitesimal  deformations are governed by partition Lie algebras.
In characteristic $0$, these higher categorical  structures are  modelled by differential graded Lie algebras, but in  characteristic $p$, they  are \mbox{more subtle.}

We give explicit models for partition Lie algebras over  general coherent rings,  both  in the 
setting  of spectral and derived algebraic \mbox{geometry.}
For the \mbox{spectral} case,  we refine operadic  Koszul duality to a functor from operads to  {divided power operads,} by taking  {`refined linear duals'} of $\Sigma_n$-representations. The derived case requires a further refinement of Koszul duality   to a more genuine setting. 
\end{abstract}

 \ \\  \\ \\ \\  \\ \\ \\  \\ \\ \\  \\ \\ \\  \\ \\ \\  \\  
\setcounter{tocdepth}{1}
\tableofcontents

  \newpage

\section{Introduction}    
Infinitesimal  deformations over a field $k$ of characteristic zero are governed by \mbox{differential} graded Lie algebras. This paradigm, which was formalised by  Lurie  \cite{lurie2011derivedX} and Pridham \cite{pridham2010unifying}, was recently generalised
 to arbitrary fields, cf.\  \cite{brantner2019deformation}.

Over $\EE_\infty$-rings,    formal moduli    are equivalent to \textit{{spectral}  {partition}  {Lie algebras}}. These  are chain complexes with extra structure, which is parametrised by 
a sifted-colimit-preserving {monad}   {$\Lie_{k,\EE_\infty}^\pi$ \hspace{-5pt}}   satisfying the following formula
 on   \mbox{coconnective  $V$:}
\begin{equation} \label{spectral} \Lie_{k, \EE_\infty}^\pi(V) = \bigoplus_r \left( \widetilde{C}^\bullet(\Sigma |\Pi_r|^\diamond, k) \otimes V^{\otimes r} \right){}^{h\Sigma_r}.\vspace{-7pt} \end{equation}
Here $\widetilde{C}^\bullet(\Sigma |\Pi_r|^\diamond, k)$ denotes\vspace{0pt} cochains on the doubly suspended $r^{th}$ \mbox{partition complex.} 

While useful for   conceptual arguments as in \cite{brantner2020purely},
this abstract  definition   can be  somewhat  elusive
  in  concrete   instances of deformation theory.
In  this work, we construct concrete models for spectral partition Lie algebras over general (coherent) rings $R$, complementing the familiar differential graded Lie algebra models in 
\mbox{characteristic zero.}
To this end,  we introduce an operad in the ordinary category of  chain \mbox{complexes:}
$${\mathbf{Lie}_{R,\EE_\infty}^{\pi}} :=  \mathbf{Lie}_R^s  \otimes \boldsymbol{\Sur}^{\boldsymbol{\vee}}_R.$$
Here $\mathbf{Lie}_R^s$ is the usual (shifted) $R$-linear Lie operad with  
$\mathbf{Lie}_R^s(r)$  concentrated in degree $1-r$, where it is 
spanned by the Lie words in $r$ letters $x_1,\ldots, x_r$  which involve each letter exactly once, modulo antisymmetry and the Jacobi identity.

The  PD surjections operad $\boldsymbol{\Sur}^{\boldsymbol{\vee}}_R$ is the operad of (nonunital)
$\EE_\infty$-$R$-algebras with divided powers which is inspired by the surjections operad of McClure--Smith \cite{mcclure2003multivariable}. 
In  homological degree $-d\leq 0$,  $\boldsymbol{\Sur}^{\boldsymbol{\vee}}_R(r)$ is given by a free $R$-module   spanned by exhaustive
sequences $(u_1,\ldots, u_{r+d})$ of elements in $\{1,\dots,r\}$  \mbox{satisfying $u_j \neq u_{j+1}$ for all $j$.}\vspace{2pt}
 
The category of ${\mathbf{Lie}_{R,\EE_\infty}^{\pi}}$-algebras in chain complexes comes with a notion of `tame weak equivalence', finer than the usual notion of a quasi-isomorphism; inverting these gives  the $\infty$-category of spectral partition Lie algebras. This is   surprising as algebras over operads are defined \mbox{using orbits,} whereas the spectral partition Lie algebra monad  involves \mbox{homotopy fixed points.} \vspace{5pt}

In the setting of simplicial commutative rings, formal moduli   are equivalent to \textit{derived partition Lie algebras}, which  are parametrised by a sifted-colimit-preserving monad $\Lie_{k,\Delta}^\pi$
satisfying a  similar formula to \eqref{spectral}, but with 
 \textit{strict} fixed points. \vspace{3pt}

Modelling derived partition Lie algebras is slightly more involved than in the spectral setting. First,
 we construct an operad ${\mathbf{Lie}_{R,\Delta}^{\pi}}$ in cosimplicial $R$-modules. The component ${\mathbf{Lie}_{R,\Delta}^{\pi}}(r)^d$ is
given by $R$-valued functions on the  
  set $P(r)_d$ of pairs $$(\sigma,S) = \left([ \sigma_0<   \ldots <   \sigma_t ]\ ,\ S_0  \subseteq  \ldots  \subseteq  S_d\right),  $$
where $[ \sigma_0<    \ldots <  \sigma_t  ]$ is a   strictly increasing\vspace{1pt} chain of partitions of   \mbox{$\underline{r}=\{1,\ldots, r\}$}   with $\sigma_0 = \Boxed{1}\hspace{-0.3 pt} \Boxed{2} \hspace{-0.3 pt} \Boxed{3} \ldots \Boxed{r}  $ and   $\sigma_t=\Boxed{123\ldots r}$, and \vspace{1pt}
$S_0  \subseteq   \ldots \subseteq  S_d = \{0,\ldots,t\}$ is an increasing chain of subsets. Here, we also allow the case $t=-1$. One can think of $(\sigma, S)$ as a levelled tree, together with a nested collection of sets of marked levels.\vspace{3pt}

To equip a   cosimplicial $R$-module $\mathfrak{g}^\bullet$ \mbox{with a restricted ${\mathbf{Lie}_{R,\Delta}^{\pi}}$-algebra structure,} we must specify an  element
$ \{a_1,  \ldots , a_r\}_{(\sigma,S)}\in \mathfrak{g}^d$
for any  tuple $\mathbf{a}=(a_1,\ldots,a_r)$ of elements in $\mathfrak{g}^d$ and all pairs 
$(\sigma,S)\in P(r)_d$. But there is more: we must also specify an element 
$\gamma_{(\sigma,S)} (a_1,  \ldots , a_r)$
with 
$$
|\Sigma_{\mathbf{a},\sigma} |\cdot  \gamma_{(\sigma,S)} (a_1,  \ldots , a_r) = \{a_1,  \ldots  , a_r\}_{\gamma(\sigma,S)},
$$
where 
$\Sigma_{\mathbf{a},\sigma}   \subset \Sigma_r$ is the group of symmetries of ${\mathbf{a}}$ fixing the chain of partitions $\sigma$. These operations satisfy  compatibility properties, which we will describe in detail.

Simplicial-cosimplicial restricted ${\mathbf{Lie}_{R,\Delta}^{\pi}}$-algebras come with a notion of \mbox{`tame weak} equivalence'; inverting these gives the $\infty$-category of \mbox{derived partition Lie algebras.}\vspace{4pt}
 
Spectral partition Lie algebras and derived partition Lie algebras do not arise as algebras over an $R$-linear $\infty$-operad. Because of this, our constructions require a  twofold refinement of operadic Koszul duality, which is of independent interest. First, we introduce a 
divided power refinement of \mbox{$\infty$-operads,} which we call `divided power (PD) $\infty$-operads'. These allow us to take `continuous duals' of the $\Sigma_n$-representations   appearing in an $\infty$-operad.
Second, we study Koszul duality for `derived $\infty$-operads'; here, the group actions are more genuine, which lets us treat structures like simplicial commutative rings.

Throughout, we rely on the formalism pro-coherent sheaves, which originated in Deligne's  \cite[Appendix]{hartshorne1966residues} and is closely related to the theory of ind-coherent sheaves \cite{gaitsgory2013indcoh}. \vspace{-2pt}

\subsection{Statement of Results}  
Before stating our main results, we will briefly recall the \mbox{formalism} of pro-coherent modules.\vspace{-3pt}
\subsection*{Pro-coherent modules}

Any finite-dimensional $k$-vector space $V$   can be \mbox{recovered} from its linear dual $ \Map_{\Vect_k}(V,k)$; \mbox{if $\dim(V) = \infty$,}  this is \mbox{no longer true.} 
However, we can  take `continuous duals' and send $V$ to the pro-finite $k$-vector space   $   \varprojlim_{\substack{W \subset V   {f.d.}} } \Map_{\Vect_k}(W,k)$; this induces an  equivalence $\Vect_k \stackrel{\sim}{\rt} \Pro (\Vect_k^{\fin})^{\op}$.\vspace{1pt}

More generally, fix a coherent  $\EE_\infty$-ring $R$ (cf.\  \cite[Proposition 7.2.4.18]{lurie2014higher}) and write $\Mod_R$ for the $\infty$-category of $R$-module spectra. We can then refine the above construction and assign to every $R$-module a \emph{pro-coherent $R$-module}. This gives a functor $\upiota\colon \Mod_R\rt \QC^\vee_R$ to the  stable $\infty$-category of pro-coherent $R$-modules which we  recall in  \Cref{def:pro-coh}. 
We state several key properties of   $\QC^\vee_R$:
\begin{enumerate}
\item $\QC^\vee_R$  admits a closed \mbox{symmetric monoidal structure;}
\item The functor $\upiota\colon \Mod_R \rightarrow \QC^\vee_R$ is a symmetric monoidal left adjoint, which is  fully faithful on   connective $R$-modules. If $R$ is eventually coconnective, then $\upiota$  is in fact fully faithful on all of $\Mod_R$, and it is an equivalence when $R$ is a discrete regular Noetherian ring;
\item The essential image of $\upiota$ is \emph{not} closed under taking duals. In fact, 
$\QC^\vee_R$ is compactly generated by all `continuous duals' $M^\vee :=  \Map_{\QC^\vee_R}(\upiota(M),R)$ of coherent \mbox{$R$-module spectra $M$.}
\end{enumerate}
If $R$ is discrete, then $\Mod_R$ can be obtained from the category of $\Ch_R$ of chain complexes by inverting quasi-isomorphisms. On the other hand, 
 $\QC^\vee_R$ is modelled by $\Ch_R$ with its 
  \textit{tame} \mbox{model structure 
 (cf.\ \cite{becker2014models}).} It has more cofibrations, but fewer weak equivalences,
than the usual model  {structure on $\Ch_R$.}  Indeed, a map of complexes $f:M\rightarrow N$ is a tame weak equivalence precisely if for all (possibly unbounded) complexes $P$ of finitely generated free $R$-modules, the induced map of complexes $\Hom(P,M) \rightarrow \Hom(P,N)$ is a quasi-isomorphism.

\subsection*{Pro-coherent symmetric sequences}
Classically, the Koszul dual $\KD(\cat{O})$ of an   augmented  $\infty$-operad $\cat{O}$ over $R$ is formed in two steps: first, we form the bar construction $\Barr({1}, \cat{O},{1})$, then we take the \vspace{2pt} $R$-linear dual to obtain   $\KD(\cat{O})$.

We refine the second step by taking `continuous duals' of symmetric sequences.
For this, 
  we introduce the $\infty$-category $\sSeq^\vee_R$ of \textit{pro-coherent symmetric \mbox{sequences}} in \Cref{procosseq}, which is the home of continuous   linear duals of    \mbox{symmetric sequences.} 

The $n^{th}$ term of a pro-coherent symmetric sequence is a pro-coherent $R[\Sigma_n]$-module. In particular, 
pro-coherent  modules are just pro-coherent symmetric sequences  concentrated in \mbox{degree $0$.}

If $R$ is an ordinary  ring, we  will  model    $\sSeq^\vee_R$ by equipping the  category $\mathbf{sSeq}_R$ of  {symmetric sequences in $\Ch_R$}
with the  {\textit{tame model structure}} in  \Cref{sseqtame}.\vspace{-2pt}

\begin{warning} 
Note that $\sSeq^\vee_R$ is usually \textit{not}  equivalent to symmetric sequences in $\QC^\vee_R$, as the spaces $B\Sigma_n$ have infinite \mbox{homological dimension.}
\end{warning}

\subsection*{Pro-coherent composition}
Let $\Op^{\aug}_R$ denote the $\infty$-category of augmented \mbox{$\infty$-operads}, i.e.\ augmented  algebra objects in  $(\sSeq_R, \circ)$, 
the monoidal \mbox{$\infty$-category} of symmetric sequences  equipped  with the\vspace{3pt}   \mbox{{usual composition product}  $\circ$.}

For $\cat{O}\in \Op^{\aug}_R$, the bar construction ${1} \circ_{\cat{O}} {1} = \Barr({1}, \cat{O},{1})$ admits a coherently  coassociative comultiplication (cf.\ \cite[Section 4.3]{lurie2011derivedX} or Section \ref{sec:refined koszul}), which is informally given by $${1} \circ_{\cat{O}} {1} \simeq {1} \circ_{\cat{O}}
\cat{O}
\circ_{\cat{O}} {1} \rightarrow {1} \circ_{\cat{O}}
{1}
\circ_{\cat{O}} {1} \simeq  ({1} \circ_{\cat{O}} {1}) \circ ({1} \circ_{\cat{O}} {1}).$$
To  continuously dualise this map, we  construct a  monoidal product \mbox{on $\sSeq^\vee_R$} and prove in Propositions \ref{prop:pro-coh composition} and \ref{prop:formula for pro-coh comp prod}:
\begin{proposition}[Pro-coherent composition product]
Let $R$ be a coherent \mbox{$\mathbb{E}_\infty$-ring.} Then $\sSeq^\vee_R$ admits a monoidal structure $\circ$, the \textit{pro-coherent composition product},  which preserves small colimits in the first and sifted colimits in the second variable.\\
If $X, Y$ are  
`continuous duals' of 
symmetric sequences that are almost perfect, we have 
\begin{equation} \label{compositionformula} X \  {\circ} \ Y  \simeq \bigoplus_{n} \Big(X_n \otimes  Y^{\otimes n}\Big){}^{h\Sigma_n}.\end{equation}
\end{proposition}

\begin{remark} If $Y\in \sSeq^\vee_R$ is concentrated in degree zero, then so is $X \ {\circ} \ Y$ for any $X\in \sSeq^\vee_R$. Hence $\QC^\vee_R$ is left-tensored over 
$(\sSeq^\vee_R, \circ\big)$; the   action 
$
\big(\sSeq^\vee_R, \circ\big)\curvearrowright \QC^\vee_R
$ \mbox{preserves sifted colimits.}
\end{remark}
For $R$ a discrete coherent ring, we will give an explicit model for this composition product:
\begin{theoremn}[\ref{thm:chain model for exotic composition} \normalfont{(Point-set model for pro-coherent ${\circ}$\hspace{3pt})}] 
The composition product 
\begin{equation}\label{explicitcompositionformula} X  \boldsymbol{\circ} Y =  \bigoplus_{n} \Big(X_n \otimes Y^{\otimes n}\Big)_{\Sigma_n}  \end{equation}
on the model category $\mathbf{sSeq}_R$ induces a monoidal structure on its $\infty$-categorical localisation. The resulting monoidal $\infty$-category is equivalent to $\sSeq^\vee_R$ with the monoidal structure $\circ$ of \eqref{compositionformula}.
\end{theoremn}

\begin{remark}
It is somewhat surprising that formula \eqref{explicitcompositionformula} agrees with formula \eqref{compositionformula} on continuous duals of 
suitably finite symmetric sequences, as  \eqref{explicitcompositionformula} involves strict orbits while \eqref{compositionformula} involves homotopy fixed points. 
This phenomenon relies on two  facts:
first, invariants and  coinvariants agree on projective $R[\Sigma_n]$-modules via the norm;
second, on bounded above complexes of finite projective $R[\Sigma_n]$-modules, invariants and homotopy fixed points are equivalent.

\end{remark}

\subsection*{Divided power operads and their algebras} We define  a new notion of \mbox{$\infty$-operad:}
\begin{definition}[Divided power operads]
Let $R$ be a coherent $\EE_\infty$-ring. A \textit{PD $\infty$-operad} is an algebra object  {in $(\sSeq^\vee_{R}, \circ)$.} Write $\Op^{\pd}_R$ and   $\Op^{\aug,\pd}_R$  for the $\infty$-categories  of  PD $\infty$-operads and augmented PD $\infty$-operads, respectively.
\end{definition}

If $R$ is discrete, consider the category $\mathbf{Op}_R$ of ordinary operads in chain complexes over $R$; these are often called dg-operads. In
 \Cref{thm:chain models for PD operads}, we  prove that inverting tame weak equivalences in $\mathbf{Op}_R$ gives the $\infty$-category $\Op^{\pd}_R$.
 
 In \Cref{thm:chain models for PD algebras}, we show that if $\mathbf{P}\in \mathbf{Op}_R$ is a dg-operad with tamely cofibrant underlying symmetric sequence, then the $\infty$-category $\Alg_{\PP}(\QC^\vee_R)$  of pro-coherent algebras over the corresponding PD $\infty$-operad $\PP$ 
can be obtained from $\mathbf{P}$-algebras in chain complexes by inverting tame weak equivalences.

\subsection*{Refined Koszul Duality}

Using that the continuous $R$-linear  
duality functor $(\sSeq_R, \circ)^{op} \rightarrow (\sSeq^\vee_R, \circ)$ is lax monoidal, we offer a refinement of the classical operadic Koszul duality construction of Ginzburg--Kapranov \cite{GinzburgKapranov:KDO}, Fresse \cite{fresse346koszul},  Salvatore \cite{salvatore1998configuration}, \mbox{and Ching \cite{ching2005bar}:}

\begin{theoremn}[\ref{thm:koszul duality for PD operads}  \normalfont{(Refined Koszul duality for operads)}]
Let $R$ be a coherent $\mathbb{E}_\infty$-ring spectrum. Then there is a commuting diagram of $\infty$-categories
$$\begin{tikzcd}
\Op^{\pd, \aug}_R\arrow[r, "\KD^{\pd}"] & \Op^{\pd, \aug, \op}_R\arrow[d, "\upupsilon"]\\
\Op^{\aug}_R\arrow[r, "\KD"]\arrow[u, "\upiota"] & \Op^{\aug, \op}_R
\end{tikzcd}$$
where the bottom functor sends an augmented $\infty$-operad to its classical Koszul dual $\infty$-operad, given by the Spanier--Whitehead dual of its bar construction.
\end{theoremn}

\begin{example}[Partition Lie algebras] Over a field $R=k$, the refined Koszul dual of the nonunital $\EE_\infty$-operad $\EE_{\infty,k}^{\mm{nu}}$ is a divided power $\infty$-operad \mbox{$\Lie_{k, \EE_\infty}^{\pi}=\KD^{\pd}(\EE_{\infty,k}^{\mm{nu}})$} which induces the spectral partition Lie algebra monad   of \cite[
Construction 1.20]{brantner2019deformation} \mbox{on $\QC_k^\vee \simeq \Mod_k$.}
In fact, our setup gives a definition of spectral partition Lie algebras over any coherent $\EE_\infty$-ring; their relation to deformation theory is the subject of future work.
\end{example}

We also offer a divided power refinement of Koszul duality for algebras:

\begin{theoremn}[\ref{thm:koszul duality for algebras} \normalfont{(Refined  Koszul duality for algebras)}] 
Let $R$ be a coherent $\mathbb{E}_\infty$-ring and $\PP$ an augmented $\infty$-operad over $R$. There is a commuting diagram  
$$\begin{tikzcd}
\Alg_{\PP}(\QC^\vee_R)\arrow[r, "\KD^{\pd}"] & \Alg_{\KD^{\pd}(\PP)}(\QC^\vee_R)^{\op}\arrow[d, "\upupsilon"]\\
\Alg_{\PP}(\Mod_R)\arrow[r, "\KD"]\arrow[u, "\upiota"] & \Alg_{\KD(\PP)}(\Mod_R)^{\op}
\end{tikzcd}$$
where the bottom functor sends a $\PP$-algebra $A$ its classical Koszul dual algebra, given by the Spanier--Whitehead dual of its bar construction.
\end{theoremn}

Given a dg-operad $\mathbf{P}$ over a coherent ring $R$, one can construct a new dg-operad $\KD(\mathbf{P})$ via the chain-level bar construction, cf.\ \Cref{chainlevelbar}.
This generalisation of quadratic duality is due to Ginzburg--Kapranov \cite{GinzburgKapranov:KDO} and studied in depth by Getzler--Jones \cite{GJ} and Fresse \cite{fresse346koszul} (see also \cite{LV} for a textbook account).\vspace{2pt}

\begin{theoremn}[\ref{chainkoszul} \normalfont{(Chain models for Koszul duality)}]
Fix a coherent ring  $R$. Let $\mathbf{P}$ be an augmented dg-operad over $R$ 
with tamely cofibrant 
  underlying symmetric sequence and let $\PP$ denote the corresponding PD $\infty$-operad. 

Then the chain-level dual operad $\KD(\mathbf{P})$ is a model for the Koszul dual PD $\infty$-operad $\KD^{\pd}(\PP)$.
Furthermore, inverting tame weak equivalences gives rise to  a  
commuting square of $\infty$-categories in which the vertical functors are equivalences
$$\begin{tikzcd}[column sep=3pc]
\dgAlg_{\mathbf{P}}[W^{-1}_\mm{tame}]\arrow[r, "\KD_{\mathbf{P}}"]\arrow[d, "\simeq"{swap}] & \dgAlg_{\KD(\mathbf{P})}[W^{-1}_{\mm{tame}}]^{\op}\arrow[d, "\simeq"]\\
\Alg_{\PP}(\QC^\vee_R)\arrow[r, "\KD^{\pd}_{\PP}"] & \Alg_{\KD^{\pd}(\PP)}(\QC^\vee_R)^{\op}.
\end{tikzcd}$$
\end{theoremn}
 
\subsection*{Explicit Models for Spectral Partition Lie Algebras}
 \Cref{chainkoszul} lets us give explicit models for spectral partition Lie algebras, using the following notation:
\begin{notation}[Nondegenerate sequences]
Given $r\geq 1$,   a \emph{nondegenerate sequence} in $\ul{r}$ is an (ordered) sequence $\mb{u}=(u_1, \dots, u_{r+d})$ of elements in $\ul{r}=\{1, \dots, r\}$ such that each $1, \dots, r$ appears in the sequence and $u_{\alpha}\neq u_{\alpha+1}$ for all $\alpha$. If $\mb{u}$ does not exhaust all of $\ul{r}$ or if $u_{\alpha}=u_{\alpha+1}$ for some $\alpha$, then $\mb{u}$ is said to be \emph{degenerate}.
\end{notation}
For the sake of exposition, we suppress signs; they are specified in \mbox{the main text:}
\begin{definitionn}[\ref{explicitone} \normalfont{(Spectral partition $L_\infty$-algebra)}]
Let $R$ be a discrete coherent ring. A \textit{spectral partition $L_\infty$-algebra} is   a chain complex of $R$-modules $\mf{g}$, together with the following algebraic structure: given $r\geq 2$ and a nondegenerate sequence $\mb{u}=(u_1,\dots, u_{r+d})$, \mbox{there is an operation}
$$\begin{tikzcd}
\{-, \dots, -\}_{\mb{u}}\colon \mf{g}^{\otimes r}\arrow[r] & \mf{g}
\end{tikzcd}$$
of homological degree $-1-d$. Furthermore, these operations satisfy:
\begin{enumerate}[label=(\alph*)]
\item \emph{Equivariance.} For every $\sigma\in \Sigma_r$, let $\sigma(\mb{u})=\big(\sigma(u_1), \dots, \sigma(u_{r+d})\big)$. Then
$$
\{x_1, \dots, x_r\}_{\sigma(\mb{u})} = \pm \{x_{\sigma^{-1}(1)}, \dots, x_{\sigma^{-1}(r)}\}_{\mb{u}}
$$ 
\item \emph{Differential.} For each nondegenerate sequence $\mb{u} =(u_1,\dots, u_{r+d}) $ in $\ul{r}$ and each tuple $x_1, \dots, x_r\in \mf{g}$, we have
\begin{align*}
\partial\{x_1, \dots, x_r\}_{\mb{u}} &= \sum_{i=1}^r  \pm \{x_1, \dots, \partial(x_i), \dots, x_r\}_{\mb{u}}\\
&+ \sum_{\alpha=1}^{r+d+1}\sum_{\substack{v=1\vspace{1pt} \\ v \neq u_{\alpha - 1}, u_{\alpha}}}^r \pm  \{x_1, \dots, x_r\}_{\mb{u}_+=(u_1, \dots, u_{\alpha-1}, v, u_{\alpha}, \dots, u_{r+d})}\vspace{-15pt}\\
&+\sum_{k=2}^{r-2}\sum_{\sigma\in \mm{UnSh}_{\mathbf{u}}(k, r-k)} \hspace{-15pt}\pm  \big\{\{x_{\sigma(1)}, \dots, x_{\sigma(k)}\}_{\mb{v}(k, \sigma)}, x_{\sigma(k+1)}, \dots, x_{\sigma(r)}\big\}{}_{\mb{w}(k, \sigma)}
\end{align*} 
In the third row, we sum  over the set $ \mm{UnSh}_{\mathbf{u}}(k, r-k)$ of $(k, r-k)$-unshuffles $\sigma$ which are compatible with $\mathbf{u}$ in the following sense: if we 
decompose the subsequence of $\mb{u}$ consisting of all $u_i\in \{\sigma(1), \dots, \sigma(k)\}$ into intervals 
$$
\mb{u}_1=\big(u_{\alpha(1)}, u_{\alpha(1)+1}, \dots, u_{\alpha(1)+\beta(1)}\big),\qquad\dots, \qquad \mb{u}_n=\big(u_{\alpha(n)}, u_{\alpha(n)+1}, \dots, u_{\alpha(n)+\beta(n)}\big)
$$
separated  in $\mb{u}$  by elements in $\{\sigma(k+1), \dots, \sigma(r)\}$, then\vspace{3pt}  \mbox{$u_{\alpha(i)+\beta(i)}= u_{\alpha(i+1)}$ for all $i$.}\\
Define $\mb{v}(k, \sigma)$ to be the sequence in $\ul{k}$ given by applying $\sigma^{-1}$ to the sequence\vspace{-2pt}
$$
\big(u_{\alpha(1)}, \dots, u_{\alpha(1)+\beta(1)-1}, u_{\alpha(2)}, \dots, u_{\alpha(i)+\beta(i)-1}, u_{\alpha(i)}, \dots, u_{\alpha(n)+\beta(n)}\big).
$$
Define $\mb{w}(k, \sigma)$ as the sequence of elements of $\ul{r-k+1}$ obtained from $\mb{u}$ by replacing each $\sigma(k+i)$ (for $i=1, \dots, r-k$) by $1+i$ and replacing each of the intervals $\mb{u}_1, \dots, \mb{u}_n$ by a single copy of $1$.

If $\mb{v}(k, \sigma)$ or $\mb{w}(k, \sigma)$ is degenerate, the corresponding \mbox{term is zero.}
\end{enumerate}
\end{definitionn}

\begin{theoremn}[\ref{chaintheoremone} \normalfont{(Chain models for spectral partition Lie algebras I)}]
Inverting tame weak equivalences on the category of spectral partition $L_\infty$-algebras gives the $\infty$-category $\Alg_{{\Lie}_{R, \EE_\infty}^\pi}(\QC^\vee_R).$
In particular, when $R=k$ is a field,  localising spectral partition $L_\infty$-algebras  at the weak equivalences gives the $\infty$-category of partition Lie algebras from \cite[Definition 5.32]{brantner2019deformation}.
\end{theoremn}

We also provide a second model as algebras over a  certain
dg-operad  $\mathbf{Lie}_{R, \EE_\infty}^\pi$  with tamely cofibrant underlying symmetric sequence.
Two ingredients are needed:
\begin{enumerate}
\item The usual (shifted) {Lie operad} $\mathbf{Lie}_{R}^s$;
\item The {PD surjections operad}   $\Surj_R^\vee$. 
\end{enumerate}

The dg-operad $\mathbf{Lie}_{R}^s$ is familiar. In weight $r$, $\mathbf{Lie}_{R}^s(r)$ sits in homological degree $1-r$, where it is generated 
 by Lie words $w(c_1,\ldots, c_r)$ in $r$ letters involving each letter exactly once, modulo Jacobi identity and
  \mbox{antisymmetry.} For example, $\mathbf{Lie}_{R}^s(3)$ is a free $R$-module   generated by the Lie words $[c_1,[c_2,c_3]], [c_3,[c_1,c_2]]$\vspace{3pt} in \mbox{degree $-2$.}

The dg-operad $\Surj_R^\vee$ constructed in Appendix \ref{app:surjections cooperad}  is an analogue of the Barratt--Eccles operad (and the surjections operad  \cite{mcclure2003multivariable})
 for (nonunital)
$\EE_\infty$-$R$-algebras with divided powers. 
The Koszul dual of $\Surj_R^\vee$ is the shifted Lie operad, see Theorem \ref{thm:KD of pd surj op}.
Here $\Surj_R^\vee(r)$ is  a coconnective chain complex, which in homological degree $-d$ is a 
free $R$-module  spanned by all nondegenerate
sequences $(u_1,\ldots, u_{r+d})$  in $\underline{r}$.

\begin{example}
Let $k$ be a field. Then each  $\Surj_R^\vee(r)$ is a chain complex of finitely generated free $k[\Sigma_r]$-modules. For any bounded above chain complex $V$, we have   \vspace{-2pt}
$$
\Surj_R^\vee\circ V = \bigoplus_r \big(\Surj_R^\vee(r)\otimes V^{\otimes r}\big)_{\Sigma_r} \cong 
\bigoplus_r \big(\Surj_R^\vee(r)\otimes V^{\otimes r}\big)^{\Sigma_r} 
\simeq \bigoplus_r (V^{\otimes r})^{h\Sigma_r} \vspace{-4pt}
$$  
The second isomorphism uses that 
the norm map is an isomorphism on 
finitely generated free $k[\Sigma_r]$-modules, and the third that $\Surj_R^\vee(r)\otimes V^{\otimes r}$ is $\Sigma_r$-fibrant.\end{example}

We then define the dg-operad $\mathbf{Lie}_{R, \EE_\infty}^\pi$ as a levelwise tensor product:
$$\mathbf{Lie}_{R, \EE_\infty}^\pi := \mathbf{Lie}_{R}^s \otimes_{\lev} \Surj_R^\vee. $$ 

We spell out the resulting structure of a $\mathbf{Lie}_{R, \EE_\infty}^\pi$-algebra in \Cref{explicittwo}, and deduce:
\begin{theorem}[Chain models for spectral partition Lie algebras II, cf.\  \Cref{spectralpliemodel}]
Inverting tame weak equivalences on the category of dg-algebras over  $\mathbf{Lie}_{R, \EE_\infty}^\pi$ gives an equivalence 
$$\dgAlg_{\mathbf{Lie}_{R, \EE_\infty}^\pi}[W^{-1}_\mm{tame}] \simeq \Alg_{{\Lie}_{R, \EE_\infty}^\pi}(\QC^\vee_R).$$
In particular, when $R=k$ is a field,  localising $\mathbf{Lie}_{k, \EE_\infty}^\pi$-algebras  at  weak equivalences gives the $\infty$-category of partition Lie algebras from \cite[Definition 5.32]{brantner2019deformation}.
\end{theorem}

\subsection*{The Derived Setting}
There is a second, more algebraic, generalisation of \mbox{classical} algebraic geometry based on simplicial commutative rings (rather than connective $\EE_\infty$-rings). Here, formal moduli are controlled by \textit{derived} \mbox{partition Lie algebras.} 
To construct point-set models for these objects, we implement the above programme in a more genuine setting. We  briefly outline our main results, but will leave detailed statements to the main text.

Let $R$ be a coherent commutative ring. Given $n\geq 0$, write $\underline{R}$ for the constant $\Sigma_n$-Mackey functor corresponding to $R$, thought of as a genuine $\Sigma_n$-spectrum, 
and consider the $\infty$-category $ \Mod_{\ulR}^{\Sigma_n}$ of $\underline{R}$-modules in $\Sp^{\Sigma_n}$.

In \Cref{def:derived sym seq}, we assemble these  into the $\infty$-category of 
 \textit{derived symmetric sequences} $\sSeq^{\gen}_{\ulR},$
which  admits a composition product $\circ$, cf.\ \Cref{con:derived products}. Passing to algebra objects leads to the    $\infty$-category $\Op^{\gen}_{\ulR}$ \mbox{of  \textit{derived $\infty$-operads}  over $R$.}
Identifying $\Mod_R$ with  derived sequences   in degree zero, we see that $\Mod_R$ is left-tensored over $\sSeq^{\gen}_{\ulR}$; for any 
derived $\infty$-operad $\cat{O}$, we obtain an $\infty$-category  of $\mathcal{O}$-algebras $\Alg_{\cat{O}}(\Mod_R)$.

The bar construction of an augmented derived $\infty$-operad $\cat{O}$ is again equipped with a comultiplication. To dualise it `continuously', we introduce the $\infty$-category $\sSeq_{\ulR}^{\gen,\vee}$ of \textit{pro-coherent  derived symmetric sequences} in \Cref{derivedProCoh}. 

   \Cref{not:restricted composition} gives  $\sSeq_{\ulR}^{\gen,\vee}$   a sifted-colimit-preserving product 
 $ \ol{\circ}$  satisfying 
$$X \ \ol{\circ} \ Y  \simeq \bigoplus_{r} \Big(X(r) \otimes Y^{\otimes r}\Big){}^{\Sigma_r}.$$

Writing
$\Op_{\ulR}^{\gen, \pd}  = \Alg(\sSeq_{\ulR}^{\gen,\vee},\ol{\circ})$ for the  $\infty$-category of \textit{derived  divided power operads}, we construct a Koszul duality functor
$$ \KD^{\pd} : \Op_{\ulR}^{\gen, \aug} \rightarrow \Op_{\ulR}^{\gen, \pd, \aug, \op}$$
in \Cref{kdforderived}, and give a version for algebras in \Cref{kdforderivedalgebras}.

Using symmetric sequences in simplicial-cosimplicial $R$-modules, we  construct point-set models for the monoidal $\infty$-categories $(\sSeq_{\ulR}^{\gen},\circ)$ and $(\sSeq_{\ulR}^{\gen,\vee},\ol{\circ}) $ in \Cref{thm:simp-cosimp model for exotic composition}. This allows us to model derived (PD) $\infty$-operads  in \Cref{thm:simp-cosimp models for derived PD operads}, and their algebras in \Cref{thm:simp-cosimp models for derived algebras} and \Cref{rem:rectification algebras non-coherent}.

\begin{example}[Derived commutative rings]
	The unique operad $\Com^{ }$ in sets with \mbox{$\Com^{ }(r) = \ast$} for all $r$ gives rise to a derived $\infty$-operad $\mathrm{Com}^{ }_R$ such that $\Alg_{\mathrm{Com}^{ }_R}(\Mod_R)_{\geq 0}$ is   the $\infty$-category of 
	simplicial commutative (i.e.\ animated) $R$-algebras.
	The  $\infty$-category  $\Alg_{\mathrm{Com}^{ }_R}(\Mod_R)$ is equivalent to the $\infty$-category  of  \textit{derived $R$-algebras}, which was also studied  by Bhatt-Mathew, Raksit \cite{raksit2020hochschild} and others, using the technique of extending monads via Goodwillie calculus from \cite[Section 3]{brantner2019deformation}.
	Our rectification result immediately implies that derived $R$-algebras are modelled by simplicial-cosimplicial $R$-algebras (cf.\ \Cref{derivedringssimplicialcosimplicial}). These have appeared in the literature much earlier, for instance in the 
	work of
	Illusie \cite[Section I-4]{Illusie} \mbox{and  Kaledin \cite[Section 3]{Kaledin}.}

\end{example}

\begin{example}
	
	The non-unital version of  $\mathrm{Com}$ is denoted $\mathrm{Com}^{\mm{nu}}$; it differs from $\mathrm{Com}^{}$  in that $\mathrm{Com}^{\mm{nu}}(0)=\emptyset$.
	For $k$   a field,  $\KD^{\pd}(\Com^{\mm{nu}}_k)$ defines derived \mbox{partition Lie algebras} in the sense of \cite[Construction 1.10]{brantner2019deformation}.
\end{example}

We then give an explicit construction of the refined Koszul duality functor for derived $\infty$-operads in  \Cref{thm:simp-cosimp KD}.  
In \Cref{def:simplicial-cosimplicial part lie}, we construct a cosimplicial restricted operad 
$$
\partLieR.
$$
that models  $\KD^{\pd}(\Com^{\mm{nu}}_k)$.
Here  $\partLieR(r)^d$ is dual to the set of \textit{nested chains of partitions} of $\underline{r}$ of length $d$, i.e.\ the \mbox{set of pairs}
$$ \big(\sigma, S\big)=\big([\hat{0}=x_0<\dots <x_t=\hat{1}], S_0\subseteq \dots\subseteq S_d\big)$$
where $\sigma$ is a nondegenerate chain of partitions of $\ul{r}$ and $S_0\subseteq \dots\subseteq S_d=\{0, \dots, t\}$ is an increasing set of subsets. We allow $d=-1$ in this definition.

The explicit description of derived Koszul duality allows us to 
construct explicit  point-set models 
for derived partition Lie algebras  in  \Cref{thm:simp-cosimp partition lie} of \mbox{the main text:}
\begin{theorem}[Simplicial-cosimplicial models for partition Lie algebras]
Inverting tame weak equivalences on the category of simplicial-cosimplicial restricted $\partLieR$-algebras   induces an equivalence of $\infty$-categories
$$\begin{tikzcd}
\mathbf{Alg}_{\partLieR}^{\simcos,\res}[W_\mm{tame}^{-1}]\arrow[r, "\simeq"] & \Alg_{\Lie^\pi_{R,\Delta}}^{\gen,\pd}(\QC^\vee_R).
\end{tikzcd}$$
Hence when $R=k$ is a field, the localisation of the category of simplicial-cosimplicial restricted algebras over $\mathbf{Lie}_{k,\Delta}^\pi$ at the weak equivalences is equivalent to the $\infty$-category of partition Lie algebras from \cite[Definition 5.47]{brantner2019deformation}.\end{theorem}
 
In \Cref{derivedPLAexplicit}, we describe simplicial-cosimplicial restricted ${\partLieR}$-algebras over a field $R=k$ as simplicial-cosimplicial modules with  explicit operations satisfying relations we specify.

\subsection{Outline}
We provide  a brief outline of the structure of the paper. In the first half, we give an $\infty$-categorical treatment of (derived) PD $\infty$-operads and their algebras; in particular, we define the (derived) PD $\infty$-operad whose algebras are spectral (derived) partition Lie algebras. The second half of the paper  provides explicit point-set models for these $\infty$-categorical objects.

We will start by collecting some results on $\infty$-categories of pro-coherent modules in \textbf{Section \ref{sec:functors on pro-coh}}. Most importantly, we show that a polynomial functor between (coherent) additive $\infty$-categories admits a natural extension to a sifted-colimit-preserving functor between the corresponding $\infty$-categories of pro-coherent modules.

In \textbf{Section \ref{sec:PD operads and koszul}}, we use this machinery to develop the theory of PD $\infty$-operads and their algebras. In particular, this leads to a refinement of the usual Koszul duality for operads (Section \ref{sec:refined koszul}). In Section \ref{sec:refined koszul} we also provide a few more details on the $\infty$-categorical bar construction, to fill in some gaps in the literature (as pointed out in \cite{dancohenhorev2022koszul}).
Section \ref{sec:derived operads} discusses the derived analogues of $\infty$-operads and PD $\infty$-operads over a simplicial commutative ring.

\textbf{Section \ref{sec:chain models}} provides chain models for PD $\infty$-operads over a discrete coherent ring. In particular, we describe the tame homotopy theory of chain complexes that is used to present the $\infty$-categories of pro-coherent modules and symmetric sequences. Using this, we give chain complex models for spectral partition Lie algebras.

Similarly, the $\infty$-categories of derived PD $\infty$-operads and their algebras admit concrete models in terms of simplicial-cosimplicial $R$-modules, which are discussed in \textbf{Section \ref{sec:simp-cosimp point-set models}}. This allows for an explicit description of derived partition Lie algebras in terms of simplicial-cosimplicial algebras with divided power operations.

Finally, \textbf{Appendix \ref{app:surjections cooperad}} gives a detailed construction of the PD surjections operad; this is used in Section \ref{sec:chain models} to provide a chain model for the spectral partition Lie PD $\infty$-operad and also to produce a cofibrant model for the Lie operad. \textbf{Appendix \ref{app:free algebras}} describes free algebras in monoidal $\infty$-categories where the tensor product does not preserve colimits in the second variable (such as the composition product).

\subsection{Acknowledgements}
The authors wish to thank Nathan Adlam, Damien Calaque, Jacob Lurie, Zhouhang Mao,  Denis Nardin,  Pelle Steffens, and Bruno Vallette for helpful conversations related to this paper, and the anonymous referee for their  valuable suggestions.

The first author would also like to thank Merton College and the Mathematical Institute at  Oxford University for their support, 
as well as the Royal Society  (University Research Fellowship   URF$\backslash$R1$\backslash$211075) and the Centre national de la recherche scientifique (CNRS).
The second author acknowledges support by the grant ANR-20-CE40-0016 HighAGT.
The third author has received funding from the European Research Council (ERC) under the European Union’s Horizon 2020 research and innovation programme (grant agreement No 768679).

 \ \\ \\

\section{Functors on pro-coherent modules}
\label{sec:functors on pro-coh}
The main goal of this section is to study sifted-colimit-preserving functors on pro-coherent modules, which will be a key ingredient for our subsequent treatment of  refined Koszul duality via divided power $\infty$-operads. 
\subsection{Pro-coherent modules}\label{sec:pro-coh} We begin by discussing  pro-coherent modules over {additive $\infty$-categories}.
This general framework will allow us to give a uniform treatment of 
several examples of interest, including pro-coherent modules over a ring and pro-coherent symmetric sequences. \vspace{4pt} First, we  recall several \mbox{preliminary definitions.}
 
An $\infty$-category $\cat{A}$ is  \textit{additive}
if it admits  finite products and coproducts and $h\cat{A}$  is an additive category, cf.\ \cite[Definition C.1.5.1]{lurie2016spectral}. This implies that products and coproducts agree; we call them `direct sums' and denoted them by $\oplus$.

\begin{definition}[Modules over additive $\infty$-categories]\label{def:module cats}
Given a small additive \mbox{$\infty$-category}
$\cat{A}$, the $\infty$-category of (left) \emph{$\cat{A}$-modules}  is the initial presentable  stable $\infty$-category 
receiving a functor from $\cat{A}$ that preserves finite direct sums:
\begin{equation}\label{eq:yoneda}\begin{tikzcd}
\cat{A}\arrow[r] & \Mod_{\cat{A}}.
\end{tikzcd}\end{equation}
Stabilising \cite[Proposition 5.3.6.2]{lurie2009higher}, we  can identify $\Mod_{\cat{A}}$ with the full subcategory of $\Fun(\cat{A}^{\op}, \Sp)$ spanned by the functors $M\colon \cat{A}^{\op}\rt \Sp$ preserving finite direct sums. The fully faithful   universal functor \eqref{eq:yoneda} 
then arises from the Yoneda embedding, using that all mapping spaces in $\cat{A}$ are grouplike $\EE_\infty$-spaces in an  essentially unique way, cf.\ \cite[Section C.1.5]{lurie2016spectral}.
\end{definition}
\begin{example}\label{mainex}
Let $R$ be a connective $\mathbb{E}_1$-ring spectrum 
in the sense of \cite[Definition 7.2.4.16]{lurie2014higher} and consider the additive $\infty$-category $\Vect_R^{\omega}$ of finitely generated  free left $R$-modules  of the form $R^{\oplus n}$. The inclusion $\Vect_R^\omega\hookrightarrow \Mod_{\Vect_R^{\omega}}$ can then be identified with the usual inclusion $\Vect_R^\omega\hookrightarrow \Mod_R$ into the $\infty$-category of left $R$-modules. For a general additive $\infty$-category $\cat{A}$, we will therefore refer to objects in $\cat{A}$ as \mbox{\emph{finitely generated  free $\cat{A}$-modules}.}
\end{example}
\begin{remark}
Our formalism is not well adapted to non-connective rings, as we only remember the mapping spaces (not spectra) between objects in $\Vect_R^{\omega}$.
\end{remark}

\begin{example}\label{ex:sums of additive cats}
Let $\cat{A}_i$ be a set of small additive categories and write $\bigoplus_i \cat{A}_i\subseteq \prod_i \cat{A}_i$ for the full subcategory spanned by tuples of objects $V_i\in \cat{A}_i$ such that almost all $V_i$ are the zero object. Then $\bigoplus_i \cat{A}_i$ is additive and $\Mod_{\bigoplus_i \cat{A}_i}\simeq \prod_i \Mod_{\cat{A}_i}$.
\end{example}

\begin{definition}[(Co)connective modules]\label{def:classical t-structure on modules}
An $\cat{A}$-module is said to be \emph{(co)connective} if the corresponding functor $\cat{A}^{\op}\rt \Sp$ takes values in (co)connective spectra. In particular, the essential image of $\cat{A}\hookrightarrow \Mod_{\cat{A}}$ consists of connective $\cat{A}$-modules.
This defines a $t$-structure $(\Mod_{\cat{A}, \geq 0}, \Mod_{\cat{A}, \leq 0})$ on $\Mod_{\cat{A}}$.
\end{definition}

\begin{example}[Opposite $\infty$-category]\label{ex:opposite}
If $\cat{A}$ is an additive $\infty$-category, then \mbox{so is  $\cat{A}^{\op}$.} 
We   write $V^\vee \in \cat{A}^{\op}$ for the object corresponding to $V \in \cat{A}$.
One can identify $\Mod_{\cat{A}^{\op}}$ with the dual stable presentable $\infty$-category of $\Mod_{\cat{A}}$, i.e.\ the full subcategory of $\Fun(\Mod_{\cat{A}}, \Sp)$ spanned by the left adjoints. We will denote the induced pairing between left and right $\cat{A}$-modules by 
$$
-\otimes_{\cat{A}}-\colon \Mod_{\cat{A}^{\op}}\times \Mod_{\cat{A}}\rt \Sp.
$$
For $V^\vee \in \cat{A}^{\op}$ and $W\in \cat{A}$, the spectrum $V^\vee\otimes_{\cat{A}} W$ is simply the spectrum corresponding to the grouplike $\mathbb{E}_\infty$-space $\Map_{\cat{A}}(V, W)$. In these terms, a right $\cat{A}$-module $M$ is connective if and only if $M\otimes_{\cat{A}} -$ is a right $t$-exact functor.
\end{example}

We   introduce several standard finiteness conditions  in this  {generalised setting:}

\begin{definition}[Finiteness conditions]
Let $\cat{A}$ be a small additive $\infty$-category. An $\cat{A}$-module $M$ is said to be:
\begin{enumerate}
\item \emph{perfect} if it is a compact object in $\Mod_{\cat{A}}$;
\item \emph{almost perfect} if for each $n$, there exists a perfect $\cat{A}$-module $N$ and a map $N\rt M$ with
$n$-connective cofibre;
\item \emph{coherent} if it is almost perfect and eventually coconnective, which means that $M$ belongs to $ \Mod_{\cat{A}, \leq N}$ for some $N\gg 0$.
\end{enumerate}
We will denote the full subcategories of $\Mod_{\cat{A}}$ spanned by the perfect, almost perfect and coherent $\cat{A}$-modules by $\Perf_{\cat{A}}$, $\APerf_{\cat{A}}$ and $\Coh_{\cat{A}}$ respectively.
\end{definition}
\begin{remark}\label{rem:connective almost perfect as realisations}
The full subcategory $\Perf_{\cat{A}, \geq 0}\subseteq \Mod_{\cat{A}}$ of connective perfect {$\cat{A}$-modules} is the smallest subcategory of $\Mod_{\cat{A}}$ that contains $\cat{A}$ and is closed under finite colimits and retracts. 
Similarly, the full subcategory $\APerf_{\cat{A}, \geq 0}\subseteq \Mod_{\cat{A}}$ of connective almost perfect $\cat{A}$-modules is the smallest subcategory that contains $\cat{A}$ and is closed under geometric realisations. In fact, every connective almost perfect $\cat{A}$-module $X$ can be obtained as the geometric realisation of a simplicial object $X_\bullet$ in $\cat{A}$, and the cofibre of the natural map $X_0 \rightarrow |X_\bullet|\simeq X$  is always $1$-connective. 
\end{remark}

\begin{definition}[Coherence]
An additive $\infty$-category $\cat{A}$ is said to be  \emph{left coherent} if the $t$-structure on $\Mod_{\cat{A}}$ restricts to a $t$-structure on  $\APerf_{\cat{A}}$. 
We will say that $\cat{A}$ is \emph{coherent} if both $\cat{A}$ and $\cat{A}^{\op}$ are left coherent.
\end{definition}
\begin{example}[Coherent $\EE_n$-rings]
If $R$ is a connective $\mathbb{E}_n$-ring spectrum as in \Cref{mainex}, then $\Vect_R^{\omega}$ is (left) coherent if and only if $R$ is a (left) coherent $\mathbb{E}_1$-ring spectrum in the sense of  \cite[Proposition 7.2.4.18]{lurie2014higher}.

\end{example}
\begin{lemma}\label{rem:detecting coherent categories}
Let $f\colon \cat{A}_0\rt \cat{A}$ be an additive functor between additive $\infty$-categories and let $f_!\colon \Mod_{\cat{A}_0}\leftrightarrows \Mod_{\cat{A}}\colon f^*$ be the induced adjoint pair. If $f^*$ detects equivalences and $f^*(\cat{A})\subseteq \APerf_{\cat{A}_0}$, then $X\in \Mod_{\cat{A}}$ is almost perfect if and only if $f^*(X)\in \Mod_{\cat{A}_0}$ is almost perfect. Since $f^*$ commutes with truncation, it then follows that $\cat{A}$ is coherent if $\cat{A}_0$ is coherent.
\end{lemma}
\begin{proof}
The functor $f^\ast$ sends almost perfect $\cat{A}$-modules to almost perfect $\cat{A}_0$-modules, because it preserves realisations and sends $\cat{A}$ into $\APerf_{\cat{A}_0}$.
On the other hand,  note that $f_!$ sends $\APerf_{\cat{A}_0}$ to $\APerf_{\cat{A}}$. If $f^\ast(X)$ belongs to $\APerf_{\cat{A}_0}$, we can write $X$ as a colimit of a simplicial diagram $\Barr_\bullet(f_!f^\ast, f_!f^\ast, X)$, which belongs to $\APerf_{\cat{A}}$.
\end{proof}
\begin{example}[Genuine equivariant spectra]\label{ex:genuine G-spectra}
Let $G$ be a finite group and write $A(G)$ for its spectral Burnside $\infty$-category, with objects given by finite $G$-sets and morphism spaces $\Map_{A(G)}(X, Y)$ given by the group completions of the $\mathbb{E}_\infty$-spaces of spans $X\lt Z\rt Y$ of $G$-sets (with disjoint union). Note that $\Mod_{A(G)}\simeq \Sp^G$ is the $\infty$-category of spectral Mackey functors \cite{barwick2017spectral}, or equivalently, genuine $G$-spectra \cite{GuillouMay2011, nardin2016}. Then $A(G)$ is a coherent additive $\infty$-category. Indeed, this follows by applying Lemma \ref{rem:detecting coherent categories} where $\cat{A}_0$ is the free additive $\infty$-category on the set of orbits $\{G/H\}$ and $f^*\colon \Mod_{A(G)}\rt \Mod_{\cat{A}_0}=\prod_{H<G} \Sp$ simply evaluates a spectral Mackey functor at $G/H$. The condition of  \Cref{rem:detecting coherent categories} follows from the fact that each $\Map_{A(G)}(X, Y)$ is an almost perfect spectrum (as it has finitely \mbox{generated homotopy groups).}
\end{example}
The following example is of key significance in  our treatment of \mbox{derived $\infty$-operads:}
\begin{notation}[Constant Mackey functors]
Given a finite group $G$ and an abelian group $A$,  let $\underline{A}\in \Sp^G$ be the Eilenberg--Mac Lane spectrum corresponding to 
the constant {Mackey functor on $A$}. Recall that this constant Mackey functor 
sends a finite $G$-set $X$ to the abelian group $\Map(X,A)^G  \cong \Map(X/G,A)$ consisting of 
of $G$-invarant functions $X\rightarrow A$; restriction maps correspond to precomposition and transfers to summation over fibres. In particular, $\underline{A}$ sends all orbits $G/H$ to the  Eilenberg--Mac Lane spectrum of $A$. This assignment sends direct sums of abelian groups to direct sums in $\Sp^{\sG}$, so taking its sifted-colimit-preserving extension provides a colimit-preserving functor $\Mod_{\ZZ, \geq 0} \rt \Sp^{\sG}; A \longmapsto  \underline{A}$ defined on  \mbox{connective $\ZZ$-module spectra.} 
\end{notation}
\begin{lemma}
The functor $\Mod_{\mathbb{Z}, \geq 0}\rt \Sp^{\sG}; A\longmapsto \ul{A}$ has a lax symmetric monoidal structure, where the symmetric monoidal structure on $\Sp^{\sG}$ is given by Day convolution.
\end{lemma}
\begin{proof}
Recall that for any $\infty$-category $\cat{C}$ and a presentable $\infty$-category $\cat{V}$, left Kan extension along the inclusion $i\colon \cat{C}\rt \mathcal{P}_{\Sigma}(\cat{C})$ into the sifted-colimit completion of $\cat{C}$ defines a fully faithful functor $i_!\colon \Fun(\cat{C}, \cat{V})\rt \Fun(\mathcal{P}_{\Sigma}(\cat{C}), \cat{V})$, whose essential image consists of those functors that preserve sifted colimits. If $\cat{C}$ is symmetric monoidal and $\cat{V}$ is closed symmetric monoidal, then $i_!$ becomes a symmetric monoidal functor with respect to Day convolution. In particular, $i_!$ preserves $\mathbb{E}_\infty$-algebras, i.e.\ if $F\colon \cat{C}\rt \cat{V}$ is a lax symmetric monoidal functor, then its sifted-colimit preserving extension $i_!(F)$ is lax symmetric monoidal as well.

Applying this to $\Mod_{\mathbb{Z}, \geq 0}=\mathcal{P}_{\Sigma}(\Vect_{\mathbb{Z}}^\omega)$, it remains to verify that $\Vect^\omega_{\mathbb{Z}}\rt \Sp^{\sG}; A\longmapsto \underline{A}$ is lax symmetric monoidal. This functor admits a factorisation
$$\begin{tikzcd}
\Vect_{\mathbb{Z}}^\omega\arrow[r] & \Sp^{\sG, \heartsuit}\arrow[r] & \Sp^{\sG}
\end{tikzcd}$$
over the heart of $\Sp^{\sG}$, i.e.\ the category of Mackey functors $A(G)\rt \cat{Ab}$ with values in (discrete) abelian groups. The inclusion $\Sp^{\sG, \heartsuit}\hookrightarrow \Sp^{\sG}$ is lax symmetric monoidal and one readily verifies that sending $A$ to the corresponding constant Mackey functor is lax symmetric monoidal.
\end{proof}
\begin{example}[Cohomological Mackey functors]
\label{ex:derived reps}
If $R$ is a connective $\mathbb{E}_1$-ring spectrum over $\ZZ$, then $\ul{R}$ defines an associative algebra in $\Sp^{\sG}$, and we let $\Mod_{\ulR}^{\sG}=\Mod_{\ulR}(\Sp^{\sG})$ denote the corresponding category of left modules. Let us point out that $\ul{R}$ \emph{differs} from the $\mathbb{E}_1$-algebra denoted $\mm{triv}_G(R)$ in \cite[Example 3.7]{patchkoriasanderswimmer}, whose modules are the ($R$-linear) derived Mackey functors $A(G)\rt \Mod_R$ of Kaledin. The $t$-structure on $\Sp^{\sG}$ induces a left and right complete $t$-structure on $\Mod_{\ulR}^{\sG}$, in which an object $M$ is (co)connective if and only if each $M(X)$ is a (co)connective spectrum for any finite $G$-set $X$. 

Write $R[\mathcal{O}_G]\subseteq \Mod_{\ulR, \geq 0}^{\sG}$ for the full (additive) subcategory spanned by the free $\ul{R}$-modules generated by finite $G$-sets $X$, i.e.\ of $\ul{R}$-modules of the form $\ul{R}\otimes \Sigma^{\infty}_+X$. The objects of $R[\mathcal{O}_G]$ are compact projective generators for $\Mod_{\ulR, \geq 0}^{\sG}$, which implies that there is an equivalence 
$$\begin{tikzcd}
\Mod_{R[\Orb_G]}\arrow[r, "\sim"] & \Mod_{\ulR}^{\sG}=\Mod_{\ulR}(\Sp^{\sG}).
\end{tikzcd}$$
We then use \Cref{rem:detecting coherent categories} (as in  Example \ref{ex:genuine G-spectra}) to show that $R[\Orb_G]$ is a coherent additive $\infty$-category if $R$ is a coherent $\mathbb{E}_1$-ring spectrum over $\ZZ$.

When $R$ is a \emph{discrete ring}, the objects in $R[\Orb_G]$ are all contained in the heart of the $t$-structure, i.e.\ they correspond to $\ul{R}$-modules in the (ordinary) category of Mackey functors $A(G)\rt \cat{Ab}$ with values in discrete abelian groups. Indeed, since all suspension spectra of finite $G$-sets are dualisable in $\Sp^{\sG}$, we have that
$$
\big(\ul{R}\otimes \Sigma^{\infty}_+X\big)(Y)\simeq \Hom_{\Sp^G}\big(\Sigma^\infty_+Y, \ul{R}\otimes \Sigma^{\infty}_+X\big)\simeq \Hom_{\Sp^G}\big(\Sigma^\infty_+(X\times Y), \ul{R}\big)\simeq \ul{R}(X\times Y)
$$
so that $\ul{R}\otimes \Sigma^\infty_+X$ corresponds to the Mackey functor $Y\mapsto \Map(X\times Y, R)^G$.

Following Yoshida \cite{yoshida1983}, the category $R[\mathcal{O}_G]$ can then be identified explicitly as follows: it is the full subcategory of the (ordinary) category $\Mod_{R[G]}^{\heartsuit}$ of discrete $R[G]$-modules spanned by the $R[G]$-modules obtained as $R$-linearisations of finite $G$-sets. We will denote such an $R$-linearisation of a $G$-set $X$ by $R[X]$.
This identification of $R[\Orb_G]$ is then induced by the functor $\ev_G\colon R[\mathcal{O}_G]\rt \Mod_{R[G]}^{\heartsuit}$ evaluating at the free $G$-set $G\in A(G)$. Indeed, this functor sends $\ul{R}\otimes \Sigma^\infty_+X$ to $R[X]$ and one readily verifies that it is fully faithful, using that
$$
\Map_{R[\Orb_G]}\big(\ul{R}\otimes \Sigma^\infty_+X, \ul{R}\otimes \Sigma^\infty_+Y\big)\simeq \Map_{\Sp^G}\big(\Sigma^\infty_+(X), \ul{R}\otimes \Sigma^\infty_+ Y\big)\simeq \ul{R}(X\times Y)=R[X\times Y]^G
$$
and likewise that $\Map_{R[G]}\big(R[X], R[Y]\big)\simeq \Map_{\cat{Set}^{G}}(X, R[Y])\simeq R[X\times Y]^G$.
\end{example}

After these recollections, we can now turn to the main topic of this subsection: 
\begin{definition}[Pro-coherent modules]\label{def:pro-coh}
Let $\cat{A}$ be a coherent additive $\infty$-category. We define the $\infty$-category of \emph{pro-coherent (left) $\cat{A}$-modules} as
$$
\QC^\vee_{\cat{A}}=\mm{Ind}\left(\Coh_{\cat{A}^{\op}}^{\op}\right).
$$
More explicitly, one can identify $\QC^\vee_{\cat{A}}$ with the $\infty$-category $\Fun_\mm{ex}(\Coh_{\cat{A}^{\op}}, \Sp)$ of exact functors $M\colon \Coh_{\cat{A}^{\op}}\rt \Sp$.
\end{definition}
Coherent modules are generally not preserved by (nonabelian) left derived functors such as tensor products. It will therefore be convenient to give a slightly different presentation of pro-coherent modules in terms of almost perfect modules.
\begin{definition}[Convergent functors]\label{def:convergence}
Let $\cat{C}$ be a stable $\infty$-category with a left complete $t$-structure.
If $\cat{V}$ is an $\infty$-category with sequential limits, a functor $F\colon \cat{C}\rt \cat{V}$ is said to be \emph{convergent} if for any object $X\in \cat{C}$, the natural map $$F(X)\xrightarrow{\simeq} \lim_m F(\tau_{\leq m} X)$$ is an equivalence. Write $\Fun_{\conv}(\cat{C}, \cat{V})\subseteq \Fun(\cat{C}, \cat{V})$ for the full subcategory spanned by the convergent functors. 
\end{definition}
\begin{remark}\label{rem:conv in terms of postnikov}
Note that a functor $F\colon \cat{C}\rt \cat{V}$ as above is convergent if and only if it preserves limits of all 
towers 
$\dots\rt X_1\rt X_0$ in $\cat{C}$ with the property that   for each $m\geq 0$, the tower $\dots\rt \tau_{\leq m}X_1\rt \tau_{\leq m}X_0$ is eventually constant. 

\end{remark} 
\begin{lemma}\label{lem:convergent functors}
Let $\cat{C}$ be a small stable $\infty$-category equipped with a left complete $t$-structure, and write $\cat{C}^+\subset \cat{C}$ for the full subcategory of eventually coconnective objects. Given another $\infty$-category 
$\cat{V}$  with small limits,   restriction determines an equivalence $\Fun_{\conv}(\cat{C}, \cat{V})\simeq \Fun(\cat{C}^+, \cat{V})$, with inverse given by right Kan extension. 
\end{lemma}
\begin{proof}
Since right Kan extension along the fully faithful inclusion $\cat{C}^+\hookrightarrow \cat{V}$ defines a fully faithful functor $\Fun(\cat{C}^+, \cat{V})\rt \Fun(\cat{C}, \cat{V})$, it suffices to verify that a functor is convergent if and only if it is right Kan extended from $\cat{C}^+$. This follows from Remark \ref{rem:conv in terms of postnikov} and the fact that for any $X\in \cat{C}$, its Postnikov tower defines a right cofinal functor $\mathbb{N}\rt \cat{C}^+_{X/}$.
\end{proof}
Since $\Coh_{\cat{A}^{\op}} \simeq \APerf_{\cat{A}^{\op}}^+$, we obtain a new characterisation of coherent modules:
\begin{corollary}\label{cor:pro-coh via aperf}
Let $\cat{A}$ be a coherent additive $\infty$-category. Then there is an equivalence $\QC^\vee_{\cat{A}}\simeq \Fun_{\mm{ex}, \conv}(\APerf_{\cat{A}^{\op}}, \Sp)$.
\end{corollary}
\begin{remark} The exact functors $\Coh_{\cat{A}^{\op}}\rt \Sp$ and $\APerf_{\cat{A}^{\op}}\rt   \Sp$ are determined by their restriction to connective objects, as all objects are eventually connective.
\end{remark}
\begin{definition}[Dually almost perfect modules]\label{daperf}
We say that a pro-coherent $\cat{A}$-module $M$ is   \emph{dually almost perfect} if the corresponding convergent exact  functor $M\colon \APerf_{\cat{A}^{\op}}\rt \Sp$ is corepresentable by an almost perfect $\cat{A}^{\op}$-module. Write $\APerf^\vee_{\cat{A}}\subseteq \QC^\vee_{\cat{A}}$ for the full subcategory spanned by these  objects, and observe that there is a (formal) equivalence of $\infty$-categories
$$\begin{tikzcd}
(-)^\vee\colon \APerf_{\cat{A}^{\op}}^{\op}\arrow[r, "\simeq"] & \APerf^\vee_{\cat{A}}.
\end{tikzcd}$$
\end{definition}
We will now describe the relation between the $\infty$-categories of $\cat{A}$-modules and pro-coherent $\cat{A}$-modules, their difference being controlled by $t$-structures. We start by endowing pro-coherent modules with a $t$-structure.
\begin{lemma}[Pro-coherent $t$-structure]\label{lem:pro-coh t-structure}
Let $\cat{A}$ be a coherent additive $\infty$-category. Then $\QC^\vee_{\cat{A}}$ carries a left complete, accessible $t$-structure such that a pro-coherent module $M$ is \emph{connective} if and only if $M\colon \Coh_{\cat{A}^{\op}}\rt \Sp$ is right $t$-exact.
\end{lemma}
\begin{proof}
The existence of the desired $t$-structure follows immediately from \cite[Proposition 1.4.4.11]{lurie2014higher}. It is left complete because the connective objects are closed under products and  the intersection $\bigcap_n \QC^\vee_{\cat{A},\geq n}$ contains only the zero object \cite[Proposition 1.2.1.19]{lurie2014higher}.
\end{proof}
\begin{remark}\label{rem:connective pro-coh}
Note that a pro-coherent module $M$ is connective if and only if the associated exact convergent functor $M\colon \APerf_{\cat{A}^{\op}}\rt \Sp$ is right $t$-exact. Indeed, for each $X\in \APerf_{\cat{A}^{\op},\geq 0}$ the spectrum $M(X)$ arises as the limit of a tower of connective spectra $M(\tau_{\leq n}X)$ with connective fibres.
\end{remark}
\begin{remark}[Relation to ind-coherent modules]
Let $R$ be a coherent commutative ring with dualising complex $\omega_R$. Then Serre duality gives an equivalence $\QC^\vee_R\simeq \mm{Ind}(\Coh_R)$. However, this equivalence does not identify the $t$-structure of Lemma \ref{lem:pro-coh t-structure} with the $t$-structure on ind-coherent sheaves from \cite[Proposition 1.2.2]{gaitsgory2017study}. Instead, the induced $t$-structure on $\mm{Ind}(\Coh_R)$ has connective objects generated by $\omega_R$ under colimits and extensions.
\end{remark}
Using the pairing $\Mod_{\cat{A}^{\op}}\times \Mod_{\cat{A}}\rt \Sp$ from Example \ref{ex:opposite}, every left $\cat{A}$-module $M$ determines an exact functor $(-)\otimes_{\cat{A}} M\colon \Coh_{\cat{A}^{\op}}\rt \Sp$. We obtain a functor $\upiota$ from $\cat{A}$-modules to pro-coherent $\cat{A}$-modules, which is part of an adjunction
$$\begin{tikzcd}
\cooliota  \colon \Mod_{\cat{A}}\arrow[r, yshift=1ex] & \QC^\vee_{\cat{A}}\colon \coolU.\arrow[l, yshift=-1ex]
\end{tikzcd}$$
Observe  that $\upiota: \Mod_{\cat{A}}\rt \QC^\vee_{\cat{A}}$ is the unique colimit-preserving extension of its restriction to $\cat{A}$. In terms of Corollary \ref{cor:pro-coh via aperf}, this restriction sends each $V$ to the convergent functor $\APerf_{\cat{A}^{\op}}\rt \Sp$ corepresented by $V$, which we view as an \mbox{object in $\cat{A}^{\op}$.}
\begin{proposition}
Let $\cat{A}$ be a coherent additive \mbox{$\infty$-category.} Then $\upiota$ exhibits $\Mod_{\cat{A}}$ as the right completion of $\QC^\vee_{\cat{A}}$.
\end{proposition} 
\begin{proof}
If $M$ is a connective $\cat{A}$-module, then $(-)\otimes_{\cat{A}} M\colon \Coh_{\cat{A}^{\op}}\rt \Sp$ is right \mbox{$t$-exact,} and so  
 $\upiota$ is a right $t$-exact functor. To verify that   $\upiota$ restricts  to an equivalence $\Mod_{\cat{A}, \geq 0}\xrightarrow{\ \simeq\ }\QC^\vee_{\cat{A},\geq 0}$, first note that we can identify
$$
\Mod_{\cat{A}, \geq 0}\subseteq \Fun(\cat{A}^{\op}, \Sp_{\geq 0}) \qquad\text{and}\qquad \QC^\vee_{\cat{A},\geq 0}\subseteq \Fun(\APerf_{\cat{A}^{\op},\geq 0}, \Sp_{\geq 0})
$$
with the full subcategories spanned by additive functors and  right exact  convergent functors, respectively, using  Remark \ref{rem:connective pro-coh}. In fact, note that every right exact functor $F\colon \APerf_{\cat{A}^{\op},\geq 0}\rt \Sp_{\geq 0}$ is automatically convergent, because the cofibre of each $F(X)\rt F(\tau_{\leq m}X)$ is the $(m+2)$-connective spectrum $F\big(\tau_{\geq m+1}X[1]\big)$. 

Unravelling the definitions, we can identify the functor  $\upiota\colon \Mod_{\cat{A}, \geq 0}\rt \QC^\vee_{\cat{A},\geq 0}$ with the functor taking left Kan extension along $\cat{A}^{\op}\rt \APerf_{\cat{A}^{\op}}$
$$\begin{tikzcd}
\Fun_{\oplus}(\cat{A}^{\op}, \Sp_{\geq 0})\arrow[r] & \Fun_{\mm{rex}}(\APerf_{\cat{A}^{\op},\geq 0}, \Sp_{\geq 0}).
\end{tikzcd}$$
In particular, $\upiota$ is fully faithful, so it only remains to check that restriction along $\cat{A}^{\op}\rt \APerf_{\cat{A}^{\op}}$ detects equivalences. This holds as any right exact functor $F\colon \APerf_{\cat{A}^{\op},\geq 0}\rt \Sp_{\geq 0}$ preserves geometric realisations as for any simplicial diagram $X_\bullet$ in $\APerf_{\cat{A}^{\op},\geq 0}$ and each $m\geq 0$, the natural map $|\mm{sk}_m F(X_\bullet)|\stackrel{\simeq}{\rt} F|\mm{sk}_m(X_\bullet)|\rt F|X_\bullet|$ has an $(m+1)$-connective cofibre.
\end{proof}
\begin{remark}[The bounded case]\label{rem:modules as subcat of pro-coh}
Let $\cat{A}$ be a coherent additive $\infty$-category such that there is an $n$ such that $\Hom_{\cat{A}}(V, W)$ is $n$-coconnective \mbox{for all $V, W\in \cat{A}$.} Then there are  inclusions of full subcategories $\cat{A}^{\op}\subseteq \Coh_{\cat{A}^{\op}}\subseteq \Mod_{\cat{A}^{\op}}$, and $\upiota$ can then be identified with the functor $\Fun_{\oplus}(\cat{A}^{\op}, \Sp)\rt \Fun_{\mm{ex}}(\Coh_{\cat{A}^{\op}}, \Sp)$ taking left Kan extension. Hence $\upiota$ is fully faithful and preserves compact objects.

\end{remark}

\begin{example}\label{ex:regular}
If $R$ is a coherent connective $\EE_1$-ring as in \Cref{mainex}, \mbox{set $\QC^\vee_R:=\QC^\vee_{\Vect_R^\omega}$.}
Then $\cooliota$ is fully faithful if and only if $R$ is eventually coconnective. One implication follows directly from Remark \ref{rem:modules as subcat of pro-coh}. For the converse, unravelling the definitions shows that for any connective $\EE_1$-ring $R$ and a left module $M\in \Mod_R$, the unit map $M\rt \coolU(\cooliota(M))$ can be identified with
$$\begin{tikzcd}
M\arrow[r] & \lim_{n\to \infty}\Big(\tau_{\leq n}R\otimes_R M\Big).
\end{tikzcd}$$
Applying this to $M=\bigoplus_{k\geq 0} R[-k]$ shows that $R$ is eventually coconnective if $\cooliota$ is fully faithful.

If $R$ is furthermore Noetherian and $\EE_\infty$, then  $\upiota$ is an equivalence if and only if $R$ is discrete regular  Noetherian, as in this case, any finitely generated $R$-module admits a finite free resolution and the inclusion $\Perf_{R^{op}} \xrightarrow{\simeq} \Coh_{R^{op}}$ is an equivalence.
\end{example}
Finally, let us mention the following condition that is somewhat dual to being connective in the $t$-structure from \Cref{lem:pro-coh t-structure}: 
\begin{definition}\label{def:tor-dimension}
We will say that a pro-coherent $\cat{A}$-module $M$ is of \emph{tor-amplitude $\leq d$} if $M\colon \Coh_{\cat{A}^{\op}}\rt \Sp$ is left $t$-exact up to a shift by $d$, i.e.\ it sends $\Coh_{\cat{A}^{\op}, \leq 0}$ to $\Coh_{\cat{A}^{\op}, \leq d}$. Let us write $\QC^\vee_{\cat{A}, \weirdleq d}$ for the full subcategory of pro-coherent modules of tor-amplitude $\leq d$.
\end{definition}
\begin{example}
Let $R$ be a coherent connective $\EE_1$-ring and $M$ a left $R$-module. Then the pro-coherent $R$-module $\upiota(M)$ is of tor-amplitude $\leq d$ if and only if $M$ is of tor-amplitude $\leq d$ in the usual sense, i.e.\ for each discrete right $R$-module $N\in \Mod_{R^{\op}}^{\heartsuit}$, the tensor product $N\otimes_R M$ is $d$-coconnective.
\end{example}
\begin{example}\label{ex:tor-dim=dual of conn}
Let $\cat{A}$ be a coherent additive $\infty$-category and $M\in \APerf(\cat{A}^{\op})$. Then the dually almost perfect module $M^\vee\in \QC^\vee_{\cat{A}}$ from Definition \Cref{daperf} has tor-amplitude $\leq d$ if and only if $M$ is $(-d)$-connective.
\end{example}

\subsection{Extended functors}\label{sec:extended functors}
We will now consider sifted-colimit-preserving 
functors $\QC^\vee_{\cat{A}}\rt \QC^\vee_{\cat{B}}$   between categories of pro-coherent modules. 
Our aim is to construct such functors as extensions of functors $\cat{A}\rt \cat{B}$,  thereby
generalising a method of the first author and Mathew \cite[Section 3.2]{brantner2019deformation}, which is  related to the   work of Illusie \cite[Section I-4]{Illusie} and
\mbox{Kaledin \cite[Section 3]{Kaledin}.}
  
\begin{notation}
If $\cat{C}, \cat{V}$ are two $\infty$-categories with sifted colimits, let $\Fun_{\Sigma}(\cat{C}, \cat{V})$ be the full subcategory of $\Fun(\cat{C}, \cat{V})$ spanned by the sifted-colimit-preserving  functors.
\end{notation}
We start by recalling that for any small  additive $\infty$-category $\cat{A}$, the objects in $\cat{A}$ form compact projective generators for $\Mod_{\cat{A}, \geq 0}$. Given another $\infty$-category $\cat{V}$ with sifted colimits, restriction along $\cat{A}\hookrightarrow \Mod_{\cat{A}, \geq 0}$ therefore defines an equivalence \cite[Proposition 5.5.8.15]{lurie2009higher}
\begin{equation}\label{diag:classical nonabelian derived}
\begin{tikzcd}
\Fun_{\Sigma}\big(\Mod_{\cat{A}, \geq 0}, \cat{V}\big)\arrow[r, "\simeq"] & \Fun\big(\cat{A}, \cat{V}\big).
\end{tikzcd}
\end{equation}
The inverse is given by left Kan extension, and sends a functor $F\colon \cat{A}\rt \cat{V}$ to its nonabelian left derived functor. \vspace{5pt}

When $\cat{A}$ is a small coherent additive $\infty$-category, there is a similar method for producing functors out of pro-coherent modules, where on the right hand side of \eqref{diag:classical nonabelian derived}, we need to enlarge $\cat{A}$ to also include some non-connective objects.\vspace{2pt}

We will use the following notion from \cite[Appendix C]{lurie2016spectral}:
\begin{definition}[op-prestable $\infty$-categories]
An $\infty$-category $\cat{C}$ is said to be \emph{op-prestable} if $\cat{C}^{\op}$ is a prestable \mbox{$\infty$-category} in the sense of
\cite[Definition C.1.2.1]{lurie2016spectral}. In other words, $\cat{C}$ is op-prestable if there is a fully faithful embedding $\iota\colon \cat{C}\hookrightarrow \cat{D}$ into a stable $\infty$-category, with essential image closed under finite limits and extensions.
 If $\iota$ is initial among such   embeddings, we call $\cat{D}$ the \emph{stable envelope} of $\cat{C}$; this is the case precisely if every object in $\cat{D}$ is an iterated suspension \mbox{of objects in $\cat{C}$.}
\begin{definition}[The $\infty$-category $\APerf^\vee_{\cat{A},\weirdleq 0 }$]\label{coconnect}
Let $\APerf^\vee_{\cat{A},\weirdleq 0 }\subseteq \APerf_{\cat{A}}^\vee$ be the full subcategory of dually almost perfect modules of tor-amplitude $\leq 0$. By \Cref{ex:tor-dim=dual of conn}, this is equivalent to the full subcategory spanned by the
 modules $M^\vee$ with $M \in \APerf_{\cat{A}^{\op}, \geq 0}$ connective. Note that $\APerf^\vee_{\cat{A},\weirdleq 0 } $ is an op-prestable $\infty$-category with stable envelope $\APerf_{\cat{A}}^\vee$.

Likewise, let  $\Perf_{\cat{A}, \weirdleq 0 }\subseteq \Perf_{\cat{A}}$ be the full subcategory of perfect $\cat{A}$-modules of tor-amplitude $\leq 0$, or equivalently, of perfect $\cat{A}$-modules with connective dual $\cat{A}^{\op}$-module. Then $\Perf_{\cat{A}}$ is the \mbox{stable envelope of $\Perf_{\cat{A}, \weirdleq 0 }$.}
\end{definition}
\begin{remark}
Note that $\APerf^\vee_{\cat{A},\weirdleq 0 }$ is generally different from $\APerf_{\cat{A}, \leq 0}$,  the full subcategory of  almost perfect  modules which are coconnective in the \mbox{$t$-structure} considered in \Cref{lem:pro-coh t-structure}. 
For example, take   \mbox{$\cat{A}=\Vect_{k[\epsilon]}^\omega$} as in \Cref{mainex}. 
The augmentation $k[\epsilon] \rightarrow k$ induces a functor $\QC^\vee_{k[\epsilon]} \rightarrow \QC^\vee_{k} \simeq \Mod_k $ which 
preserves (dually) almost perfect objects. In $\QC^\vee_{k} \simeq \Mod_k$, these are just complexes bounded below (above) with finite-dimensional homotopy groups.
The discrete $k[\epsilon]$-module $k$ is connective  almost perfect, and so $k^{\vee}$ belongs to $\APerf_{\cat{A}, \weirdleq 0}^{\vee}$.
However, $k^\vee$ is not almost perfect, as $k$ is not  dually almost perfect since  $k \otimes_{k[\epsilon]} k$ does not have \mbox{bounded above homotopy.}

Note also that $\APerf^\vee_{\cat{A},\weirdleq 0 }$ can be  different from $\APerf^\vee_{\cat{A}, \leq 0}$,  the full subcategory of \textit{dually} almost perfect  modules which are coconnective. For example, take $k[\epsilon_1]$ the trivial square-zero extension of $k$ by a class in  degree $1$.
Then $k[\epsilon_1]^\vee = k[\epsilon_1]$ belongs to  $\APerf^\vee_{\cat{A},\weirdleq 0 }$, but  is not coconnective as \mbox{there is a nonzero map $\Sigma k[\epsilon_1] \rightarrow k[\epsilon_1]$.}

\end{remark}

\end{definition}
Recall that a simplicial object in an $\infty$-category is called $m$-skeletal if it is the left Kan extension of its restriction to $\Delta^{\op}_{\leq m}$.
\begin{notation}[Finite stable geometric realisations]
If $\cat{C}$ is an op-prestable $\infty$-category and $\cat{V}$ admits geometric realisations, then a functor $F\colon \cat{C}\rt \cat{V}$ is said to \emph{preserve finite stable geometric realisations} if the following condition holds:
if $X_\bullet$ is a simplicial object  in $\cat{C}$ such that the image in the stable envelope of $\cat{C}$ is $m$-skeletal for some $m$ and has its geometric realisation contained in $\cat{C}$, then the  natural map $|F(X_\bullet)|\rt F(|X_\bullet|)$ is an equivalence. We  write $\Fun_{\sigma}(\cat{C}, \cat{V})\subseteq \Fun(\cat{C}, \cat{V})$ for the full subcategory spanned by the functors preserving finite stable geometric realisations.

\end{notation}
\begin{remark}If $\cat{C}$ is already stable, we will also refer to finite stable geometric realisations  as finite geometric realisations.
\end{remark}
\begin{definition}[Regular functors]
Let $\cat{A}$ be a coherent additive $\infty$-category. 
If $\cat{V}$ is an $\infty$-category with sequential colimits, then a functor $F\colon \APerf_{\cat{A}}^\vee\rt \cat{V}$ is said to be \emph{regular} if the composite $\APerf_{\cat{A}^{\op}} \xrightarrow{\simeq} 
(\APerf_{\cat{A}}^\vee)^{\op} \rt \cat{V}^{\op}$, $V\mapsto F(V^\vee)$ is convergent in the sense of \Cref{def:convergence}. Write $\Fun_{\reg}\big(\APerf_{\cat{A}}^\vee, \cat{V}\big)\subseteq \Fun\big(\APerf_{\cat{A}}^\vee, \cat{V}\big)$ for the full subcategory spanned by the regular functors. 
\end{definition}
We begin by restricting from pro-coherent to dually almost perfect modules:
\begin{proposition}\label{prop:left extension from unbounded almost perfect}
Let $\cat{A}$ be a coherent additive $\infty$-category and $\cat{V}$ a presentable $\infty$-category. Restriction determines an equivalence of $\infty$-categories
$$\begin{tikzcd}
\Fun_\Sigma\big(\QC^\vee_{\cat{A}}, \cat{V}\big)\arrow[r, "\simeq"] & \Fun_{\sigma, \reg}\big(\APerf_{\cat{A}}^\vee, \cat{V}\big),
\end{tikzcd}$$
the inverse of which is given by left Kan extension.
\end{proposition}
\begin{proof}
Recall that each $\QC^\vee_{\cat{A}}$ is compactly generated by $\Coh_{\cat{A}^{\op}}^{\op}$. The proof of \cite[Proposition 3.9]{brantner2019deformation} then shows that restriction and left Kan extension determine an adjoint equivalence $\Fun_{\sigma}\big(\Coh_{\cat{A}^{\op}}^{\op}, \cat{V}\big)\leftrightarrows \Fun_{\Sigma}\big(\QC^\vee_{\cat{A}}, \cat{V}\big)$, where the domain is the full subcategory of functors preserving finite geometric realisations. Hence, it suffices to verify that restriction and left Kan extension determine an adjoint equivalence $\Fun_{\sigma}\big(\Coh_{\cat{A}^{\op}}^{\op}, \cat{V}\big)\leftrightarrows\Fun_{\sigma, \reg}\big(\APerf_{\cat{A}}^\vee, \cat{V}\big)$. This follows from (the opposite of) Lemma \ref{lem:convergent functors} and the fact that given an $m$-skeletal simplicial diagram $X_\bullet^\vee$ in $\APerf_{\cat{A}}^\vee$, there is a sequence of $m$-skeletal diagrams $(\tau_{\leq n} X_{\bullet})^\vee$ in $\Coh_{\cat{A}^{\op}}^{\op}$ \mbox{with colimit $X_\bullet^\vee$.}
\end{proof}
In a second step, we restrict even further from $\APerf_{\cat{A}}^\vee$ to the $\infty$-category $\APerf^\vee_{\cat{A},\weirdleq 0 }$ of dually almost perfect modules of tor-amplitude $\leq 0$ (cf.\ \Cref{coconnect}):
 
\begin{proposition}\label{prop:left extension from almost perfect}
Let $\cat{A}$ be an additive $\infty$-category and $\cat{V}$ an $\infty$-category with sifted colimits. Then restriction determines an equivalence of $\infty$-categories
$$\begin{tikzcd}
\Fun_\Sigma\big(\Mod_{\cat{A}}, \cat{V}\big)\arrow[r, "\simeq"] & \Fun_\sigma\big(\Perf_{\cat{A},{\weirdleq 0}}, \cat{V}\big)
\end{tikzcd}$$
with inverse given by left Kan extension. If $\cat{A}$ is furthermore coherent, this extends to a commuting square
$$\begin{tikzcd}
\Fun_{\Sigma}\big(\QC^\vee_{\cat{A}}, \cat{V}\big)\arrow[r, "\simeq"]\arrow[d, "\upiota^*"{swap}] & \Fun_{\sigma, \reg}\big( \APerf^\vee_{\cat{A},\weirdleq 0 }, \cat{V}\big)\arrow[d]\\
\Fun_\Sigma\big(\Mod_{\cat{A}}, \cat{V}\big)\arrow[r, "\simeq"{swap}] &  \Fun_\sigma\big(\Perf_{\cat{A},{\weirdleq 0}}, \cat{V}\big)
\end{tikzcd}$$
where the horizontal equivalences are given by restriction, with inverses given by left Kan extension. 
\end{proposition}

The assertion about sifted-colimit preserving functors out of $\Mod_{\cat{A}}$ follows from exactly the same argument as \cite[Proposition 3.14]{brantner2019deformation}. For the second assertion in Proposition \ref{prop:left extension from almost perfect}, we will need two auxiliary observations.
\begin{lemma}\label{lem:cechresolution}
Let $\cat{A}$ be an additive $\infty$-category and $M\in\APerf_{\cat{A},\geq n}$. There is a right cofinal  functor $\Delta\rt (\APerf_{\cat{A},\geq n+1})_{M/}$ such that the underlying cosimplicial diagram   is $1$-coskeletal   and the above diagram exhibits $M$ as its limit.
\end{lemma}
\begin{proof}
Since $M$ is almost perfect and $n$-connective, there exists a cofibre sequence $V[n]\rt M\rt M^0$ of $\cat{A}$-modules with $V\in \cat{A}$ and $M^0\in\APerf_{\cat{A},\geq n+1}$, cf.\ \Cref{rem:connective almost perfect as realisations}. Let $M^\bullet$ be the `\v{C}ech conerve' of the map $M\rt M^0$. This determines a functor $\phi\colon \Delta\rt (\APerf_{\cat{A},\geq n+1})_{M/}$ with the desired two properties. It remains to verify that $\phi$ is right cofinal. To this end, let $N$ be an $(n+1)$-connective almost perfect module equipped with a map $f\colon M\rt N$. We have to show that the over-category $\Delta_{/f}$ is contractible. Note that the projection $\Delta_{/f}\rt \Delta$ is the right fibration classifying the simplicial space
$$\begin{tikzcd}
\Delta^{\op}\arrow[r] & \An; & {[n]}\arrow[r, mapsto] & \Map(M^n, N)\times_{\Map(M, N)} \{f\}.
\end{tikzcd}$$
We have to check that the geometric realisation of this simplicial space is \mbox{contractible,} for which it suffices to show that the natural map $|\Map(M^\bullet, N)|\rt \Map(M, N)$ is an equivalence. Since $M^\bullet$ is the \v{C}ech conerve of $M\rt M^0$, the above diagram is the \v{C}ech nerve of the map $\Map(M^0, N)\rt \Map(M, N)$. It therefore suffices to verify that this map of spaces induces a surjection on $\pi_0$. In other words, for any map $g\colon M\rt N$, we need to provide a null-homotopy of the composition $V[n]\rt M\rt N$.
This follows immediately from the assumption that $N$ was $(n+1)$-connective and that $\Hom(V, -)\colon \Mod_{\cat{A}}\rt \Sp$ is $t$-exact for all $V\in \cat{A}$.
\end{proof}
\begin{lemma}\label{lem:extendfromcoconn}
Let $\cat{A}$ be a coherent additive  $\infty$-category and $\cat{V}$ an $\infty$-category with sifted colimits. For any functor $F\colon \APerf_{\cat{A}}^\vee\rt \cat{V}$, the following are equivalent:
\begin{enumerate}
\item $F$ preserves finite geometric realisations.
\item $F$ is left Kan extended from its restriction to $ \APerf^\vee_{\cat{A},\weirdleq 0 }$, which preserves finite stable geometric realisations.
\end{enumerate}
\end{lemma}
\begin{proof}
Set $\cat{X}=\APerf_{\cat{A}}^\vee$. Given $m\geq 0$, write $\cat{X}_m =\APerf_{\geq -m, \cat{A}}^{\vee}\subseteq \cat{X}$ for the full subcategory spanned by those dually almost perfect modules with $(-m)$-connective duals. Note that each $\cat{X}_m$ is op-prestable, with stable envelope $\cat{X}$, and there is a colimit sequence of $\infty$-categories
$\begin{tikzcd}[cramped]
\cat{X}_0\ar[r] & \cat{X}_1\ar[r] & \mc{X}_2\arrow[r] & \dots\ar[r] & \cat{X}.
\end{tikzcd}$

Arguing as in \cite[Proposition 3.11]{brantner2019deformation}, it suffices to verify inductively that for all $m$, the functor  $F\big|_{\cat{X}_m}$ preserves finite stable geometric realisations if and only if it is right Kan extended from $\cat{X}_{ m-1}$ and the restriction $F\big|_{\cat{X}_{m-1}}$ preserves finite stable geometric realisations. First, if $F\big|_{\cat{X}_{m}}$ preserves finite stable geometric realisations, we have to prove that for every $M\in \cat{X}_{ m}$, the map
$\colim\limits_{M_\alpha\in (\cat{X}_{m-1})_{/M}} F(M_\alpha)\rightarrow F(M)$
is an equivalence in $\cat{V}$. Using the opposite of Lemma \ref{lem:cechresolution}, we can replace the colimit in the domain by a finite stable geometric realisation; the result then follows from the assumption that $F$ preserves such geometric realisations.

For the converse, let $M_\bullet$ be a finite simplicial diagram in $\cat{X}_{m}$ with $M_{-1}=|M_\bullet|$ contained in $\cat{X}_{m}$ as well. There exists a fibre sequence in $\cat{X}$ of the form $M_{-1, 0}\rt M_{-1}\rt V[m]$, with $V\in \cat{A}$ and $M_{-1, 0}\in \cat{X}_{m-1}$. For all $p\geq 0$, the composite $M_p\rt M_{-1}\rt V[m]$ has fibre $M_{p, 0}$ in $\cat{X}_{m-1}$ as well. Let $M_{\bullet, \bullet}$ be the bisimplicial diagram arising as the \v{C}ech nerve of the natural transformation $M_{\bullet, 0}\rt M_{\bullet}$, so that $M_p=\colim_q M_{p, q}$. We then have a commuting square
$$\begin{tikzcd}
F(M_{-1})=F(\colim_{p, q} M_{p, q})\arrow[r]\arrow[d] & \colim_{p} F(\colim_q M_{p, q})=\colim_p F(M_p)\arrow[d]\\
\colim_q F(\colim_q M_{p, q})\arrow[r] & \colim_{p, q} F(M_{p, q}).
\end{tikzcd}$$
Assuming that $F\big|\cat{X}_m$ is left Kan extended from $\cat{X}_{m-1}$, the two vertical maps are equivalences by the opposite of (the proof of) Lemma \ref{lem:cechresolution}. For each $q\geq 0$, $|M_{\bullet, q}|$ is the finite stable geometric realisation of a simplicial diagram in $\cat{X}_{m-1}$. Since $F\big|_{\cat{X}_{m-1}}$ preserves such geometric realisations by assumption, the bottom horizontal map is an equivalence. This implies that the top horizontal map is an equivalence, i.e.\ $F\big|_{\cat{X}_{m}}$ preserves finite stable geometric realisations.
\end{proof}
\begin{proof}[Proof (of Proposition \ref{prop:left extension from almost perfect})]
The first equivalence follows in exactly the same way as \cite[Proposition 3.14]{brantner2019deformation}. Alternatively, it follows by substituting perfect $\cat{A}$-modules for almost perfect $\cat{A}$-modules in \Cref{lem:cechresolution} and \Cref{lem:extendfromcoconn}.  
For the top equivalence when $\cat{A}$ is coherent, it suffices by Proposition \ref{prop:left extension from unbounded almost perfect} to verify that restriction and left Kan extension define an adjoint equivalence 
$$\begin{tikzcd}
\Fun_{\sigma, \reg}\big(\APerf_{\cat{A}}^\vee, \cat{V}\big)\arrow[r, "\simeq"] & \Fun_{\sigma, \reg}\big(\APerf^\vee_{\cat{A},\weirdleq 0 }, \cat{V}\big).
\end{tikzcd}$$ It suffices to verify that $F\colon \APerf_{\cat{A}}^\vee\rt \cat{V}$ is regular and preserves finite geometric realisations if and only if it is left Kan extended from $\APerf^\vee_{\cat{A},\weirdleq 0 }$ and its restriction to $\APerf^\vee_{\cat{A},\weirdleq 0 }$ is regular and preserves finite stable geometric realisations. This follows from Lemma \ref{lem:extendfromcoconn} by unravelling the regularity conditions in $(1)$ and $(2)$.
\end{proof}
 
We will now use Proposition \ref{prop:left extension from almost perfect} to construct functors $\QC^\vee_{\cat{A}}\rt \cat{V}$, respectively $\Mod_{\cat{A}}\rt \cat{V}$ from functors $\cat{A}\rt \cat{V}$, thereby generalising \cite[Section 3.2]{brantner2019deformation} to coherent rings:
\begin{definition}[Right extendable functors]\label{def:right extendable}
Let $\cat{A}$ be an additive \mbox{$\infty$-category} and $\cat{V}$ an $\infty$-category with small limits and colimits. 
\mbox{A functor $F\colon \cat{A}\rt \cat{V}$ is:} 
\begin{enumerate}
\item  
\emph{right extendable} if its right Kan extension $F^R\colon  \Perf_{\cat{A},\weirdleq 0}\rt \cat{V}$ along the inclusion $ \cat{A}\hookrightarrow \Perf_{\cat{A},\weirdleq 0}$ preserves finite stable geometric realisations. 
\item if $\cat{A}$ is coherent: \emph{coherently right extendable} if its right Kan extension $F^R\colon  \APerf^\vee_{\cat{A},\weirdleq 0 }\rt \cat{V}$ along the inclusion $ \cat{A}\hookrightarrow \APerf^\vee_{\cat{A},\weirdleq 0 }$ is regular and preserves finite stable geometric realisations. 
\end{enumerate}
\end{definition}
\begin{cons}[Right-left extension]
Given a right-extendable functor $F$ as in \Cref{def:right extendable} (1), 
the \emph{right-left} \emph{derived functor} of $F$ is given by 
the sifted-colimit-preserving functor 
$F^{RL}\colon \QC_{\cat{A}}\rt \cat{V}$ \mbox{provided by \Cref{prop:left extension from almost perfect}.}

If $\cat{A}$ is coherent and $F$ is coherently right-extendable,  
the \emph{coherent right-left} \emph{derived functor} of $F$ is 
the sifted-colimit-preserving functor $
F^{RL}\colon \QC^\vee_{\cat{A}} \rightarrow \cat{V}$ provided by \Cref{prop:left extension from almost perfect}.
\end{cons} 

Note that a coherently right extendable functor is in particular right extendable.

\begin{remark}\label{rem:left right extended functors}
In the setting of Definition \ref{def:right extendable}, 
restriction along $\cat{A}\rt \Mod_{\cat{A}}$ defines an equivalence between the full subcategory of $\Fun(\cat{A}, \cat{V})$ spanned by the right extendable functors and the full subcategory of $\Fun(\Mod_{\cat{A}}, \cat{V})$ on those functors that preserve sifted colimits and also finite totalisations of  diagrams in $\cat{A}$, by Proposition \ref{prop:left extension from almost perfect} and Remark \ref{rem:connective almost perfect as realisations}.

Likewise, restriction along $\cat{A}\rt \QC^\vee_{\cat{A}}$ defines an equivalence between the full subcategory of $\Fun(\cat{A}, \cat{V})$ on the coherently right extendable functors and the full subcategory of $\Fun(\QC^\vee_{\cat{A}}, \cat{V})$ on those functors that preserve all sifted colimits and also totalisations of cosimplicial diagrams in $\cat{A}$.
\end{remark}

Generalising \cite[Theorem 3.27]{brantner2019deformation}, our main source of examples comes from functors of finite degree \cite{eilenberg1954groups}, or a  mild generalisation thereof:

\begin{proposition}\label{prop:finite degree implies extendable}
Let $\cat{A}$ be a (coherent) additive $\infty$-category and $\cat{V}$ a stable $\infty$-category with small limits and colimits, equipped with a right complete $t$-structure such that $\cat{V}_{\leq 0}$ is closed under countable direct sums. Then:
\begin{enumerate}
\item Let $F\colon \cat{A}\rt \cat{V}$ be a functor of finite degree with values in $\cat{V}_{\leq 0}$. Then $F$ is (coherently) right extendable.
\item More generally, let $F_1\rt F_2\rt \dots$ be a countable sequence of functors $F_i\colon \cat{A}\rt \cat{V}$ as in (1). Then  $F:=\colim_i F_i$ is (coherently) right extendable and the natural map $\colim_i \big(F_i^{RL}\big)\rt F^{RL}$ is an equivalence.
\end{enumerate}
\end{proposition}

\begin{proof}
We will only deal with the coherent case, the non-coherent case being similar but easier.

For (1), consider the \emph{opposite} functor $F^{\op}\colon \cat{A}^{\op}\rightarrow \cat{V}^{\op}$, which is also a functor of finite degree $r$. As in \cite[Proposition 3.35]{brantner2019deformation}, the left Kan extension $F'\colon \Mod_{\cat{A}^{\op},\geq 0}\rightarrow \cat{V}^{\op}$ preserves sifted colimits and is $r$-excisive.  Theorem 3.36 and Proposition 3.37  in  \cite{brantner2019deformation} together imply that 
 $F'$ preserves finite stable totalisations. It also preserves limits of Postnikov towers: indeed, since $F^{\op}\colon \cat{A}^{\op}\rightarrow \cat{V}^{\op}$ takes values in connective objects, each map $F'(M)\rightarrow F'(\tau_{\leq n}M)$ has $(n+1)$-connective fibres, so that the tower for $M$ converges by left completeness of the $t$-structure on $\cat{V}^{\op}$. We conclude that the restriction of $F'$ to $\APerf_{\cat{A}^{\op}}$ is convergent and preserves finite stable totalisations. Passing to opposite categories, we see that $F^R\colon \APerf_{\cat{A}}^\vee\rightarrow \cat{V}$ is regular and preserves finite stable geometric realisations.

For (2), the previous argument gives a functor $F':=\lim_i F'_i\colon \APerf_{\cat{A}^{\op},\geq 0}\rightarrow \cat{V}^{\op}$ which preserves finite stable totalisations and limits of Postnikov towers. We claim that $F'$ is the left Kan extension of $F^{\op}\colon \cat{A}^{\op}\rightarrow \cat{V}^{\op}$; dually, this means that $F^R\simeq \colim_i F_i^R$,  which implies  assertion (2). For the claim, note that $F'$ agrees with $F^{\op}$ on $\cat{A}^{\op}$, so that it suffices to show that $F'$ preserves geometric realisations. For any simplicial diagram $M_\bullet$ in $\APerf_{\cat{A}}$, we have
$$\begin{tikzcd}
{|F'(M_\bullet)|=\big|\lim_i F'_i(M_\bullet)\big|}\arrow[r, "\simeq"] & \lim_i \big|F'_i(M_\bullet)\big|\arrow[r, "\simeq"] & \lim_i F'_i\big(|M_\bullet|\big)=F'\big(|M_\bullet|\big).
\end{tikzcd}$$
The first equivalence uses that geometric realisations commute with limits of towers of connective objects in $\cat{V}^{\op}$; this in turn follows from the fact that geometric realisations commute with countable products, since $\cat{V}^{\op}$ is left complete and connective objects are closed under products. The second equivalence follows because each $F'_i$ preserves sifted colimits by construction. 
Passing to opposite categories, we deduce that the functors $F_i^R, F^R\colon \APerf_{\cat{A}}^\vee\rightarrow \cat{V}$ have the desired properties.
\end{proof} 

 Proposition \ref{prop:finite degree implies extendable} lets us  extend certain functors between  additive $\infty$-categories:
\begin{definition}\label{def:polynomial functor}
Let $\cat{A}, \cat{B}$ be additive $\infty$-categories. A functor $F\colon \cat{A}\rightarrow \cat{B}$ is called \emph{locally polynomial} if it arises as the colimit of a sequence $F_1\rightarrow F_2\rightarrow \dots$ of functors from $\cat{A}$ to $\cat{B}$, such that:
\begin{enumerate}
\item For each $X\in \cat{A}$, the sequence $F_1(X)\rightarrow F_2(X)\rightarrow \dots$ is \mbox{eventually constant.}

\item Each $F_i\colon \cat{A}\rightarrow \cat{B}\rightarrow \Perf(\cat{B})$ is a functor of finite degree, i.e.\ there exists an $r\geq 0$ such that the cross-effect $\mm{cr}_{r+1}\colon \cat{A}^{\times r+1}\rightarrow \Perf(\cat{B})$ vanishes.
\end{enumerate}

The composition of two locally polynomial functors is again locally polynomial.  
\begin{notation}\label{addpoly}
Write $\mathscr{A}\mm{dd}^\mm{coh,poly}$ for the (non-full) subcategory of $\mathscr{C}\mm{at}_\infty$ spanned by   coherent additive $\infty$-categories and locally polynomial functors between them.\end{notation}
\end{definition}
\begin{corollary}\label{cor:polynomial functors between coherent}
Let $F\colon \cat{A}\rightarrow \cat{B}$ be a locally polynomial functor between coherent additive $\infty$-categories. Then the following diagram admits a unique extension as indicated
\begin{equation}\label{diag:pro-coh extension}\begin{tikzcd}
\cat{A}\arrow[d, "F"{swap}]\arrow[r, hook] & \Mod_{\cat{A}}\arrow[r, "\upiota"]\arrow[d, "F'", dotted] & \QC^\vee_{\cat{A}}\arrow[d, "F''", dotted]\\
\cat{B}\arrow[r, hook] & \Mod_{\cat{B}}\arrow[r, "\upiota"{swap}] &  \QC^\vee_{\cat{B}}
\end{tikzcd}\end{equation}
such that $F'$ preserves sifted colimits and finite totalisations
and $F''$ 
 preserves sifted colimits and all totalisations of cosimplicial diagrams in $\cat{A}$. 
\end{corollary}

\begin{proof}
Uniqueness  follows immediately from Remark \ref{rem:left right extended functors}. For existence, we may assume without restriction that  $F$ has finite degree $r$. Indeed,  if $F$ is a sequential colimit of finite degree functors $F_i$, we can simply take the sequential colimit of the extensions $F'_i$ and $F''_i$, which  has the desired properties by Proposition \ref{prop:finite degree implies extendable}.

For the existence of $F''$, we  will apply Proposition \ref{prop:finite degree implies extendable}. Indeed, note that   $\QC^\vee_{\cat{B}}=\mm{Ind}\big(\Coh_{\cat{B}^{\op}}^{\op}\big)$ admits a \emph{second, right complete} $t$-structure, such that $F$ takes values in coconnective objects: the connective part of this $t$-structure is the ind-completion of $(\Coh_{\cat{B}^{\op},\leq 0})^{\op}$, and its coconnective part is the ind-completion of $(\Coh_{\cat{B}^{\op},\geq 0})^{\op}$. We will not use this second $t$-structure elsewhere.

For the existence of $F'$, it suffices to show that $F''$ maps $\Perf_{\cat{A},\weirdleq 0}$ into $\Perf_{\cat{B},\weirdleq 0}$; the desired extension $F'$ is then the left Kan extension of the following composite:
$$\begin{tikzcd}
\Perf_{\cat{A},\weirdleq 0}\arrow[r, "F''"] &  \Perf_{\cat{B}, \weirdleq 0}\arrow[r] & \Mod_{\cat{B}}.
\end{tikzcd}$$
To see that $F''$ preserves the duals of perfect connective objects, we observe that the functor  $F''\colon \APerf^\vee_{\cat{A},\weirdleq 0 }\rightarrow \APerf^\vee_{\cat{B},\weirdleq 0 }$ is  opposite to an $r$-excisive functor $\APerf_{\cat{A}^{\op},\geq 0}\rightarrow \APerf_{\cat{B}^{\op},\geq 0}$ sending $\cat{A}^{\op}$ to $\cat{B}^{\op}$.

Such functors preserve perfect objects. Indeed, write $\cat{E} \subset \APerf_{\cat{A}^{\op},\geq 0}$ for the full subcategory of all $M$ for which  $F''(M)$ perfect. Given a cofibre sequence $X \rightarrow V \rightarrow C$ with $X\in \cat{E}$ and $V\in \cat{A}^{\op}$, we form the strongly coCartesian cube with `initial legs'   $X \rightarrow V$. Its colimit $V \oplus_X   \ldots \oplus_X V$ can be identified with $V \oplus C \oplus \ldots \oplus C$. Hence $F(C)$ perfect, as it is a  retract of the perfect module $F(V \oplus C \oplus \ldots \oplus C)$.
 As $\cat{E}$ also contains $0$ and is closed under retracts, we conclude  $\Perf_{\cat{A}^{\op},\geq 0}  \subset \cat{E}$.
\end{proof}
\begin{example}[Divided orbits]\label{dG1}
Given a coherent $\EE_1$-ring $R$ and a finite group $G$, we will write $R[G]=R \otimes \Sigma^\infty_+ G$ for the associated group ring, which is again coherent. 
 The evident map $R[G] \rightarrow R$ induces a limit-preserving functor $\Mod_{R} \rightarrow \Mod_{R[G]}$. Restricting its left adjoint  induces an additive functor  $$(-)_{hG}: \Perf_{R[G]} \rightarrow \Perf_{R},$$ which on underlying spectra takes homotopy orbits.
We right-left extend using \Cref{cor:polynomial functors between coherent} to obtain a  functor
$$(-)_{dG}: \QC^\vee_{R[G]}\longrightarrow \QC^\vee_{R}. $$
The functor $(-)_{dG}$ behaves like a mix between homotopy orbits and fixed points. 

Indeed, if $V\in \APerf_{R[G], \geq 0}$ is almost perfect and connective, we can write $V = |V_\bullet|$ as a realisation of a simplicial diagram  in $\Vect_{R[G]}$ and compute $$V_{dG} \simeq |(V_{\bullet})_{hG}| \simeq V_{hG}.$$

If, on the other hand, we have $V \in \APerf_{R[G], \weirdleq 0}^\vee$, then we can find  a cosimplicial diagram $V^\bullet$ in $\Vect_{R[G]}$ with $V \simeq \Tot(V^\bullet) $.
As the norm is an equivalence on objects in $\Vect_{R[G]}$, we 
compute
 $V_{d\sG} \simeq    \Tot(V^\bullet_{h\sG}) \simeq  \Tot((V^\bullet)^{h\sG}) \simeq  \Tot(V^\bullet)^{h\sG} \simeq V^{h\sG}. $
In fact, since the functors $(-)_{dG}$, $(-)_{hG}$, and $(-)^{hG}$ are exact, we obtain identifications
$V_{dG} \simeq V_{hG}$ for all $V\in \APerf_{R[G]}$  and 
$V_{dG} \simeq V^{hG}$ for all $V\in \APerf^\vee_{R[G]}$.
\end{example}
\begin{example}[Derived orbits and genuine fixed points]\label{ex:derived fixed points}
For a coherent $\mathbb{E}_1$-ring spectrum $R$ over $\mathbb{Z}$ and a finite group $G$, let $\Mod^{\sG}_{\ulR}=\Mod_{\ulR}(\Sp^{\sG})$ be the $\infty$-category from Example \ref{ex:derived reps}. Recall that taking genuine $G$-fixed points defines a functor $(-)^{\sG}\colon \Mod^{\sG}_{\ulR}\rt \Mod_R$ that preserves both limits and colimits; in terms of spectral Mackey functors, this simply evaluates at the trivial $G$-orbit $G/G\in A(G)$. Its left adjoint is the functor
$$\begin{tikzcd}
\mm{triv}_{\sG}\colon \Mod_R\arrow[r] & \Mod_{\ulR}^{\sG}; & M\arrow[r] & \ul{R}\otimes_R M
\end{tikzcd}$$
where $\ul{R}\otimes_R M$ denotes the spectral Mackey functor given by $(\ul{R}\otimes_R M)(X)=\ul{R}(X)\otimes_R M$. One can think of this as an $\ul{R}$-linearised version of endowing $M$ with the trivial $G$-action.

The constant spectral Mackey functor $\ul{R}$ has the rather special feature that $\mm{triv}_G$ also preserves limits: indeed, since the genuine $H$-fixed points jointly detect limits, this follows from the fact that $\mm{triv}_G(M)^H=\ul{R}(G/H)\otimes_R M\simeq M$. We will denote by $(-)_G\colon \Mod_{\ulR}^{\sG}\rt \Mod_R$ the corresponding left adjoint to $\mm{triv}_G$.

These three functors restrict to (adjoint) functors between finitely generated free objects
$$(-)_{\sG}\colon R[\Orb_G]\rt \Vect_R^\omega,\quad\qquad \mm{triv}_G\colon \Vect_R^\omega\rt R[\Orb_G],\quad\qquad (-)^{\sG}\colon R[\Orb_G]\rt \Vect_R^\omega.$$
Indeed, one readily verifies that $(\ul{R}\otimes\Sigma^{\infty}_+(G/H))_G$ and $(\ul{R}\otimes\Sigma^{\infty}_+(G/H))^G$ are both equivalent to $R$. When $R$ is a discrete ring, Example \ref{ex:derived reps} identifies $R[\Orb_G]$ with the ordinary category of $R[G]$-modules $R[X]$ induced by finite $G$-sets, and the above three functors coincide with taking (strict) $G$-coinvariants, trivial $G$-modules and $G$-invariants, respectively.

Using \Cref{cor:polynomial functors between coherent}, we then obtain colimit-preserving functors
$$(-)_{G}: \QC^\vee_{R[\Orb_G]}\longrightarrow \QC^\vee_{R} \qquad \qquad (-)^{G}: \QC^\vee_{R[\Orb_G]}\longrightarrow \QC^\vee_{R} $$
taking \emph{derived orbits} and \emph{derived genuine fixed points}. Note that the derived genuine fixed points functor $(-)^G$ behaves as expected on dually almost perfect objects: for any cosimplicial diagram $V^\bullet$ in $R[\Orb_G]$ one has that $\Tot(V^\bullet)^G\simeq \Tot((V^\bullet)^G)$.
\end{example}

\subsection{Monoidal structures}
Corollary \ref{cor:polynomial functors between coherent}  provides the  main source of functors between categories of pro-coherent modules for us. To express the functoriality of these derived functors, let $\StPrSigma\subset \mathscr{C}\mm{at}_\infty$ be the (non-full) subcategory of (large) $\infty$-categories on the stable presentable $\infty$-categories with sifted-colimit-preserving functors between them. Both $\StPrSigma$ and $\mathscr{A}\mm{dd}^\mm{coh,poly}$  are closed \mbox{under finite products.}

\begin{theorem}\label{thm:derived functors}
There is a natural transformation of symmetric monoidal functors
$$
\begin{tikzcd}[column sep=3pc]
\mathscr{A}\mm{dd}^\mm{coh,poly} \arrow[r, bend left=30, "\Mod", ""{below, name=U}]\arrow[r, bend right=30, "\QC^\vee"{swap}, ""{name=D}] & \StPrSigma  \arrow[Rightarrow, from=U, to=D, end anchor={[xshift=1pt]}, "\upiota"]
\end{tikzcd}$$
sending each coherent additive $\infty$-category $\cat{A}$ to $\QC^\vee_{\cat{A}}$ and each locally polynomial functor to its (right-left) derived functor.
\end{theorem}
\begin{proof}
The entire diagram can be described as a single product-preserving functor $F\colon \mathscr{A}\mm{dd}^\mm{coh,poly}\rightarrow \mm{Ar}\big(\StPrSigma\big)$ to the arrow category. To construct $F$, consider the subcategory $\cat{X}\subseteq \Fun(\Delta[2], \mathscr{C}\mm{at}_\infty)$ consisting of sequences of the form $\cat{A}\rightarrow \Mod_{\cat{A}}\rightarrow \QC^\vee_{\cat{A}}$, with maps between them given by natural diagrams as in \eqref{diag:pro-coh extension}, where $F$ is polynomial and $F'$ (respectively $F''$) preserves sifted colimits and totalisations of coskeletal (respectively all) cosimplicial diagrams in $\cat{A}$. The functor $F$ then arises from the  zig-zag $\mathscr{A}\mm{dd}^\mm{coh,poly} \stackrel{\simeq}{\lt} \cat{X}\rightarrow \mm{Ar}\big(\StPrSigma\big)$, where the left functor is an equivalence by Corollary \ref{cor:polynomial functors between coherent}.

To see that $F$ preserves finite products, it suffices to verify that the natural maps $\cat{A}\times \cat{B}\rightarrow \Mod_{\cat{A}}\times \Mod_{\cat{B}}$ and $\cat{A}\times \cat{B}\rightarrow \QC^\vee_{\cat{A}}\times \QC^\vee_{\cat{B}}$ extend to equivalences
$$
\Mod_{\cat{A}\times \cat{B}}\simeq \Mod_{\cat{A}}\times \Mod_{\cat{B}} \qquad\quad \QC^\vee_{\cat{A}\times \cat{B}}\simeq \QC^\vee_{\cat{A}}\times \QC^\vee_{\cat{B}}.
$$
The first equivalence follows from $\Mod_{\cat{A}\times \cat{B}, \geq 0}\simeq \Mod_{\cat{A}, \geq 0}\times \Mod_{\cat{B},\geq 0}$, which holds because both $\infty$-categories have $\cat{A}\times \cat{B}$ as compact projective generators. This implies  that the natural map $\Coh_{\cat{A}^{\op}\times\cat{B}^{\op}}\rightarrow \cat{Coh}_{\cat{A}^{\op}}\times\Coh_{\cat{B}^{\op}}$ is an equivalence as well, and the second equivalence follows by ind-completing.
\end{proof}

\begin{remark}\label{rem:extending linear}
If $F\colon \bigoplus_i \cat{A}_i\rightarrow \cat{B}$ is a polynomial functor which is \emph{additive} in the $k$-th variable, then its extension $F^{RL}\colon \prod_i\QC^\vee(\cat{A}_i)\rightarrow \QC^\vee_{\cat{B}}$ preserves sifted colimits in each variable and small colimits in the $k$-th variable.
\end{remark}
\begin{example}[Monoidal structure on pro-coherent modules]\label{ex:tensor product}
Let $R$ be a coherent $\mathbb{E}_{n+1}$-algebra. Then the $\mathbb{E}_n$-monoidal structure $\otimes_R$ on $\Mod_R$ restricts to a tensor product on the additive $\infty$-category $\Vect_R$. Since this is linear in each variable, this determines an $\mathbb{E}_n$-algebra in $\mathscr{A}\mm{dd}^\mm{coh,poly}$. By \Cref{thm:derived functors}, $\QC^\vee_{R}$ inherits an $\mathbb{E}_n$-monoidal structure, which  preserves colimits in each variable by   \Cref{rem:extending linear}, and the functor $\upiota\colon \Mod_R\rightarrow \QC^\vee_{R}$ is $\mathbb{E}_n$-monoidal. More explicitly,    \Cref{prop:left extension from almost perfect} can be used to realise $\QC^\vee_{R}$ as an $\mathbb{E}_n$-monoidal localisation of $\Fun(\APerf_{R^{\op},\weirdleq 0}, \Sp)$, equipped with the Day convolution product.
\end{example} 
\noindent In the presence of a symmetric monoidal structure on $\QC^\vee_{\cat{A}}$ satisfying \mbox{mild conditions,}  dually almost perfect modules and almost perfect modules are  related by duality:

\begin{proposition}\label{prop:dual is indeed dual}
Let $\cat{A}$ be a coherent additive $\infty$-category equipped with a nonunital symmetric monoidal structure $\otimes$ which preserves finite sums in each variable, and moreover satisfies the following conditions:
\begin{enumerate}
\item The nonunital closed monoidal structure on $\QC^\vee_{\cat{A}}$, constructed as in \Cref{ex:tensor product}, admits a unit $\mb{1}$, which is eventually connective. 
\item Every object in $\cat{A}$ is dualisable, with dual contained in $\cat{A}$.
\end{enumerate}
Then taking duals determines an equivalence
$$\begin{tikzcd}
(-)^\vee:=\Hom(-, \mb{1})\colon \APerf_{\cat{A}}\arrow[r, "\simeq"]& (\APerf_{\cat{A}}^\vee)^{\op}.
\end{tikzcd}$$
\end{proposition}
\begin{proof}
It suffices to show that $(-)^\vee\colon \QC^\vee_{\cat{A}}\rightarrow (\QC^\vee_{\cat{A}})^{\op}$ preserves totalisations of cosimplicial objects in $\cat{A}$. This will imply the result because 
 $(-)^\vee$  restricts to an equivalence on $\cat{A}$ by (1)  and preserves  small colimits.

So let $M^\bullet$ be a cosimplicial diagram in $\cat{A}$.
As $\QC^\vee_{\cat{A}}\subset\Fun((\APerf^\vee_{\cat{A},\weirdleq 0 })^{\op}, \Sp)$ is a reflective subcategory, it is enough to prove that for any object $N\in \APerf^\vee_{\cat{A},\weirdleq 0 }$, the natural map 
$\begin{tikzcd}[column sep=1pc]
\big|\Hom\big(N, (M^\bullet)^\vee\big)\big|\arrow[r] & \Hom\big(N, \mm{Tot}(M^\bullet)^\vee\big)
\end{tikzcd}$
is an equivalence; this implies that  $\mm{Tot}(M^\bullet)^\vee$ is the geometric realisation of   $(M^\bullet)^\vee$.
We can identify the above map with the  composite map
$$\begin{tikzcd}
\big|\Hom\big(N\otimes M^\bullet, \textbf{1}\big)\big|\arrow[r] & \Hom\big(\mm{Tot}(N\otimes M^\bullet), \textbf{1}\big)\arrow[r] & \Hom\big(N\otimes \mm{Tot}(M^\bullet), \textbf{1}\big).
\end{tikzcd}$$
The second map is an equivalence since $\otimes$ preserves totalisations of cosimplicial diagrams in $\cat{A}$  by \Cref{cor:polynomial functors between coherent}, and the first map is an equivalence because $\Hom(-, \textbf{1})\colon \APerf^\vee_{\cat{A},\weirdleq 0 }\rightarrow \Sp$ is right $t$-exact up to a shift, by our assumption that $\textbf{1}$ is eventually connective (cf.\ \Cref{rem:connective pro-coh}).
\end{proof}
\begin{remark}\label{rem:ext poly noncoherent case}
Of course, \Cref{cor:polynomial functors between coherent} and \Cref{thm:derived functors} have analogues for additive $\infty$-categories that are not coherent (with the same proofs): there is a symmetric monoidal functor $\Mod\colon \mathscr{A}\mm{dd}^\mm{poly}\rt \StPrSigma$ sending each additive $\infty$-category $\cat{A}$ to $\Mod_{\cat{A}}$ and each locally polynomial functor $F$ to its right-left derived functor. 
\end{remark}

\subsection{$\mathscr{O}$-monoidal structures}
Given a (coloured) $\infty$-operad $\cat{O}^{\otimes}  \rightarrow \N(\Fin_\ast)$, which we   informally also call $\cat{O}$,  \Cref{thm:derived functors} shows that $\cat{O}$-monoidal structures on $\cat{A}$ and $\cat{B}$ induce canonical $\cat{O}$-monoidal structures on $\QC^\vee_{\cat{A}}$ and $\QC^\vee_{\cat{B}}$, respectively. We have seen that the right-left extension of a strong monoidal polynomial functor $\cat{A}\rightarrow \cat{B}$ is again  strong monoidal, and  will  now establish a refinement to (op)lax
 $\cat{O}$-monoidal functors, which is needed  for our  treatment of \mbox{PD operads:}
\begin{proposition}\label{prop:derived functor of lax monoidal}
Let $\cat{A}, \cat{B}$ be $\mc{O}$-algebras in $\mathscr{A}\mm{dd}^\mm{coh,poly}$ and let $F\colon \cat{A}\rightarrow \cat{B}$ be an (op)lax $\mc{O}$-monoidal functor with  $F_x\colon \cat{A}_x\rightarrow \cat{B}_x$ of finite degree  for each colour $x\in \mc{O}$. \vspace{3pt}
Then $F^{RL}\colon \QC^\vee_{\cat{A}}\rightarrow \QC^\vee_{\cat{B}}$ admits a \mbox{natural (op)lax $\mc{O}$-monoidal structure.}
\end{proposition}
The proof relies on two observations concerning Kan extensions of lax $\mc{O}$-monoidal functors along $\mc{O}$-monoidal functors: namely, there is a canonical lax $\mc{O}$-monoidal structure on the right Kan extension, and in good cases also on the \mbox{left Kan extension.}
\begin{lemma}\label{lem:right kan extension lax monoidal}
Let  $\phi\colon \cat{C}_0\rightarrow \cat{C}$ be an \mbox{$\mc{O}$-monoidal} functor. If $F\colon \cat{C}_0\rightarrow \cat{D}$ is a lax $\mc{O}$-monoidal functor, then the following two assertions are equivalent:
\begin{enumerate}
\item For every colour $x\in\mc{O}$, there exists a functor $G_x\colon \cat{C}_x\rightarrow \cat{D}_x$ and a natural transformation $F_x\rightarrow G_x\circ \phi_x$ exhibiting $G_x$ as the right Kan extension of $F_x\colon \cat{C}_{0, x}\rightarrow \cat{D}_x$ along $\phi_x\colon \cat{C}_{0, x}\rightarrow\cat{C}_x$.
\item There exists a lax $\mc{O}$-monoidal functor $G\colon \cat{C}^{\otimes}\rightarrow \cat{D}^{\otimes}$ and a natural transformation $F\rightarrow G\circ \phi$ over $\mc{O}^{\otimes}$ exhibiting $G$ as the right Kan extension (relative to $\mc{O}^{\otimes}$) of $F\colon \cat{C}_0^{\otimes}\rightarrow \cat{D}^{\otimes}$ along $\phi\colon \cat{C}^{\otimes}_0\rightarrow \cat{C}^{\otimes}$.
\end{enumerate}
In this case, the fibre of the natural transformation $G\rightarrow F\circ \phi$ over a colour $x\in\mc{O}$ exhibits a right Kan extension of $F_x$ along $\phi_x$.
\end{lemma}
We  make use of the Day convolution product, cf.\ \cite{glasman2016day} and   {\cite[Section 2.2.6]{lurie2014higher}:} recall that for any small $\mc{O}$-monoidal $\infty$-category $\cat{C}$ and any presentably $\mc{O}$-monoidal $\infty$-category $\cat{D} \in \Alg_{\cat{O}}(\mm{\cat{P}r^L})$, there is another presentably $\mc{O}$-monoidal $\infty$-category $\Fun(\cat{C}, \cat{D})$ such that $\cat{O}$-algebras in $\Fun(\cat{C}, \cat{D})$ are lax $\cat{O}$-monoidal functors $\cat{C}\rightarrow \cat{D}$, with  fibre over  $x\in \mc{O}$  given by $\Fun(\cat{C}, \cat{D})_x=\Fun(\cat{C}_x, \cat{D}_x)$.

 We will apply this in particular when the target is spaces,  and  the $\mc{O}$-monoidal structure arises from the \mbox{cartesian product.}

\begin{proof}
Unraveling the definitions, we have to verify that the map
$$\begin{tikzcd}
\Alg_{\cat{C}^\otimes}(\cat{D})\times_{\Alg_{\cat{C}_0^{\otimes}}(\cat{D})}\Alg_{\cat{C}_0^{\otimes}}(\cat{D})_{F/}\arrow[r] & \prod\limits_{x\in \mc{O}}\Fun(\cat{C}_x, \cat{D}_x)\times_{\Fun(\cat{C}_{0, x},\cat{D}_x)}\Fun(\cat{C}_{0, x}, \cat{D}_x)_{F_x/}
\end{tikzcd}$$
preserves and detects terminal objects. As the Yoneda embedding $\cat{D}\rightarrow \cat{P}(\cat{D})$ is $\mc{O}$-monoidal for the Day convolution product on $\cat{P}(\cat{D})$, the above map is the pullback of the same map with $\cat{D}$ replaced by $\cat{P}(\cat{D})$, along a fully faithful functor. Since the Yoneda embedding preserves limits, it then suffices to verify that the corresponding map for $\cat{P}(\cat{D})$ preserves and detects terminal objects. Consequently, we may assume that $\cat{D}\in \Alg_{\cat{O}}(\mm{\cat{P}r^L})$ is an $\mc{O}$-monoidal presentable $\infty$-category.
In this situation, consider the $\mc{O}$-monoidal categories $\Fun(\cat{C}_0, \cat{D})$ and $\Fun(\cat{C}, \cat{D})$ given by Day convolution. Since $F$ is an $\mc{O}$-algebra in $\Fun(\cat{C}_0, \cat{D})$, we can form the pullback of $\mc{O}$-monoidal $\infty$-categories
$$
\Fun(\cat{C}, \cat{D})_{F/}^{\otimes}=\Fun(\cat{C}, \cat{D})^{\otimes}\times_{\Fun(\cat{C}_0, \cat{D})^{\otimes}} \Fun(\cat{C}_0, \cat{D})^{\otimes}_{F/}
$$
where $\Fun(\cat{C}_0, \cat{D})^{\otimes}_{F/}$ is the $\mc{O}$-monoidal $\infty$-category from \cite[Theorem 2.2.2.4]{lurie2014higher}. The above map can then be identified with the map
$$
\Alg_{\mc{O}}\big(\Fun(\cat{C}, \cat{D})_{F/}\big)\rightarrow \prod_{x\in \mc{O}}\big(\Fun(\cat{C}, \cat{D})_{F/}\big)_x.
$$
This map preserves and detects terminal objects by \cite[Corollary 3.2.2.3]{lurie2014higher}.
\end{proof}
\begin{lemma}\label{lem:left kan extension lax monoidal}
Let $\mc{O}$ be an operad and let $\phi\colon \cat{C}_0\rightarrow \cat{C}$ be an $\mc{O}$-monoidal functor. Let $F\colon \cat{C}_0\rightarrow \cat{D}$ be a lax $\mc{O}$-monoidal functor between $\cat{O}$-monoidal $\infty$-categories with the following property: for every colour $x\in \mc{O}$ and every $c\in \cat{C}_{x}$, the diagram
\begin{equation}\label{diag:kan extension colimit}\begin{tikzcd}
(\cat{C}_{0, x})_{/c}=\cat{C}_{0, x}\times_{\cat{C}_x} (\cat{C}_x)_{/c}\arrow[r] & \cat{C}_{0, x}\arrow[r, "F"] & \cat{D}_x
\end{tikzcd}\end{equation}
admits a colimit, which is preserved by each  $\psi(d_1, \dots, d_n, -)\colon \cat{D}_x\rightarrow \cat{D}_y$ for $\psi\in \mc{O}(y_1, \dots, y_n, x; y)$ and $d_i\in \cat{D}_{y_i}$. In this case, the left Kan extension of $F$ along $\phi$ exists and is   lax $\mc{O}$-monoidal.
\end{lemma}
\begin{proof}
This is essentially a consequence of \cite[Proposition 3.1.3.3]{lurie2014higher}; we include an argument for the reader's convenience. We can endow $\cat{P}(\cat{D})$ with the Day convolution $\mc{O}$-monoidal structure and let $\cat{V}$ be the left Bousfield localisation of $\cat{P}(\cat{D})$ at the natural maps from the colimits of \eqref{diag:kan extension colimit}, computed in $\cat{P}(\cat{D})$, to the representable presheaf on their colimit in $\cat{D}$. Our assumptions imply that $\cat{V}$ is an $\mc{O}$-monoidal localisation of $\cat{P}(\cat{D})$ and that the Yoneda embedding $h\colon \cat{D}\hookrightarrow \cat{V}$ is a fully faithful $\mc{O}$-monoidal functor preserving the colimits \eqref{diag:kan extension colimit}. As $\cat{V}\in \Alg_{\cat{O}}(\mm{\cat{P}r^L})$, we can equip $\Fun(\cat{C}, \cat{V})$ with the Day convolution $\mc{O}$-monoidal structure, so that left Kan extension defines an $\mc{O}$-monoidal functor
$
\mm{Lan}_\phi\colon \Fun(\cat{C}_0, \cat{V})\rightarrow \Fun(\cat{C}, \cat{V}).
$

As $\mc{O}$-algebras for the Day convolution product can be identified with lax  {$\mc{O}$-monoidal} functors, it follows that the left Kan extension of $h\circ F\colon \cat{C}_0\rightarrow \cat{V}$ along $\phi$ carries a canonical lax $\mc{O}$-monoidal structure. Since the Yoneda embedding $h\colon \cat{D}\hookrightarrow \cat{V}$ is $\mc{O}$-monoidal and preserves the colimits \eqref{diag:kan extension colimit}, we have that $\mm{Lan}_\phi(h\circ F)\simeq h\circ \mm{Lan}_\phi(F)$, so that $\mm{Lan}_\phi(F)$ inherits a lax $\mc{O}$-monoidal structure.
\end{proof}
\begin{proof}[Proof (of Proposition \ref{prop:derived functor of lax monoidal})]
We  first treat  the case where $F\colon \cat{A}\rightarrow \cat{B}$ is lax $\mc{O}$-monoidal. The construction of the lax monoidal structure on $F^{RL}\colon \QC^\vee_{\cat{A}}\rightarrow \QC^\vee_{\cat{B}}$ then proceeds in two steps: first taking a right Kan extension and then a left Kan extension, we obtain a diagram
$$\begin{tikzcd}
\cat{A}\arrow[d, "F"{swap}]\arrow[r] & \APerf^\vee_{\cat{A},\weirdleq 0 }\arrow[d, dotted, "{F^R}"]\arrow[r] & \QC^\vee_{\cat{A}}\arrow[d, dotted, "{F^{RL}}"]\\
\cat{B}\arrow[r] & \APerf^\vee_{\cat{B},\weirdleq 0 }\arrow[r] & \QC^\vee_{\cat{B}}.
\end{tikzcd}$$
The horizontal functors are all (strong) $\mc{O}$-monoidal by \Cref{thm:derived functors}.  Lemma \ref{lem:right kan extension lax monoidal}   implies that the right Kan extension $F^R$ is lax $\mc{O}$-monoidal. Next, we note that for every $M\in \QC^\vee_{\cat{A}}$, the over-category $(\APerf^\vee_{\cat{A},\weirdleq 0 })_{/M}$ is sifted (it admits finite sums) and that the $\mc{O}$-monoidal structure on $\QC^\vee_{\cat{B}}$ preserves sifted colimits in each variable. Lemma \ref{lem:left kan extension lax monoidal} shows that the left Kan extension $F^{RL}$ of $F^R$ \mbox{is lax $\mc{O}$-monoidal.}

For the oplax $\mc{O}$-monoidal case, one instead uses (the opposite of) Lemma \ref{lem:left kan extension lax monoidal} to show that $F^R$ is oplax $\mc{O}$-monoidal, using that for any $M\in \APerf_{\cat{A},\weirdleq 0}^\vee$, the under-category $\cat{A}_{M/}$ admits a right cofinal  
functor from $\Delta$ and that the tensor product on $\APerf_{\cat{B}}^\vee$ preserves totalisations. Next, (the opposite of) Lemma \ref{lem:right kan extension lax monoidal} shows that the left Kan extension $F^{RL}$ of $F^R$ is oplax $\mc{O}$-monoidal.
\end{proof}
\begin{remark}
The exact same proof shows that if $F\colon \cat{A}\rt \cat{B}$ is a lax $\mc{O}$-monoidal functor between additive $\infty$-categories, the left-right derived functor $F^{RL}\colon \Mod_{\cat{A}}\rightarrow \Mod_{\cat{B}}$ is (op)lax $\mc{O}$-monoidal. If $\cat{A}$ and $\cat{B}$ are coherent, then there is an equivalence $F^{RL}\circ \upiota\simeq \upiota\circ F^{RL}$ of lax $\mc{O}$-monoidal functors $\Mod_{\cat{A}}\rightarrow \QC^\vee_{\cat{B}}$.
\end{remark}

\  
\section{PD operads and refined Koszul duality}
\label{sec:PD operads and koszul}
We proceed to the  main abstract contribution of this paper: a twofold  refinement of classical operadic Koszul duality. First, we   show that the Koszul dual of an augmented $\infty$-operad $\cat{O}$ is not just an $\infty$-operad, but a divided power (`PD')  {$\infty$-operad} $\KD^{\pd}(\cat{O})$, which controls   Koszul duals of $\cat{O}$-algebras. \mbox{In a second, orthogonal,}  step, we replace $\infty$-operads $\cat{O}$ by \textit{derived $\infty$-operads}: here, the group actions are `more genuine', which means that 
derived operads   can parametrise structures  like simplicial commutative rings. We then   set up a refined Koszul duality in this  setting. 

\subsection{A reminder on $\infty$-operads}\label{sec:operads} 
To set the stage, let us recall a
definition of  $\infty$-operads as algebras in the category of symmetric sequences with the composition product. We follow the discussion in \cite[Section 4.1.2]{brantnerthesis}, which generalises a $1$-categorical construction of Kelly \cite{kelly2005ontheoperads} and Trimble \cite{trimblenotes} to the higher categorical setting. An alternative approach has been proposed by \mbox{Haugseng in \cite{haugseng2017infty}.}
\begin{notation}
Recall that the $\infty$-category $\PrL$ of presentable $\infty$-categories and colimit-preserving functors 
admits   the structure of a closed symmetric monoidal $\infty$-category by \cite[Proposition 4.8.1.15]{lurie2014higher}. The $\infty$-category of \emph{presentably symmetric monoidal $\infty$-categories} is given by $\CAlg(\PrL)$, and can be identified with the $\infty$-category  \mbox{of commutative algebras in $\PrL$.} 

Explicitly, a presentably symmetric monoidal $\infty$-category $\cat{V}$ is a symmetric monoidal $\infty$-category with  presentable  underlying $\infty$-category and a product $\otimes$  which distributes \mbox{over colimits.}

Let $R$ be an $\mathbb{E}_\infty$-ring   and consider the presentably symmetric monoidal \mbox{$\infty$-category} $\big(\Mod_R, \otimes_R\big)$ of $R$-modules. The $\infty$-category $\CAlg_R(\PrL)$ of presentably symmetric monoidal $R$-linear $\infty$-categories is given by the under-category $\CAlg(\PrL)_{(\Mod_R, \otimes_R)/}$.
\end{notation}
\begin{definition}[Symmetric sequences]\label{def:sym seq}
Let $R$ be an $\mathbb{E}_\infty$-ring. The \mbox{$\infty$-category} $\sSeq_R$ of $R$-linear \emph{symmetric sequences} is the free symmetric monoidal $R$-linear $\infty$-category generated by an object $\textbf{1}$. The universal symmetric monoidal structure on $\sSeq_R$ will be denoted by $\otimes$.
\end{definition}

The $\infty$-category $\sSeq_R$ can be described more explicitly as follows (cf.\ \cite[Section 4.1.2]{brantnerthesis}). Write $B\Sigma$ for the (nerve of the) category $\Fin^{\cong}$ of finite sets and bijections. The disjoint union of finite sets makes $(B\Sigma, \sqcup)$ the free symmetric monoidal $\infty$-category generated by the object $\ul{1}$ (cf.\  \cite[Proposition 2.2.4.9]{lurie2014higher}). We can then identify $(\sSeq_R, \otimes)$ with $\Fun(B\Sigma, \Mod_R)$, equipped with the Day convolution product of $\otimes_R$ and $\sqcup$ (cf.\  \cite[Corollary 4.8.1.12]{lurie2014higher}).

\begin{notation}[Symmetric sequences in arity $r$]
	For each $r$, let $B\Sigma_r$ denote the groupoid of finite sets of cardinality $r$ and bijections between them; up to equivalence, it has one object with automorphism group $\Sigma_r$. There is an adjoint pair $\iota_r\colon \Fun(B\Sigma_r, \Mod_R)\leftrightarrows \sSeq_R\colon \ev_r$ given by restriction and left Kan extension, respectively. The left adjoint $\iota_r$ is fully faithful and a symmetric sequence is said to be concentrated in \emph{arity} $r$ if it is contained in its essential image. The above adjunction then induces an equivalence between symmetric sequences concentrated in arity $r$ and modules over the group ring $R[\Sigma_r]$. Under this identification, the symmetric sequence $\mb{1}^{\otimes r}$ corresponds to the free $R[\Sigma_r]$-module of rank $1$. If $X$ is a symmetric sequence, we will denote its arity $r$ piece by $X(r)$.\vspace{-3pt}
\end{notation}

The universal property of $\sSeq_R$ asserts that for any  $\cat{V}\in \CAlg_R(\PrL)$, evaluation at $\textbf{1}$ defines an equivalence
$$
\mm{ev}_{\textbf{1}}\colon \Fun^{\mm{L}, \otimes}_R(\sSeq_R, \cat{V})\xrightarrow{\   \simeq \ } \cat{V},
$$
where the domain is the $\infty$-category of symmetric monoidal $R$-linear colimit-preserving functors $\sSeq_R\rightarrow \cat{V}$. Setting $\cat{V}=\sSeq_R$  gives an equivalence $$\End^{\mm{L}, \otimes}_R(\sSeq_R)\xrightarrow{\   \simeq \ }  \sSeq_R,$$ which categorifies the well-known identity $\Map_{\Rings}(\ZZ[t],\ZZ[t]) \cong \ZZ[t]$.

\begin{definition}[Composition product]\label{def:comp prod}
The \emph{composition product} $\circ$ on $\sSeq_R$
 is the monoidal structure   corresponding to the \emph{opposite} of the
evident monoidal structure on $\End^{\mm{L}, \otimes}_R(\sSeq_R)$ 
under the above equivalence. The unit of $\circ$  is the object $\textbf{1}$.
\end{definition}
\begin{remark}\label{def:prop of comp prod}
The definition of the composition product implies that the inverse of $$\End^{\mm{L}, \otimes}_R(\sSeq_R)\xrightarrow{\   \simeq \ } \sSeq_R$$ sends $Y$ to $(-)\circ Y$. In particular,   the composition product preserves colimits in the first variable. Similarly, for any symmetric sequence $X$, there are functors 
$$\begin{tikzcd}
\End^{\mm{L}, \otimes}_R(\sSeq_R)\arrow[r] & \End^{\mm{L}}_R(\sSeq_R)\arrow[r, "\mm{ev}_X"] & \sSeq_R
\end{tikzcd}$$
preserving sifted colimits and finite sifted limits.
Here the first functor forgets the monoidal structure and the second evaluates at $X$. This implies that  $\circ$ preserves sifted colimits and finite sifted limits in its second variable.
\end{remark}

\begin{remark}[Explicit formula for composition product]\label{rem:composition product formula}
Let $X$ and $Y$ be symmetric sequences. Unraveling the definitions, one sees that for each $r$, there is an $R$-linear left adjoint functor
$$\begin{tikzcd}
\Mod_{R[\Sigma_r]}\simeq \Fun(B\Sigma_r, \Mod_R)\arrow[r, hook, "\iota_r"] & \sSeq_R\arrow[r, "(-)\circ Y"] & \sSeq_R
\end{tikzcd}$$ 
sending the generating object $R[\Sigma_r]$ to the $r$-fold Day convolution product $\mb{1}^{\otimes r}\circ Y\simeq Y^{\otimes r}$. This implies that for each $r$, there is a natural equivalence $X(r)\circ Y\simeq X(r)\otimes_{\Sigma_r} Y^{\otimes r}$. Since every symmetric sequence $X$ can naturally be decomposed as $X\simeq \bigoplus_r X(r)$, we then obtain \vspace{-3pt}
\begin{equation}\label{eq:composition product formula}
X\circ Y\simeq \bigoplus_r X(r)\otimes_{\Sigma_r} Y^{\otimes r},
\end{equation}
which reproduces the classical 1-categorical formula for the composition product of symmetric sequences (see for instance \cite[Section 2.2.2]{fresse2009modules}).
\end{remark}
For any symmetric sequence $X$, the functor $X\circ (-)$ preserves symmetric sequences concentrated in arity $0$. Consequently, there is a (left) action
$$\begin{tikzcd}
\circ\colon \sSeq_R\times \Mod_R\arrow[r] & \Mod_R\vspace{-5pt}
\end{tikzcd}$$
of $\big(\sSeq_R, \circ\big)$ on $\Mod_R$, preserving sifted colimits and finite totalisations.
\begin{definition}[$\infty$-operads and cooperads]  \label{opcoop}
An  \emph{$\infty$-operad} $\PP$ over an $\EE_\infty$-ring $R$ is an associative algebra object in $\sSeq_R$ with respect to the composition product $\circ$. A \emph{$\PP$-algebra} is a left $\PP$-module in $\Mod_R$, equipped with the $\sSeq_R$-tensored structure described above.  We will write $\Op_R$ for the $\infty$-category of $R$-linear $\infty$-operads and $\Alg_{\PP}$ for the $\infty$-category of $\PP$-algebras. 

Dually, an \emph{$\infty$-cooperad} $\CC$ is a coassociative coalgebra in $\sSeq_R$ with respect to the composition product, and a (conilpotent) 
$\CC$-coalgebra is a left $\CC$-comodule in $\Mod_R$. We will write $\cat{Coop}_R$ for the $\infty$-category of $\infty$-cooperads.
\end{definition}
\begin{remark}
The $\infty$-operads in \Cref{opcoop} are often referred to as $\infty$-operads with one colour. Note that $\Op_R$ is   compactly generated   by Theorem \ref{thm:free algebra}.
\end{remark}
\subsection{The levelwise tensor product}
The category of symmetric sequences can be equipped  with yet another symmetric monoidal structure $\otimes_{\lev}$ -- the \emph{levelwise tensor product}. Its unit is the constant symmetric sequence on $R$, i.e.\ the $R$-linearisation of the $\mathbb{E}_\infty$-operad. This tensor product is compatible with the composition product in the following sense:
\begin{proposition}\label{prop:levelwise tensor lax monoidal}
The functor $\otimes_{\lev}\colon (\sSeq_R\times \sSeq_R,\  \circ \times \circ ) \rightarrow (\sSeq_R,\circ)$ has both a natural lax and oplax monoidal structure with respect to the composition product. In particular, for all $A,B,C,D \in \sSeq_R$, there are natural morphisms 
$$(A\circ B)\otimes_{\lev} (C\circ D) \rightarrow (A\otimes_{\lev} C)\circ (B\otimes_{\lev} D)$$
$$ (A\otimes_{\lev} C)\circ (B\otimes_{\lev} D) \rightarrow (A\circ B)\otimes_{\lev} (C\circ D).$$
\end{proposition}
In particular, this implies that the levelwise tensor product of two $\infty$-operads is again an $\infty$-operad, and a similar statement holds for $\infty$-cooperads. The proof requires a preliminary observation:
\begin{lemma}\label{lem:conjugating by adjunction}
Let $F\colon \cat{C}\leftrightarrows \cat{D}\colon G$ be an adjunction between $\infty$-categories. Then the induced  functor between endomorphism categories $\End(\cat{D})\rightarrow \End(\cat{C}); \, T\mapsto GTF$ inherits a lax monoidal structure.
\end{lemma}
\begin{proof}
Let $\pi\colon \cat{M}\rightarrow \Delta^1$ denote the correspondence classifying the adjoint pair $(F, G)$ \cite[Section 5.2.2]{lurie2009higher} and let $\End_{/\Delta^1}(\cat{M})$ denote the category of endofunctors of $\cat{M}$ compatible with the projection to $\Delta^1$. Restricting such endofunctors to the fibre $\cat{C}$ over $0$, respectively $\cat{D}$ over $1$, defines monoidal functors $\End(\cat{C})\lt \End_{/\Delta^1}(\cat{M})\rightarrow \End(\cat{D})$ with respect to composition. The right functor admits a right adjoint, given by relative right Kan extension over $\Delta^1$. Since right adjoints to monoidal functors are lax monoidal \cite[Corollary 7.3.2.7]{lurie2014higher}, we obtain a composite lax monoidal functor $\End(\cat{D})\rightarrow \End_{/\Delta^1}(\cat{M})\rightarrow \End(\cat{C})$.

To see that this functor indeed sends $T$ to $G\circ T\circ F$, let $X\in \cat{C}\subseteq \cat{M}$. By (the opposite of) \cite[Proposition 4.3.1.9]{lurie2009higher}, the relative right Kan extension of $T\colon \cat{D}\rightarrow \cat{D}\rightarrow \cat{M}$, restricted to $\cat{C}$, can be computed by the right Kan extension of 
$$\begin{tikzcd}
\cat{C}\times_{\Fun(\{0\}, \cat{M})}\Fun(\Delta^1, \cat{M})\times_{\Fun(\{1\}, \cat{M})}\cat{D}\arrow[r, "\pi"] & \cat{D}\arrow[r, "T"] & \cat{D}\arrow[r, "G"] & \cat{C}
\end{tikzcd}$$
along the projection $q\colon \cat{C}\times_{\Fun(\{0\}, \cat{M})}\Fun(\Delta^1, \cat{M})\times_{\Fun(\{1\}, \cat{M})}\cat{D}$. Note that $q$ is a cartesian fibration; its fibre over $X\in \cat{C}$ is given by $\cat{M}_{X/}\times_{\cat{M}}\cat{D}$. Each of these fibres has an initial object, given by the coCartesian arrow $u_X\colon X\rightarrow F(X)$, so that $q$ admits a left adjoint section sending $X$ to $u_X$. The right Kan extension along $q$ is then equivalent to the restriction along this left adjoint; this is precisely $GTF$, as desired (this argument also shows that the relative Kan extension  exists).
\end{proof}
\begin{proof}[Proof of Proposition \ref{prop:levelwise tensor lax monoidal}]
Write  $\cat{bisSeq}_R$ for the free $R$-linear symmetric monoidal $\infty$-category on \emph{two} objects $\mb{1}_L$ and $\mb{1}_R$. Explicitly, this is the $\infty$-category of functors $B\Sigma\times B\Sigma\rightarrow \Mod_R$, with the Day convolution product. There are three natural fully faithful $R$-linear symmetric monoidal functors $\iota_L, \iota_R, \Delta_!\colon \sSeq_R\hookrightarrow \cat{bisSeq}_R$, determined by $\iota_L(\textbf{1})=\textbf{1}_L$, $\iota_R(\textbf{1})=\textbf{1}_R$ and $\Delta_!(\textbf{1})=\textbf{1}_L\otimes \textbf{1}_R$. Write $\Delta^*$ for the functor restricting along the diagonal $\Delta\colon B\Sigma\rightarrow B\Sigma^{\times 2}$. Then $\Delta^*$ is right adjoint to $\Delta_!$, and also  left adjoint (via the norm \cite[Proposition 6.1.6.12]{lurie2014higher}).

Let $\cat{E}\subseteq \End^{\otimes, L}_R(\cat{bisSeq}_R)$ be the full monoidal subcategory of symmetric monoidal $R$-linear endofunctors which furthermore preserve the essential images of $\iota_L$ and $\iota_R$. Evaluation at $\textbf{1}_L$ and $\textbf{1}_R$ then determines an equivalence $\cat{E}\simeq \sSeq_R^{\times 2}$. The inverse sends a pair $(X, Y)$ to the endofunctor of bisymmetric sequences 
\begin{equation}\label{eq:bisymmetric sequence action}
Z\longmapsto \bigoplus_{p, q} Z(p, q)\otimes_{\Sigma_p\times \Sigma_q} \iota_L(X)^{\otimes p}\otimes \iota_R(Y)^{\otimes q}.
\end{equation}
Note that the equivalence $\cat{E}\simeq \sSeq_R^{\times 2}$ identifies composition in $\cat{E}$ with the opposite of the composition product on each of the factors of $\sSeq_R^{\times 2}$. We now consider $\Delta^*$ as the right adjoint to $\Delta_!$; since $\Delta_!$ is symmetric monoidal, $\Delta^*$ inherits a lax symmetric monoidal structure, so that conjugation by $\Delta_!$ and $\Delta^*$ sends symmetric monoidal functors to lax symmetric monoidal functors:
$$\begin{tikzcd}[column sep=4pc]
\sSeq_R^{\times 2}\simeq \cat{E}\arrow[r, "T\mapsto \Delta^*T\Delta_!"] & \End^{\mm{lax}-\otimes}_R(\sSeq_R).
\end{tikzcd}$$
Using Equation \eqref{eq:bisymmetric sequence action}, one sees that the above functor sends $(X, Y)$ simply to the endofunctor $(-)\circ (X\otimes_{\lev} Y)$. In particular, it takes values in the full subcategory $\End^{\otimes, L}_R(\sSeq_R)$ of strong monoidal endofunctors. Applying Lemma \ref{lem:conjugating by adjunction} shows that the above functor is lax monoidal with respect to composition, so that $\otimes_{\lev}$ is indeed lax monoidal for the composition product. Viewing $\Delta^*$ instead as the left adjoint to $\Delta_!$, the opposite of Lemma \ref{lem:conjugating by adjunction} provides the desired oplax monoidal structure. 
\end{proof}

We will now consider symmetric sequences and the composition product in the context of extended functors, cf.\ Section \ref{sec:extended functors}.
\begin{definition}\label{def:finite symmetric set}
A symmetric sequence of sets $X\colon B\Sigma\rightarrow \cat{Set}$ is said to be \emph{finite} if each $X(r)$ is a finite set, which is empty for all but finitely many $r\geq 0$. It is said to be $\Sigma$-free if each $X(r)$ is a (possibly empty) free $\Sigma_r$-set.
\end{definition}

\begin{definition}[Finitely generated free symmetric sequences] \label{ffree}
Let $R$ be a connective $\mathbb{E}_\infty$-ring spectrum. An $R$-linear symmetric sequence is said to be \emph{finitely generated free} if it arises as the $R$-linearisation of a finite $\Sigma$-free sequence of sets. Write  $R[\Sigma]\subseteq \sSeq_R$ for the full subcategory of finitely generated free symmetric sequences. One can identify $R[\Sigma]$ with the smallest full subcategory of $\sSeq_R$ which is closed under finite direct sums and contains all objects $\textbf{1}^{\otimes r}$.
\end{definition}
\begin{remark}
The additive $\infty$-category $R[\Sigma]$ can be identified with the direct sum $\bigoplus_{r\geq 0} \Vect_{R[\Sigma_r]}$ of the additive categories of finitely generated free $R[\Sigma_r]$-modules.
\end{remark}
The objects in $R[\Sigma]$ are compact generators of $\sSeq_R$. Consequently, the fully faithful inclusion $R[\Sigma]\rightarrow \sSeq_R$ induces an equivalence $\Mod_{R[\Sigma]}\simeq \sSeq_R$.
\begin{remark}[Almost perfect symmetric sequences]
The equivalence $\Mod_{R[\Sigma]}\simeq \sSeq_R$ identifies almost perfect $R[\Sigma]$-modules with symmetric sequences $X$ that are almost perfect in the sense that each $X(r)$ is an almost perfect $R[\Sigma_r]$-module and for each $m\geq 0$, $\tau_{\leq m}X(r)$ is trivial for all but finitely many $r$.
\end{remark}
\begin{lemma}\label{lem:products extended}
The full subcategory $R[\Sigma]\hookrightarrow \sSeq_R$ is closed under the Day tensor product $\otimes$, the levelwise tensor product $\otimes_{\lev}$ and the composition product $\circ$. Furthermore, $\otimes$ and $\otimes_{\lev}$ are additive in each variable and $\circ\colon R[\Sigma]\times R[\Sigma]\rightarrow R[\Sigma]$ is locally polynomial and additive in the first variable.
\end{lemma}
Note that $R[\Sigma]$ contains the monoidal unit for $\otimes$ and $\circ$, but not for $\otimes_{\lev}$.
\begin{proof}
It is clear that $R[\Sigma]$ is closed under $\otimes$ and $\otimes_{\lev}$, and since both tensor products preserve colimits in each variable, their restrictions are additive in each variable. Equation \eqref{eq:composition product formula} then implies that $R[\Sigma]$ is closed under the composition product as well. Furthermore, we can write the composition product functor $\circ$ as a filtered colimit $X\circ Y=\colim_n F_n(X, Y)$, where $F_n(X, Y)=\bigoplus_{r\leq n} X(r)\otimes_{\Sigma_r} Y^{\otimes r}$. 

Since each functor $Y\mapsto Y^{\otimes r}$ is of degree $r$, it follows that $F_n$ is of degree $n$. On the other hand, the sequence of $F_n(X, Y)$ stabilizes since every $X\in R[\Sigma]$ is concentrated in finitely many arities.
\end{proof}
Combining \Cref{lem:products extended} and \Cref{cor:polynomial functors between coherent}, we can deduce:
\begin{corollary}\label{cor:sseq products as extensions}
All three tensor products $\otimes, \otimes_{\lev}, \circ$ are the right-left extension of their restriction to $R[\Sigma]$.
\end{corollary}

\subsection{Pro-coherent symmetric sequences and PD operads}\label{sec:PD operads}
Using \Cref{ffree},  we can introduce a refined version of symmetric sequences; linear duals of ordinary symmetric sequences are naturally equipped with this structure.
\begin{definition}[Pro-coherent symmetric sequences] \label{procosseq}
Let $R$ be a coherent (connective) $\mathbb{E}_\infty$-ring spectrum. A \emph{pro-coherent symmetric sequence} over $R$ is a pro-coherent module over the coherent additive $\infty$-category $R[\Sigma]$. We will write $\sSeq^\vee_R$ for the $\infty$-category of pro-coherent symmetric sequences over $R$.
\end{definition}
\begin{proposition}\label{prop:pro-coh composition}
Let $R$ be a coherent $\mathbb{E}_\infty$-ring spectrum. Then the $\infty$-category of pro-coherent symmetric sequences can be equipped with
\begin{enumerate}
\item a closed symmetric monoidal structure $\otimes$;
\item a composition product $\circ$  preserving sifted colimits in each variable and small colimits in the first variable;
\item a sifted-colimit-preserving action $\circ\colon \sSeq^\vee_R\times \QC^\vee_R\rightarrow \QC^\vee_R$ of $\big(\sSeq^\vee_R, \circ\big)$;
\item a closed symmetric monoidal structure $\otimes_{\lev}$, together with a lax and oplax monoidal structure on $\otimes_{\lev}\colon \sSeq^\vee_R\times\sSeq^\vee_R\rightarrow \sSeq^\vee_R$ with respect to the composition product;
\end{enumerate}
which are 
  right-left extended from the corresponding functors on  the \mbox{$\infty$-category $R[\Sigma]$.} Furthermore, the natural functors $  \sSeq_R\rightarrow \sSeq_R^\vee$ and $ \Mod_R\rightarrow \QC^\vee_R$ intertwine all of these monoidal structures.
\end{proposition}
\begin{proof}
Almost all assertions follow from \Cref{thm:derived functors} and Lemma \ref{lem:products extended}. To see that $\otimes_{\lev}$ is (op)lax monoidal with respect with the composition product, we use Proposition \ref{prop:derived functor of lax monoidal} and Proposition \ref{prop:levelwise tensor lax monoidal}. Finally, note that $\otimes_{\lev}$ a priori only defines a nonunital symmetric monoidal structure on $\sSeq^\vee_R$ (because it does not have a monoidal unit contained in $R[\Sigma]$). However, the image $\upiota(\mathbb{E}_{\infty,R})$ is easily seen to provide a (connective) unit: indeed, $\upiota(\mathbb{E}_{\infty,R})\otimes_{\lev} (-)$ is the right-left extended functor of its restriction to $R[\Sigma]$, which is the identity since $\mathbb{E}_{\infty,R}$ serves as the unit for $\otimes_{\lev}$ in the $\infty$-category of symmetric sequences.
\end{proof}
The composition product $\circ$ on $\sSeq_R^\vee$ coincides with the usual composition product on ordinary symmetric sequences. Surprisingly, there are many other pro-coherent symmetric sequences on which $\circ$  acts like a \emph{restricted composition product}:

\begin{proposition}\label{prop:formula for pro-coh comp prod} 
Given $X, Y \in \sSeq_R^\vee$, there is a natural map
$$\begin{tikzcd}
\nu\colon X\circ Y \arrow[r] & \bigoplus_{r\geq 0} \big(X(r)\otimes Y^{\otimes r}\big)^{h\Sigma_r}
\end{tikzcd}$$
which is an equivalence whenever $X$ and $Y$ are dually almost perfect (cf.\ \Cref{daperf} for $\cat{A}=R[\Sigma])$. If $R$ is eventually coconnective, it is furthermore an equivalence when both $X$ and $Y$ are the colimits of filtered diagrams in $\APerf^\vee_{R[\Sigma], \weirdleq m} := (\APerf_{R[\Sigma], \geq m})^\vee$, 
 for some $m$.
\end{proposition}
\begin{proof}
Since the composition product is obtained by right-left extension, it suffices to describe $\nu$ when $X$ and $Y$ are contained in $\APerf^\vee_{R[\Sigma], \weirdleq 0}$. In turn, the domain and codomain of $\nu$ are both functors that are right Kan extended from $R[\Sigma]$ to $\APerf^\vee_{R[\Sigma], \weirdleq 0}$. It therefore remains to describe $\nu$ when $X$ and $Y$ are finitely generated free. In this case, the norm map provides a natural equivalence $\nu\colon X\circ Y\xrightarrow{\   \simeq \ } \bigoplus_{r\geq 0} \big(X(r)\otimes Y^{\otimes r}\big){}^{h\Sigma_r}$, because $X(r)$ is $\Sigma_r$-free. 

In particular, this implies that the resulting map $\nu$ is an equivalence for all $X, Y$ in $\APerf^\vee_{R[\Sigma],  \weirdleq 0}$. As both domain and codomain   preserve geometric realisations of skeletal diagrams, $\nu$ is an equivalence whenever $X$ and $Y$ are dually \mbox{almost perfect.}

Finally, suppose that $X, Y$ are colimits of filtered diagrams in $\APerf^\vee_{R[\Sigma], \weirdleq m}$. Then each $X(r)\otimes Y^{\otimes r}$ is a filtered colimit of objects in $\APerf^\vee_{R[\Sigma], \weirdleq m}$ as well. Under the assumption that $R$ is $n$-coconnective, $R[\Sigma]$ is also $n$-coconnective and one finds that $X(r)\otimes Y^{\otimes r}$ is a filtered colimit of $n'$-coconnective objects for some $n'$. Taking homotopy fixed points commutes with such filtered colimits of $n'$-coconnective objects, so that $\nu$ is an equivalence for $X$ and $Y$ as well.
\end{proof}
\begin{definition}[PD $\infty$-operads]\label{def:PD operads}
Let $R$ be a coherent $\mathbb{E}_\infty$-ring spectrum. A \emph{PD $\infty$-operad} $\PP$ over $R$ is an associative algebra in the $\infty$-category of pro-coherent symmetric sequences, with respect to the composition product. We will write $\Op^{\pd}_R$ for the $\infty$-category of PD $\infty$-operads.

An \emph{algebra} over a PD $\infty$-operad $\PP$ is a pro-coherent $R$-module $A$ equipped with a left $\PP$-module structure with respect to the composition action. We will denote the $\infty$-category of (pro-coherent) $\PP$-algebras by $\Alg_{\PP}(\QC^\vee_R)$.
\end{definition}
\begin{example}[Underlying operads]
Every ordinary $\infty$-operad gives rise to a PD $\infty$-operad via the functor $\upiota\colon \sSeq_R\rt \sSeq_R^\vee$. Conversely, every PD $\infty$-operad has an underlying $\infty$-operad, via the right adjoint $\upupsilon\colon \sSeq^\vee_R\rt \sSeq_R$.
\end{example}
\begin{remark}
The action of pro-coherent symmetric sequences on $\QC^\vee_R$ defines a sifted-colimit-preserving monoidal functor $\sSeq^\vee_R\rightarrow \End_{\Sigma}(\QC^\vee_R)$ with respect to the composition product. Using Lemma \ref{lem:conjugating by adjunction},  conjugating by the adjoint pair $\upiota\colon \Mod_R\leftrightarrows \QC^\vee_R\colon \upupsilon$ yields a lax monoidal functor $\sSeq^\vee_R\rightarrow \End(\Mod_R)$. In particular, every PD $\infty$-operad $\PP$ determines a monad $T_{\PP}$ on $\Mod_R$. This monad \emph{differs} from the monad induced by the underlying $\infty$-operad of $\PP$.
When $R$ is eventually coconnective, $\upupsilon$ preserves colimits (Remark \ref{rem:modules as subcat of pro-coh}) so that $T_{\PP}$ preserves sifted colimits.  
\end{remark}
\begin{example}
Let $k$ be a field, so that $\QC^\vee_k\simeq \Mod_k$ (Example \ref{ex:regular}). Suppose that $\PP$ is a PD $\infty$-operad over $k$ which is dually almost perfect. By Proposition \ref{prop:formula for pro-coh comp prod}, $\PP$ determines a monad on $\Mod_k$ which preserves sifted colimits and is given on eventually coconnective $k$-modules by
$$
\mm{Free}_{\PP}(V) = \bigoplus_{r\geq 0} \big(\PP(r)\otimes V^{\otimes r}\big)^{h\Sigma_r}.
$$
We will produce examples of these kinds of $\infty$-operads by Koszul duality.
\end{example}

\subsection{Refined Koszul duality}
\label{sec:refined koszul}
We will now discuss a refinement of the classical operadic Koszul duality functor \cite{GinzburgKapranov:KDO, fresse346koszul, salvatore1998configuration, ching2005bar} to the setting of PD $\infty$-operads. Recall that the classical Koszul duality functor is defined in two steps. First, every augmented $\infty$-operad gives rise to an $\infty$-cooperad by the bar construction. One then takes the Spanier--Whitehead dual of the bar construction to obtain an $\infty$-operad, usually referred to as the (classical) \emph{Koszul dual $\infty$-operad}. We will refine each of these two steps to the setting of PD $\infty$-operads.

\subsubsection*{The bar construction for PD operads}
Recall that for any monoidal $\infty$-category $\cat{C}$ with geometric realisations and totalisations, there is an adjoint pair
$$\begin{tikzcd}
\Barr\colon \Alg^{\aug}(\cat{C})\arrow[r, yshift=1ex] & \Coalg^{\aug}(\cat{C})=\Alg^{\aug}(\cat{C}^{\op})^{\op}\colon \Cobar\arrow[l, yshift=-1ex]
\end{tikzcd}$$ 
given by the $\infty$-categorical bar construction and cobar construction \cite[Section 5.2.2]{lurie2014higher}. If $A$ is an augmented algebra in $\cat{C}$, the underlying object of $\Barr(A)$ can be identified with the relative tensor product $\mb{1}\otimes_A \mb{1}$, computed as the realisation of the two-sided simplicial bar construction $\Barr_\bullet(\mb{1}, A, \mb{1})$. We will give a more rigorous account of the bar construction below (see in particular Corollary \ref{cor:luriebar}), including a few arguments that are not completely worked out in \cite[Section 5.2.2]{lurie2014higher} (as pointed out in \cite{dancohenhorev2022koszul}). For now, specialising to the case where $\cat{C}$ is the $\infty$-category of pro-coherent symmetric sequences, we obtain:
\begin{definition}[Bar construction for PD operads]
Let $R$ be a coherent $\mathbb{E}_\infty$-ring spectrum. We will write $\Barr\colon \Op^{\pd, \aug}_R\leftrightarrows \Coop^{\pd, \aug}_R\colon \Cobar$ for the $\infty$-categorical bar and cobar construction in the $\infty$-category $\sSeq^\vee_R$ of pro-coherent symmetric sequences, with respect to the composition product $\circ$.
\end{definition}
Our next goal will be to relate algebras over an augmented PD $\infty$-operad $\PP$ to coalgebras over its bar construction $\Barr(\PP)$. To do this, we will need a variant of the $\infty$-categorical bar construction of \cite[Section 5.2.2]{lurie2014higher} for left modules and left comodules. 
\begin{notation}[Bimodule $\infty$-categories]
Recall that a \emph{bimodule $\infty$-category} is a triple $(\Cl , \Cm, \Cr)$ consisting of monoidal $\infty$-categories $\Cl$ and $\Cr$, together with commuting left and right actions $\Cl\curvearrowright \Cm\curvearrowleft \Cr$; more precisely, it is an algebra in $\mathscr{C}\mm{at}_\infty$ over the coloured operad $\mathscr{BM}$ from \cite[Definition 4.3.1.1]{lurie2014higher}.

A \emph{left module} in $\Cm$ is given by a tuple $(A, M)$ consisting of an associative algebra $A\in \Alg(\Cl)$ and a left $A$-module $M\in \LMod_A(\Cm)$. We will write $\LMod(\Cm)$ for the $\infty$-category of left modules in $\Cm$. The canonical projection $\pi\colon \LMod(\Cm)\rt \Alg(\Cl)$ is a cartesian fibration \cite[Corollary 4.2.3.2]{lurie2014higher}. Furthermore, $\pi$ is a map of right $\Cr$-module categories, where $\Cr$ acts trivially on $\Alg(\Cl)$ \cite[Proposition 4.3.2.5, Proposition 4.3.2.6]{lurie2014higher}. On  underlying objects, the tensoring of a left $A$-module $M$ in $\Cm$ with $X\in \Cr$ is given by the left $A$-module $M\otimes X$.
\end{notation}
\begin{theorem}\label{thm:bar for modules}
Let $(\Cl, \Cm, \Cr)$ be a bimodule $\infty$-category such that $\Cl, \Cm$ and $\Cr$ all admit geometric realisations and totalisations, and such that the units $\mb{1}_{\Cl}$ and $\mb{1}_{\Cr}$ are both terminal and initial. Then there is a commuting diagram
$$\begin{tikzcd}[column sep=3pc, row sep=2pc]
\LMod(\Cm)\arrow[d, "\pi"{swap}]\arrow[r, yshift=1ex, "\Barr"] & \cat{LComod}(\Cm)=\cat{LMod}(\Cm^{\op})^{\op}\arrow[d, xshift=-7ex, "\pi"]\arrow[l, yshift=-1ex, "\Cobar"]\\
\Alg(\Cl)\arrow[r, yshift=1ex, "\Barr"] & \Coalg(\Cl)=\Alg(\Cl^{\op})^{\op}\arrow[l, yshift=-1ex, "\Cobar"]
\end{tikzcd}$$
where the rows are adjunctions. Furthermore, the following assertions hold:
\begin{enumerate}
\item The functor $\Barr\colon \LMod(\Cm)\rt \cat{LComod}(\Cm)$ preserves coCartesian arrows.

\item If the action $\Cm\times \Cr\rt \Cm$ preserves geometric realisations in the first variable, then \mbox{$\Barr\colon \LMod(\Cm)\rt \cat{LComod}(\Cm)$} is a right $\Cr$-linear functor.
\end{enumerate}
\end{theorem}
\begin{remark}
Given an associative algebra  $A$ in $\Cl$, we obtain a functor $\Barr_A\colon \LMod_A(\Cm)\rt \LComod_{\Barr(A)}(\Cm)$ between fibres. This admits a right adjoint, which first applies the functor $\Cobar\colon \LComod_{\Barr(A)}(\Cm)\rt \LMod_{\Cobar(\Barr(A))}(\Cm)$ and then restricts scalars along the unit map $A\rt \Cobar(\Barr(A))$.
\end{remark}
We postpone the proof of \Cref{thm:bar for modules} to the end of this section and first discuss some applications.
To start, suppose that $A$ is an associative algebra in $\Cl$ with augmentation $\epsilon\colon A\rt \mb{1}$. Restriction and induction along $\epsilon$ define an adjoint pair (cf.\ the proof of \cite[Proposition 5.2.2.5]{lurie2014higher})
$$\begin{tikzcd}
\epsilon_!\colon \LMod_A(\Cm)\arrow[r, yshift=1ex] &  \LMod_{\mb{1}}(\Cm)\simeq \Cm\colon \epsilon^*.\arrow[l, yshift=-1ex]
\end{tikzcd}$$
Explicitly, $\epsilon_!$ sends each left $A$-module $M$ to the geometric realisation of the simplicial bar construction $\Barr_\bullet(\mb{1}, A, M)$. Considering \Cref{thm:bar for modules} in the case where $\Cr=\ast$ then yields:
\begin{proposition}\label{prop:comonads from bar construction}
Let $(\Cl, \Cm)$ be a left module $\infty$-category as in \Cref{thm:bar for modules} and suppose that the left action $\Cl\times \Cm\rt \Cm$ preserves geometric realisations in the first variable. For any associative algebra $A\in \Alg(\Cl)$, there is an equivalence of comonads on $\Cm$
$$
\epsilon_!\circ\epsilon^*\simeq \Barr(A)\otimes (-).
$$
\end{proposition}
\begin{proof}
Part (1) of \Cref{thm:bar for modules} provides a commuting triangle of left adjoints
$$\begin{tikzcd}[column sep=1.7pc, row sep=1pc]
\LMod_A(\Cm)\arrow[rd, "\epsilon_!"{swap}]\arrow[rr, "\Barr_A"] & & \LComod_{\Barr(A)}(\Cm)\arrow[ld, "\mm{forget}"]\\
& \Cm. & 
\end{tikzcd}$$
By \cite[Corollary 5.8, Corollary 8.9]{haugseng2020monads}, this induces a natural map of comonads $\mu\colon \epsilon_!\circ \epsilon^*\rt \Barr(A)\otimes (-)$. It remains to verify that the underlying map of endofunctors of $\Cm$ is an equivalence. For each $M\in \Cm$, the map $\mu$ can be identified with the natural map $|\Barr_\bullet(\mb{1}, A, \epsilon^*M)|\rt |\Barr_\bullet(\mb{1}, A, \mb{1})|\otimes M$. This map is an equivalence by the assumption that the left action of $\Cl$ on $\Cm$ preserves geometric realisations in $\Cl$.
\end{proof}
\begin{notation}[Trivial algebras and cotangent fibre]
Let $\PP$ be an augmented PD $\infty$-operad with augmentation $\epsilon\colon \PP\rt \mb{1}$. We denote by $\cot_{\PP}\colon \Alg_{\PP}\leftrightarrows \Mod_R\colon \triv_{\PP}$ the adjoint pair induced by the augmentation $\epsilon$. We will refer to these functors as taking \emph{cotangent fibre}, respectively  \emph{trivial} $\PP$-algebra. 
\end{notation}
\begin{corollary}\label{cor:bar construction for algebras}
Let $\PP$ be an augmented PD $\infty$-operad. Then there is a commuting diagram of left adjoint functors
$$\begin{tikzcd}
\Alg_{\PP}\arrow[rr, "\Barr_{\PP}"]\arrow[rd, "\cot_{\PP}"{swap}] & & \Coalg_{\Barr(\PP)}\arrow[ld, "\mm{forget}"]\\
& \Mod_R
\end{tikzcd}$$
and an equivalence of comonads $\mm{cot}_{\PP}\circ \mm{triv}\simeq \Barr(\PP)$.
\end{corollary}
\begin{proof}
Apply Proposition \ref{prop:comonads from bar construction} where $\Cr=\ast$ and $\Cl$ is the $\infty$-category $\sSeq_{1//1}$ of augmented symmetric sequences with the composition product, acting from the \mbox{left on $\Mod_R$.}
\end{proof}
As another application of \Cref{thm:bar for modules}, we shall give another possible definition of the bar construction of an associative algebra, due to Lurie \cite{Unstablechrom}; it is more convenient for later purposes.
\begin{definition}[Coendomorphisms object]
Let $\Cl$ be a monoidal $\infty$-category, $\Cm$ a left $\Cl$-module $\infty$-category and $M$ an object in $\Cm$. Consider an object $X\in \Cl$ together with a map $\lambda\colon M\rt X\otimes M$ in $\Cm$. We will say that $\lambda$ exhibits $X$ as a \emph{coendomorphism object} of $M$ if for every object $Y$ in $\Cl$, the natural map
$$\begin{tikzcd}
\Map_{\Cl}(X, Y)\arrow[r] & \Map_{\Cm}(X\otimes M, Y\otimes M)\arrow[r, "\lambda^*"] & \Map_{\Cm}(M, Y\otimes M)
\end{tikzcd}$$
is an equivalence. Similarly, let $C\in \cat{Coalg}(\Cl)$ be a coassociative coalgebra in $\Cl$ and denote by $\LComod_C(M)=\LComod_C(\Cm)\times_{\Cm} \{M\}$ the space of left $C$-comodule structures on $M$. Then $\lambda\in \LComod_C(M)$ exhibits $C$ as a \emph{coendomorphism coalgebra} of $M$ if for each coalgebra $D$, the natural map which corestricts the coaction of $C$ on $M$ to a coaction of $D$ on $M$ is an equivalence:
$$\begin{tikzcd}
\Map_{\Coalg(\Cl)}(C, D)\arrow[r, "\sim"] & \LComod_D(M); & f\arrow[r, mapsto] & (f\otimes \mm{id})\circ \lambda
\end{tikzcd}$$
\end{definition}
\begin{lemma}[{cf.\ \cite[Proposition 7]{Unstablechrom}}]\label{lem:coendomorphism objects}
Let $(\Cl, \Cm)$ be a left module $\infty$-category and $M\in \Cm$. Then the following assertions hold:
\begin{enumerate}
\item Let $\lambda\in \LComod_C(M)$. Then $\lambda$ exhibits $C$ as a coendomorphism coalgebra of $M$ if and only if the underlying map $M\rt C\otimes M$ exhibits $C$ as a coendomorphism object of $M$.

\item Suppose that $M$ admits a coendomorphism object $X$. Then $M$ admits a coendomorphism coalgebra.
\end{enumerate}
\end{lemma}
\begin{proof}
By the dual of \cite[Theorem 4.7.1.34]{lurie2014higher}, there exists a monoidal $\infty$-category $\Cl[M]$ with  objects given by tuples of objects $X\in \Cl$ and maps $M\rt X\otimes M$ in $\Cm$, such that $\Coalg(\Cl[M])$ is equivalent to the $\infty$-category of coalgebras together with a left comodule structure on $M$. By definition, a coendomorphism object of $M$ is an initial object of $\Cl[M]$, while a coendomorphism coalgebra of $M$ is an initial object of $\Coalg(\Cl[M])$. The assertions then follow from the fact that the forgetful functor $\Coalg(\Cl[M])\rt \Cl[M]$ preserves and detects initial objects.
\end{proof}
\begin{cons}[Koszul complex]\label{cons:koszul complex}
Let $\cat{C}$ be a monoidal $\infty$-category such that the tensor product preserves geometric realisations in the first variable and the monoidal unit $\mb{1}$ is both initial and terminal. We consider $\cat{C}$ as a bimodule $\infty$-category over itself. If $A$ is an associative algebra in $\cat{C}$, then \Cref{thm:bar for modules} provides a right $\cat{C}$-linear functor $\Barr_A\colon \LMod_A(\cat{C})\rt \LComod_{\Barr(A)}(\cat{C})$. Write $K(A)$ for the value of this functor on the free left $A$-module $A$; it follows from part (1) of \Cref{thm:bar for modules} this is simply given by the trivial comodule $\mb{1}$. Since the free $A$-module $A$ has a commuting right $A$-module structure, we obtain a natural object
$$
K(A)\in \RMod_A\big(\LComod_{\Barr(A)}(\cat{C})\big)\simeq \LComod_{\Barr(A)}\big(\RMod_A(\cat{C})\big)
$$
such that the underlying right $A$-module is the terminal object $\mb{1}$. We will refer to $K(A)$ as the \emph{Koszul complex} of $A$.
\end{cons}
\begin{proposition}\label{prop:bar construction as coendomorphisms}
Let $\cat{C}$ be a monoidal $\infty$-category such that the tensor product preserves geometric realisations in the first variable and the monoidal unit $\mb{1}$ is both initial and terminal. Then the left $\Barr(A)$-comodule structure on the Koszul complex $K(A)\simeq \mb{1}$ exhibits $\Barr(A)$ as the coendomorphism coalgebra of $\mb{1}\in \RMod_A(\cat{C})$.
\end{proposition}
\begin{proof}
By Lemma \ref{lem:coendomorphism objects}, it suffices to verify that the right $A$-linear map $\mb{1}\rt \Barr(A)\otimes \mb{1}$ exhibits $\Barr(A)$ as an endomorphism object of the trivial right $A$-module $\mb{1}$, i.e.\ for every object $Y\in \cat{C}$, the map
$$\begin{tikzcd}
\Map_{\cat{C}}\big(\Barr(A), Y\big)\arrow[r] & \Map_{\RMod_A(\cat{C})}\big(\Barr(A)\otimes \mb{1}, Y\otimes \mb{1}\big)\arrow[r] & \Map_{\RMod_A(\cat{C})}\big(\mb{1}, Y\otimes \mb{1}\big)
\end{tikzcd}$$
is an equivalence. To see this, using Proposition \ref{prop:comonads from bar construction} and writing $\epsilon\colon A\rt \mb{1}$ for the augmentation, the above map can be identified with the composite
$$\begin{tikzcd}
\Map_{\cat{C}}\big(\epsilon_!\epsilon^*(\mb{1}), Y\big)\arrow[r] & \Map_{\RMod_A(\cat{C})}\big(\epsilon^*\epsilon_!\epsilon^*(\mb{1}), \epsilon^*Y\big)\arrow[r] & \Map_{\RMod_A(\cat{C})}\big(\epsilon^*(\mb{1}), \epsilon^*Y\big)
\end{tikzcd}$$
where the first map applies $\epsilon^*$ and the second map restricts along the unit $\epsilon^*(\mb{1})\rt \epsilon^*\epsilon_!\epsilon^*(\mb{1})$. This composite is an equivalence since $(\epsilon_!, \epsilon^*)$ is an adjoint pair.
\end{proof}
Finally, we turn to the proof of \Cref{thm:bar for modules}. The argument, which we learned from  Lurie, is a direct modification of the construction of the bar-cobar adjunction in \cite[Section 5.2.2]{lurie2014higher}. We start by recalling some terminology from loc.\ cit.
\begin{notation}[Pairings]\label{not:pairing}
Recall that a \emph{pairing} of $\infty$-categories $\cat{C}$ and $\cat{D}$ is a right fibration $\lambda\colon \cat{M}\rt \cat{C}\times \cat{D}$. An object $M\in \cat{M}$ with image $(C, D)$ is called \emph{left universal} if it is terminal in $\{C\}\times_{\cat{C}} \cat{M}$. The pairing $\lambda$ is called \emph{left representable} if every $C\in \cat{C}$ is the image of a left universal object. We denote by $\cat{CPair}\subseteq \Fun(\Lambda^0[2], \mathscr{C}\mm{at}_\infty)$ the full subcategory spanned by the pairings of $\infty$-categories, and by $\cat{CPair}^\mm{L}$ the subcategory of $\cat{CPair}$ on the left representable pairings and maps of pairings preserving left universal objects. Both $\cat{CPair}$ and $\cat{CPair}^\mm{L}$ are closed under the cartesian product in $\Fun(\Lambda^0[2], \mathscr{C}\mm{at}_\infty)$.

There is an equivalence of $\infty$-categories $\cat{Pair}^\mm{L}\simeq \Fun([1], \mathscr{C}\mm{at}_\infty)$ \cite[Proposition 2.2]{li2015stack}. By \cite[Construction 5.2.1.9]{lurie2014higher}, this equivalence sends a left representable pairing $\lambda\colon \cat{M}\rt \cat{C}\times\cat{D}$
to the unique functor $F_\lambda\colon \cat{C}\rt \cat{D}^{\op}$ which admits a natural equivalence
$$
\lambda^{-1}(C, D)\simeq \Map_{\cat{D}^{\op}}\big(F_\lambda(C), D\big).
$$
\end{notation}
\begin{example}\label{ex:twisted arrows of bimodule cat}
Let $\cat{C}$ be a bimodule $\infty$-category, given by $\Cl\curvearrowright \Cm \curvearrowleft \Cr$. Because taking twisted arrow $\infty$-categories preserves products, one obtains a bimodule object $\Tw(\cat{C})$ in the $\infty$-category of left (and right) representable pairings of the form
\vspace{-50pt}
$$\begin{tikzcd}[column sep=1.2pc]
\Tw(\Cl)\arrow[d]\arrow[r, bend left=90, looseness=2, start anchor={[xshift=6ex, yshift=-2.5ex]}, end anchor={[xshift=-6ex, yshift=-2.5ex]}] & \Tw(\Cm)\arrow[d] \arrow[r, bend left=90, , looseness=2, start anchor={[xshift=6ex, yshift=-2.5ex]}, end anchor={[xshift=-6ex, yshift=-2.5ex]}, leftarrow] & \Tw(\Cr)\arrow[d]\\
\Cl\times \Cl^{\op}\arrow[r, bend left=90, , looseness=2, start anchor={[xshift=6ex, yshift=-2.5ex]}, end anchor={[xshift=-6ex, yshift=-2.5ex]}] & \Cm\times\Cm^{\op} \arrow[r, bend left=90, , looseness=2, start anchor={[xshift=6ex, yshift=-2.5ex]}, end anchor={[xshift=-6ex, yshift=-2.5ex]}, leftarrow] & \Cr\times\Cr^{\op}.
\end{tikzcd}$$
Taking left module objects, we then obtain a map of pairings $\LMod(\Tw(\Cm))\rt \Alg(\Tw(\Cl))$, together with a fibrewise right action of $\Tw(\Cr)$ on $\LMod(\Tw(\Cm))$ \cite[Proposition 4.3.2.5, Proposition 4.3.2.6]{lurie2014higher}.
\end{example}
\begin{lemma}\label{lem:bimodule pairing}
Let $\lambda\colon (\Ml, \Mm, \Mr)\rt (\Cl, \Cm, \Cr)\times (\Dl, \Dm, \Dr)$ be a bimodule object in the $\infty$-category $\cat{Pair}$ of pairings and let 
$$
A_\mb{1}\in \Alg(\Ml)\times_{\Alg(\Dl)} \{\mb{1}\} \qquad\qquad \text{and}\qquad\qquad B_{\mb{1}}\in \Alg(\Mr)\times_{\Alg(\Dr)} \{\mb{1}\}
$$
be two algebras with images $(A, \mb{1})$ and $(B, \mb{1})$ in $\Alg(\Cl)\times \Alg(\Dl)$ and $\Alg(\Cr)\times \Alg(\Dr)$. Consider the induced pairing between categories of bimodules 
$$
\lambda_{A, B}\colon {}_{A_\mb{1}}\cat{BMod}_{B_{\mb{1}}}(\Mm)\rt {}_A\cat{BMod}_B(\Cm)\times \Dm
$$
(where we identify ${}_{\mb{1}}\cat{BMod}_\mb{1}(\Dm)\simeq \Dm$). Then the following assertions hold:
\begin{enumerate}
\item If $\lambda$ is left representable and the $\infty$-category $\Dm$ admits totalisations of cosimplicial objects, then $\lambda_{A, B}$ is left representable.

\item Suppose that the pairing $\lambda\colon \Mm\rt \Cm\times \Dm$ is right representable and that there exist augmentations $A_{\mb{1}}\rt \mb{1}$ and $B_{\mb{1}}\rt \mb{1}$ in $\Alg(\Ml)$ and $\Alg(\Mr)$. Then $\lambda_{A, B}$ is right representable and the associated functor can be identified with
$$\begin{tikzcd}
G_{\lambda_{A, B}}\colon \Dm\arrow[r, "G_{\lambda}"] & \Cm\arrow[r, "\triv"] & {}_A\cat{BMod}_B(\Cm)
\end{tikzcd}$$
where the second functor is the restriction along the induced augmentations $A\rt \mb{1}$ and $B\rt \mb{1}$ in $\Alg(\Cl)$ and $\Alg(\Cr)$.
\end{enumerate}
\end{lemma}
\begin{remark}\label{rem:bimodules fiber}
Suppose we are in the situation of Lemma \ref{lem:bimodule pairing} and fix an object $D\in \Dm$. Then the actions $\Ml\curvearrowright \Mm\curvearrowleft \Mr$ restrict to actions $\Ml\times_{\Dl}\{\mb{1}\}\curvearrowright \Mm\times_{\Dm}\{D\}\curvearrowleft \Mr\times_{\Dr}\{\mb{1}\}$ and there is a natural equivalence
$$\begin{tikzcd}
{}_{A_\mb{1}}\cat{BMod}_{B_{\mb{1}}}\big(\Mm\times_{\Dm} \{D\}\big)\arrow[r, "\sim"] & {}_{A_\mb{1}}\cat{BMod}_{B_{\mb{1}}}(\Mm)\times_{\Dm} \{D\}.
\end{tikzcd}$$
\end{remark}
\begin{proof}
Part (1) follows from the following adaptation of \cite[Lemma 5.2.2.40]{lurie2014higher}. For a bimodule $M\in {}_A\cat{BMod}_B(\Cm)$, we have to show that induced right fibration 
$$\begin{tikzcd}
\cat{E}_M = \{M\}\times_{{}_A\cat{BMod}_B(\Cm)} {}_{A_\mb{1}}\cat{BMod}_{B_{\mb{1}}}(\Mm)\arrow[r] & \Dm
\end{tikzcd}$$
is representable (i.e.\ $\cat{E}_M$ admits a terminal object). To do this, we will proceed in two steps.

First, let us suppose that $M=A\otimes V\otimes B$ is the free $A$-$B$-bimodule on an object $V\in \Cm$. In this case, the right fibration $\cat{E}_M\rt \Dm$ is representable by the same argument as \cite[Lemma 5.2.2.32]{lurie2014higher}: taking a left representable object $\tilde{V}\in \{V\}\times_{\Cm}\Mm$, the free bimodule $A_{\mb{1}}\otimes \tilde{V}\otimes B_{\mb{1}}$ is a terminal object in $\cat{E}_{A\otimes V\otimes B}$.

For a general bimodule $M$, let $M_\bullet=A^{\otimes \bullet+1}\otimes M\otimes B^{\otimes 1+\bullet}$ be its bar construction, so that $M=|M_\bullet|$. Let $\chi_M\colon \Dm^{\op}\rt \An$ be the presheaf classified by the right fibration $\cat{E}_M\rt \Dm$. We claim that $\chi_M\simeq \lim \chi_{M_\bullet}$. Assuming this, it follows that $\chi_M$ is representable, because it is a totalisation of representable presheaves and $\Dm$ admits totalisations.

It suffices to verify the claim at each point $D\in \Dm$. To do this, consider the commuting diagram
$$\begin{tikzcd}
\cat{F}_D={}_{A_\mb{1}}\cat{BMod}_{B_{\mb{1}}}(\Mm)\times_{\Dm} \{D\}\arrow[r, "g"]\arrow[d, "q"{swap}] & \Mm\times_{\Dm} \{D\}\arrow[d, "p"]\\
{}_A\cat{BMod}_B(\Cm)\arrow[r, "g'"] & \Cm
\end{tikzcd}$$
where the horizontal functors forget the bimodule structure. By \cite[Corollary 5.2.2.39]{lurie2014higher}, it now suffices to check that for every simplicial object $N_\bullet\colon \Del^{\op}\rt \cat{F}_D$ lifting $M_\bullet$, there exists a geometric realisation in $\cat{F}_D$ that is preserved by $q$. For any such $N_\bullet$, the image $g(N_\bullet)$ is a lift of the image $g'(M_\bullet)$ of the simplicial bar construction, which is split. By \cite[Corollary 4.7.2.11]{lurie2014higher}, $g(N_\bullet)$ is a split simplicial object as well. Remark \ref{rem:bimodules fiber} now implies that $N_\bullet$ admits a realisation in $\cat{F}_B$ (by monadicity). To see that this realisation is preserved by $q$, let $N^+_\bullet\colon \Del_+^{\op}\rt \cat{F}_B$ denote the resulting colimiting cocone. Then $g(N^+_\bullet)$ is split and hence $q(N^+_\bullet)$ is sent to a split augmented simplicial object by the forgetful functor $g'$. Again by monadicity, this implies that $q(N^+_\bullet)$ is a colimiting cocone in ${}_A\cat{BMod}_B(\Cm)$, as desired.

\medskip

Part (2) follows from follows from the fact that restriction along the augmentations $A_{\mb{1}}\rt \mb{1}$ and $B_{\mb{1}}\rt \mb{1}$ yields a functor
\begin{equation}\label{diag:bimod right rep}\begin{tikzcd}
\Mm\simeq {}_{\mb{1}}\cat{BMod}_{\mb{1}}(\Mm)\arrow[r] & {}_{A_\mb{1}}\cat{BMod}_{B_{\mb{1}}}(\Mm)
\end{tikzcd}\end{equation}
that preserves right representable objects (by \cite[Proposition 5.2.1.17]{lurie2014higher}). By definition, such representable objects are the terminal objects in the fibres over each $D\in \Dm$. Remark \ref{rem:bimodules fiber} now identifies the functor between fibres with the functor
$$\begin{tikzcd}
\Mm\times_{\Dm}\{M\}\arrow[r] & {}_{A_{\mb{1}}}\cat{BMod}_{B_{\mb{1}}}\big(\Mm\times_{\Dm}\{M\}\big)
\end{tikzcd}$$
restricting along the augmentations of $A_{\mb{1}}$ and $B_{\mb{1}}$. This functor preserves terminal objects.
\end{proof}
The pairing from Lemma \ref{lem:bimodule pairing} has some additional structure in the case where each bimodule category arises from the natural two-sided action of a monoidal category on itself:
\begin{notation}
If $\cat{C}$ is a monoidal $\infty$-category and $A\in \Alg(\cat{C})$ is an associative algebra in $\cat{C}$, then the category ${}_A\cat{BMod}_A(\cat{C})$ is the underlying category of a nonsymmetric $\infty$-operad (see \cite[Theorem 3.3.3.9, Theorem 4.4.1.28]{lurie2014higher})
$$\begin{tikzcd}
\cat{Mod}^{\mm{Assoc}}_A(\cat{C})^{\otimes}\arrow[r] & \mm{Assoc}^{\otimes}
\end{tikzcd}$$
with the property that $\Alg(\cat{Mod}^{\mm{Assoc}}_A\big(\cat{C})\big)\simeq \Alg(\cat{C})_{A/}$.
\end{notation}
Now suppose that $\lambda\colon \cat{M}\rt \cat{C}\times\cat{D}$ is a pairing of monoidal $\infty$-categories and let $A_{\mb{1}}\in \Alg(\cat{M})$ be an algebra with image $(A, \mb{1})$ in $\Alg(\cat{C})\times\Alg(\cat{D})$. We obtain a pairing of nonsymmetric $\infty$-operads
$$\begin{tikzcd}
\lambda^{\otimes}_A\colon \cat{Mod}^{\mm{Assoc}}_{A_\mb{1}}(\cat{M})^{\otimes}\arrow[r] & \cat{Mod}^{\mm{Assoc}}_{A}(\cat{C})^{\otimes}\times_{\mm{Assoc}^{\otimes}} \cat{D}^{\otimes}
\end{tikzcd}$$
where we identify $\cat{Mod}^{\mm{Assoc}}_{\mb{1}}(\cat{D})^{\otimes}\simeq \cat{D}^{\otimes}$ \cite[Proposition 3.4.2.1]{lurie2014higher}. Since $A$ can be considered as an associative algebra in $\cat{Mod}^{\mm{Assoc}}_{A}(\cat{C})^{\otimes}$, we can consider the nonsymmetric $\infty$-operad $\cat{E}_A^{\otimes}$ defined as the fibre product
$$\begin{tikzcd}
\cat{E}_A^{\otimes}\arrow[r] \arrow[d] & \cat{Mod}^{\mm{Assoc}}_{A_\mb{1}}(\cat{M})^{\otimes}\arrow[d]\\
\mm{Assoc}^{\otimes}\arrow[r, "A"] & \cat{Mod}^{\mm{Assoc}}_{A}(\cat{C})^{\otimes}.
\end{tikzcd}$$
\begin{lemma}
In the above situation, suppose that $\lambda\colon \cat{M}\rt \cat{C}\times\cat{D}$ is left representable and that $\cat{D}$ admits totalisations. Then $\cat{E}_A^{\otimes}$ is a lax monoidal $\infty$-category, i.e.\ the map $\cat{E}_A^{\otimes}\rt \mm{Assoc}^{\otimes}$ is a locally cocartesian fibration.
\end{lemma}
\begin{proof}
Consider the map of correspondences
$$\begin{tikzcd}
\cat{C}\arrow[d, hookrightarrow] & \cat{M}\arrow[l]\arrow[r]\arrow[d, hookrightarrow] & \cat{D}\arrow[d, hookrightarrow]\\
\cat{P}(\cat{C}) & \cat{P}(\cat{M})\arrow[l]\arrow[r] & \cat{P}(\cat{D}).
\end{tikzcd}$$ 
Endowing all presheaf categories with the monoidal structure given by Day convolution \cite[Corollary 4.8.1.12]{lurie2014higher}, this gives a diagram of monoidal $\infty$-categories and monoidal functors.

Considering $A$ and $A_{\mb{1}}$ as associative algebras in $\cat{P}(\cat{C})$ and $\cat{P}(\cat{M})$ under the Yoneda embedding, we can form a similar nonsymmetric $\infty$-operad
$$
\widehat{\cat{E}}_A^{\otimes}=\mm{Assoc}^{\otimes}\times_{\cat{Mod}^{\mm{Assoc}}_{A}(\cat{P}(\cat{C}))^{\otimes}} \cat{Mod}^{\mm{Assoc}}_{A_\mb{1}}\big(\cat{P}(\cat{M})\big)^{\otimes}.
$$
Since the Day convolution product preserves colimits in each argument, this is a fibre product of monoidal $\infty$-categories and monoidal functors between them \cite[Theorem 3.4.4.2]{lurie2014higher}.

We now obtain a diagram
$$\begin{tikzcd}
\cat{E}_A^{\otimes}\arrow[rr]\arrow[rd, "p"{swap}] & & \widehat{\cat{E}}_A^{\otimes}\arrow[ld, "q"]\\
& \mm{Assoc}^{\otimes}
\end{tikzcd}$$
where the top map is fully faithful and $q$ is a cocartesian fibration. We will use this to prove that $p$ is a locally cocartesian fibration.

To this end, note that the map of $\infty$-operads $\cat{E}_A^{\otimes}\rt \mm{Assoc}^{\otimes}$ is a locally cocartesian fibration if each active morphism $\alpha\colon \langle n\rangle \rt \langle 1\rangle$ in $\mm{Assoc}^{\otimes}$ (there are $n!$ of these) admits locally cocartesian lifts. Let us therefore pick an active morphism $\alpha\colon \langle n\rangle \rt \langle 1\rangle$ in $\mm{Assoc}^{\otimes}$ and $n$ objects $M_1, \dots, M_n$ in 
$$
\big(\cat{E}_A^{\otimes}\big)_{\langle 1\rangle}\simeq \{A\}\times_{{}_{A}\cat{BMod}_{A}(\cat{C})} {}_{A_{\mb{1}}}\cat{BMod}_{A_{\mb{1}}}(\cat{M}).
$$
There exists a $q$-cocartesian lift of $\alpha$ in $\widehat{\cat{E}}_A^{\otimes}$ of the form (up to a permutation of the $M_i$)
$$\begin{tikzcd}
\tilde{\alpha}\colon \big(M_1, \dots, M_n\big)\arrow[r] & M_1\otimes_{\widehat{\cat{E}}_A} \dots \otimes_{\widehat{\cat{E}}_A} M_n
\end{tikzcd}$$
where the target denotes the tensor product of $\widehat{\cat{E}}_A$. To see that $\alpha$ admits a locally $p$-cocartesian lift, it suffices to verify that there is an initial object in the full subcategory $(\cat{E}_A^{\otimes})_{\langle 1\rangle}\subseteq (\widehat{\cat{E}}_A^{\otimes})_{\langle 1\rangle}$ that receives a map from the tensor product $M_1\otimes_{\widehat{\cat{E}}_A} \dots \otimes_{\widehat{\cat{E}}_A} M_n$.

To prove this, let us consider the image of $M_1\otimes_{\widehat{\cat{E}}_A} \dots \otimes_{\widehat{\cat{E}}_A} M_n$ in $\cat{BMod}_{A_{\mb{1}}}(\cat{P}(\cat{M}))$, which is given by the $n$-fold relative tensor product
$$
M_1\widehat{\otimes}_{A_{\mb{1}}}\dots \widehat{\otimes}_{A_\mb{1}} M_n.
$$
Here we abuse notation by identifying each $M_i$ with its image in $\cat{BMod}_{A_{\mb{1}}}(\cat{M})\subseteq \cat{BMod}_{A_{\mb{1}}}(\cat{P}(\cat{M}))$, and we write $\widehat{\otimes}_{A_{\mb{1}}}$ to indicate that the relative tensor product is computed at the presheaf level. Explicitly, this relative tensor product is given by the geometric realisation in $\cat{BMod}_{A_{\mb{1}}}(\cat{P}(\cat{M}))$ of the simplicial diagram
$$\begin{tikzcd}[row sep=0.1pc]
Q_\bullet\colon \Del^{\op}\arrow[r] &  \cat{BMod}_{A_{\mb{1}}}(\cat{M});\\
{[k]}\arrow[r, mapsto] & M_1\otimes A_{\mb{1}}^{\otimes k} \otimes M_2\otimes A_{\mb{1}}^{\otimes k}\otimes \dots \otimes A_{\mb{1}}^{\otimes k}\otimes M_n.
\end{tikzcd}$$
We will show that this diagram also has a colimit in $\cat{BMod}_{A_{\mb{1}}}(\cat{M})$. First, note that the image of $Q_\bullet$ in $\cat{D}\simeq \cat{BMod}_{\mb{1}}(\cat{D})$ is essentially constant on a certain object $D\in \cat{D}$. It therefore lifts to a diagram of the form
$$\begin{tikzcd}
\Del^{\op}\arrow[r, "Q_\bullet"] & \cat{F}_D=\cat{BMod}{}_{A_\mb{1}}(\cat{M})\times_{\cat{D}} \{D\}\arrow[r, "g"]\arrow[d] & \cat{M}\times_{\cat{D}} \{D\}\arrow[d]\\
& \cat{BMod}_A(\cat{C})\arrow[r, "g'"] & \cat{C}
\end{tikzcd}$$
Next, note that the image of $Q_\bullet$ in $\cat{BMod}_{A}(\cat{C})$ is the diagram $A\otimes A^{\otimes k}\otimes \dots\otimes A^{\otimes k} \otimes A$. This diagram has $A$ as a colimit and the embedding $\cat{BMod}_{A}(\cat{C})\hookrightarrow \cat{BMod}_{A}(\cat{P}(\cat{C}))$ preserves the colimit: this follows from the fact that the underlying diagram in $\cat{C}$ is split. 

We can now repeat the argument used in the proof of Lemma \ref{lem:bimodule pairing}. The image $g(Q_\bullet)$ lifts a split object in $\cat{C}$ and is therefore split, so that $Q_\bullet$ admits a colimit in $\cat{F}_D$ by Remark \ref{rem:bimodules fiber}. Because $\cat{BMod}_{A_{\mb{1}}}(\cat{M})\rt \cat{D}$ is a cartesian fibration, the fibre inclusion $\cat{F}_D\hookrightarrow \cat{BMod}_{A_{\mb{1}}}(\cat{M})$ preserves colimits. We conclude that $Q_\bullet$ admits a colimit in $\cat{BMod}{}_{A_\mb{1}}(\cat{M})$, which is furthermore preserved by the functors
$$\begin{tikzcd}
\cat{BMod}{}_{A_\mb{1}}(\cat{M})\arrow[r] & \cat{BMod}_A(\cat{C})\arrow[r] & \cat{BMod}_A(\cat{P}(\cat{C})).
\end{tikzcd}$$
All in all, we therefore obtain a canonical map in $\cat{BMod}_{A_{\mb{1}}}(\cat{P}(\cat{M}))$ of the form
$$\begin{tikzcd}
\beta\colon M_1\widehat{\otimes}_{A_{\mb{1}}}\dots \widehat{\otimes}_{A_\mb{1}} M_n=\big|Q_\bullet\big|_{\cat{P}(\cat{M})}\arrow[r] & \big|Q_\bullet\big|_{\cat{M}}
\end{tikzcd}$$
from the geometric realisation computed in $\cat{BMod}_{A_{\mb{1}}}(\cat{P}(\cat{M}))$ to the geometric realisation computed in $\cat{BMod}_{A_{\mb{1}}}(\cat{M})$. By definition, the map $\beta$ exhibits $|Q_\bullet|_{\cat{M}}$ as the initial object in $\cat{M}$ which receives a map from $M_1\widehat{\otimes}_{A_{\mb{1}}}\dots \widehat{\otimes}_{A_\mb{1}} M_n$. 

Since the image of $\beta$ in $\cat{BMod}_A(\cat{P}(\cat{C}))$ is (equivalent to) the identity on $A$, the map $\beta$ determines an arrow in the fibre $(\widehat{\cat{E}}_A^{\otimes})_{\langle 1\rangle}$  of the form
$$
\beta_A \colon M_1\otimes_{\widehat{\cat{E}}_A} \dots \otimes_{\widehat{\cat{E}}_A} M_n\rt X
$$ 
Now recall that the fully faithful inclusion $\big(\cat{E}_A^{\otimes}\big)_{\langle 1\rangle}\subseteq \big(\widehat{\cat{E}}_A^{\otimes}\big)_{\langle 1\rangle}$ is the base change of the inclusion $\cat{BMod}_{A_{\mb{1}}}(\cat{M})\subseteq \cat{BMod}_{A_{\mb{1}}}(\cat{P}(\cat{M}))$. Consequently, the map $\beta_A$ exhibits $X$ as the initial object in $(\cat{E}_A^{\otimes})_{\langle 1\rangle}$ receiving a map from $M_1\otimes_{\widehat{\cat{E}}_A} \dots \otimes_{\widehat{\cat{E}}_A} M_n$. The composite $\beta_A\circ \tilde{\alpha}\colon (M_1, \dots, M_n)\rt X$ then provides the desired locally $p$-cocartesian lift of $\alpha$, as desired.
\end{proof}

\begin{corollary}[{\cite[Proposition  5.2.2.27]{lurie2014higher}}]\label{cor:luriebar}
Let $\lambda\colon \cat{M}\rt \cat{C}\times \cat{D}$ be a pairing of monoidal $\infty$-categories such that the following conditions hold:
\begin{enumerate}
\item The unit $\mb{1}\in \cat{D}$ is an initial object and the functor $\cat{M}\times_{\cat{D}}\{\mb{1}\}\rt \cat{C}$ is an equivalence. 
\item The pairing $\lambda$ is left representable. 
\item The $\infty$-category $\cat{D}$ admits totalisations of cosimplicial objects.
\end{enumerate}
Then the induced pairing $\Alg(\lambda)\colon \Alg(\cat{M})\rt \Alg(\cat{C})\times \Alg(\cat{D})$ is left representable.
\end{corollary}
\begin{proof}
Given $A\in \Alg(\cat{C})$, we have to show that the fibre $\{A\}\times_{\Alg(\cat{C})} \Alg(\cat{M})$ admits a terminal object. To this end, let us start by noting that there exists a unique lift $A_{\mb{1}}\in \{A\}\times_{\Alg(\cat{C})} \Alg(\cat{M})\times_{\Alg(\cat{D})} \{\mb{1}\}$; this follows from the monoidal equivalence $\cat{M}\times_{\cat{D}}\{\mb{1}\}\simeq \cat{C}$. Since the functor
$$\begin{tikzcd}
\{A\}\times_{\Alg(\cat{C})_{A/}} \Alg(\cat{M})_{A_{\mb{1}}/}\arrow[r] & \{A\}\times_{\Alg(\cat{C})} \Alg(\cat{M})
\end{tikzcd}$$
preserves terminal objects \cite[Proposition 5.2.2.30]{lurie2014higher}, it suffices to verify that the domain has a terminal object. But now we can identify
$$
\{A\}\times_{\Alg(\cat{C})_{A/}} \Alg(\cat{M})_{A_{\mb{1}}/}=\{A\}\times_{\Alg(\Mod^{\mm{Assoc}}_A(\cat{C}))} \Alg\big(\Mod^{\mm{Assoc}}_{A_\mb{1}}(\cat{M})\big)\simeq \Alg(\cat{E}^\otimes_A).
$$
Because $\cat{E}^{\otimes}_A$ is a lax monoidal $\infty$-category with a terminal object, its category of associative algebras admits a terminal object as well by \cite[Proposition 3.2.2.1]{lurie2014higher}: indeed, if $\cat{A}$ is an $\infty$-category with a lax monoidal structure, then a terminal object in $\cat{A}$ also determines a $p$-terminal object for $p\colon \cat{A}^{\otimes}\rt \mm{Assoc}^{\otimes}$.
\end{proof}
The technical heart of \Cref{thm:bar for modules} is the following analogue of Corollary \ref{cor:luriebar}:
\begin{proposition}\label{prop:representability for modules}
Let $\lambda\colon (\Ml, \Mm)\rt (\Cl, \Cm)\times (\Dl, \Dm)$ be a left module object in the $\infty$-category $\cat{Pair}$ of pairings, such that the following conditions hold:
\begin{enumerate}
\item The unit $\mb{1}\in \Dl$ is an initial object and the functor $\Ml\times_{\Dl}\{\mb{1}\}\rt \Cl$ is an equivalence. 
\item The pairings $\Ml\rt \Cl\times\Dl$ and $\Mm\rt \Cm\times \Dm$ are both left representable. 
\item The $\infty$-categories $\Dl$ and $\Dm$ both admit totalisations of cosimplicial objects.
\end{enumerate}
In this case, the pairing $\LMod(\Mm)\rt \LMod(\Cm)\times\LMod(\Dm)$ is left representable and the forgetful functor $\LMod(\Mm)\rt \Alg(\Ml)$ preserves left representable objects. 
\end{proposition}
\begin{proof}
The proof of this result follows the lines of the proof of \cite[Proposition 5.2.2.27]{lurie2014higher}. Let $(A, M)\in \LMod(\Cm)$ be a tuple of an associative algebra $A$ in $\Cl$ and a left $A$-module $M$ in $\Cm$. Consider the cartesian fibration taking underlying algebras
$$\begin{tikzcd}
\pi\colon \LMod(\Mm)\times_{\LMod(\Cm)} \{(A, M)\} \arrow[r] & \Alg(\Ml)\times_{\Alg(\Cl)} \{A\}=\Alg(\Ml)_{A}
\end{tikzcd}$$
The target admits a terminal (i.e.\ left universal) object $A_L$ by \cite[Proposition 5.2.2.27]{lurie2014higher} or Corollary \ref{cor:luriebar}. We have to prove that the domain admits a terminal object of the form $(A_L, M_L)$. We will do this by proving that each fibre of $\pi$ admits a terminal object and that the change-of-fibre functors preserve these terminal objects. Then $\pi$ admits a fully faithful right adjoint; its  value on $A_L$ is the desired $(A_L, M_L)$.

To see this, note that $\Alg(\Ml)_A \rt \Alg(\Dl)$ is a right fibration represented by $\Barr(A)$ \cite[Proposition 5.2.2.27]{lurie2014higher}. In particular, $\Alg(\Ml)_A$ admits an initial object $A_{\mb{1}}$; its image in $\Alg(\Dl)$ is the initial object $\mb{1}$. Let $A'\in \Alg(\Ml)_A$ be any other lift of $A$ and let $f\colon A_{\mb{1}}\rt A'$ denote the unique map in $\Alg(\Ml)_A$. By Lemma \ref{lem:adjoint to restriction} below, restriction of modules along $f$ defines a right adjoint functor
$$\begin{tikzcd}
f^*\colon \LMod_{A'}(\Mm)\times_{\LMod_A(\Cm)} \{M\}\arrow[r] & \LMod_{A_{\mb{1}}}(\Mm)\times_{\LMod_A(\Cm)} \{M\}.
\end{tikzcd}$$
In particular, this implies that the functor $f^*$ preserves and detects terminal objects. The codomain $\LMod_{A_{\mb{1}}}(\Mm)\times_{\LMod_A(\Cm)} \{M\}$ admits a terminal object by a similar, but easier argument as in Lemma \ref{lem:bimodule pairing} or \cite[Lemma 5.2.2.40]{lurie2014higher}. Consequently, each fibre of $\pi$ admits a terminal object, which is preserved by all change-of-fibre functors. 
\end{proof}
\begin{lemma}\label{lem:adjoint to restriction}
Consider the setting of Proposition \ref{prop:representability for modules} and let $A\in \Alg(\Cl)$ be an associative algebra. Let $f\colon A_{\mb{1}}\rt A'$ be a map in the fibre $\Alg(\Ml)\times_{\Alg(\Cl)}\{A\}$ with domain given by the initial object. Then the restriction functor $f^*$ between $\infty$-categories of modules is the right adjoint in a relative adjunction
$$\begin{tikzcd}[column sep=1.5pc]
\LMod_{A_{\mb{1}}}(\Mm)\arrow[rr, yshift=1ex, "f_!"]\arrow[rd] & & \LMod_{A'}(\Mm)\arrow[ld]\arrow[ll, "f^*", yshift=-1ex]\\
& \LMod_A(\Cm).
\end{tikzcd}$$
\end{lemma}
\begin{proof}
Let $N\in \LMod_{A_{\mb{1}}}(\Mm)$ and consider the bar construction $\Barr_\bullet(A', A_{\mb{1}}, N)$. If the geometric realisation of $\Barr_\bullet(A', A_{\mb{1}}, N)$ in $\LMod_{A'}(\Mm)$ exists, then it computes the value of the putative left adjoint $f_!$ on $N$. It therefore suffices to verify that each $\Barr_\bullet(A', A_{\mb{1}}, N)$ admits a geometric realisation in $\LMod_{A'}(\Mm)$. Since $\LMod_{A'}(\Mm)$ is monadic over $\Mm$, it suffices to verify that the underlying simplicial diagram in $\Mm$ is split.

To see this, note that the image of $\Barr_\bullet(A', A_{\mb{1}}, N)$ under the projection $p\colon \Mm\rt \Dm$ is given by $\Barr_\bullet(p(A'), \mb{1}, p(N))$. This simplicial diagram is constant on $X=p(A')\otimes p(N)$, so that we can think of $\Barr_\bullet(A', A_{\mb{1}}, N)$ as a simplicial diagram in the fibre $\Mm\times_{\Dm} X$. Now consider the right fibration $q\colon \Mm\times_{\Dm} X\rt \Cm$. By \cite[Corollary 4.7.2.11]{lurie2014higher}, it suffices to verify that image of $\Barr_\bullet(A', A_{\mb{1}}, N)$ in $\Cm$ is a split simplicial diagram. This image is simply the split simplicial diagram $\Barr_\bullet(A, A, q(N))$. We conclude that $f^*$ indeed admits a left adjoint $f_!$. Furthermore, the image of $f_!(N)$ in $\LMod_A(\Cm)$ agrees with $|\Barr_\bullet(A, A, q(N))|\simeq q(N)$, so that $f_!$ and $f^*$ form a relative adjunction over $\LMod_A(\Cm)$.
\end{proof}
\begin{proof}[Proof of \Cref{thm:bar for modules}]
The commuting square of adjunctions is an immediate consequence of Proposition \ref{prop:representability for modules}, applied to the pairing $\Tw(\Cl)\rt \Cl\times\Cl^{\op}$ and $\Tw(\Cm)\rt \Cm\times\Cm^{\op}$ (and its opposite for the cobar functors). 

For assertion (1), an inspection of the proof of Proposition \ref{prop:representability for modules} and \cite[Lemma 5.2.2.40]{lurie2014higher} shows that $\Barr\colon \LMod(\Cm)\rt \LComod(\Cm)$ sends $(A, M)$ to a tuple of a coalgebra and a comodule, with underlying objects given by the inductions $\epsilon_!(\mb{1})$ and $\epsilon_!(M)$ along the augmentation $\epsilon\colon A\rt \mb{1}$. Consider a coCartesian arrow in $\LMod(\Cm)$ of the form $(A, M)\rt (B, f_!(M))$, where $f\colon A\rt B$ is a map of algebras. Denoting the augmentation maps of $A$ and $B$ by $\epsilon_A$ and $\epsilon_B$ respectively, the image of this arrow under the bar construction is given on underlying objects by the natural map 
$$\begin{tikzcd}
\big(\epsilon_{A!}(\mb{1}), \epsilon_{A!}(M)\big)\arrow[r] &  \big(\epsilon_{B!}(\mb{1}), \epsilon_{B!}(f_!(M))\big).
\end{tikzcd}$$
This map is coCartesian as soon as the natural map $\epsilon_{A!}(M)\rt \epsilon_{B!}(f_!(M))$ is an equivalence, which follows from transitivity of extension of scalars and the fact that $\epsilon_A\simeq\epsilon_B\circ f$.

For assertion (2) about right $\Cr$-linearity of the bar construction, it suffices to verify that the right action $\LMod(\Tw(\Cm))\times \Tw(\Cr)\rt \LMod(\Tw(\Cm))$ preserves left universal arrows (see Notation \ref{not:pairing}). Since $\Tw(\Cr)$ classifies the identity functor on $\Cr$, this comes down to the assertion that for any object $X\in \Cr$, the natural map $\Barr(A, M\otimes X)\rt \Barr(A, M)\otimes X$ is an equivalence. This map can be identified with the canonical map $\big(\epsilon_!(\mb{1}), \epsilon_!(M\otimes X)\big)\rt \big(\epsilon_!(\mb{1}), \epsilon_!(M)\otimes X\big)$. The map 
$$\begin{tikzcd}
\epsilon_!(M\otimes X)=\big|\Barr_\bullet(\mb{1}, A, M)\otimes X\big|\arrow[r] & \big|\Barr_\bullet(\mb{1}, A, M)\big|\otimes X =\epsilon_!(M)\otimes X
\end{tikzcd}$$
is now an equivalence because $(-)\otimes X$ preserves geometric realisations.
\end{proof}

\subsubsection*{Refined Koszul Duality}
The \emph{Koszul dual} of an augmented PD $\infty$-operad $\PP$ now arises from the bar construction $\Barr(\PP)$ by linear duality.
\begin{notation}[Linear dual symmetric sequences]\label{not:linear duality}
Let $R$ be a coherent $\mathbb{E}_\infty$-ring spectrum and recall that the $\infty$-category $\sSeq^\vee_R$ comes equipped with the \emph{levelwise} tensor product $\otimes_{\mm{lev}}$; its unit is the $\mathbb{E}_\infty$-operad. If $X$ is a pro-coherent symmetric sequence, we will write $X^\vee$ for its dual with respect to the levelwise tensor product and refer to it as the \emph{linear dual} of $X$.
\end{notation}
\begin{remark}[Refined linear duality]\label{rem:refined linear duality}
The pro-coherent linear duality described above refines the usual operation of taking $R$-linear dual symmetric sequences, in the sense that there is a commuting diagram
$$\begin{tikzcd}
\sSeq_R^{\vee, \op}\arrow[r, "(-)^\vee"] & \sSeq_R^\vee\arrow[d, "\upupsilon"]\\
\sSeq^{\op}_R\arrow[r, "(-)^\vee"]\arrow[u, "\upiota"] & \sSeq_R.
\end{tikzcd}$$
The bottom functor is the usual functor taking the levelwise linear dual $X^\vee(r)=\Hom_R(X(r), R)$. If $X$ is an ordinary symmetric sequence, then its dual in $\sSeq^\vee_R$ crucially need \emph{not} arise from a symmetric sequence, i.e.\ need not be contained in the essential image of $\upiota$.
\end{remark}
Informally, $X^\vee$ is the pro-coherent symmetric sequence given in each arity $r$ by the continous $R$-linear dual of $X(r)$. This is substantiated by the following observation:
\begin{proposition}\label{prop:dual symmetric sequence}
Let $R$ be a coherent $\mathbb{E}_\infty$-ring spectrum. Then $(-)^\vee\colon \sSeq^\vee_R\rt \sSeq^{\vee, \op}_R$ is the right-left extension of the functor $R[\Sigma]\rt R[\Sigma]^{\op}$ sending a finite type free symmetric sequence $X$ to the $R$-linear dual symmetric sequence $X^\vee(r)=\Hom_R(X(r), R)$. Furthermore, it restricts to an equivalence
$$\begin{tikzcd}
(-)^\vee\colon \APerf^\vee_{R[\Sigma]}\arrow[r, "\simeq"] & \APerf_{R[\Sigma]}^{\op}.
\end{tikzcd}$$
\end{proposition}
\begin{proof}
When $X$ is finitely generated free, $X^\vee(r)=\Hom_R(X(r), R)$ indeed defines an object in $R[\Sigma]$. Furthermore, there are canonical maps in the $\infty$-category of symmetric sequences $\mathbb{E}_\infty\rt X\otimes_{\lev} X^\vee$ and $X^\vee\otimes_{\lev} X\rt \mathbb{E}_\infty$ exhibiting $X^\vee$ as the dual of $X$. Consequently, they remain dual in pro-coherent symmetric sequences as well. The result now follows from (the proof of) Proposition \ref{prop:dual is indeed dual}.
\end{proof}
\begin{proposition}\label{prop:dual symmetric sequence lax monoidal}
Let $R$ be a coherent $\mathbb{E}_\infty$-ring spectrum. Then $(-)^\vee\colon \sSeq_R^{\vee, \op}\rt \sSeq^{\vee}_R$ is lax monoidal with respect to the extended composition product $\circ$. Furthermore, it restricts to a (strong) monoidal equivalence $\APerf^\vee_{R[\Sigma]}\simeq \APerf_{R[\Sigma]}^{\op}$.
\end{proposition}
In particular, (continuous) linear duality sends PD $\infty$-cooperads to PD $\infty$-operads and restricts to an equivalence between the $\infty$-categories of almost perfect PD $\infty$-cooperads and dually almost perfect PD $\infty$-operads.
\begin{proof}
Consider the functor $F\colon \sSeq_R^{\vee, \op}\rt \Fun(\sSeq_R^{\vee, \op}, \An_\mm{large})$ sending each pro-coherent symmetric sequence $X$ to the (large) presheaf $\Map(X\otimes_{\lev} -, \mathbb{E}_\infty)$. This takes values in the essential image of the Yoneda embedding, and the corresponding functor precisely sends $X\mapsto X^\vee$. We endow the (large) presheaf category $\Fun\big(\sSeq_R^{\vee, \op}, \An_\mm{large}\big)$ with the Day convolution product with respect to $\circ$. Since the Yoneda embedding is a fully faithful monoidal functor, it suffices to endow the functor $F$ with a lax monoidal structure. By the universal property of Day convolution \cite[Section 2.2.6]{lurie2014higher}, such a lax monoidal structure is equivalent to a lax monoidal structure on the functor adjoint to $F$
$$\begin{tikzcd}[column sep=4pc]
\sSeq_R^{\vee, \op}\times \sSeq_R^{\vee, \op}\arrow[r, "\otimes_{\lev}"] & \sSeq_R^{\vee, \op}\arrow[r, "{\Map(-, \mathbb{E}_\infty)}"] & \An.
\end{tikzcd}$$
The first functor is lax monoidal by Proposition \ref{prop:pro-coh composition} and the second functor is lax monoidal since $\mathbb{E}_\infty$ is an algebra with respect to $\circ$ (and the Yoneda embedding $\cat{C}\rt \mm{Fun}(\cat{C}^{\op}, \An)$ is monoidal for the Day convolution product).

For the final assertion, we have to verify that the natural map $\mu\colon X^\vee\circ Y^\vee\rt (X\circ Y)^\vee$ is an equivalence when $X$ and $Y$ are dually almost perfect. Since $\circ$ and $(-)^\vee$ preserve sifted colimits and finite totalisations, we may assume that $X, Y\in \APerf^\vee_{R[\Sigma], \weirdleq 0}$. Using furthermore that both functors preserve totalisations of diagrams in $R[\Sigma]$, we can reduce to the case where $X$ and $Y$ are contained in $R[\Sigma]$. In this case the result follows by inspection.
\end{proof}
\begin{definition}[Koszul dual PD operad]\label{def:koszul dual PD operad}
Let $R$ be a coherent $\mathbb{E}_\infty$-ring. If $\PP$ is an augmented PD $\infty$-operad, we define its \emph{Koszul dual} PD $\infty$-operad $\KD^{\pd}(\PP)=\Barr(\PP)^\vee$ to be the linear dual of the bar construction.
\end{definition}
\begin{theorem}[Refined Koszul duality for operads]\label{thm:koszul duality for PD operads}
Let $R$ be a coherent $\mathbb{E}_\infty$-ring spectrum. Then there is a commuting diagram of $\infty$-categories
$$\begin{tikzcd}
\Op^{\pd, \aug}_R\arrow[r, "\KD^{\pd}"] & \Op^{\pd, \aug, \op}_R\arrow[d, "\upupsilon"]\\
\Op^{\aug}_R\arrow[r, "\KD"]\arrow[u, "\upiota"] & \Op^{\aug, \op}_R
\end{tikzcd}$$
where the bottom functor sends an augmented $\infty$-operad to its classical Koszul dual $\infty$-operad, given by the Spanier--Whitehead dual of its bar construction.
\end{theorem}
\begin{proof}
The functor $\upiota$ is monoidal with respect to the composition product and preserves geometric realisations. Consequently, it commutes with the bar construction. The result then follows from the fact that linear duality in pro-coherent symmetric sequences provides a lax monoidal refinement of Spanier--Whitehead duality of ordinary symmetric sequences (Remark \ref{rem:refined linear duality}).
\end{proof}
Since linear duality preserves pro-coherent symmetric sequences concentrated in arity $0$, it furthermore sends coalgebras over $\infty$-cooperads to algebras over $\infty$-operads:
\begin{definition}[Koszul dual algebra]
Let $R$ be a coherent $\mathbb{E}_\infty$-ring and $\PP$ an augmented $R$-linear PD $\infty$-operad, with Koszul dual $\KD^{\pd}(\PP)$. If $A$ is a $\PP$-algebra, we define its \emph{Koszul dual} $\KD(\PP)$-algebra  to be the linear dual of its bar construction $\KD^{\pd}(A)=\Barr_{\PP}(A)^\vee$ (Corollary \ref{cor:bar construction for algebras}).
\end{definition}
\begin{theorem}\label{thm:koszul duality for algebras}
Let $R$ be a coherent $\mathbb{E}_\infty$-ring spectrum and $\PP$ an augmented $\infty$-operad over $R$. Then there is a commuting diagram of $\infty$-categories
$$\begin{tikzcd}
\Alg_{\PP}(\QC^\vee_R)\arrow[r, "\KD^{\pd}"] & \Alg_{\KD^{\pd}(\PP)}(\QC^\vee_R)^{\op}\arrow[d, "\upupsilon"]\\
\Alg_{\PP}(\Mod_R)\arrow[r, "\KD"]\arrow[u, "\upiota"] & \Alg_{\KD(\PP)}(\Mod_R)^{\op}
\end{tikzcd}$$
where the bottom functor sends a $\PP$-algebra $A$ its classical Koszul dual algebra, given by the Spanier--Whitehead dual of its bar construction.
\end{theorem}
\begin{proof}
The proof of \Cref{thm:koszul duality for PD operads} carries over, using instead that $\upiota$ commutes with the bar construction of \Cref{thm:bar for modules}.
\end{proof}
We will now illustrate how \Cref{thm:koszul duality for PD operads} and \Cref{thm:koszul duality for algebras} concretely refine the usual Koszul duality for $\infty$-operads and their algebras.
\begin{definition}[Almost finitely presented $\infty$-operads]
Let $R$ be a connective $\mathbb{E}_\infty$-ring spectrum. An augmented $\infty$-operad $\PP$ is said to be \emph{connective} if each $\PP(r)$ is a connective spectrum. A connective augmented operad $\PP$ is  \emph{almost of finite presentation} if it defines an almost compact object in the compactly generated $\infty$-category $\Op^{\aug}_{R, \geq 0}$ of connective augmented operads, in the sense of \cite[Definition 7.2.4.8]{lurie2014higher}: this means $\tau_{\leq m}\PP$ is a compact object in the $\infty$-category $\Op^{\aug}_{R, \geq 0, \leq m}$ of augmented $\infty$-operads that are connective and $m$-coconnective, for each $m\geq 0$.
\end{definition}
\begin{proposition}\label{prop:PD koszul duality in aft case}
Let $R$ be a coherent $\mathbb{E}_\infty$-ring spectrum and let $\PP$ be a connective augmented $\infty$-operad over $R$ which is almost of finite presentation. Then the Koszul dual PD $\infty$-operad $\KD^{\pd}(\PP)$ is dually almost perfect. The induced monad $\KD^{\pd}(\PP)\colon \QC^\vee_R\rt \QC^\vee_R$ preserves sifted colimits and dually almost perfect objects, and the resulting monad on $\APerf^\vee_R$ can be identified with
$$
\KD^{\pd}(\PP)(V)\simeq \bigoplus_r \big(\KD(\PP)^{\pd}(r)\otimes V^{\otimes r}\big)^{h\Sigma_r} \simeq \big(\cot_{\PP}\circ \triv_{\PP}(V^\vee)\big)^\vee.
$$
\end{proposition}
Notice that the above differs from the free algebra over the classical Koszul dual $\KD(\PP)$, even for finitely generated free $R$-modules.
\begin{proof}
Note that the bar construction restricts to a functor 
$\Barr\colon \Op^{\aug}_{R, \geq 0}\rt \sSeq_{R, \geq 0}$ from connective augmented operads to connective symmetric sequences (indeed, this is simply the bar construction for augmented algebras in the monoidal $\infty$-category $\sSeq_{R, \geq 0}$). This functor preserves colimits and sends a free augmented $\infty$-operad $\mm{Free}(X)$ to $X$, so that it preserves almost compact objects for formal reasons. 

It follows that the bar construction of an almost finitely presented connective augmented $\infty$-operad $\PP$ is almost perfect as a symmetric sequence, so that its Koszul dual is dually almost perfect. The dually almost perfect symmetric sequences are closed under the composition product, so that the free $\KD^{\pd}(\PP)$-algebra functor preserves dually almost perfect objects. Since linear duality is an equivalence on dually almost perfect symmetric sequences (Proposition \ref{prop:dual symmetric sequence lax monoidal}), there is an equivalence of monads $\KD^{\pd}(\PP)\simeq (-)^\vee\circ \Barr(\PP)\circ (-)^\vee$. The formulas for the monad $\KD^{\pd}(\PP)$ then follow from Proposition \ref{prop:formula for pro-coh comp prod} and Corollary \ref{cor:bar construction for algebras}.
\end{proof}
One can verify that the nonunital $\mathbb{E}_\infty$-operad is almost finitely presented. In particular, its bar construction is the symmetric sequence $$\Barr(\mathbb{E}_{\infty,R}^\mm{nu})(r)\simeq R\wedge \Sigma|\Pi_r|^\diamond$$ of reduced-unreduced suspensions of the partition complex; this is indeed an almost perfect symmetric sequence (which is all we need). We will write $\Hom(\Sigma|\Pi_r|^\diamond, R)\in \QC^\vee_{R[\Sigma_r]}$ for its pro-coherent $R$-linear dual.
\begin{definition}[The spectral Lie PD operad]\label{def:spectral partition lie}
The \emph{spectral partition Lie PD $\infty$-operad} is the PD Koszul dual $\mm{Lie}^\pi_{R, \mathbb{E}_\infty}=\KD^{\pd}(\mathbb{E}^\mm{nu}_\infty)$. 
\end{definition} 
\begin{corollary}
The monad associated to the spectral partition Lie PD $\infty$-operad agrees (over a field $k$) with the spectral partition Lie monad from \cite[Definition 5.32]{brantner2019deformation} and is given on dually almost perfect objects by 
$$
\mm{Lie}^\pi_{R, \mathbb{E}_\infty}(V)=\bigoplus_r \big(\Hom(\Sigma|\Pi_r|^\diamond, R)\otimes V^{\otimes r}\big){}^{h\Sigma_r}.
$$
\end{corollary}

\subsection{Derived operads and derived PD operads}
\label{sec:derived operads}
In this section, we will describe a derived refinement of the notion of $\infty$-operad and PD $\infty$-operad over a coherent simplicial commutative ring, which also accounts (in a rather strict way) for the genuine equivariant homotopy theory of the symmetric group actions. We will first discuss the derived version of classical $\infty$-operads and then turn to the pro-coherent setting.

\subsubsection*{Derived operads}
Recall that the $\infty$-category of symmetric sequences over $R$ is generated by free $\Sigma_r$-modules, for various $r$. We will now introduce an $\infty$-category of derived symmetric sequences over $R$ that will be generated by $\Sigma_r$-orbits.
\begin{definition}[Derived symmetric sequences]\label{def:derived sym seq}
Let $R$ be a simplicial commutative ring. We define the $\infty$-category of \emph{derived symmetric sequences} to be
$$
\sSeq^{\gen}_{\ulR} \simeq  \prod_{r\geq 0} \Mod_{\ulR}^{\Sigma_r}
$$
where $\Mod^{\Sigma_r}_{\ulR}$ is the $\infty$-category from Example \ref{ex:derived reps}. In other words, a derived symmetric sequence over $R$ has an arity $r$ component given by a module over the constant cohomological Mackey functor $\underline{R}$ in the $\infty$-category of genuine $\Sigma_r$-spectra. 
\end{definition}
Let us point out that the definition of the $\infty$-category $\sSeq^{\gen}_{\ulR}$ also makes sense when $R$ is an $\mathbb{E}_\infty$-algebra (or even an $\mathbb{E}_1$-algebra) over $\mathbb{Z}$. However, we will only be interested in the case where $R$ is a simplicial commutative ring, because in that case we can endow $\sSeq^{\gen}_{\ulR}$ with a strict version of the composition product (see Construction \ref{con:derived products}, which proceeds by induction from the case of a discrete ring $R$).
\begin{notation}
Recall from Example \ref{ex:derived reps} that each $\Mod^{\Sigma_r}_{\ulR}\simeq \Mod_{R[\Orb_{\Sigma_r}]}$ can be obtained as the $\infty$-category of modules over the full additive subcategory $R[\Orb_{\Sigma_r}]\subseteq \Mod^{\Sigma_r}_{\ulR}$ spanned by the free $\ul{R}$-modules on finite $\Sigma_r$-sets.

Write $\ROS:=\bigoplus_{r\geq 0} R[\Orb_{\Sigma_r}]$ for the sum of all of these additive $\infty$-categories (Example \ref{ex:sums of additive cats}). One can identify $\ROS\subseteq \sSeq^{\gen}_{\ulR}$ with the full subcategory spanned by the $\ul{R}$-linearisations $R[X]$ of finite symmetric sequences of sets (Definition \ref{def:finite symmetric set}) and Example \ref{ex:sums of additive cats} shows that
$$
\sSeq^{\gen}_{\ulR}\simeq \Mod_{\ROS}
$$
coincides with the $\infty$-category of modules over the additive $\infty$-category $\Mod_{\ROS}$. Example \ref{ex:derived reps} shows that $\ROS$ is coherent if $R$ is a coherent simplicial ring.
\end{notation}
\begin{example}
Let $R$ be a discrete commutative ring and let $\sSeq_R^\heartsuit$ denote the (ordinary) category of symmetric sequences of discrete $R$-modules. This can also be identified with the heart of the $t$-structure on $\sSeq_R$ provided by Definition \ref{def:classical t-structure on modules}. It follows from Example \ref{ex:derived reps} that $\ROS$ can be identified with the full subcategory of the  category $\sSeq^{\heartsuit}_R$ of symmetric sequences of discrete $R$-modules, spanned by the symmetric sequences $R[X]$ with $X$ a finite symmetric sequence of sets.
\end{example}
\begin{lemma}\label{lem:extension to SCR}
There is a natural sifted-colimit-preserving functor $\SCR\rt \mathscr{A}\mm{dd}$ sending $R\mapsto \ROS$.
\end{lemma}
\begin{proof}
It suffices to show that in each individual arity $r$, the assignment $R\mapsto R[\Orb_{\Sigma_r}]$ extends to a functor $\SCR\rt \mathscr{A}\mm{dd}$ preserving sifted colimits. To see this, notice that $\Mod_{\ulR}^{\Sigma_r}=\Mod_{\ulR}(\Sp^{\Sigma_r})$ depends functorially on the simplicial commutative ring $R$, via
$$\begin{tikzcd}[column sep=2.5pc]
F\colon \SCR\arrow[r] & \CAlg(\Mod_{\mathbb{Z}, \geq 0})\arrow[r, "R\mapsto \ul{R}"] & \CAlg(\Sp^{\Sigma_r})\arrow[r, "\Mod"] & \cat{Pr^L}.
\end{tikzcd}$$
Here the first functor sends a simplicial commutative rings to the corresponding $\mathbb{E}_\infty$-ring spectrum over $\ZZ$, the second functor sends this to the corresponding constant Mackey functor and the last functor sends $A\in \CAlg(\Sp^{\Sigma_r})$ to the (stable) presentable $\infty$-category $\Mod_A(\Sp^{\Sigma_r})$. The first two functors manifestly preserve sifted colimits (which are computed on the underlying object) and the last functor preserves sifted colimits by \cite[Corollary 4.8.5.13]{lurie2014higher}.

For any map of simplicial rings $R\rt R'$, the induced left adjoint functor $\Mod_{\ulR}^{\Sigma_r}\rt \Mod_{\ul{R}'}^{\Sigma_r}$ simply induces along $\ul{R}\rt \ul{R}'$. In particular, this sends the full subcategory $R[\Orb_{\Sigma_r}]$ to $R'[\Orb_{\Sigma_r}]$. One then obtains $\SCR\rt \mathscr{A}\mm{dd}; R\mapsto R[\Orb_{\Sigma_r}]$ as a diagram of full subcategories. This preserves sifted colimits because the functor $\mathscr{A}\mm{dd}\rt \cat{Pr^L}; \cat{A}\mapsto \Mod_{\cat{A}}$ preserves colimits and detects equivalences, by the universal property discussed in Definition \ref{def:module cats}.
\end{proof}

\begin{example}[Borel derived symmetric sequences]\label{ex:borel derived sym seq}
For any simplicial commutative ring $R$, there is a fully faithful inclusion $R[\Sigma]\hookrightarrow \ROS$ with essential image given by $R$-linearised finite $\Sigma$-free symmetric sequences of sets. This induces a fully faithful inclusion $\sSeq_R\rt \sSeq^{\gen}_{\ulR}$. We will refer to the essential image of the inclusion as the \emph{Borel derived symmetric sequences}.
\end{example}
\begin{example}\label{ex:linearised symmetric sequences}
There is a functor $\sSeq(\cat{Set})\rt \sSeq^{\gen}_{\ulR}$ sending a set-valued symmetric sequence $X$ to its $R$-linearisation $R[X]$. This functor is uniquely characterized by the fact that it preserves coproducts and sends a symmetric sequence in arity $r$ of the form $\Sigma_r/H$ to the object $R[\Sigma_r/H]$ in $\ROS$.
\end{example}
\begin{example}[Discrete symmetric sequences]\label{ex:discrete sym seq}
If $R$ is a discrete commutative ring, every discrete symmetric sequence $X$ determines a derived symmetric sequence $X^\mm{gen}$, given by the additive functor $\ROS^{\op}\rt\Sp$ sending each $R[S]$ to the discrete abelian group $\smash{\Hom_{\sSeq^\heartsuit_R}}(R[S], X)$. This encodes the data of all fixed points $X(r){}^H$ with $H<\Sigma_r$.
\end{example}
Just like on symmetric sequences, there is a plethora of monoidal structures on derived symmetric sequences.
\begin{cons}[Monoidal structures on derived symmetric sequences]\label{con:derived products}
If $R$ is a discrete ring, then the full subcategory $\ROS\subseteq \sSeq_R$ is closed under the monoidal structures $\circ$, $\otimes$ and $\otimes_{\mm{lev}}$ from Section \ref{sec:operads}. Note that for every map of rings $f\colon R\rt S$, the induced functor $\ROS\rt S[\Orb_\Sigma]$ preserves these monoidal structures, as well as all the compatibilities between them, e.g.\ the natural transformation exhibiting $\otimes_{\lev}$ as (op)lax monoidal with respect to $\circ$. We can then define all of these structures for a simplicial commutative ring $R$ as well, using functoriality over polynomial rings and extending by sifted colimits (using Lemma \ref{lem:extension to SCR}). 

As in Lemma \ref{lem:products extended}, all of these monoidal structures are given by locally polynomial functors. If $R$ is a coherent simplicial commutative ring, \Cref{thm:derived functors} shows that they extend to monoidal structures $\circ, \otimes, \otimes_{\mm{lev}}$ on $\sSeq^{\gen}_{\ulR}$, which preserve sifted colimits and finite totalisations. Furthermore, all of these monoidal structures preserve all colimits in the first variable. 
\end{cons}
\begin{definition}[Derived $\infty$-(co)operads]
Let $R$ be a coherent simplicial commutative ring. We define a \emph{derived $\infty$-operad} over $R$ to be an associative algebra in $\sSeq^{\gen}_{\ulR}$ with respect to the derived composition product $\circ$. Likewise, a \emph{derived $\infty$-cooperad} is a coalgebra in $\sSeq^{\gen}_{\ulR}$. We will write $\Op^{\gen}_{\ulR}$ and $\coOp^{\gen}_{\ulR}$ for the $\infty$-categories of derived $\infty$-operads and $\infty$-cooperads, respectively.
\end{definition}
\begin{example}
The inclusion $R[\Sigma]\hookrightarrow \ROS$ is preserves the composition product. Consequently, the inclusion $\sSeq_R\rt \sSeq^{\gen}_{\ulR}$ of the Borel derived symmetric sequences preserves the composition product and its right adjoint is lax monoidal for the composition product. It follows that there is an adjoint pair $\Op_R\leftrightarrows \Op_{\ulR}^{\gen}$ where the left adjoint includes $R$-linear $\infty$-operads into the derived $\infty$-operads and the right adjoint sends each derived $\infty$-operad to the underlying `Borel operad'.
\end{example}
\begin{example}[Algebraic operads]\label{ex:algebraic operads}
Let $R$ be a discrete coherent ring. Then there is an adjoint pair $F\colon \sSeq^{\gen}_{\underline{R}, \geq 0}\leftrightarrows \sSeq_R^\heartsuit\colon (-)^\mm{gen}$ where the right adjoint is as in Example \ref{ex:discrete sym seq}. The left adjoint is the nonabelian derived functor of the inclusion $\ROS\rt \sSeq_R^\heartsuit$. This functor is monoidal for the composition product, so that $(-)^\mm{gen}$ is lax symmetric monoidal. Consequently, every classical $R$-linear operad $\PP$ determines a derived $\infty$-operad $\PP^\mm{gen}$. This construction can be understood more concretely in terms of our point-set models, see Remark \ref{rem:connective simplicial-cosimplicial}. 
\end{example}
\begin{example}[Derived commutative operad]\label{ex:derived commutative operad}
Applying the previous example to the commutative operad gives a derived $\infty$-operad that we will denote by $\mm{Com}$. Unraveling the definitions, $\mm{Com}$ is the derived symmetric sequence given in each arity $r$ by the free $R$-module on the point, equipped with the trivial $\Sigma_r$-action. Using this, one easily sees that the derived symmetric sequence underlying $\mm{Com}$ is the unit for the levelwise tensor product.
\end{example}
\begin{remark}[Formula for derived composition product]\label{rem:derived action via derived orbits}
For a coherent simplicial ring $R$, the $r$-fold Day convolution product of a derived symmetric sequence $Y$ admits a genuine $\Sigma_r$-equivariant structure. More precisely, there is a functor
$$\begin{tikzcd}
\sSeq^{\gen}_{\ulR}\arrow[r] & \prod_{q\geq 0} \Mod^{\Sigma_r\times \Sigma_q}_{\ulR}=\Mod_{R[\Orb_{\Sigma_r\times\Sigma}]}; & Y\arrow[r] & Y^{\otimes r}
\end{tikzcd}$$
obtained by left-right extending a polynomial functor $T_r\colon \ROS\rt R[\Orb_{\Sigma_r\times\Sigma}]$: when $R$ is a discrete ring, $T_r$ simply sends the linearisation $R[K]$ of a finite symmetric sequence of sets to $R[K^{\otimes r}]$ and one extends to general simplicial commutative rings by sifted colimits (Lemma \ref{lem:extension to SCR}).
The composition product on $\sSeq^{\gen}_{\ulR}$ can then be identified with
$$
X\circ Y \cong \bigoplus_{r\geq 0} \big(X(r)\otimes Y^{\otimes r}\big)_{\Sigma_r}
$$
where we take the tensor product of $X(r)$ and $Y^{\otimes r}$ in $\prod_q \Mod^{\Sigma_r\times \Sigma_q}_{\ulR}$ and then take genuine $\Sigma_r$-orbits (Example \ref{ex:derived fixed points}). Indeed, both functors are obtained as left-right extensions and coincide on $\ROS$.
\end{remark}
The derived composition product $\circ$ restricts to an action $\circ\colon \sSeq^{\gen}_{\ulR}\times \Mod_R\rt \Mod_R$, where we identify $\Mod_R$ with the full subcategory of derived symmetric \mbox{sequences concentrated in arity $0$.}
\begin{example}
Let $X$ be a symmetric sequence of sets and let $R[X]$ be the associated derived symmetric sequence (Example \ref{ex:linearised symmetric sequences}). The induced endofunctor of $\Mod_R$ is the right-left extended functor of the functor sending a finitely generated free $R$-module $V$ to
$\bigoplus_{r\geq 0} \big(X(r)_+\wedge V^{\otimes r}\big){}_{\Sigma_r}.$
\end{example}
\begin{definition}[Algebras over derived operads]
Let $R$ be a simplicial commutative ring and $\PP$ a derived $\infty$-operad. We define a $\PP$-algebra to be a left $\PP$-module in $\Mod_R$ with respect to the composition product. Likewise, a coalgebra over a derived $\infty$-cooperad is a left comodule in $\Mod_R$ with respect to the composition product. We will write $\Alg_{\PP}(\Mod_R)$ and $\Coalg_{\CC}(\Mod_R)$ for the $\infty$-categories of (co)algebras.
\end{definition}
\begin{example}[Derived commutative algebras]\label{ex:derived commutative}
The monad associated to the derived commutative $\infty$-operad is the right-left extension of the functor sending a finitely generated free $R$-module $V$ to the symmetric algebra $\mm{Sym}_R(V)=\bigoplus_{n\geq 0} (V^{\otimes n})_{\Sigma_n}$.
In particular, the $\infty$-category of \emph{connective} algebras over the derived $\infty$-operad $\mm{Com}$ is the $\infty$-category of simplicial commutative (i.e.\ animated) $R$-algebras. 
The $\infty$-category of all algebras for this monad is the $\infty$-category of derived rings, as studied  by Bhatt-Mathew, Raksit \cite{raksit2020hochschild} and others. 
\end{example}

\subsubsection*{Derived PD operads}
We will now discuss a version of derived $\infty$-operads with divided powers, following the discussion in Section \ref{sec:PD operads}. 
\begin{definition}[Pro-coherent derived symmetric sequences] \label{derivedProCoh}
Let $R$ be a coherent simplicial commutative ring. A \emph{derived pro-coherent symmetric sequence} over $R$ is a pro-coherent module over the additive $\infty$-category $\ROS$. We will denote the $\infty$-category of pro-coherent derived symmetric sequences by $\sSeq^{\gen, \vee}_{\ulR}$.
\end{definition}
\begin{remark}
The fully faithful inclusion $\sSeq_R\hookrightarrow \sSeq^{\gen}_{\ulR}$ of the Borel derived symmetric sequences extends to a fully faithful inclusion $\sSeq_R^{\vee}\hookrightarrow \sSeq_{\ulR}^{\gen, \vee}$ between pro-coherent objects.
\end{remark}
Recall that the $\infty$-category $\sSeq^{\gen, \vee}_{\ulR}$ is (often) a further enlargement of the $\infty$-category of derived symmetric sequences, which also contains the continuous $R$-linear duals of almost perfect derived symmetric sequences. 
We start by studying the operation of taking $R$-linear dual pro-coherent derived symmetric sequences. To this end, note that we can use \Cref{thm:derived functors} to endow $\sSeq^{\gen, \vee}_{\ulR}$ with the levelwise tensor product. The unit for this tensor product is the derived symmetric sequence $\mm{Com}$, given in each arity by the trivial $\Sigma_r$-representation on $R$.
\begin{notation}[Linear dual pro-coherent derived symmetric sequences]
Let $R$ be a coherent simplicial commutative ring. If $X$ is a pro-coherent derived symmetric sequence, we will write $X^\vee$ for its dual with respect to the levelwise tensor product and refer to it as the $R$-linear dual of $X$.
\end{notation}
We have the following analogue of Proposition \ref{prop:dual symmetric sequence}.
\begin{proposition}\label{prop:dual derived symmetric sequences}
Let $R$ be a coherent simplicial commutative ring and consider the functor $(-)^\vee\colon \sSeq^{\gen, \vee}_{\ulR}\rt \sSeq^{\gen, \vee, \op}_{\ulR}$ taking $R$-linear duals. This functor is the right-left extension of the equivalence $\ROS\rt \ROS^{\op}$ sending each $R$-linearised finite symmetric sequence of sets $R[X]$ to the $R$-linear dual symmetric sequence $R[X]^\vee(r)=\Hom_R\big(R[X](r), R\big)$. Furthermore, it restricts to an equivalence
$$\begin{tikzcd}
(-)^\vee\colon \APerf^\vee_{\ROS}\arrow[r, "\simeq"] & \APerf_{\ROS}^{\op}.
\end{tikzcd}$$
\end{proposition} 
In particular, the image of $\upiota\colon \sSeq^{\gen}_{\ulR}\rt \sSeq^{\gen, \vee}_{\ulR}$ is typically not closed under duality.
\begin{proof}
Notice that the symmetric sequence $R[X]^\vee$ indeed defines an object in $\ROS$ (isomorphic to $R[X]$ itself) and that the resulting functor $(-)^\vee\colon \ROS\rt \ROS^{\op}$ is an equivalence. We claim that $R[X]^\vee$ is indeed the dual of $R[X]$ with respect to the levelwise tensor product. For each arity $r$, there are canonical maps in $R[\Orb_{\Sigma}]$ of the form $\mm{Com}(r)\rt X(r)\otimes_{\lev} X(r)^\vee$ and $X(r)^\vee\otimes_{\mm{lev}} X(r)\rt \mm{Com}(r)$ exhibiting $X^\vee$ as the dual of $X$. The result now follows from Proposition \ref{prop:dual is indeed dual}.
\end{proof}
By \Cref{thm:derived functors}, the composition product on derived symmetric sequences (Construction \ref{con:derived products}) extends to a composition product on pro-coherent derived symmetric sequences. For applications to Koszul duality, we will be more interested in a version of the composition product based on strict invariants, rather than strict orbits:
\begin{definition}[Restricted composition product]\label{not:restricted composition}
Let $R$ be a simplicial commutative ring. Conjugating the composition product on $\ROS$ by the self-equivalence $(-)^\vee\colon \ROS\rt \ROS^{\op}$ yields another monoidal structure, usually referred to as the \emph{restricted composition product}. Explicitly, this monoidal structure on $\ROS$ can be identified with
$$
X\bcirc Y = (X^\vee\circ Y^\vee){}^\vee\cong \bigoplus_r \big(X(r)\otimes Y^{\otimes r}\big){}^{\Sigma_r}.
$$
Here $X(r)\otimes Y^{\otimes r}$ defines an object in $R[\Orb_{\Sigma_r\times \Sigma}]$, as in Remark \ref{rem:derived action via derived orbits} and $(-)^{\Sigma_r}$ takes genuine $\Sigma_r$-fixed points.

This has properties analogous to the usual composition product; for example, the levelwise tensor product $\otimes_{\lev}$ is both lax and oplax monoidal with respect to $\bcirc$. The norm maps $\big(X(r)\otimes Y^{\otimes r}\big){}_{\Sigma_r}\rt \big(X(r)\otimes Y^{\otimes r}\big){}^{\Sigma_r}$ determine a natural map 
$$
\mm{Nm}\colon X\circ Y\rt X\bcirc Y.
$$
This endows the identity functor with the structure of a lax monoidal functor $\big(\ROS, \circ\big)\rt \big(\ROS, \bcirc\big)$ (see e.g.\ \cite{fresse2000homotopy}). Note that the norm map is an equivalence if $X$ is $\Sigma$-free or $Y$ is concentrated in arity $\geq 1$. All of these properties and structures are verified directly when $R$ is a discrete ring and hold for simplicial commutative rings by taking sifted colimits (as in Construction \ref{con:derived products}).
\end{definition}
Using the results from Section \ref{sec:extended functors}, the various products considered above can now be extended to pro-coherent derived symmetric sequences:
\begin{proposition}\label{prop:derived pro-coh composition}
Let $R$ be a coherent ring. Then the monoidal structures $\circ, \bcirc, \otimes$ and $\otimes_{\lev}$ all admit right-left extensions to monoidal structures on the categories $\sSeq^{\gen}_{\ulR}$ and $\sSeq^{\gen, \vee}_{\ulR}$. Furthermore, these monoidal structures have the following properties:
\begin{enumerate}
\item Each of the four monoidal structures $\circ, \bcirc, \otimes$ and $\otimes_{\lev}$ preserves sifted colimits and all colimits in the first variable.

\item There is a commuting square of left adjoint functors
$$\begin{tikzcd}
\sSeq_R\arrow[r]\arrow[d, "\upiota"{swap}]\arrow[r, hook] & \sSeq_{\ulR}^{\gen}\arrow[d, "\upiota"]\\
\sSeq_R^\vee\arrow[r, hook] & \sSeq^{\gen, \vee}_{\ulR}
\end{tikzcd}$$
where the horizontal functors include the Borel (pro-coherent) derived symmetric sequences. All of these functors are (symmetric) monoidal with respect to $\circ, \bcirc, \otimes$ and $\otimes_{\lev}$. Here we identify $\bcirc=\circ$ on symmetric sequences and pro-coherent symmetric sequences.

\item The functor $\otimes_{\lev}$ is both lax and oplax monoidal with respect to $\circ$ and $\bcirc$.

\item There is a natural norm map $\mm{Nm}\colon X\circ Y\rt X\bcirc Y$ that endows the identity functor with the structure of a lax monoidal functor. The norm map is an equivalence if $X$ is a pro-coherent Borel derived symmetric sequence or if $Y$ is concentrated in arity $\geq 1$.
\end{enumerate}
\end{proposition}
\begin{proof}
The proof of Proposition \ref{prop:pro-coh composition} carries over mutatis mutandis. Note that the horizontal fully faithful inclusions in (2) are induced by the fully faithful inclusion $R[\Sigma]\rt \ROS$. The full subcategory $R[\Sigma]$ is closed under each of the four tensor products and furthermore the two composition products $\circ$ and $\bcirc$ coincide on $R[\Sigma]$ (since the norm map is an equivalence on Borel derived symmetric sequences). \Cref{thm:derived functors} then implies that the functors in the diagram are (symmetric) monoidal for each of the four products.
\end{proof}
\begin{example}\label{ex:strict invariants}
When $R$ is a coherent simplicial ring, the monoidal structure $\bcirc$ on $\sSeq^{\gen, \vee}_{\ulR}$ is given by the formula
$$
X\bcirc Y = \bigoplus_r \big(X(r)\otimes Y^{\otimes r}\big){}^{\Sigma_r}.
$$
Here $X(r)\otimes Y^{\otimes r}$ defines an object in $\QC^\vee_{R[\Orb_{\Sigma_r\times \Sigma}]}$ (using Remark \ref{rem:derived action via derived orbits}) and $(-)^{\Sigma_r}$ is the derived genuine fixed points functor from Example \ref{ex:derived fixed points}.

In particular, when $X\simeq \mm{Tot}(X^\bullet)$ and $Y\simeq \mm{Tot}(Y^\bullet)$ arise as totalisations of cosimplicial diagrams in $\ROS$, the value is given by the derived strict invariants
$$
X\bcirc Y \simeq \mm{Tot}\Big(\bigoplus_{r\geq 0}\big(X^\bullet(r)\otimes (Y^{\bullet})^{\otimes r}\big){}^{\Sigma_r}\Big).
$$
\end{example}
\begin{definition}[Derived PD operads]\label{def:derived PD operad}
Let $R$ be a coherent simplicial commutative ring. A \emph{derived PD $\infty$-operad} over $R$ is defined to be an associative algebra in $\sSeq^{\gen, \vee}_{\ulR}$ with respect to the restricted composition product $\bcirc$. We will denote the $\infty$-category of derived PD $\infty$-operads by $\Op_{\ulR}^{\gen,\pd}$.
\end{definition}
\begin{remark}
By part (4) of Proposition \ref{prop:derived pro-coh composition}, there is a forgetful functor from derived PD $\infty$-operads to algebras in $\sSeq^{\gen, \vee}_{\ulR}$ with respect to the composition product $\circ$. This forgetful functor is an equivalence for derived PD $\infty$-operads without operations in arity $0$.
\end{remark}
The restricted composition product $\bcirc$ induces an action of $\sSeq^{\gen, \vee}_{\ulR}$ on $\QC^\vee_R$.
\begin{definition}[Algebras over derived PD operads]
An algebra over a derived PD $\infty$-operad $\PP$ is a left $\PP$-module in $\QC^\vee_R$ with respect to the $\bcirc$-action. We will write $\Alg^{\gen}_{\PP}(\QC^\vee_R)$ for the $\infty$-category of $\PP$-algebras.
\end{definition}
\begin{example}[Divided power algebras]
The derived commutative $\infty$-operad $\mm{Com}$ from Example \ref{ex:derived commutative} admits a nonunital version $\mm{Com}^\mm{nu}$. Since $\mm{Com}^\mm{nu}$ is trivial in arity $0$, its image in $\sSeq^{\gen, \vee}_{\ulR}$ has the structure of a derived PD $\infty$-operad. The corresponding monad on $\QC^\vee_R$ is the right-left extension of the functor sending a finitely generated free $R$-module $V$ to the divided power algebra $\Gamma_R(V)=\bigoplus_{r\geq 1} (V^{\otimes r})^{\Sigma_r}$.
\end{example}

\subsection{Refined Koszul duality for derived PD operads}
Finally, we shall discuss a refinement of the classical Koszul duality for $\infty$-operads to the setting of derived $\infty$-operads. As a first step, the $\infty$-categorical bar construction yields a functor from augmented derived $\infty$-operads to derived $\infty$-cooperads and from derived algebras to derived coalgebras
$$\begin{tikzcd}[column sep=1.6pc]
\Barr\colon \Op_{\ulR}^{\gen,\aug} \arrow[r] & \Coop^{\gen, \aug}_{\ulR} & \Barr\colon \Alg^{\gen}_{\PP}(\Mod_R)\arrow[r] & \Coalg_{\Barr(\PP)}(\Mod_R).
\end{tikzcd}$$
If $\PP$ is a derived $\infty$-operad, then we define its (refined) Koszul dual to be the pro-coherent $R$-linear dual of its bar construction. This carries the structure of a derived PD $\infty$-operad by the following observation:
\begin{proposition}\label{prop:dualisation derived lax monoidal}
Linear duality induces oplax monoidal functors 
$$\begin{tikzcd}
\big(\sSeq^{\gen, \vee}_{\ulR}, \circ\big)\arrow[r, "(-)^\vee"] & \big(\sSeq^{\gen, \vee, \op}_{\ulR}, \bcirc\big) & \big(\sSeq^{\gen, \vee}_{\ulR}, \bcirc\big)\arrow[r, "(-)^\vee"] & \big(\sSeq^{\gen, \vee, \op}_{\ulR}, \circ\big)
\end{tikzcd}$$
restricting to (strong) monoidal equivalences between almost perfect \mbox{(dually almost perfect) objects.}
\end{proposition}
\begin{proof}
Let us only treat the first case and write $F\colon \ROS \rt \ROS^{\op}\subseteq \sSeq^{\gen, \vee, \op}_{\ulR}$ for the functor taking $R$-linear duals. By Proposition \ref{prop:dual derived symmetric sequences}, linear duality is the right-left extension of $F$. The construction of the restricted composition product (\Cref{not:restricted composition}) implies that $F$ is a strong monoidal functor. Since $\bcirc$ preserves sifted colimits in each variable, the right extension $F^R\colon \APerf^\vee_{\ROS, \weirdleq 0}\rt \sSeq_{\ulR}^{\gen, \vee, \op}$ remains strong monoidal. By (the opposite of) Lemma \ref{lem:right kan extension lax monoidal}, the left extension $F^{LR}$ then inherits an oplax monoidal structure. It is strong monoidal on dually almost perfect objects by the construction and the fact that $\bcirc$ and $\circ$ both preserve finite geometric realisations and totalisations in each variable.
\end{proof}
\begin{definition}[Koszul duality for derived operads]\label{kdforderived}
Let $R$ be a coherent simplicial commutative ring and $\PP$ an augmented derived $\infty$-operad over $R$. We define the \emph{Koszul dual derived PD $\infty$-operad} of $\PP$ to be the pro-coherent $R$-linear dual of the bar construction $\KD^{\pd}(\PP)=\Barr(\PP)^\vee$.
\end{definition}
This refines the Koszul duality of \Cref{thm:koszul duality for PD operads}: Koszul duality fits into a commuting square where the vertical arrows include the Borel derived (PD) $\infty$-operads
$$\begin{tikzcd}[column sep=2.2pc]
\Op_R^{\aug}\arrow[r, "\KD^{\pd}"] \arrow[d, hook] & \big(\Op_R^{\pd, \aug}\big)^{\op}\arrow[d, hook]\\
\Op^{\gen,\aug}_{\ulR}\arrow[r, "\KD^{\pd}"] & \big(\Op^{\gen, \pd, \aug}_{\ulR}\big)^{\op}.
\end{tikzcd}$$
\begin{definition}\label{kdforderivedalgebras}
Let $R$ be a coherent simplicial commutative ring and $\PP$ an augmented derived $\infty$-operad over $R$. If $A$ is a $\PP$-algebra, then we define its \emph{Koszul dual $\KD^{\pd}(\PP)$-algebra} to be the pro-coherent $R$-linear dual of the bar construction $\KD^{\pd}(A)=\Barr_{\PP}(A)^\vee$.
\end{definition}
Let us give a more explicit description of the monad associated to the Koszul dual of a derived $\infty$-operad satisfying some finiteness conditions:
\begin{definition}[Almost finitely presented derived $\infty$-operads]
Let $R$ be a simplicial commutative ring. An augmented derived $\infty$-operad $\PP$ over $R$ is said to be \emph{connective} if its underlying derived symmetric sequence is connective. A connective augmented derived $\infty$-operad $\PP$ is said to be \emph{almost finitely presented} if it defines an almost compact object in the $\infty$-category $\Op^{\gen, \aug}_{R, \geq 0}$ of connective augmented derived $\infty$-operads, in the sense of \cite[Definition 7.2.4.8]{lurie2014higher}.
\end{definition}
We can then describe the monad induced by $\KD^{\pd}(\PP)$ in terms of the adjunction $$\cot_{\PP}\colon \Alg^{\gen}_{\PP}(\QC^\vee_R)\leftrightarrows \QC^\vee_R\colon \triv_{\PP}$$ arising from the augmentation map of derived $\infty$-operads $\PP\rt \mb{1}$:
\begin{proposition}\label{prop:derived PD koszul duality in aft case}
Let $R$ be a coherent simplicial commutative ring and $\PP$ an almost finitely presented augmented derived $\infty$-operad over $R$. Then $\KD^{\pd}(\PP)$ is dually almost perfect. The induced monad $\KD^{\pd}(\PP)\colon \QC^\vee_R\rt \QC^\vee_R$ preserves sifted colimits and dually almost perfect objects, and the resulting monad on $\APerf^\vee_R$ can be identified with
$$
\KD^{\pd}(\PP)(V)\simeq \bigoplus_r \big(\KD(\PP)^{\pd}(r)\otimes V^{\otimes r}\big)^{\Sigma_r} \simeq \big(\cot_{\PP}\circ \triv_{\PP}(V^\vee)\big)^\vee.
$$
\end{proposition}
\begin{proof}
The proof of Proposition \ref{prop:PD koszul duality in aft case} carries over verbatim. The first formula for the monad $\KD^{\pd}(\PP)$ follows from Example \ref{ex:strict invariants}. The second equivalence follows from Proposition \ref{prop:comonads from bar construction} (applied to $\Cl=(\sSeq^{\gen, \vee}_{\ulR})_{1//1}$ acting on $\Cm=\QC^\vee_R$) and the fact that linear duality gives a monoidal equivalence between dually almost perfect and almost perfect pro-coherent derived symmetric sequences.
\end{proof}
One can verify that the nonunital commutative derived $\infty$-operad is almost finitely presented. In particular, its bar construction is the derived symmetric sequence $\Barr(\Com^{\mm{nu}})(r)\simeq R[\Sigma|\Pi_r|^\diamond]$ of reduced-unreduced suspensions of the nerve of the partition complex; this is indeed an almost perfect derived symmetric sequence (which is in fact all that we need), arising as the geometric realisation of a simplicial diagram of derived symmetric sequences as in Example \ref{ex:linearised symmetric sequences}. Let us write $\Hom(\Sigma|\Pi_r|^\diamond, \underline{R})$ for its (pro-coherent) $R$-linear dual in $\QC^\vee_{R[\Orb_{\Sigma_r}]}$; in terms of the simplicial-cosimplicial models from Section \ref{sec:simp-cosimp point-set models}, this can be described by the cosimplicial $R$-module of $R$-valued functions on the simplicial set $\Sigma|\Pi_r|^\diamond$.
\begin{definition}[The derived partition Lie PD operad]\label{def:derived partition lie operad}
The \emph{(derived) partition Lie PD $\infty$-operad} is the Koszul dual $\mm{Lie}^\pi_{R, \Delta}=\KD^{\pd}(\mm{Com}^\mm{nu})$. 
\end{definition} 
\begin{corollary}\label{cor:derived partition lie monad}
The monad associated to the derived partition Lie PD $\infty$-operad agrees (over a field $k$) with the partition Lie monad from \cite[Definition 5.47]{brantner2019deformation} and is given on dually almost perfect objects by 
$$
\mm{Lie}^\pi_{R, \Delta}(V)=\bigoplus_r \big(\Hom(\Sigma|\Pi_r|^\diamond, \underline{R})\otimes V^{\otimes r}\big){}^{\Sigma_r}.
$$
\end{corollary}

\ \\ \\
\section{Chain models for PD operads}
\label{sec:chain models}
In the previous section we have given an $\infty$-categorical discussion of PD $\infty$-operads over coherent $\mathbb{E}_\infty$-ring spectra. Every PD $\infty$-operad determines a sifted-colimit-preserving monad on the $\infty$-category of pro-coherent $R$-modules, which can be constructed as a right-left extended functor and can be described by a formula involving the \emph{divided orbits} of Example \ref{dG1}.

The purpose of this section is to provide explicit point-set models for these $\infty$-categorical constructions in the case where $R$ is a discrete coherent ring. In particular, we give a presentation of PD $\infty$-operads and their algebras in terms of chain complexes of $R$-modules. As a motivation for all the constructions appearing in this section, we shall give the following example:
\begin{example}\label{ex:barratt-eccles}
Let $k$ be a field and let $\ope{E}^\mm{nu}$ denote the (nonunital) Barratt--Eccles operad, given by $\ope{E}^\mm{nu}(r)=C_*(E\Sigma_r)$. In particular, each $\ope{E}^\mm{nu}(r)$ is given by a chain complex of finitely generated free $k[\Sigma_r]$-modules, in nonnegative degrees. For any chain complex $V$, the composition product $\ope{E}^\mm{nu}(r)$ then computes the free nonunital $\mathbb{E}_\infty$-algebra
$$
\ope{E}^\mm{nu}\circ V = \bigoplus_{r>0} \big(\ope{E}^\mm{nu}(r)\otimes V^{\otimes r}\big)_{\Sigma_r} \simeq \bigoplus_{r>0} V^{\otimes r}_{h\Sigma_r}.
$$
The last equivalence uses that $\ope{E}^\mm{nu}(r)$ is a projective resolution of the trivial $\Sigma_r$-representation.

On the other hand, consider the linear dual $\ope{E}^{\mm{nu}, \vee}$ of the Barratt--Eccles operad. This does not admit an obvious operad structure, but in Appendix \ref{app:surjections cooperad}, we construct a dg-operad $\mathbf{\Sur}^\vee$ whose underlying symmetric sequence is chain homotopic to $\ope{E}^{\mm{nu}, \vee}$. Leaving this issue aside, note that $\ope{E}^{\mm{nu},\vee}\circ V$ does \emph{not} compute the free $\mathbb{E}_\infty$-algebra $\bigoplus_{r>0} V^{\otimes r}_{h\Sigma_r}$ on $V$, even though $\ope{E}^{\mm{nu}, \vee}$ is quasi-isomorphic to $\ope{E}^\mm{nu}$. Indeed, even though $\ope{E}^{\mm{nu}, \vee}(r)$ is a chain complex of finitely generated free $k[\Sigma_r]$-modules, it is not a projective resolution of the trivial $\Sigma_r$-representation; instead it is an injective resolution. Consequently, for any bounded above complex $V$ we now have that
$$
\ope{E}^{\mm{nu}, \vee}\circ V = \bigoplus_{r>0} \big(\ope{E}^{\mm{nu}, \vee}(r)\otimes V^{\otimes r}\big)_{\Sigma_r} \cong \bigoplus_{r>0} \big(\ope{E}^{\mm{nu}, \vee}(r)\otimes V^{\otimes r}\big)^{\Sigma_r} \simeq \bigoplus_{r>0} (V^{\otimes r})^{h\Sigma_r}
$$
computes a free nonunital $\mathbb{E}_\infty$-algebra with divided powers. Here the second isomorphism uses that $\ope{E}^{\mm{nu}, \vee}(r)$ is a complex of finitely generated free $k[\Sigma_r]$-modules, so that the norm map is an isomorphism (which also holds for $\ope{E}^\mm{nu}$) and the last equivalence uses crucially that $\ope{E}^{\mm{nu}, \vee}(r)$ is an injective resolution of the trivial $\Sigma_r$-representation.
\end{example}

\begin{remark}
Notice that the above computation is not in conflict with the standard homotopy theory for operads from e.g.\ \cite{HinichHomotopy, bergermoerdijk2003}: even though the symmetric sequence $\ope{E}^{\mm{nu}, \vee}$ consists of complexes of free $\Sigma_r$-modules, it is not $\Sigma$-cofibrant in the usual sense and is hence usually excluded from considerations.\end{remark}

We will show that the $\infty$-category of PD $\infty$-operads over $R$ can be described by a homotopy theory of dg-operads in which more dg-operads are $\Sigma$-cofibrant, and hence  fewer dg-operads are weakly equivalent to one another. In particular, in Section \ref{sec:spectral part lie} we will describe a dg-operad controlling the theory of spectral partition Lie algebras, which is $\Sigma$-cofibrant (only) in this more liberal sense.

We start by discussing a chain model for the $\infty$-category of pro-coherent $R[G]$-modules and the divided orbits functor for a finite group $G$; in fact, for our later description of the $\infty$-category of derived $\infty$-operads, we will simultaneously treat pro-coherent modules over the additive $\infty$-category $R[\Orb_G]$ from Example \ref{ex:derived reps}.

\subsection{Chain models for pro-coherent modules}
\label{sec:equivariant models}
Throughout this section, we fix a discrete (commutative) ring $R$, a finite group $G$ and a full subcategory $\mc{F}\subseteq \Orb_G$ of the orbit category. We will only make use of the extreme cases where $\mc{F}$ contains only the trivial subgroup (later in this section) and where $\mc{F}=\Orb_G$ (in Section \ref{sec:simp-cosimp point-set models}).
\begin{definition}\label{def:admissible subgroup}
A $G$-set is said to be $\mc{F}$-\emph{admissible}, or briefly admissible, if each orbit is contained in $\mc{F}$, and a subgroup $H<G$ is said to be admissible if $G/H$ is an admissible orbit. 
If $R$ is a ring, then an $R$-linear $G$-representation $V$ is said to be a (\emph{finite}) \emph{$\mc{F}$-admissible representation} if $V\cong R[S]$ is the $R$-linearisation of an (finite) $\mc{F}$-admissible $G$-set. We will denote by $R[\mc{F}]$ the full subcategory of the category of (discrete) $R$-linear $G$-representations spanned by the finite $\mc{F}$-admissible representations.
\end{definition}
Similar to  \Cref{ex:derived reps}, $R[\mc{F}]$ is an additive category, which is coherent when $R$ is a coherent ring. By the formalism of Section \ref{sec:pro-coh}, we therefore obtain an $\infty$-category $\Mod_{R[\mc{F}]}$ of $R[\mc{F}]$-modules, as well as a fully faithful functor of $\infty$-categories in the situation where $R$ is coherent
$$
\upiota\colon \Mod_{R[\mc{F}]}\rt \QC^\vee\big(R[\mc{F}]\big).
$$
We will give model-categorical presentations of these $\infty$-categories in terms of chain complexes.
\begin{notation}
We denote the category of chain complexes of (discrete) $R[G]$-modules by $\Ch_{R[G]}$. This category is naturally enriched and tensored over the category $\Ch_R$ of chain complexes of $R$-modules. We will write $\Hom_{R[G]}(X, Y)$ for the mapping complex. If $X$ is a chain complex, denote the $n$-fold suspension by $X[n]$ and the cone of the $n$-fold suspension by $X[n, n+1]$.
\end{notation}
\begin{definition}\label{def:tame equivalence}
A complex $P$ of $R[G]$-modules is said to be \emph{$\mc{F}$-quasifree} if it is given in each degree by an $\mc{F}$-admissible $G$-representation. It is \emph{$\mc{F}$-quasiprojective} if it is the retract of an $\mc{F}$-quasifree complex of $R[G]$-modules.

A map of complexes of $R[G]$-modules $X\rt Y$ is is said to be an \emph{$\mc{F}$-tame weak equivalence} if the induced map on mapping complexes $\Hom_{R[G]}(P, X)\rt \Hom_{R[G]}(P, Y)$ is a quasi-isomorphism for every $\mc{F}$-quasiprojective object $P$.
\end{definition}
Taking $P=R[G/H]$, one sees that every $\mc{F}$-tame weak equivalence induces quasi-isomorphisms on $H$-fixed points for all admissible $H<G$.
\begin{proposition}\label{prop:chain model structures on G-reps}
Let $R$ be a ring, $G$ a finite group and $\mc{F}\subseteq \Orb_G$ a full subcategory. Then the category $\Ch_{R[G]}$ can be endowed with the following two combinatorial model structures:
\begin{enumerate}
\item the \emph{$\mc{F}$-projective} model structure, the weak equivalences of which are maps inducing quasi-isomorphisms on $H$-fixed points and fibrations are maps inducing surjections on $H$-fixed points, for all $\mc{F}$-admissible subgroups $H<G$.
\item the \emph{$\mc{F}$-tame} model structure, in which the weak equivalences are the $\mc{F}$-tame weak equivalences, the cofibrations are degreewise split monomorphisms with an $\mc{F}$-quasiprojective \mbox{cokernel} and the fibrations are maps inducing surjections on $H$-fixed points, for all  $\mc{F}$-admissible subgroups $H<G$. 
\end{enumerate}
Furthermore, both model structures are naturally enriched over $\Ch_R$, equipped with the projective model structure.
\end{proposition}
\begin{warning}
Since these model structures are enriched over $\Ch_R$, their associated $\infty$-categories are stable. However, unlike for many of the usual model structures on chain complexes, a short exact sequence in $\Ch_{R[G]}$ need not define a cofibre sequence in the associated stable $\infty$-category (because $H$-fixed points are not exact).
\end{warning}
\begin{proof}
We will only prove part (2), following the argument in \cite{Nuiten} (see also \cite{becker2014models}); part (1) follows a similar, but more classical proof. We first observe that for a map $p\colon X\rt Y$, the following four properties are equivalent:
\begin{enumerate}[label={(\alph*)}]
\item $p$ has the right lifting property against the cofibrations.
\item the map $\Hom_{R[G]}(P, X)\rt \Hom_{R[G]}(P, Y)$ is an acyclic fibration for all $\mc{F}$-quasifree complexes $P$.

\item the map $\Hom_{R[G]}(T, X)\rt \Hom_{R[G]}(T, Y)$ is an acyclic fibration for all bounded above complexes $T$ of finite $\mc{F}$-admissible representations.

\item $p$ is both a fibration and an $\mc{F}$-tame weak equivalence. 
\end{enumerate}
The equivalences between (a), (b) and (d) are formal. We write $\cat{T}$ for the set of complexes $T$ appearing in (c). The fact that (c) implies the stronger condition (b) relies on an inductive argument on the $G$-sets of $R$-linear generators of $P$, using that for every generator $x\in P$ there exists a subcomplex $x\in T\subseteq P$ with $T\in \cat{T}$ (see \cite[Lemma 8.6]{Nuiten} for more details). 

We then define the following sets of generating cofibrations and trivial cofibrations:
$$
I=\big\{T\rt T[0, 1]\colon T\in \mc{T}\big\}\qquad \qquad J=\big\{0\rt R[S][n, n+1]\colon S\in \mc{F}\big\}.
$$
By construction, a map has the right lifting property against $J$ if and only if it is a fibration and $I$ generates the class of cofibrations. It then remains to verify that a transfinite composition of pushouts of maps in $J$ is a cofibration and an $\mc{F}$-tame equivalence: this is clear, since such maps are summand inclusions $X\rt X\oplus Y$ where $Y$ is chain homotopic to zero.
\end{proof}
\begin{remark}\label{rem:O-tame we detected by finite stuff}
As a consequence of the proof, a map $X\rt Y$ is an $\mc{F}$-tame weak equivalence if and only if $\Hom_{R[G]}(T, X)\rt \Hom_{R[G]}(T, Y)$ is a quasi-isomorphism for every bounded above chain complex of finite $\mc{F}$-admissible $G$-representations.
\end{remark}
\begin{example}
Suppose that $G$ is the trivial group. If $R=k$ is a field, then the projective and tame model structures on $\Ch_k$ are easily seen to coincide (this holds more generally when $R$ is a regular Noetherian ring, by Example \ref{ex:regular} and Corollary \ref{cor:chain models for pro-coh G-reps} below). In general, these model structures are different because of complexes of projective modules in negative degrees. The classical example is given by the two complexes of modules over $k[\epsilon]/\epsilon^2$
$$\begin{tikzcd}[column sep=1pc, row sep=0.6pc]
\dots\arrow[r] & 0\arrow[r] & 0\arrow[r] & k[\epsilon]/\epsilon^2\arrow[r, "\epsilon"] & k[\epsilon]/\epsilon^2\arrow[r, "\epsilon"] &  k[\epsilon]/\epsilon^2\arrow[r, "\epsilon"] & \dots\\
\dots\arrow[r, "\epsilon"] & k[\epsilon]/\epsilon^2\arrow[r, "\epsilon"] & k[\epsilon]/\epsilon^2\arrow[r, "\epsilon"] & k[\epsilon]/\epsilon^2\arrow[r] & 0\arrow[r] &  0\arrow[r] & \dots
\end{tikzcd}$$
These are quasi-isomorphic, but not tamely equivalent: indeed, induction along $k[\epsilon]/\epsilon^2\rt k$ defines a left Quillen functor for the tame model structure on $\Ch_{k[\epsilon]}$, but sends the above two (tamely cofibrant) complexes to complexes of $k$-vector spaces that are not quasi-isomorphic.
\end{example}
\begin{example}[Divided orbits]\label{ex:chain divided orbits}
Let $R$ be a coherent ring, $G$ a finite group and let $\mc{F}\subseteq \Orb_G$ contain only the trivial subgroup. Taking $G$-orbits gives a left Quillen functor $(-)_G\colon \Ch_{R[G]}^{\mc{F}\mm{-tame}}\rt \Ch_{R}^{\mm{tame}}$. It will follow from Remark \ref{rem:derived functor} that the left derived functor models the divided orbits functor $(-)_{dG}$ of Example \ref{dG1}.

Concretely, note that on complexes of $R[G]$-modules that are projectively cofibrant, the left derived functor simply computes the homotopy orbits. This is in particular the case for a \emph{bounded below} chain complex of projective $R[G]$-modules. However, a \emph{bounded above} complex $X$ of finitely generated projective $R[G]$-modules need not be projectively cofibrant. Instead it is \emph{fibrant} in the the model structure on $G$-objects in $\Ch_R$ given as follows: considering the projective model structure on $\Ch_R$ a cofibration (resp. weak equivalence) of $G$-objects is a cofibration (resp. weak equivalence) on the underlying object in $\Ch_R$. Consequently, its $G$-orbits coincide with its homotopy fixed points:
$$\begin{tikzcd}
X_G\arrow[r, "\mm{Nm}", "\cong"{swap}] & X^G\arrow[r, "\sim"] & X^{hG}.
\end{tikzcd}$$
For example, let $C_*(EG; R)\stackrel{\simeq}{\rt} R$ be the standard resolution of the trivial module by finite free $R[G]$-modules. Then both $C_*(EG; R)$ and the $R$-linear dual $C^*(EG; R)$ are quasi-isomorphic to $R$, but the left Quillen functor $(-)_G$ sends $C_*(EG; R)$ to the group homology and $C^*(EG; R)$ to the group cohomology of $G$. In particular, the composite quasi-isomorphism $C_*(EG; R)\rt R\rt C^*(EG; R)$ is not a tame weak equivalence. 
\end{example}
\begin{example}[Derived orbits and fixed points]\label{ex:chain derived fixed points}
In the case where $\mc{F}=\Orb_G$, taking $G$-orbits gives a left Quillen functor $(-)_G\colon \Ch_{R[G]}^{\mc{F}\mm{-tame}}\rt \Ch_{R}^{\mm{tame}}$. Remark \ref{rem:derived functor} will show that the left derived functor models the derived orbits functor $(-)_{G}$ of Example \ref{ex:derived fixed points}.

In addition, consider the functor $(-)^G\colon \Ch_{R[G]}\rt \Ch_{R}$. This functor is not a left adjoint, but it does have a left derived functor: indeed, it preserves $\mc{F}$-tame cofibrations and $\mc{F}$-tame trivial cofibrations and hence restricts to a functor between cofibrant objects sending $\mc{F}$-tame weak equivalences to tame weak equivalences. Since $(-)^G$ preserves pushouts along cofibrations (without differentials, cofibrations are summand inclusions) and infinite direct sums, the associated functor $\LL(-)^G$ of stable $\infty$-categories preserves colimits. We will later identify $\LL(-)^G$ with the derived fixed points (Remark \ref{rem:derived functor}).
\end{example}
\begin{notation}
We will denote the stable $\infty$-categories associated to the model categories of Proposition \ref{prop:chain model structures on G-reps} by $\cat{D}_\mc{F}(R[G])$ and $\cat{D}_\mc{F}^\mm{tame}(R[G])$. Because the $\mc{F}$-projective model structure is a right Bousfield localisation of the $\mc{F}$-tame model structure, there is a fully faithful left adjoint $\cat{D}_\mc{F}(R[G])\hookrightarrow \cat{D}_\mc{F}^\mm{tame}(R[G])$.
\end{notation}
Our goal will be to show that $\cat{D}_\mc{F}(R[G])$ and $\cat{D}^{\mm{tame}}_\mc{F}(R[G])$ model the $\infty$-categories of (pro-coherent) modules over the additive category $R[\mc{F}]$ from Definition \ref{def:admissible subgroup}. We start by endowing both $\infty$-categories with a $t$-structure that will correspond to the $t$-structure from Lemma \ref{lem:pro-coh t-structure}.
\begin{lemma}\label{lem:t-structure chain}
Let $R$ be a ring and $X$ a complex of $R[G]$-modules.
\begin{enumerate}
\item Then $X$ is weakly equivalent to a nonnegatively graded complex of $R[G]$-modules in the $\mc{F}$-projective model structure if and only if $\pi_*(X^H)=0$ for each $*<0$ and each $\mc{F}$-admissible subgroup $H<G$.
\item Then $X$ is weakly equivalent to a nonnegatively graded complex of $R[G]$-modules in the $\mc{F}$-tame model structure if and only if $\pi_*\Hom_{R[G]}(T, X)=0$ for each $*<0$ and each nonpositively graded chain complex of finite $\mc{F}$-admissible $G$-representations $T$.
\end{enumerate}
\end{lemma}
\begin{proof}
Note that the two conditions are invariant under ($\mc{F}$-projective, resp.\ $\mc{F}$-tame) weak equivalences, by definition and by Remark \ref{rem:O-tame we detected by finite stuff}. Furthermore, they are clearly satisfied by every nonnegatively graded complex, so that the conditions are indeed necessary.

To see that the conditions are sufficient, we may assume that $X$ is $\mc{F}$-projectively cofibrant (for (1)) or is $\mc{F}$-tamely cofibrant (for (2)). Let us write $\mc{K}$ for the class of $\mc{F}$-tamely cofibrant objects $Y\in \Ch_{R[G]}$ with the property that $\pi_*\Hom_{R[G]}(Y, X)=0$ for $\ast\leq 0$. Because $\Hom_{R[G]}(-, X)$ sends short exact sequences $0\to Y'\to Y\to Y''\to 0$ of $\mc{F}$-tamely cofibrant objects to short exact sequences of chain complexes, it follows that $Y\in\mc{K}$ whenever $Y'\in\mc{K}$ and $Y''\in \mc{K}$. Likewise, suppose that $Y_0\hookrightarrow Y_1\hookrightarrow \dots$ is a (transfinite) sequence of inclusions between objects in $\mc{K}$, such that each $Y_{\alpha+1}/Y_\alpha\in \mc{K}$ and $Y_{\beta}=\colim_{\alpha<\beta} Y_\alpha$ for each limit ordinal. Then $\colim_\alpha Y_\alpha$ is an object in $\mc{K}$ as well.

For (1), we have that $R[G/H][n]\in \mc{K}$ for each admissible $H<G$ and each $n<0$. The small object argument then implies that $\mc{K}$ contains all $\mc{F}$-projectively cofibrant complexes concentrated in negative degrees. In particular, the stupid truncation $X_{\leq -1}$ is contained in $\mc{K}$, so that there exists a null-homotopy $h$ of the inclusion $i\colon X_{\leq -1}\hookrightarrow X$. Such a null-homotopy endows the inclusion $\tau_{\geq 0}X\hookrightarrow X$ with the structure of a deformation retract, using that $X_0$ admits a direct sum decomposition $Z_0(X)\oplus \ker(1-hd)$. In particular, $X$ is chain homotopy equivalent to a nonnegatively graded complex.

For (2), we have that $\mc{K}$ contains all negatively graded chain complexes of finite $\mc{F}$-admissible $G$-representations, so that the argument from \cite[Lemma 8.6]{Nuiten} shows that $\mc{K}$ contains all $\mc{F}$-projectively cofibrant complexes concentrated in negative degrees. In particular, $X_{\leq -1}\in \mc{K}$, so that the same argument shows that $X$ is chain homotopy equivalent to the nonnegatively graded complex $\tau_{\geq 0}X$.
\end{proof}
\begin{lemma}\label{lem:t-structure on complexes of G-reps}
The $\infty$-categories $\cat{D}_\mc{F}(R[G])$ and $\cat{D}^{\mm{tame}}_\mc{F}(R[G])$ both come equipped with a left complete $t$-structure, in which an object is connective if and only if it is weakly equivalent to a chain complex of $R[G]$-modules concentrated in degrees $\geq 0$. Furthermore, the fully faithful functor $\iota\colon \cat{D}_\mc{F}(R[G])\hookrightarrow \cat{D}^{\mm{tame}}_\mc{F}(R[G])$ exhibits the domain as the right completion of the target. 
\end{lemma}
\begin{warning}
It is \emph{not} true (unless only the trivial subgroup is admissible) that an object is coconnective if and only if it is weakly equivalent to a complex in degrees $\leq 0$.
\end{warning}
\begin{proof}
Let us write $\cat{D}_\mc{F}(R[G])_{\geq 0}\subseteq \cat{D}_\mc{F}(R[G])$ and $\cat{D}^{\mm{tame}}_\mc{F}(R[G])_{\geq 0}\subseteq \cat{D}^{\mm{tame}}_\mc{F}(R[G])$ for the full subcategories of connective objects as defined in \Cref{lem:t-structure on complexes of G-reps}. 
Note that $\cat{D}_\mc{F}(R[G])_{\geq 0}$ is the smallest subcategory of $\cat{D}_\mc{F}(R[G])$ that is closed under colimits and contains the objects $R[G/H]$ for each admissible $H<G$. It is also closed under extensions in $\cat{D}_\mc{F}(R[G])$: this follows from the characterisation of the connective objects from Lemma \ref{lem:t-structure chain}, because each short exact sequence $X'\to X\to X''$ of complexes of $R[G]$-modules with $X\to X''$ a fibration induces a short exact sequence as well. Using \cite[Proposition 1.4.4.11]{lurie2014higher}, we then find that $\cat{D}_\mc{F}(R[G])_{\geq 0}$ is the subcategory of connective objects for a certain $t$-structure. Unravelling the definitions, an object $X$ is (co)connective in this $t$-structure if and only if $X^H$ is (co)connective for each admissible $H<G$. This implies that the $t$-structure is left complete.

For the tame case, we note that the fully faithful inclusion $\iota\colon \cat{D}_\mc{F}(R[G])\hookrightarrow \cat{D}^{\mm{tame}}_\mc{F}(R[G])$ restricts to an equivalence $\cat{D}_\mc{F}(R[G])_{\geq 0}\simeq \cat{D}^{\mm{tame}}_\mc{F}(R[G])_{\geq 0}$. This follows from the fact that an $\mc{F}$-projective weak equivalence $f\colon X\rt Y$ between bounded below chain complexes is also an $\mc{F}$-tame weak equivalence: for any $T$ as in Remark \ref{rem:O-tame we detected by finite stuff}, the map $\Hom(T, X)\rt \Hom(T, Y)$ is isomorphic to the map of bounded below complexes $\big(T^\vee\otimes_R X\big){}^{G}\rt \big(T^\vee\otimes_R Y\big){}^{G}$, which is easily seen to be a quasi-isomorphism by a filtration argument. 

In particular, the connective objects in $\cat{D}^{\mm{tame}}_\mc{F}(R[G])$ are closed under colimits and extensions and hence form the connective part of a $t$-structure. Furthermore, $\iota$ exhibits $\cat{D}_\mc{F}(R[G])$ as its right completion, since it restricts to an equivalence between connective objects. Finally, the connective objects in $\cat{D}_\mc{F}^{\mm{tame}}(R[G])$ are closed under products and every $\infty$-connective object is contractible, since this was already the case in $\cat{D}_\mc{F}(R[G])$. This implies that the $t$-structure is left complete.
\end{proof}
\begin{proposition}\label{prop:chain model for genuine G-reps}
Let $R$ be a ring, $G$ a finite group and $\mc{F}\subseteq \Orb_G$ a full subcategory. Then the natural functor $R[\mc{F}]\rt \cat{D}_\mc{F}(R[G])$ induces a $t$-exact equivalence of stable $\infty$-categories
$$\begin{tikzcd}
F\colon \Mod_{R[\mc{F}]}\arrow[r, "\simeq"] & \cat{D}_\mc{F}(R[G]).
\end{tikzcd}$$
\end{proposition}
\begin{proof}
The functor $R[\mc{F}]\rt \cat{D}_\mc{F}(R[G])$ sends each finite $\mc{F}$-admissible $G$-representation $V$ to itself, viewed as a complex concentrated in degree $0$. Note that each such $V$ is cofibrant in the $\mc{F}$-projective model structure on $\Ch_{R[G]}$. Since $\Hom_{R[G]}(V, -)$ preserves direct sums, it follows that the objects $V$ form a set of compact connective generators for $\cat{D}_\mc{F}(R[G])$.

The universal property of $\Mod_{R[\mc{F}]}$ (Definition \ref{def:module cats}) now gives rise to a sifted-colimit-preserving functor $F$. Since $F$ maps the compact generators $R[\mc{F}]$ of $\Mod_{R[\mc{F}]}$ to compact generators of $\cat{D}_\mc{F}(R[G])$, it is an equivalence. It identifies the $t$-structures because in both categories, an object $Y$ is (co)connective if and only if the spectrum of maps $R[G/H]\rt Y$ is connective for each admissible subgroup $H<G$.
\end{proof}
Next, we will provide a set of compact generators for the $\infty$-category $\cat{D}^\mm{tame}_\mc{F}(R[G])$, following the argument in \cite[Proposition 7.14]{Neeman2008}, \cite[Proposition 8.8]{Nuiten}.
\begin{proposition}\label{prop:compact generators for genuine reps}
Let $R$ be a coherent ring, $G$ a finite group and $\mc{F}\subseteq \Orb_G$. We write $\mc{K}$ for the set of complexes of $R[G]$-modules $Q$ satisfying the following conditions:
\begin{enumerate}
\item $Q$ is a bounded above complex of finite $\mc{F}$-admissible $G$-representations.
\item the $R$-linear dual complex $Q^\vee$ is $m$-coconnective for some $m$: its $H$-fixed points have vanishing homology in degrees $>m$, for all admissible $H<G$.
\end{enumerate}
Then $\mc{K}$ provides a set of compact generators for $\cat{D}^\mm{tame}_\mc{F}(R[G])$.
\end{proposition}
\begin{proof}
We first verify that every object $Q\in \mc{K}$ is compact. To this end, let $Y_\alpha$ be a set of $\mc{F}$-quasifree complexes of $R[G]$-modules and let $Y_\infty=\bigoplus Y_\alpha$ be their direct sum. Furthermore, let $Y^{(n)}_\alpha$ denote the quotient of $Y_\alpha$ by its subcomplex in degree $<n$, so that $Y_\alpha$ is the limit of $Y^{(n)}_\alpha$ as $n\to -\infty$. Now consider the diagram of abelian groups of homotopy classes of maps
\begin{equation}\label{diag:convergence of homotopy classes of maps}\begin{tikzcd}
\bigoplus_\alpha\big[Q, Y_\alpha\big]\arrow[r]\arrow[d, "\phi"{swap}] & \dots\arrow[r] & \bigoplus_\alpha\big[Q, Y^{(n)}_\alpha\big]\arrow[r]\arrow[d, "\phi_{n}"] & \bigoplus_\alpha\big[Q, Y^{(n+1)}_\alpha\big]\arrow[d, "\phi_{n+1}"] \\
\big[Q, Y_\infty\big]\arrow[r] & \dots\arrow[r] & \big[Q, Y^{(n)}_\infty\big]\arrow[r] & \big[Q, Y^{(n+1)}_\infty\big].
\end{tikzcd}\end{equation}
We have to prove that $\phi$ is a bijection. First, observe that for each $n$ and $\alpha$ (allowing $\alpha=\infty$), there is an isomorphism of bounded below complexes
$$
\Hom_{R[G]}\big(Q, Y^{(n)}_\alpha\big)\cong \big(Q^\vee\otimes_R Y^{(n)}_\alpha\big)^G.
$$
In particular, this implies that $[Q, Y^{(n)}_\alpha]=0$ for all $n>0$. Furthermore, the fibre $Z^{(n)}_\alpha$ of $Y^{(n)}_\alpha\rt Y^{(n+1)}_\alpha$ is an $\mc{F}$-admissible $G$-representation, concentrated in a single degree $n$. Since $(Q^\vee)^H$ has vanishing homology in a range $[0, m]$ for all admissible subgroups $H$, the fibre
$$
\big(Q^\vee\otimes_R Z_\alpha^{(n)}\big)^{G}\rt \Hom_{R[G]}\big(Q, Y^{(n)}_\alpha\big)\rt \Hom_{R[G]}\big(Q, Y^{(n+1)}_\alpha\big)
$$
then has homology groups in the range $[n, n+m]$. Consequently, the horizontal towers in \eqref{diag:convergence of homotopy classes of maps} stabilise for very negative $n$ and converge.
It therefore suffices to prove by (descending) induction that each map $\phi_n$ is bijective. This follows because the induced map on mapping fibres identifies with the bijection
$$
\bigoplus_\alpha \big(Q^\vee\otimes_R Z_\alpha^{(n)}\big)^{G}\rt \big(Q^\vee\otimes_R \bigoplus_\alpha Z_\alpha^{(n)}\big)^{G}.
$$
Next, consider the class of objects generated by $\mc{K}$ under colimits and desuspensions. By Remark \ref{rem:O-tame we detected by finite stuff}, it suffices to show that this class contains any bounded above complex $T$ of finite $\mc{F}$-admissible $G$-representations. For such $T$, the $R$-linear dual $T^\vee$ is a nonnegatively graded chain complex of finite $\mc{F}$-admissible $G$-representations. In particular, $T^\vee$ is a cofibrant object with respect to the $\mc{F}$-projective model structure. By Proposition \ref{prop:chain model for genuine G-reps}, $T^\vee$ can be considered as an object in the $\infty$-category $\Mod_{R[\mc{F}]}$; in this sense, it is an almost perfect module over $R[\mc{F}]$.

We will inductively define a chain model for the Postnikov tower of $T^\vee$ with respect to the $t$-structure on $\cat{D}_{\mc{F}}(R[G])$ of Lemma \ref{lem:t-structure on complexes of G-reps}
$$
T^\vee\rt \dots\rt P_n\rt P_{n-1}\rt \dots \rt P_0.
$$
To do this, we proceed as follows: let $F_n$ be the fibre of $\tau_{\leq n}(T^\vee)\rt \tau_{\leq n}(T^\vee)$ in the $\infty$-category $\Mod_{R[\mc{F}]}$. Since $R[\mc{F}]$ is coherent, $F_n$ is an $n$-connective, almost perfect module over $R[\mc{F}]$. This implies that $F_n$ can be modelled at the chain level by a complex $Q_n$ of finite $\mc{F}$-admissible $G$-representations, concentrated in degrees $\geq n$. Finally, one can then model each $\tau_{\leq n}(T^\vee)$ by $P_n=\bigoplus_{i=0}^n Q_i$, with a certain differential.

The upshot of this is the following: each $P_n$ in the above tower is a connective, $n$-coconnective chain complex of $\mc{F}$-admissible $G$-representations. Furthermore, the tower of $P_n$ stabilises in each degree, so that $T^\vee\rt \lim P_n=\bigoplus_{i\geq 0} Q_i$ is a weak equivalence between cofibrant objects and hence a chain homotopy equivalence. Dualizing, we then obtain that $\colim P^\vee_n\rt T^{\vee\vee}\cong T$ is a chain homotopy equivalence as well. Furthermore, each $P_n^\vee\rt P_{n+1}^\vee$ is an $\mc{F}$-tame cofibration between objects in $\mc{K}$, so that the colimit agrees with the homotopy colimit. It follows that $T$ can be realised as a (filtered) homotopy colimit of objects in $\mc{K}$, as desired.
\end{proof}
\begin{corollary}\label{cor:chain models for pro-coh G-reps}
Let $R$ be a coherent ring, $G$ a finite group and $\mc{F}\subseteq \Orb_G$ a full subcategory. Then there are natural equivalences, compatible with $t$-structures
$$\begin{tikzcd}
\Mod_{R[\mc{F}]}\arrow[d, hookrightarrow, "\upiota"{left}]\arrow[r, "\simeq"] & \cat{D}_{\mc{F}}(R[G])\arrow[d, hookrightarrow]\\
\QC^\vee_{R[\mc{F}]}\arrow[r, "\simeq"] & \cat{D}^\mm{tame}_{\mc{F}}(R[G]).
\end{tikzcd}$$
\end{corollary}
\begin{proof}
The top equivalence is Proposition \ref{prop:chain model for genuine G-reps}. For the bottom equivalence, consider the full subcategory $\cat{D}^\mm{tame}_{\mc{F}}(R[G])^\omega$ of compact generators. By Proposition \ref{prop:compact generators for genuine reps}, $R$-linear duality provides a fully faithful functor $\cat{D}^\mm{tame}_{\mc{F}}(R[G])^\omega\hookrightarrow \cat{D}_{\mc{F}}(R[G])^{\op}$. Using Proposition \ref{prop:chain model for genuine G-reps}, its essential image can be identified with the full subcategory $\Coh(R[\mc{F}])^{\op}$ of coherent $R[\mc{F}]$-modules. This induces the desired equivalence $\QC^\vee_{R[\mc{F}]}\simeq \cat{D}^\mm{tame}_{\mc{F}}(R[G])$. 

Unraveling the definitions, the composite $R[\mc{F}]\hookrightarrow \QC^\vee_{R[\mc{F}]}\hookrightarrow \cat{D}^\mm{tame}_{\mc{F}}(R[G])$ is simply the natural inclusion sending a finite $\mc{F}$-admissible $G$-representation to itself, viewed as a complex in degree $0$. This yields the desired commuting square. Since the vertical functors are equivalences on connective objects, it follows that the bottom equivalence identifies connective objects and hence preserves the $t$-structures.
\end{proof}
\begin{remark}\label{rem:derived functor}
Consider a map of coherent rings $f\colon R\rt S$ and a map $\phi\colon G\rt H$ such that induction maps $\mc{F}_G$ to $\mc{F}_H$. This determines a left Quillen functor 
$$\begin{tikzcd}
F\colon \Ch_{R[G]}^{\mc{F}_G\mm{-tame}}\arrow[r] & \Ch_{S[H]}^{\mc{F}_H\mm{-tame}}; & X\arrow[r, mapsto] & S[H]\otimes_{R[G]} X.
\end{tikzcd}$$
In particular, it restricts to a functor $F\colon R[\mc{F}_G]\rt S[\mc{F}_H]$ (viewed as complexes in degree $0$). 
The associated left derived functor $\LL F\colon \cat{D}_{\mc{F}_G}^\mm{tame}(R[G])\rt \cat{D}_{\mc{F}_H}^\mm{tame}(S[H])$ preserves colimits and totalisations of cosimplicial objects in $R[\mc{F}_G]$ (which can simply be computed as total complexes). Under the equivalence of Corollary \ref{cor:chain models for pro-coh G-reps}, this means that $\LL F$ presents the pro-coherent right-left extension of the functor $F\colon R[\mc{F}_G]\rt S[\mc{F}_H]$. For example, the divided orbits and derived orbits functors arise in this way (Example \ref{ex:chain divided orbits}, \ref{ex:chain derived fixed points}). The same argument applies to the derived functor of $G$-fixed points discussed in Example \ref{ex:chain derived fixed points}, even though it is not left Quillen.
\end{remark}
\begin{remark}\label{rem:connectivity in chain terms}
 
Suppose that $R$ is a coherent ring. Corollary \ref{cor:chain models for pro-coh G-reps} then provides the following alternative characterisation of the connective objects with respect to the $t$-structure on $\cat{D}^\mm{tame}_{\mc{F}}(R[G])$ from \Cref{lem:t-structure on complexes of G-reps}: an object $X\in \cat{D}^\mm{tame}_{\mc{F}}(R[G])$ is connective if and only if 
$$
\pi_*\Hom_{R[G]}(Q, X)=0\qquad \text{for all }*<0
$$
for each of the compact generators $Q$ from Proposition \ref{prop:compact generators for genuine reps}. 
\end{remark}
\begin{remark}\label{rem:chain models for daperf}
In the $\mc{F}$-tame model structure on $\Ch_{R[G]}$, geometric realisations of simplicial objects can be computed by taking normalised chains in the simplicial direction and then taking total complexes, using direct sums. Using this, the $\infty$-category $\APerf_{R[\mc{F}]}$ simply arises from the dg-category of bounded below complexes of finite $\mc{F}$-admissible representations. Dually, totalisations of cosimplicial objects can be computed by taking normalised chains and then taking total complexes using direct products. Consequently, the $\infty$-category $\APerf^\vee_{R[\mc{F}]}$ arises from the dg-category of bounded above complexes of finite $\mc{F}$-admissible $G$-representations. Note that $R$-linear duality identifies these two subcategories.
\end{remark}

\subsection{Explicit PD operads and their algebras}
Using the homological algebra from the previous section, we will now provide explicit chain models for PD $\infty$-operads. We begin by giving a description of the $\infty$-category of pro-coherent symmetric sequences.

\subsubsection*{Explicit pro-coherent symmetric sequences}
Consider the model categories of Proposition \ref{prop:chain model structures on G-reps} for all symmetric groups, using only the case where $\mc{F}\subseteq\Orb_{\Sigma_n}$ consists of the trivial subgroup. This yields a model-categorical presentation of the $\infty$-category of pro-coherent symmetric sequences over $R$.
\begin{definition}[The tame model structure on symmetric sequences]\label{sseqtame}
Let $R$ be a ring and let $\mathbf{sSeq}_R := \Ch_{R[\Sigma]}$ denote the category of symmetric sequences of chain complexes of $R$-modules. The \emph{tame model structure} on $\mathbf{sSeq}_R$ is the cofibrantly generated model structure whose fibrations are the surjections and whose cofibrations are injections whose cokernel is given in each arity $r$ by a complex of projective $R[\Sigma_r]$-modules.
\end{definition}
The standard projective model structure on symmetric sequences, whose weak equivalences are the quasi-isomorphisms, is a right Bousfield localisation of the tame model structure. The results from the previous section can now be summarised as follows:
\begin{corollary}
Let $R$ be a commutative ring. Then the underlying $\infty$-category of the projective model structure on symmetric sequences is equivalent to the $\infty$-category $\sSeq_R$ from \Cref{def:sym seq}, i.e.\ $\mathbf{sSeq}_R[W_{\mm{proj}}^{-1}]\simeq \sSeq_R$.

If $R$ is coherent, the fully faithful left adjoint of $\infty$-categories $\mathbf{sSeq}_R[W_{\mm{proj}}^{-1}]\hookrightarrow \mathbf{sSeq}_R[W_{\mm{tame}}^{-1}]$ is naturally equivalent to the fully faithful functor $\upiota\colon \sSeq_R\hookrightarrow \sSeq_R^\vee$. Furthermore, a map between bounded below symmetric sequences is a tame weak equivalence if and only if it is a quasi-isomorphism.
\end{corollary}
Our next goal will be to give a model-categorical description of the various monoidal structures on $\sSeq^\vee_R$, as described in Proposition \ref{prop:pro-coh composition}.
\begin{lemma}\label{lem:chain model tensor product pro-coh}
Let $R$ be a ring. Then the tame and projective model structures on $\mathbf{sSeq}_R$ satisfy the pushout-product axiom with respect to the Day convolution product $\otimes$ and levelwise tensor product $\otimes_{\lev}$ of symmetric sequences of chain complexes.
The induced closed symmetric monoidal structures on the projective localisation $\sSeq_R$ and (if $R$ is coherent) tame localisation $\sSeq^\vee_R$ coincide with those from \Cref{cor:sseq products as extensions} and Proposition \ref{prop:pro-coh composition}.
\end{lemma}
\begin{proof}
We will only treat the tame case, the projective case is proven in the same way. The pushout-product axiom is readily verified for both tensor products. Write $\otimes^{\sLL}$ and $\otimes_{\lev}^{\sLL}$ for the induced closed monoidal structures on $\sSeq^\vee_R$. Using Remark \ref{rem:chain models for daperf}, we see that the restriction of $\otimes^{\sLL}$ (and likewise $\otimes_{\lev}^{\sLL}$) to dually almost perfect objects can be identified with the composite
$$\begin{tikzcd}
(\APerf^\vee_{R[\Sigma]})^{\times 2}\arrow[r, "(-)^\vee", "\simeq"{swap}] &  \big(\APerf_{R[\Sigma]}^{\times 2}\big)^{\op}\arrow[r, "\otimes^{\sLL}"] & \APerf_{R[\Sigma]}^{\op}\arrow[r, "(-)^\vee", "\simeq"{swap}] & \APerf^\vee_{R[\Sigma]}.
\end{tikzcd}$$
Since each step preserves totalisations of cosimplicial objects, Remark \ref{rem:left right extended functors} implies that $\otimes^{\sLL}$ and $\otimes_{\lev}$ are obtained by right-left extension from their restriction to $R[\Sigma]$. The result follows from the fact that both coincide with the usual Day convolution and levelwise tensor product on the full subcategory $R[\Sigma]\hookrightarrow \sSeq^\vee_R$.
\end{proof} 
We want to carry out a similar analysis for the composition product on pro-coherent symmetric sequences. 
\begin{proposition}\label{prop:composition product compatible with tame}
Let $R$ be a ring, let $\circ$ denote the usual composition product on $\mathbf{sSeq}_R$ and let $X$ be a tamely cofibrant symmetric sequence. Then the following assertions hold:
\begin{enumerate}
\item The functor $(-)\circ X\colon \mathbf{sSeq}_R\rt \mathbf{sSeq}_R$ is a left Quillen functor for the tame model structure.
\item The functor $X\circ(-)\colon \mathbf{sSeq}_R\rt \mathbf{sSeq}_R$ preserves tamely cofibrant objects, as well as tame cofibrations and tame weak equivalences between tamely cofibrant objects. 
\item The induced functor of $\infty$-categories $X\circ^{\LL}(-)\colon \mathbf{sSeq}_R[W_{\mm{tame}}^{-1}]\rt \mathbf{sSeq}_R[W_{\mm{tame}}^{-1}]$ preserves sifted colimits.
\end{enumerate}
\end{proposition}
\begin{remark}
If $X$ is tamely cofibrant, then the derived functor of $X\circ (-)$ sends a tamely cofibrant symmetric sequence $Y$ to $\bigoplus_r (X(r)\otimes Y^{\otimes r})_{d\Sigma_r}$, where $(-)_{d\Sigma_r}$ is the divided orbits functor (Example \ref{ex:chain divided orbits}). 
\end{remark}
The third assertion requires some preliminary observations. First, note that it can be reduced to a purely model-categorical assertion as follows: 
\begin{lemma}\label{lem:model categorical sifted-colimit-preserving}
Let $F\colon \cat{M}\rt \cat{N}$ be a functor between combinatorial model categories preserving cofibrant objects and weak equivalences between them. Suppose that $F$ preserves all sifted colimits. Then the induced functor of $\infty$-categories preserves all sifted colimits if the following condition is satisfied: for every category $\cat{I}$ with finite coproducts, the induced functor $\Fun(\cat{I}, \cat{M})\rt \Fun(\cat{I}, \cat{N})$ preserves projectively cofibrant objects.
\end{lemma}
\begin{proof}
A functor between $\infty$-categories preserves sifted colimits if and only if it preserves colimits of diagrams indexed by ordinary categories $\cat{I}$ with finite coproducts \cite[Appendix A]{hnp2019}. Since $\cat{M}$ is a combinatorial model category, every $\cat{I}$-diagram in its associated $\infty$-category can be rectified to an $\cat{I}$-diagram in $\cat{M}$ itself \cite[Proposition 1.3.4.25]{lurie2014higher}. It therefore suffices to verify that $F$ preserves homotopy colimits of $\cat{I}$-diagrams. This follows from the fact that $F$ preserves sifted colimits and projectively cofibrant $\cat{I}$-diagrams.
\end{proof}
\begin{example}\label{ex:iterated tensor preserves cofib}
Let $\cat{M}$ be a combinatorial monoidal model category. If $\cat{I}$ has finite coproducts, then the projective model structure on $\Fun(\cat{I}, \cat{M})$ satisfies the pushout-product axiom for the levelwise tensor product on $\cat{M}$. Consequently, the functor $\Fun(\cat{I}, \cat{M})\rt \Fun(\cat{I}, \cat{M})$ sending $Y\mapsto Y^{\otimes p}$ preserves projectively cofibrant objects, as well as projective cofibrations and weak equivalences between them.
\end{example}
\begin{example}\label{ex:tensoring with tame G-rep}
Suppose that $X\in \mathbf{sSeq}_R$ comes with a $G$-action such that each $X(q)$ is a chain complex of projective $R[G\times \Sigma_q]$-modules. If $Y\rt Z$ is a map of $G$-equivariant symmetric sequences which is a tame cofibration without $G$-action, then $X\otimes_G Y\rt X\otimes_G Z$ is again a tame cofibration of symmetric sequences. In other words, the left adjoint functor $X\otimes_{G} (-)\colon \mathbf{sSeq}_R^{G}\rt \mathbf{sSeq}_R$ preserves tame cofibrations (ignoring the $G$-action in the domain). Consequently, for any category $\cat{I}$ the left adjoint functor $X\otimes_G(-)\colon \Fun(\cat{I}, \mathbf{sSeq}_R^{G})\rt \Fun(\cat{I},\mathbf{sSeq}_R)$ preserves projective cofibrations with respect to the tame model structure on $\mathbf{sSeq}_R$ (ignoring the $G$-action in the domain).
\end{example}
\begin{proof}[Proof (of Proposition \ref{prop:composition product compatible with tame})]
Part (1) is easily verified. For (2), a combination of Examples \ref{ex:iterated tensor preserves cofib} and \ref{ex:tensoring with tame G-rep} shows that each functor $Y\mapsto X(p)\otimes_{\Sigma_p} Y^{\otimes p}$ preserves tamely cofibrant objects and cofibrations between them. Taking the direct sum over $p$ then shows that $X\circ(-)$ preserves tamely cofibrant objects and cofibrations between them. To prove that it preserves tame weak equivalences between tamely cofibrant objects, it suffices to verify that it preserves trivial cofibrations between tamely cofibrant objects. 

Up to retracts, every such trivial cofibration is a transfinite composition of maps $Y\rt Y\oplus Z$, where $Y$ is tamely cofibrant and $Z=R[\Sigma_r][n, n+1]$ is a contractible complex in some arity $r$. Since $X\circ(-)$ preserves transfinite compositions, it suffices to verify that each $
X\circ Y\rt X\circ \big(Y\oplus Z\big)$ is a trivial cofibration. For each $p$, consider the $\Sigma_p$-equivariant symmetric sequence 
$$
L(p)= X((-)+p)\circ Y.
$$
The map $X\circ Y\rt X\circ \big(Y\oplus Z\big)$ is then obtained as a transfinite composition of inclusions whose cokernels are given by the symmetric sequences
$L(p)\otimes_{\Sigma_p} Z^{\otimes p}$.
We have to prove that these cokernels are contractible. 

Because $X$ and $Y$ are tamely cofibrant, $L(p)$ is given in each arity $q$ by a complex of projective $R[\Sigma_p\times \Sigma_q]$-modules. Examples \ref{ex:iterated tensor preserves cofib} and \ref{ex:tensoring with tame G-rep} show that each $L(p)\otimes_{\Sigma_p} Z^{\otimes p}$ is tamely cofibrant. It remains to verify that it is also tamely weakly contractible. To see this, write $L(p)=\lim_n F^{\geq -n}L(p)$ as the limit of its brutal truncations, keeping everything in degrees $\geq -n$. Since $Z$ is bounded above, we then have that
$$
L(p)\otimes_{\Sigma_p} Z^{\otimes p}\cong \lim_n\Big(\big(F^{\geq -n}L(p)\big)\otimes_{\Sigma_p} Z^{\otimes p}\Big).
$$
But now notice that $F^{\geq -n}L(p)$ is a bounded below object that is tamely cofibrant; this implies that it is \emph{projectively} cofibrant as well, so that $F^{\geq -n}L(p)\otimes_{\Sigma_p}(-)$ sends all non-equivariant tame (trivial) cofibrations of symmetric sequences to tame (trivial) cofibrations. This implies that each $F^{\geq -n}L(p)\otimes_{\Sigma_p} Z^{\otimes p}$ is tamely weakly contractible in $\mathbf{sSeq}_R$, so that the (homotopy) limit for $n\to \infty$ is tamely weakly contractible as well.

Having established (1) and (2), it follows that the composition product can be derived in each variable. Assertion (3) is then a consequence of Lemma \ref{lem:model categorical sifted-colimit-preserving} and Examples \ref{ex:iterated tensor preserves cofib} and \ref{ex:tensoring with tame G-rep}.
\end{proof}
Let us also record the following variant of Proposition \ref{prop:composition product compatible with tame} for the projective model structure:
\begin{proposition}\label{prop:composition product compatible with projective}
Suppose that $X$ and $Y$ are two projectively cofibrant symmetric sequences. Then $X\circ Y$ is projectively cofibrant as well. Consequently, the composition product induces a functor of $\infty$-categories $\circ^{\LL}\colon \mathbf{sSeq}_R[W_{\mm{proj}}^{-1}] \times \mathbf{sSeq}_R[W_{\mm{proj}}^{-1}]\rt \mathbf{sSeq}_R[W_{\mm{proj}}^{-1}]$ preserving colimits in the first variable and sifted colimits in the second variable.
\end{proposition}
\begin{proof}
For the second part, recall that the projective model structure is a right Bousfield localisation of the tame model structure. If the composition product preserves projectively cofibrant objects, the full subcategory $\mathbf{sSeq}_R[W_{\mm{proj}}^{-1}] \subseteq \mathbf{sSeq}_R[W_{\mm{tame}}^{-1}]$ is therefore closed under the composition product from Proposition \ref{prop:composition product compatible with tame}. Since this full subcategory is closed under colimits, it follows that $\circ^{\LL}\colon \mathbf{sSeq}_R[W_{\mm{proj}}^{-1}] \times \mathbf{sSeq}_R[W_{\mm{proj}}^{-1}]\rt \mathbf{sSeq}_R[W_{\mm{proj}}^{-1}]$ preserves colimits in the first variable and sifted colimits in the second variable.

To see that $X\circ Y$ is projectively cofibrant, let us assume that $Y$ is projectively cofibrant and write $\mc{K}$ for the class of tamely cofibrant objects $X$ such that $X\circ Y$ is projectively cofibrant. Since the projective model structure is a right Bousfield localisation of the tame model structure, a tamely cofibrant object is projectively cofibrant if and only if it is tamely equivalent to a projectively cofibrant object, and the projectively cofibrant objects are closed under homotopy colimits (in the tame model structure). Using that $(-)\circ Y$ preserves homotopy colimits in the tame model structure by Proposition \ref{prop:composition product compatible with tame}, it follows that the class $\mc{K}$ is closed under tame weak equivalences and homotopy colimits (in the tame model structure). Using the small object argument, any projectively cofibrant object can be obtained using homotopy colimits from the projectively cofibrant symmetric sequences $R[\Sigma_r][n]$, for some $r\geq 0$ and $n\in\mathbb{Z}$. It therefore suffices to verify that $R[\Sigma_r][n]\circ Y$ is projectively cofibrant. But $R[\Sigma_r][n]\circ Y\simeq Y^{\otimes r}[n]$ is projectively cofibrant by Lemma \ref{lem:chain model tensor product pro-coh}.
\end{proof}
\begin{theorem}[Chain models for pro-coherent composition]\label{thm:chain model for exotic composition}
Let $R$ be a coherent ring. The composition product on the tame model category $\mathbf{sSeq}_R$ induces a monoidal structure on its $\infty$-categorical localisation. The resulting monoidal $\infty$-category is equivalent to $\sSeq^\vee_R$ with the monoidal structure $\circ$ of Proposition \ref{prop:pro-coh composition}.
\end{theorem}
\begin{proof}
By (1) and (2) of Proposition \ref{prop:composition product compatible with tame}, the composition product restricts to a monoidal product on the full subcategory of tamely cofibrant symmetric sequences, which preserves weak equivalences in each variable. By part (1) and (3) of Proposition \ref{prop:composition product compatible with tame}, the resulting monoidal structure $\circ^{\sLL}$ on the $\infty$-category $\sSeq^\vee$ preserves sifted colimits. 

By Remark \ref{rem:left right extended functors}, $\circ^{\sLL}$ is the right-left extension of its restriction to $R[\Sigma]$ if it preserves totalisations of cosimplicial diagrams in $R[\Sigma]$. This follows from the same argument as in Lemma \ref{lem:chain model tensor product pro-coh}: using Remark \ref{rem:chain models for daperf}, the restriction of $\circ^{\sLL}$ to dually almost perfect objects can be identified with the functor sending $(X, Y)\mapsto \big(X^\vee\circ^{\sLL} Y^\vee\big)^\vee$, where $(-)^\vee$ takes $R$-linear dual symmetric sequences. Since $(-)^\vee$ is an equivalence between dually almost perfect objects and almost perfect objects, it follows that the restriction of $\circ^{\sLL}$ to $\APerf^\vee_{R[\Sigma]}$ preserves totalisations.
The result now follows the fact that $\circ^{\sLL}$ restricts to the usual composition product on $R[\Sigma]$.
\end{proof}
For an arbitrary ring, the same proof applies to the projective model structure, where we obtain a composition product on $\sSeq_R\simeq \mathbf{sSeq}_R[W_{\mm{proj}}^{-1}]$ by Proposition \ref{prop:composition product compatible with projective}. Using Corollary \ref{cor:sseq products as extensions} to describe the composition product from Definition \ref{def:comp prod} as a left-right derived functor, we therefore obtain the following:
\begin{corollary}
Let $R$ be a ring. The composition product on the projective model category $\mathbf{sSeq}_R$ induces a monoidal structure on its $\infty$-categorical localisation. The resulting monoidal $\infty$-category is equivalent to $\sSeq_R$ with the monoidal structure $\circ$ from Definition \ref{def:comp prod}. 
\end{corollary}

\subsubsection*{Rectification of PD operads and their algebras}
Write $\mathbf{sSeq}_R^c$ for the full subcategory of $\mathbf{sSeq}_R$ spanned by those symmetric sequences that are tamely cofibrant. \Cref{thm:chain model for exotic composition} implies that there is a zig-zag of monoidal functors
$$\begin{tikzcd}
\big(\mathbf{sSeq}_R, \circ\big) & 
\big(\mathbf{sSeq}_R^{c}, \circ\big)\arrow[r] \arrow[l] &  \big(\sSeq^\vee_R, \circ\big)
\end{tikzcd}$$
which exhibits the $\infty$-category $\sSeq^\vee_R$ as a monoidal localisation of the category of symmetric sequences of chain complexes of $R$-modules at the tame weak equivalences.
In particular, any dg-operad $\mathbf{P}$ over $R$ defines a PD $\infty$-operad $\PP$, i.e.\ an associative algebra in $\sSeq^{\vee}_R$, and every dg-algebra over such a dg-operad $\mathbf{P}$ defines a pro-coherent algebra over the corresponding PD $\infty$-operad. The goal of this section is to prove that all PD $\infty$-operads and algebras can be rectified in this way, or more precisely, that the tame homotopy theory of dg-operads presents the $\infty$-category of PD $\infty$-operads. 

We begin by describing the tame homotopy theory of dg-operads and their algebras in more detail. As the existence of model structures on categories of algebras is typically a subtle issue \cite{bergermoerdijk2003, batanin2017homotopy}, we will use semi-model structures, which first appeared in work of Hovey \cite{hovey1998monoidal}:
\begin{reminder}[Semi-model categories]
Recall that a (left) semi-model structure on a presentable category $\cat{M}$ consists of classes of weak equivalences, fibrations and cofibrations satisfying the usual axioms of a Quillen model category, with the following exceptions (see e.g.\ \cite[Ch.\ 12]{fresse2009modules} or \cite{SpitzweckOperads}): fibrations are only required to have the right lifting property against trivial cofibrations \emph{with cofibrant domain}, and only maps with cofibrant domain factor into a trivial cofibration, followed by a fibration. 

We will only deal with cofibrantly generated semi-model structures, where the generating trivial cofibrations have cofibrant domains. Essentially all model categorical results have an obvious analogue in this setting. In fact, all such `tractable' semi-model structures are Quillen equivalent to combinatorial model categories (by a version of Dugger's theorem \cite{dugger2001combinatorial}); one can use this to carry over any result that is invariant under Quillen equivalence. Notably, for any small category $\cat{I}$ there an equivalence of $\infty$-categories $\Fun(\cat{I}, \cat{M})[W^{-1}]\rt \Fun(\cat{I}, \cat{M}[W^{-1}]\big)$ \cite[Proposition 1.3.4.25]{lurie2014higher} and Lemma \ref{lem:model categorical sifted-colimit-preserving} applies in this setting as well.
\end{reminder}
\begin{proposition}\label{prop:model cats of dg-operads}
The following categories carry cofibrantly generated semi-model structures whose weak equivalences and fibrations are (pointwise) tame weak equivalences and fibrations on the underlying objects:
\begin{enumerate}
\item the category $\dgOp_{R}$ of $R$-linear dg-operads.
\item the category $\dgAlg_{\mathbf{P}}$ of $R$-linear dg-algebras over a dg-operad $\mathbf{P}$ which is tamely \emph{$\Sigma$-cofibrant}, i.e.\ whose underlying symmetric sequence is tamely cofibrant.

\item for any category $\cat{I}$, the category of $\cat{I}$-diagrams in $\dgOp_R$ or $\dgAlg_{\mathbf{P}}$.
\end{enumerate}
\end{proposition}
For the analogous result using the projective model structure instead of the tame one, see for example \cite[Theorem 12.2.A, Theorem 12.3.A]{fresse2009modules} and \cite{SpitzweckOperads}.
\begin{proof}
Part (3) is formal. The existence of the (cofibrantly generated) semi-model structures (1) and (2) follows from the transfer theorem for semi-model structures  \cite[Theorem 12.1.4]{fresse2009modules}: one has to verify that for any map with a cofibrant domain that is a pushout of a generating trivial cofibration, the map between the underlying symmetric sequences or complexes is also a trivial cofibration. 

For (1), this means that $\mathbf{P}\rt \mathbf{P}\amalg \mm{Free}_{\Op}(X)$ is a trivial cofibration of symmetric sequences whenever $\mathbf{P}$ is cofibrant and $X$ is tamely cofibrant and contractible. Using the small object argument to write $\mathbf{P}$ as the retract of an iterated pushout of cell attachments, it will suffice to verify this assertion in the case where $\mathbf{P}=\mm{Free}_{\Op}(Y)$ is the free dg-operad on a tamely cofibrant symmetric sequence. In this case, note that the free operad $\mm{Free}_{\Op}(Y)=\colim_n T^{(n)}(Y)$ can be written as the colimit of the sequence of maps (cf.\  \Cref{thm:free algebra})
\begin{equation}\label{eq:free operad}\begin{tikzcd}[column sep=4pc]
i_n\colon T^{(n-1)}(Y)=1\oplus \big(Y\circ T^{(n-2)}(Y)\big)\arrow[r, "\mm{id}\oplus (Y\circ i_{n-1})"] & 1\oplus \big(Y\circ T^{(n-1)}(Y)\big)=T^{(n)}(Y).
\end{tikzcd}\end{equation}
The result then follows by induction, using that the trivial cofibration $Y\rt Y\oplus X$ induces trivial cofibrations on iterated composition products (Proposition \ref{prop:composition product compatible with tame}).

For (2), we need to verify that $A\rt A\amalg (\mathbf{P}\circ X)$ is a trivial cofibration of complexes whenever $A$ is a cofibrant algebra and $X$ is a cofibrant contractible complex. Again, one can use the small object argument to reduce to $A=\mathbf{P}\circ Y$ being free on a complex of projective $R$-modules. Then $\mathbf{P}\circ Y\rt \mathbf{P}\circ (Y\oplus X)$ is a trivial cofibration of complexes by Proposition \ref{prop:composition product compatible with tame}.
\end{proof}
\begin{proposition}\label{prop:forgetting preserves cofibrancy}
Let $\cat{I}$ be a small category with finite coproducts. Then the forgetful functors
$$
\Fun\big(\cat{I}, \dgOp_R\big)\rt \Fun\big(\cat{I}, \mathbf{sSeq}_R\big)\qquad\qquad \Fun\big(\cat{I}, \dgAlg_{\mathbf{P}}\big)\rt \Fun\big(\cat{I}, \Ch_{R}\big)
$$
preserve cofibrations with cofibrant domain (for the semi-model structures as in Proposition \ref{prop:model cats of dg-operads}).
\end{proposition}
\begin{proof}
In the operad case,  say that a map $f\colon \mathbf{P}\rt \mathbf{Q}$ in $\Fun(\cat{I}, \dgOp_R)$ is \emph{good} if it is a cofibration and for each cofibrant $\cat{I}$-diagram of symmetric sequences $X$:
\begin{enumerate}
\item $\mathbf{P}\amalg \mm{Free}_{\Op}(X)$ is a cofibrant $\cat{I}$-diagram of symmetric sequences.
\item $\mathbf{P}\amalg \mm{Free}_{\Op}(X)\rt \mathbf{Q}\amalg \mm{Free}_{\Op}(X)$ is a cofibration of $\cat{I}$-diagrams of symmetric sequences.
\end{enumerate}
We have to verify that every cofibration with cofibrant domain is good. To see this, note that good maps are closed under transfinite compositions and retracts. Furthermore, consider a map $\mathbf{P}\rt \mathbf{P}\amalg_{\mm{Free}_{\Op}(M)}\mm{Free}_{\Op}(N)=\mathbf{Q}$, where $\mathbf{P}$ satisfies condition (1) above and $M\rt N$ is a cofibration of $\cat{I}$-diagrams of symmetric sequences. For every $X$, there is then a sequence of monomorphisms 
$$
\mathbf{P}\amalg \mm{Free}_{\Op}(X)=F(0)\hookrightarrow F(1)\hookrightarrow F(2)\hookrightarrow \dots \rt \colim F(n)=\mathbf{Q}\amalg\mm{Free}_{\Op}(X)
$$
whose associated graded can be identified with $\mathbf{P}\amalg \mm{Free}_{\Op}(X\oplus N/M)$, with grading given by word length in $N/M$ (cf.\ \cite[Section 5]{bergermoerdijk2003}). By assumption (1), the associated graded is cofibrant, so that the above sequence consists of cofibrations and condition (2) is verified as well. Consequently, every cofibration whose domain satisfies (1) is good.

Finally, note that the initial operad satisfies condition (1), so that all cofibrations with cofibrant domain are good. Indeed, this follows from the formula for the free operad $\mm{Free}_{\Op}(Y)$ as the colimit over a sequence of maps $i_n\colon T^{(n-1)}(Y)\rt T^{(n)}(Y)$ as in \eqref{eq:free operad}. Proposition \ref{prop:composition product compatible with tame} (or its proof) then shows that each of these maps is a cofibration between cofibrant $\cat{I}$-diagrams of symmetric sequences, so that the colimit of the sequence is a cofibrant $\cat{I}$-diagram as well. 

In the case of algebras over a dg-operad $\mathbf{P}$ whose underlying symmetric sequence is tamely cofibrant, we proceed in exactly the same way, using that for any $\mathbf{P}$-algebra $A$ and cofibration $M\rt N$, there is a filtration on $A\amalg_{\mathbf{P}\circ M} \mathbf{P}\circ N$ with associated graded $A\amalg \mathbf{P}(N/M)$. In the last step, one has to prove that $\mathbf{P}\circ (-)$ preserves cofibrations of $\cat{I}$-diagrams of complexes of $R$-modules; this follows from Proposition \ref{prop:composition product compatible with tame}.
\end{proof}
\begin{theorem}[Rectification of PD $\infty$-operads]\label{thm:chain models for PD operads}
Let $R$ be a coherent ring. Then the underlying $\infty$-category of the tame semi-model structure on dg-operads over $R$ is equivalent to the $\infty$-category $\Op^{\pd}_R$ of PD $\infty$-operads over $R$. More precisely, there is a commuting square
$$\begin{tikzcd}
\dgOp_R[W_{\mm{tame}}^{-1}]\arrow[d]\arrow[r, "\simeq"] & \Op^{\pd}_R\arrow[d]\\
\mathbf{sSeq}_R[W_{\mm{tame}}^{-1}]\arrow[r, "\simeq"] & \sSeq^\vee_R.
\end{tikzcd}$$
\end{theorem}
\begin{proof}
\Cref{thm:chain model for exotic composition} exhibits $\sSeq^\vee_R$ as the monoidal localisation of $\mathbf{sSeq}_R$ at the tame weak equivalences, with respect to the composition product. This gives rise to the above square. To see that the functor $\Phi\colon \dgOp_R[W_{\mm{tame}}^{-1}]\rt \Op^{\pd}_R$ is an equivalence, notice that both vertical functors are monadic right adjoints: for the left functor, this follows from Lemma \ref{lem:model categorical sifted-colimit-preserving} and Proposition \ref{prop:forgetting preserves cofibrancy} and for the right functor, this follows from \Cref{thm:free algebra}. It follows that $\Phi$ is a right adjoint detecting equivalences; to see that it is an equivalence, it suffices to show that it induces an equivalence between the two monads. 

By \Cref{thm:free algebra}, the left monad is the (left) derived functor of the functor sending a symmetric sequence $X$ of chain complexes to the free operad, given by the colimit of the sequence \eqref{eq:free operad}. On the other hand, \Cref{thm:free algebra} shows that the right monad takes the free algebra with respect to $\circ$, which is given by the same construction \eqref{eq:free operad} at the $\infty$-categorical level. \Cref{thm:chain model for exotic composition} then implies that $\Phi$ induces an equivalence between these two monads.
\end{proof}
\begin{remark}[Chain models for $R$-linear $\infty$-operads]\label{rem:chain models simplicial}
The category $\dgOp_R$ also admits the standard (projective) semi-model structure whose weak equivalences are the quasi-isomorphisms. This is a right Bousfield localisation of the tame model structure, whose associated $\infty$-category is equivalent to the full subcategory $\Op_R\subseteq \Op_R^{\pd}$ spanned by the $R$-linear $\infty$-operads as in Definition \ref{opcoop}. The standard semi-model structure on $\dgOp_R$ also models the $\infty$-category $\Op_R$ of $R$-linear $\infty$-operads when $R$ is not coherent: indeed, the proofs of Proposition \ref{prop:forgetting preserves cofibrancy} and Theorem \ref{thm:chain models for PD operads} carry over verbatim from the tame to the projective setting. 
\end{remark}
\begin{theorem}[Chain models for algebras over PD operads]\label{thm:chain models for PD algebras}
Let $R$ be a coherent ring and $\mathbf{P}$ a dg-operad over $R$ whose underlying symmetric sequence is tamely cofibrant. Then the underlying $\infty$-category of the tame model structure on $\dgAlg_{\mathbf{P}}$ is equivalent to the $\infty$-category $\Alg_{\PP}(\QC^\vee_R)$ of pro-coherent algebras over the associated PD $\infty$-operad $\PP$. In other words, there is a commuting square
$$\begin{tikzcd}
\dgAlg_{\mathbf{P}}[W_\mm{tame}^{-1}]\arrow[d]\arrow[r, "\simeq"] & \Alg_{\PP}(\QC^\vee_R)\arrow[d]\\
\Ch_{R}[W_\mm{tame}^{-1}]\arrow[r, "\simeq"] & \QC^\vee_R.
\end{tikzcd}$$
\end{theorem}
\begin{proof}
The proof is similar to \Cref{thm:chain models for PD operads}: \Cref{thm:chain model for exotic composition} provides the desired square of $\infty$-categories and shows that the bottom arrow is a (monoidal) equivalence. The vertical functors are both monadic right adjoints (for the left, this follows from Lemma \ref{lem:model categorical sifted-colimit-preserving} and Proposition \ref{prop:forgetting preserves cofibrancy}). It suffices to verify that the top functor induces an equivalence between the two monads. The left monad is the derived functor of the functor sending a complex of $R$-modules $M$ to $\mathbf{P}\circ M$. \Cref{thm:chain model for exotic composition} implies that this derived functor is indeed naturally equivalent to the right monad $\PP\circ (-)\colon \QC^\vee_R\rt \QC^\vee_R$.
\end{proof}
\begin{remark}
The equivalence from \Cref{thm:chain models for PD algebras} is natural in the dg-operad $\mathbf{P}$. In particular, this implies that any tame weak equivalence $\mathbf{P}\rt \mathbf{Q}$ between dg-operads whose underlying symmetric sequence is  tamely cofibrant over $R$ induces a Quillen equivalence $\dgAlg_{\mathbf{P}}\rightleftarrows \dgAlg_{\mathbf{Q}}$. 
\end{remark}
\begin{remark}
Let $\ope{P}$ be a tamely $\Sigma$-cofibrant dg-operad. The category $\dgAlg_{\ope{P}}$ also admits a more standard model structure whose weak equivalences are the quasi-isomorphisms. This is a right Bousfield localisation of the tame model structure, whose associated $\infty$-category is equivalent to the full subcategory of $\Alg_{\PP}(\QC^\vee_R)$ generated under colimits by free $\PP$-algebras on all desuspensions of $R$. 

Note that this is typically \emph{not} equivalent to an $\infty$-category of algebras over an operad in $\Mod_R$. In particular, a quasi-isomorphism between two tamely $\Sigma$-cofibrant dg-operads need not induce a Quillen equivalence between their categories of algebras, with the standard semi-model structure.
\end{remark}

\subsection{Explicit Koszul duality}
Finally, we will present the refined Koszul duality functor $\KD^{\pd}\colon \Op_R^{\pd, \op}\rt \Op_R^{\pd}$ by the classical bar dual operad \mbox{defined by \mbox{Ginzburg--Kapranov} \cite{GinzburgKapranov:KDO}.}
\begin{cons}[Chain-level bar construction]\label{chainlevelbar}
Let $R$ be a ring and let $\dgOp^{\aug}_R$ denote the category of augmented dg-operads over $R$. If $\mathbf{P}$ is an augmented dg-operad, we will denote its \emph{chain-level bar construction} by $\dgBar(\mathbf{P})$ \cite[Section 2]{GJ}. 

Recall that $\dgBar(\mathbf{P})$ is a coaugmented dg-cooperad, which can be described explicitly as follows (see e.g.\ \cite[Section 6.5]{LV} for a textbook account). It is the cofree conilpotent cooperad $\mm{Cofree}_{\Coop^\mm{conil}}(\ol{\mathbf{P}}[1])$ generated by the suspension of the augmentation ideal of $\mathbf{P}\rt {1}$, whose underlying symmetric sequence is given by complexes of rooted trees with vertices labelled by elements of $\ol{\mathbf{P}}[1]$. The differential is given by the sum $\partial=\partial_{\mathbf{P}}+\partial_{\mm{Bar}}$, where $\partial_{\mathbf{P}}$ is the differential induced by the differential on $\mathbf{P}$ and the bar differential $\partial_{\mm{Bar}}$ is given by contracting internal edges of trees and multiplying the adjacent elements in $\ol{\mathbf{P}}[1]$. The \emph{chain-level bar dual operad} is defined to be the $R$-linear dual augmented dg-operad $\KD(\mathbf{P})=\dgBar(\mathbf{P})^\vee$.

If $A$ is a dg-algebra over $\mathbf{P}$, then its \emph{chain level bar construction} $\dgBar_{\mathbf{P}}(A)$ is the dg-coalgebra over $\dgBar(\mathbf{P})$ defined as follows (see e.g.\ \cite[Section 11.2]{LV} for a textbook account). Consider the cofree coalgebra $\dgBar(\mathbf{P})\circ A$, whose underlying chain complex consists of trees with vertices labelled by $\ol{\mathbf{P}}[1]$ and leaves labelled by $A$. This is endowed with the differential $\partial=\partial_A+\partial_{\mm{Bar}}$, where $\partial_A$ is the differential induced by the differentials on $A$ and $\dgBar(\mathbf{P})$, while $\partial_\mm{Bar}$ is given by removing leaf vertices and applying the corresponding element of $\ol{\mathbf{P}}[1]$ to the elements in $A$ labelling the leaves. The \emph{chain-level bar dual algebra} is defined to be the $R$-linear dual $\KD_{\mathbf{P}}(A)=\dgBar_{\mathbf{P}}(A)^\vee$, which is an algebra over $\KD(\mathbf{P})$.
\end{cons}

\begin{theorem}[Chain models for Koszul duality]\label{chainkoszul}  
Fix a coherent ring  $R$. Let $\mathbf{P}$ be an augmented dg-operad over $R$ 
with tamely cofibrant 
  underlying symmetric sequence and let $\PP$ denote the corresponding PD $\infty$-operad. 

Then the chain-level dual operad $\KD(\mathbf{P})$ is a model for the Koszul dual PD $\infty$-operad $\KD^{\pd}(\PP)$.
Furthermore, there is a commuting square of $\infty$-categories in which the vertical functors are equivalences
$$\begin{tikzcd}[column sep=3pc]
\dgAlg_{\mathbf{P}}[W^{-1}_\mm{tame}]\arrow[r, "\KD_{\mathbf{P}}"]\arrow[d, "\simeq"{swap}] & \dgAlg_{\KD(\mathbf{P})}[W^{-1}_{\mm{tame}}]^{\op}\arrow[d, "\simeq"]\\
\Alg_{\PP}(\QC^\vee_R)\arrow[r, "\KD^{\pd}_{\PP}"] & \Alg_{\KD^{\pd}(\PP)}(\QC^\vee_R)^{\op}.
\end{tikzcd}$$
\end{theorem}

\begin{proof}
Lemma \ref{lem:chain model tensor product pro-coh} implies that the derived functor of $R$-linear duals on $\mathbf{sSeq}_R$ presents the $R$-linear dual of pro-coherent symmetric sequences. Our main task will therefore be to prove that the chain-level bar construction $\dgBar(\mathbf{P})$ presents the $\infty$-categorical bar construction.

To this end, recall that the \emph{chain level Koszul complex} $K(\mathbf{P})$ is a symmetric sequence of the form $K(\mathbf{P})=\dgBar(\mathbf{P})\circ \mathbf{P}$, whose elements are given by trees with non-leaf vertices labelled by $\ol{\mathbf{P}}[1]$ and leaf vertices labelled by $\mathbf{P}$. The differential contracts (internal) edges and multiplies the labels of the adjacent vertices. Then $K(\mathbf{P})$ becomes a left comodule over $\dgBar(\mathbf{P})$ and a right module over $\mathbf{P}$. Now observe that the category of right $\mathbf{P}$-modules in symmetric sequences admits a model structure in which the weak equivalences and fibrations are detected on the underlying object. The canonical map $\pi\colon K(\mathbf{P})\rt \mb{1}$ is then a tame weak equivalence and exhibits $K(\mathbf{P})$ as a cofibrant resolution of $\mb{1}$ as a right $\mathbf{P}$-module (see e.g.\ \cite[Proposition 4.1.4]{fresse346koszul}, whose proof carries over to the present context). This implies that the natural map
$$\begin{tikzcd}
\mb{1}\circ^{\sLL}_{\mathbf{P}} \mb{1} \simeq K(\mathbf{P})\circ_{\mathbf{P}} \mb{1} \arrow[r, "\mm{coact}"] & \big(\dgBar(\mathbf{P})\circ K(\mathbf{P})\big)\circ_{\mathbf{P}} \mb{1}\simeq \dgBar(\mathbf{P})\circ \big(K(\mathbf{P})\circ_{\PP} \mb{1}\big)\arrow[r, "\mm{id}\otimes \pi"] & \dgBar(\mathbf{P})
\end{tikzcd}$$
is an equivalence, where the first map applies the coaction of $\dgBar(\mathbf{P})$ and the second maps applies the natural augmentation to $\mb{1}$ on the second factor. Lemma \ref{lem:coendomorphism objects} then shows that $\dgBar(\mathbf{P})$ represents the coendomorphism coalgebra of the right $\PP$-module $\mb{1}$; by Proposition \ref{prop:bar construction as coendomorphisms}, this means that $\dgBar(\mathbf{P})$ represents the $\infty$-categorical bar construction $\Barr(\PP)$. 

Furthermore, the chain level Koszul complex $K(\mathbf{P})$ represents the Koszul complex $K(\PP)$ of Constrution \ref{cons:koszul complex}. By \Cref{thm:chain model for exotic composition} and Lemma \ref{lem:chain model tensor product pro-coh}, we thus obtain a commuting diagram of functors
$$\begin{tikzcd}[column sep=2.5pc]
\dgAlg_{\mathbf{P}}^c\arrow[rr, "K(\mathbf{P})\circ_{\mathbf{P}} (-)"]\arrow[d] & & \dgCoalg^c_{\dgBar(\mathbf{P})}\arrow[d]\arrow[r, "(-)^\vee"] & \dgAlg_{\KD(\mathbf{P})}^{\op}\arrow[d]\\
\Alg_{\PP}(\QC^\vee_R)\arrow[rr, "K(\PP)\circ_{\PP}(-)"] & & \Coalg_{\Barr(\PP)}(\QC^\vee_R)\arrow[r, "(-)^\vee"] & \Alg_{\KD^{\pd}(\PP)}(\QC^\vee_R)^{\op}.
\end{tikzcd}$$
Here $\dgAlg_{\mathbf{P}}^c$ and $\dgCoalg^c_{\dgBar(\mathbf{P})}$ denote the categories of dg-algebras and coalgebras whose underlying complex of $R$-modules is tamely cofibrant. The top functors then preserve tame weak equivalences and the two left horizontal functors can be identified with the chain-level bar construction $\dgBar_{\mathbf{P}}(-)$ and the $\infty$-categorical bar construction $\Barr_{\PP}(-)$ of an algebra. Inverting the tame weak equivalences then gives the desired square from the theorem.
\end{proof}

\subsection{Spectral partition Lie algebras}\label{sec:spectral part lie}
We will now describe various chain models for the PD $\infty$-operad $\cat{Lie}^\pi_{R, \mathbb{E}_\infty}$ whose algebras (when $R$ is a field) coincide with the spectral partition Lie algebras from \cite{brantner2019deformation}. Since $\cat{Lie}^\pi_{R, \mathbb{E}_\infty}$ arises as the PD Koszul dual $\infty$-operad of the nonunital $\mathbb{E}_\infty$-operad (Definition \ref{def:spectral partition lie}), we have the following:
\begin{proposition}\label{prop:koszul dual of barratt-eccles}
Let $\ope{E}^\mm{nu}_R$ be a dg-operad modeling the $R$-linear nonunital $\mathbb{E}_\infty$-operad, such that $\ope{E}^\mm{nu}_R(0)=0$, $\ope{E}^\mm{nu}_R(1)=R\cdot 1$ and each $\ope{E}^\mm{nu}_R(r)$ is a complex of finitely generated projective $R[\Sigma_r]$-modules. Then the Koszul dual dg-operad
$$
\KD(\ope{E}^\mm{nu}_R)=\dgBar(\ope{E}^\mm{nu}_R)^\vee
$$
is a cofibrant object for the tame model structure on $\dgOp_R$, which models the spectral partition Lie PD $\infty$-operad $\cat{Lie}^\pi_{R, \mathbb{E}_\infty}$. In particular, there is an equivalence of $\infty$-categories
$$\begin{tikzcd}
\dgAlg_{\KD(\ope{E}^\mm{nu}_R)}[W_\mm{tame}^{-1}]\arrow[r, "\sim"] & \Alg_{\cat{Lie}^\pi_{R, \mathbb{E}_\infty}}(\QC^\vee_R).
\end{tikzcd}$$  
When $R=k$ is a field, this means that the $\infty$-category $\Alg_{\cat{Lie}^\pi_{R, \mathbb{E}_\infty}}(\Mod_k)$ arises as the localisation of $\dgAlg_{\KD(\ope{E}_{k})}$ at the quasi-isomorphisms.
\end{proposition}
For example, one can take $\ope{E}^\mm{nu}_R$ to be the (nonunital) chains on the Barratt--Eccles operad, or the surjections operad from \cite{mcclure2003multivariable, bergerfresse}.
\begin{notation}[Chain level cobar construction]
If $\ope{C}$ is a coaugmented dg-cooperad, we will denote its (chain-level) \emph{cobar construction} by $\Omega(\ope{C})$. Recall that this is the free dg-operad generated by the desuspension $\ol{\ope{C}}[-1]$ of the coaugmentation ideal, with differential $\partial=\partial_{\ope{C}}+\partial_{\mm{cobar}}$; here $\partial_{\ope{C}}$ is the differential induced by the differential on $\ope{C}$ and $\partial_{\mm{cobar}}$ is induced by the partial cocomposition of $\ope{C}$ (see e.g.\ \cite[\S 6.5]{LV} for more details).
\end{notation}
\begin{proof}
By \Cref{chainkoszul}, $\KD(\dgob{E}^\mm{nu}_R)$ is a dg-operad model for the PD $\infty$-operad $\cat{Lie}^\pi_{R, \mathbb{E}_\infty}=\KD(\mathbb{E}_\infty^\mm{nu})$. Furthermore, the conditions on $\ope{E}^\mm{nu}_R$ imply that there is an isomorphism to the chain-level cobar construction
$$
\KD(\dgob{E}_R^\mm{nu})\cong \Omega((\dgob{E}^\mm{nu}_R)^{\vee}).
$$
The cobar construction $\Omega((\dgob{E}^\mm{nu}_R)^{\vee})$ defines a cofibrant object in the tame model structure on $\dgOp_R$: indeed, it can be obtained by a sequence of cell attachments, where in each stage one attaches generators from the tamely cofibrant complex of $R[\Sigma_r]$-modules $(\dgob{E}^\mm{nu}_R)^{\vee}(r)[-1]$ (cf.\ \cite[\S 6]{HinichHomotopy}). In particular, it follows that $\KD(\dgob{E}_R^\mm{nu})$ is tamely cofibrant as a symmetric sequence (\Cref{prop:forgetting preserves cofibrancy}), so that $\dgAlg_{\KD(\ope{E}^\mm{nu}_R)}$ carries a semi-model structure and \Cref{thm:chain models for PD algebras} applies.
\end{proof}
We will now apply Proposition \ref{prop:koszul dual of barratt-eccles} to the surjections operad from \cite{mcclure2003multivariable, bergerfresse} and obtain a  combinatorial presentation of spectral partition Lie algebras. 
\begin{notation}[Nondegenerate sequences]  \label{nondegseq}
Given $r\geq 1$,   a \emph{nondegenerate sequence} in $\ul{r}$ is an (ordered) sequence $\mb{u}=(u_1, \dots, u_{r+d})$ of elements in $\ul{r}=\{1, \dots, r\}$ such that each $1, \dots, r$ appears in the sequence and $u_{\alpha}\neq u_{\alpha+1}$ for all $\alpha$. If $\mb{u}$ does not exhaust all of $\ul{r}$ or if $u_{\alpha}=u_{\alpha+1}$ for some $\alpha$, then $\mb{u}$ is said to be \emph{degenerate}.
\end{notation}

\begin{definition}[Spectral partition $L_\infty$-algebras]\label{explicitone}
Let $R$ be a \mbox{discrete coherent ring.} A \textit{spectral partition $L_\infty$-algebra} is   a chain complex of $R$-modules $\mf{g}$, together with the following algebraic structure: for every $r\geq 2$ and every nondegenerate sequence $\mb{u}=(u_1,\dots, u_{r+d})$, there is an operation
$$\begin{tikzcd}
\{-, \dots, -\}_{\mb{u}}\colon \mf{g}^{\otimes r}\arrow[r] & \mf{g}
\end{tikzcd}$$
of homological degree $-1-d$. Furthermore, these operations satisfy:
\begin{enumerate}[label=(\alph*)]
\item \emph{Equivariance.} For every $\sigma\in \Sigma_r$, let $\sigma(\mb{u})=\big(\sigma(u_1), \dots, \sigma(u_{r+d})\big)$. Then
$$
\{x_1, \dots, x_r\}_{\sigma(\mb{u})} = \pm_{(\sigma, x)}\{x_{\sigma^{-1}(1)}, \dots, x_{\sigma^{-1}(r)}\}_{\mb{u}}
$$
where $\pm_{(\sigma, x)}$ is the Koszul sign associated to the permutation $\sigma$ of $x_1, \dots, x_r$.
\item\label{it:differential} \emph{Differential.} For each nondegenerate sequence $\mb{u} =(u_1,\dots, u_{r+d}) $ in $\ul{r}$ and each tuple $x_1, \dots, x_r\in \mf{g}$, we have
\begin{align*}
\partial\{x_1, \dots, x_r\}_{\mb{u}} &= \sum_{i=1}^r (-1)^{|x_1|+\dots +|x_{i-1}|} \{x_1, \dots, \partial(x_i), \dots, x_r\}_{\mb{u}}\\
&+ \sum_{\alpha=1}^{r+d+1}\sum_{\substack{v=1\vspace{1pt} \\ v \neq u_{\alpha - 1}, u_{\alpha}}}^r \pm_{(\mb{u}_+, \alpha)} \{x_1, \dots, x_r\}_{\mb{u}_+=(u_1, \dots, u_{\alpha-1}, v, u_{\alpha}, \dots, u_{r+d})}\\
&+\sum_{k=1}^{r-1}\sum_{\sigma\in \mm{UnSh}_{\mathbf{u}}(k, r-k)} \hspace{-15pt}\pm_{(\sigma, x)}\pmc \big\{\{x_{\sigma(1)}, \dots, x_{\sigma(k)}\}_{\mb{v}(k, \sigma)}, x_{\sigma(k+1)}, \dots, x_{\sigma(r)}\big\}{}_{\mb{w}(k, \sigma)}
\end{align*}
The sign $\pm_{(\mb{u}_+, \alpha)}$ is  associated to the element $v$ in $\mb{u}_+$ as  in    \Cref{sgn:differential sur}.

In the third row, we sum  over the set $ \mm{UnSh}_{\mathbf{u}}(k, r-k)$ of $(k, r-k)$-unshuffles $\sigma$ which are compatible with $\mathbf{u}$, in the following sense: if we 
decompose the subsequence of $\mb{u}$ consisting of all $u_i\in \{\sigma(1), \dots, \sigma(k)\}$ into intervals 
$$
\mb{u}_1=\big(u_{\alpha(1)}, u_{\alpha(1)+1}, \dots, u_{\alpha(1)+\beta(1)}\big),\qquad\dots, \qquad \mb{u}_n=\big(u_{\alpha(n)}, u_{\alpha(n)+1}, \dots, u_{\alpha(n)+\beta(n)}\big)
$$
separated  in $\mb{u}$  by elements in $\{\sigma(k+1), \dots, \sigma(r)\}$, thenwe have \vspace{3pt}  \mbox{$u_{\alpha(i)+\beta(i)}= u_{\alpha(i+1)}$ for all $i$.}\\
We then define $\mb{v}(k, \sigma)$ to be the sequence in $\ul{k}$ given by applying $\sigma^{-1}$ to the sequence\vspace{-2pt}
\begin{equation}\label{eq:concatenated sequence}
\big(u_{\alpha(1)}, \dots, u_{\alpha(1)+\beta(1)-1}, u_{\alpha(2)}, \dots, u_{\alpha(i)+\beta(i)-1}, u_{\alpha(i)}, \dots, u_{\alpha(n)+\beta(n)}\big).
\end{equation}
Define $\mb{w}(k, \sigma)$ as the sequence of elements of $\ul{r-k+1}$ obtained from $\mb{u}$ by replacing each $\sigma(k+i)$ (for $i=1, \dots, r-k$) by $1+i$ and replacing each of the intervals $\mb{u}_1, \dots, \mb{u}_n$ by a single copy of $1$.

If either of these sequences is degenerate or of length $1$, the corresponding term is zero. Otherwise,  
the sign $\pmc$ is dictated by \Cref{rem:koszul sign for caesuras}, as follows. There is a unique (non-ordered) bijection $$\begin{tikzcd}
\phi\colon \Caes{\mb{w}(k, \sigma)}\star \Caes{\mb{v}(k, \sigma)}\arrow[r] & \Caes{\mb{u}}
\end{tikzcd}$$ 
between the the concatenation of the linear orders of caesuras (Definition \ref{def:caesura}) in $\mb{w}(k, \sigma)$ and $\mb{v}(k, \sigma)$ and the caesuras in $\mb{u}$, with the following properties: $\phi$ sends a caesura in $\mb{v}(k, \sigma)$ to the corresponding caesura in the subsequence \eqref{eq:concatenated sequence} of $\mb{u}$, and restricts to an order-preserving map on $\Caes{\mb{w}(k, \sigma)}$. Then $\pmc$ is the sign of the bijection $\phi$.
\end{enumerate}
\end{definition}

\begin{theorem}[Chain models for spectral partition Lie algebras I] \label{chaintheoremone}
Inverting tame weak equivalences on the category of spectral partition $L_\infty$-algebras gives the $\infty$-category $\Alg_{{\Lie}_{R, \EE_\infty}^\pi}(\QC^\vee_R).$
In particular, when $R=k$ is a field,  localising spectral partition $L_\infty$-algebras  at  the quasi-isomorphisms gives the $\infty$-category of partition Lie algebras from \cite[Definition 5.32]{brantner2019deformation}.
\end{theorem}

\begin{proof}
Write $\mathbf{C}$ for the cooperad given by the linear dual of the surjections operad described in \cite{bergerfresse}. 
Proposition \ref{prop:koszul dual of barratt-eccles} shows that spectral partition Lie algebras over $R$ can be described by algebras over the cobar construction $\Omega(\mathbf{C})$. Without differential, $\Omega(\mathbf{C})$ is the free operad generated by the symmetric sequence underlying the coaugmentation ideal $\overline{\mathbf{C}}[-1]$: this symmetric sequence is spanned in arity $r$ and degree $-1-d$ by the $\Sigma_r$-set of nondegenerate \mbox{sequences $\mb{u}=(u_1, \dots, u_{r+d})$ in $\ul{r}$.}

The equation in \ref{it:differential} then simply asserts that the action of $\Omega(\mathbf{C})$ by operations $\{-, \dots, -\}_{\mb{u}}$ is compatible with the differential. Indeed, note that the differential of $\mb{u}$ in $\Omega(\mathbf{C})$ takes the form $\partial_{\mathbf{C}}(\mb{u})+\partial_{\mm{cobar}}(\mb{u})$, where the first term is simply the differential of $\mb{u}$ in $\mathbf{C}$ and the second term uses the partial cocomposition of $\mathbf{C}$.
By \cite[\S 1.2.3]{bergerfresse} (see also Appendix \ref{app:surjections cooperad}) $\partial_{\mathbf{C}}(\mb{u})$ is the sum of all sequences $(u_1, \dots, u_{\alpha}, v, u_{\alpha+1}, \dots , u_{r+d})$ with a certain sign. This accounts for the second line in the above equation. The third line corresponds to the action by $\partial_{\mm{cobar}}(\mb{u})$, using the description of the cocomposition in $\mathbf{C}$ dual to the formula for the composition in the surjections operad from \cite[\S 1.2.4]{bergerfresse}.
\end{proof}
In the remainder of this section, we will introduce another model for the spectral partition Lie PD $\infty$-operad, which is  smaller than the Koszul dual of the Barratt--Eccles operad and closer to the classical (shifted) Lie operad.
\begin{notation}[Shifted Lie operad]
We will denote by $\ope{Lie}^s_R$ the dg-operad whose algebras are shifted dg-Lie algebras, i.e.\ complexes $\mf{g}$ such that $\mf{g}[-1]$ is a dg-Lie algebra. Likewise, write 
$$
\ope{Lie}^s_{\infty, R}=\KD(\Com^{\mm{nu}}_R)=\Omega(\Com^{\mm{nu}, \vee}_R)
$$
for the dg-operad defining (shifted) $L_\infty$-algebras over $R$. There is a natural map $\ope{Lie}^s_{\infty, R}\rt \ope{Lie}^{s}_R$ taking the quotient by all generating operations in arity $\geq 3$.
\end{notation}
We will model the PD $\infty$-operad $\mm{Lie}_{R,\mathbb E_\infty}^\pi$ by a modification of the standard (shifted) Lie operad that also encorporates divided power operations, using the \emph{PD surjections operad} constructed in detail in Appendix \ref{app:surjections cooperad}:
\begin{definition}[PD surjections operad, see Appendix \ref{app:surjections cooperad}]
Let $\Surj_R$ denote the symmetric sequence underlying the $R$-linear surjections operad from \cite{mcclure2003multivariable, bergerfresse}. Explicitly, $\Surj_R(r)$ is a chain complex given in degree $d$ by the free $R$-module on the set of nondegenerate sequences $\mb{u}=(u_1, \dots, u_{r+d})$ in $\ul{r}$ (cf.\ \Cref{nondegseq}). By \Cref{thm:surjections cooperad}, this symmetric sequence admits the structure of a cooperad such that the canonical map $\Surj_R\rt \ope{coCom}^{\mm{nu}}$ to the cocommutative cooperad is a quasi-isomorphism.

We define the \emph{PD surjections operad} $\Surj^\vee_R$ to be the $R$-linear dual of this cooperad. See \Cref{subsec:App PD surj operad} for more details, including a description of the differential and \mbox{composition in $\Surj^\vee_R$.}
\end{definition}
\begin{remark}
Note the substantial difference between the PD surjections operad $\Surj^\vee_R$ and the standard models for the $\EE_\infty$-operad, such as the Barratt--Eccles operad or the surjections operad: the latter are given in each arity by a projective resolution of the trivial $\Sigma_r$-module, while $\Surj^\vee_R(r)$ is an injective resolution of the trivial $\Sigma_r$-module. In particular, $\Surj^\vee_R$ is a tamely $\Sigma$-cofibrant dg-operad such that $\Surj^\vee_R\circ V\simeq \bigoplus_{r\geq 0} (V^{\otimes r})^{h\Sigma_r}$ for every bounded above complex $V$ of projective $R$-modules (cf.\ Example \ref{ex:barratt-eccles}). This implies that the canonical map $\Com^{\mm{nu}}\rt \Surj^\vee_R$ cannot be a tame weak equivalence, although it is a quasi-isomorphism.
\end{remark}
Recall that given two dg-operads $\ope{P}$ and $\ope{Q}$, their levelwise (or Hadamard) tensor product $\ope{P} \otimes_{\lev} \ope{Q}$ has a natural operad structure. The commutative operad is the unit for this tensor product.
\begin{definition}[Spectral partition Lie dg-operad]\label{Def:PlieOperad}
	Let $R$ be a coherent ring. We define the \emph{spectral partition Lie dg-operad} $\mathbf{Lie}_{R,\EE_\infty}^{\pi}$ to be the tensor product 
	\[
	\mathbf{Lie}_{R,\EE_\infty}^{\pi} = \ope{Lie}^s_R\otimes_{\lev} \Surj^\vee_R.
	\]
\end{definition}
\begin{remark}\label{rem:spectral lie operad}
Each $\mathbf{Lie}_{R,\EE_\infty}^{\pi}(r)$ provides a resolution of the $R[\Sigma_r]$-module $\ope{Lie}^s_R(r)$ by a bounded above complex of finitely generated projective $R[\Sigma_r]$-modules. In particular, $\ope{Lie}^s_R(r)\rt \mathbf{Lie}_{R,\EE_\infty}^{\pi}(r)$ provides an injective resolution of the $\Sigma_r$-action, so that
\begin{align*}
\mathbf{Lie}_{R,\EE_\infty}^{\pi}\circ V &\cong \bigoplus_{r} \big(\mathbf{Lie}_{R,\EE_\infty}^{\pi}(r)\otimes V^{\otimes r}\big)_{\Sigma_r}\\
&\cong \bigoplus_{r} \big(\mathbf{Lie}_{R,\EE_\infty}^{\pi}(r)\otimes V^{\otimes r}\big)^{\Sigma_r}\simeq \bigoplus_r \big(\ope{Lie}^s_{R}(r)\otimes V^{\otimes r}\big)^{h\Sigma_r}
\end{align*}
for any bounded above complex $V$ of projective $R$-modules. In particular, this implies that for any algebra $\mf{g}$ over $\mathbf{Lie}_{R,\EE_\infty}^{\pi}$, the homotopy groups $\pi_*(\mf{g})$ have the structure of a graded restricted Lie algebra.
\end{remark}
\begin{theorem}[Chain models for spectral partition Lie algebras] \label{spectralpliemodel}
The dg-operad $\mathbf{Lie}_{R,\EE_\infty}^{\pi}$ is a tamely $\Sigma$-cofibrant model for the spectral partition PD $\infty$-operad $\mm{Lie}^\pi_{R, \EE_\infty}$. Consequently, $\mathbf{Alg}_{\mathbf{Lie}_{R,\EE_\infty}^\pi}(\mathbf{Ch}_R)$ admits a semi-model structure whose underlying $\infty$-category is equivalent to the $\infty$-category of spectral partition Lie algebras
$$\begin{tikzcd}
\dgAlg_{\ope{Lie}_{R,\EE_\infty}^\pi}[W_{\mm{tame}}^{-1}]\arrow[r, "\sim"] & \Alg_{{{\Lie}_{R,\EE_\infty}^\pi}}(\QC^\vee_R).
\end{tikzcd}$$
\end{theorem}
This follows immediately from Proposition \ref{prop:koszul dual of barratt-eccles} and the following result:

\begin{proposition}
Let $\ope{E}^\mm{nu}_R$ be a dg-operad modeling the nonunital $\mathbb{E}_\infty$-operad as in Proposition \ref{prop:koszul dual of barratt-eccles}, for example the nonunital Barratt--Eccles operad. There exists a commuting diagram of dg-operads
	\[\begin{tikzcd}
		\KD(\Surj^\vee_R)\arrow[d] \arrow["\simeq",r] &
		\dgob{Lie}^s_{\infty, R} \arrow[d]\arrow[ld]\arrow["\simeq",r]&
		\dgob{Lie}^s_R\arrow[d]\\
		\KD(\ope{E}^\mm{nu}_R) \arrow["\simeq"',"f",r]	& \dgob{Lie}^s_{\infty, R} \otimes_{\lev} \Surj^\vee_R \arrow["\simeq"',r] &\mathbf{Lie}_{\EE_\infty, R}^{\pi} 
	\end{tikzcd}
	\]	
in which all horizontal arrows are tame weak equivalences and the vertical arrows are quasi-isomorphisms.	
\end{proposition}
\begin{proof}
We shall start by describing the top row. The first map arises from the quasi-isomorphism of cooperads $\Surj_R\rt \mb{coCom}^\mm{nu}_R$ by taking the cobar construction. Since both of these cooperads consist of projectively cofibrant $R$-modules (ignoring the $\Sigma_r$-action), the induced map between the cobar constructions is a quasi-isomorphism. Furthermore, the map $\dgob{Lie}^s_{\infty, R}\rt \dgob{Lie}^s_{R}$ is a quasi-isomorphism \cite[Theorem 6.8]{fresse346koszul}. Since the top row consists of dg-operads in nonnegative degrees, these two quasi-isomorphisms are also tame weak equivalences.

The right square is obtained by taking the levelwise tensor product with $\ope{Com}^\mm{nu}_R\rt \ope{Sur}^\vee_R$. Note that the map $\dgob{Lie}^s_{\infty, R} \otimes_{\lev} \Surj^\vee_R\rt \mathbf{Lie}_{\EE_\infty, R}^{\pi} $ can be identified in arity $r$ with the map between mapping complexes
$$\begin{tikzcd}
\Hom_R\big(\Surj_R(r), \dgob{Lie}^s_{\infty, R}\big)\arrow[r] & \Hom_R\big(\Surj_R(r), \dgob{Lie}^s_{R}\big).
\end{tikzcd}$$
Since each $\Surj_R(r)$ is tamely cofibrant, each of these maps is a tame weak equivalence.

The map $\dgob{Lie}^s_{\infty, R}=\KD(\Com^{\mm{nu}}_R)\rt \KD(\ope{E}^\mm{nu}_R)$ can be identified with the Koszul dual of quasi-isomorphism of dg-operads $\ope{E}^\mm{nu}_R\rt \Com^{\mm{nu}}_R$.
It then remains to produce the tame weak equivalence $f$ making the triangle commute. For this, it will be convenient to consider the linear dual situation and instead produce a map $\phi$ of conilpotent dg-cooperads
\begin{equation}\label{eq:cooperad lifting}\begin{tikzcd}
& \dgBar(\ope{E}^\mm{nu}_R)\arrow[d]\\
\dgBar(\Com^{\mm{nu}}_R)\otimes_{\mm{lev}} \Surj_R\arrow[r]\arrow[ru, dotted, "\phi"] & \dgBar(\Com^{\mm{nu}}_R)
\end{tikzcd}\end{equation}
where $\Surj_R$ is the cooperad constructed in Appendix \ref{app:surjections cooperad}. Given such a map $\phi$, we simply take $f$ to be its $R$-linear dual. The resulting map $f$ is then indeed a tame weak equivalence: indeed, note that both solid maps in \eqref{eq:cooperad lifting} are quasi-isomorphisms, so that $\phi$ is a quasi-isomorphism as well. Since both the domain and codomain of $\phi$ are projectively cofibrant symmetric sequences, it follows that $\phi$ is a tame weak equivalence as well and its linear dual $f$ remains a tame weak equivalence.

To produce the lift $\phi$, we will proceed by induction: for each $n\geq 0$, let $\ope{E}^\mm{nu, \leq n}_R$ denote the linear quotient of $\ope{E}^\mm{nu}_R$ by all operations of arity $>n$ that are contained in the kernel of the map $\ope{E}^\mm{nu}_R\rt \Com^{\mm{nu}}_R$. These form a tower of dg-operads such that $\ope{E}_R^{\mm{nu}, \leq 1}\cong \Com^{\mm{nu}}_R$ and $\ope{E}^\mm{nu}_R\cong \lim_n \ope{E}_R^{\mm{nu}, \leq n}$. 

Now recall that the chain-level bar construction of an augmented dg-operad $\ope{P}$ is given by the cofree conilpotent cooperad on the suspension $\ol{\ope{P}}[1]$ of the augmentation ideal, together with a certain differential on it (see \Cref{chainlevelbar}). This implies in particular that $\dgBar(\ope{E}^\mm{nu}_R)\cong \lim_n \dgBar(\ope{E}_R^{\mm{nu}, \leq n})$. It therefore suffices to inductively construct a compatible family of maps $$\phi_n\colon \dgBar(\Com^{\mm{nu}}_R)\otimes_{\mm{lev}} \Surj_R\rt \dgBar(\ope{E}_R^{\mm{nu}, \leq n}).$$ The map $\phi_1$ is just the bottom map in Diagram \eqref{eq:cooperad lifting}.
For the inductive step, note that 
$$\begin{tikzcd}
I^{(n+1)}\arrow[r] & \ope{E}_R^{\mm{nu}, \leq n+1}\arrow[r, two heads] & \ope{E}_R^{\mm{nu}, \leq n}
\end{tikzcd}$$
is a square zero extension of dg-operads with kernel $I^{(n+1)}$. This implies that the cooperad $\dgBar(\ope{E}_R^{\mm{nu}, \leq n+1})$ is obtained from $\dgBar(\ope{E}_R^{\mm{nu}, \leq n})$ by the dual of a cell attachment, adding cogenerators in arity $n+1$. More precisely, for each $n$ there is a pullback square of conilpotent dg-cooperads
$$
\begin{tikzcd}
\dgBar(\ope{E}_R^{\mm{nu}, \leq n+1})\arrow[d]\arrow[r] & \mm{Cofree}_{\mm{Coop}^\mm{conil}}(I^{(n+1)}[1, 2])\arrow[d]\\
\dgBar(\ope{E}_R^{\mm{nu}, \leq n})\arrow[r] & \mm{Cofree}_{\mm{Coop}^\mm{conil}}(I^{(n+1)}[2]).
\end{tikzcd}$$
To find an extension of $\phi_n\colon \dgBar(\Com^{\mm{nu}}_R)\otimes_{\mm{lev}} \Surj_R\rt \dgBar(\ope{E}_R^{\mm{nu}, \leq n})$, it then suffices to find a lift in the following diagram of symmetric sequences
$$\begin{tikzcd}
& & I^{(n+1)}[1, 2]\arrow[d]\\
\dgBar(\Com^{\mm{nu}}_R)\otimes_{\mm{lev}} \Surj_R\arrow[r, "\phi_n"]\arrow[rru, dotted] &  \dgBar(\ope{E}_R^{\mm{nu}, \leq n})\arrow[r] & I^{(n+1)}[2]
\end{tikzcd}$$
But now notice that $\dgBar(\Com^{\mm{nu}}_R)\otimes_{\mm{lev}} \Surj_R$ is a (projectively) $\Sigma$-cofibrant symmetric sequence and that $I^{(n+1)}$ is the part of the kernel of the acyclic fibration $\ope{E}^\mm{nu}_R\rt \Com^{\mm{nu}}_R$ concentrated in arity $n+1$. It follows that $I^{(n+1)}[1, 2]\rt I^{(n+1)}[2]$  is an acyclic fibration, so the desired lift exists.
\end{proof}

\begin{remark}
	There cannot exist a model $\ope{P} \stackrel{\simeq}{\rt} \mathbf{Lie}_{\EE_\infty}^{\pi}$ which is bounded below. Indeed, any such model would be quasi-isomorphic and therefore tamely weak equivalent to the Lie operad, but $\mathbf{Lie}_{\EE_\infty}^{\pi}$ is not tamely weak equivalent to the Lie operad.
\end{remark}
The explicit description of the PD surjections operad in \Cref{subsec:App PD surj operad} leads to the following alternative chain level description of spectral partition  Lie algebras:
\begin{corollary}[Explicit spectral partition Lie algebras -- chain model II]\label{explicittwo}
Let $R$ be a coherent ring. Then a spectral partition Lie algebra over $R$ can be described by a chain complex of $R$-modules $\mf{g}$, together with the following algebraic structure: for every $r\geq 2$, every operation $\lambda\in \ope{Lie}^s(r)$, and every nondegenerate sequence $\mb{u}=(u_1,\dots, u_{r+d})$ in $\ul{r}$, there is an operation
$$\begin{tikzcd}
\{-, \dots, -\}_{\lambda, \mb{u}}\colon \mf{g}^{\otimes r}\arrow[r] & \mf{g}
\end{tikzcd}$$
of homological degree $1-r-d$. Furthermore, these operations satisfy:
\begin{enumerate}[label=(\alph*)]
\item \emph{Equivariance.} For every $\sigma\in \Sigma_r$ and all $x_1, \dots, x_r\in \mf{g}$, 
$$
\{x_1, \dots, x_r\}_{\sigma(\lambda), \sigma(\mb{u})} = \pm_{(\sigma, x)} \{x_{\sigma(1)}, \dots, x_{\sigma(r)}\}_{\lambda, \mb{u}}
$$
where $\pm_{(\sigma, x)}$ is the Koszul sign associated to the permutation $\sigma$ of $x_1, \dots, x_r$.
\item \emph{Differential.} For each $\lambda\in \ope{Lie}^s(r)$ and a nondegenerate sequence $\mb{u}$ in $\ul{r}$ and $x_1, \dots, x_r\in \mf{g}$, one has
\begin{align*}
\partial\{x_1, \dots, x_r\}_{\lambda, \mb{u}} &= \sum_{i=1}^r (-1)^{|x_1|+\dots+|x_{i-1}|} \{x_{1}, \dots, \partial(x_i), \dots , x_{r}\}_{\lambda, \mb{u}}\\
&+ \sum_{\alpha=1}^{r+d+1}\sum_{v=1}^r \pm_{(\mb{u}_+, \alpha)} \{x_1, \dots, x_r\}_{\lambda, \mb{u}_+=(u_1, \dots, u_{\alpha-1}, v, u_{\alpha}, \dots, u_{r+d})}.
\end{align*}
Each term where $\mb{u}_+$ is a degenerate sequence is zero and if $\mb{u}_+$ is nondegenerate the sign is as in \Cref{sgn:differential sur}.

\item \emph{Composition.} Let $r, s\geq 2$ and take $\lambda\in \ope{Lie}^s(r)$ and $\mb{u}=(u_1, \dots, u_{r+d})$ a nondegenerate sequence in $\ul{r}$, as well as $\mu\in \ope{Lie}(s)$ and $\mb{v}=(v_1, \dots, v_{s+e})$ a nondegenerate sequence in $\ul{s}$. For each $1\leq k\leq r$, we then have
\begin{align*}
\{x_1, \dots, \{x_k, \dots, x_{k+s-1}\}_{\mu, \mb{v}}, \dots, x_{r+s-1}\}_{\lambda, \mb{u}} &= \sum_{\mb{w}} \pmc \{x_1, \dots, x_{r+s-1}\}_{\lambda\circ_k \mu, \mb{w}}
\end{align*}
where $\lambda\circ_k \mu$ is the partial composition of $\lambda$ and $\mu$ in the Lie operad. Here the sum runs over all nondegenerate sequences $\mb{w}=(w_1, \dots, w_{r+s-1+d+e})$ in $\ul{r+s-1}$ with the following properties:
\begin{itemize}[label={-}]
\item The subsequence of $\mb{w}$ with values in $\{k, \dots, k+s-1\}$ has the form $\big(w_{\alpha(1)}, \dots, w_{\alpha(i+s+e-1)}\big)$ for some $i$, where
$$
w_{\alpha(i)}=v_1+(k-1),\quad  w_{\alpha(i+1)}=v_2+(k-1), \quad \dots, \quad w_{\alpha(i+s+e-1)}=v_{s+e}+(k-1).
$$
\item Consider the sequence $\mb{w}'$ with values in $\{1, \dots, k-1, k, k+s, \dots, r+s-1\}$ obtained from $\mb{w}$ as follows: remove all elements $w_{\alpha(i+1)}, \dots, w_{\alpha(i+s+e-1)}$ appearing above and replace all elements $w_{\alpha(1)}, \dots, w_{\alpha(i)}$ in the sequence above by $k$. Then the resulting sequence $\mb{w}'$ (of length $r+d$) coincides with the sequence $\mb{u}$ under the obvious order-preserving bijection 
$$
\{1, \dots, k-1, k, k+s, \dots, r+s-1\}\cong \{1, \dots, r\}.
$$ 
\end{itemize}
Furthermore, the sign $\pmc$ arises from \Cref{rem:koszul sign for caesuras}, as follows: there is a (non-ordered) bijection $\Caes{\mb{w}}\cong \Caes{\mb{u}}\star \Caes{\mb{v}}$ between the linearly ordered sets of caesuras (Definition \ref{def:caesura}) of the sequence $\mb{w}$ and the concatenation of the linear orders of caesuras in $\mb{u}$ and $\mb{v}$. Then $\pmc$ is the sign of this bijection.
\end{enumerate}
\end{corollary}

\begin{remark}
Suppose that $\ope{P}$ is an $R$-linear dg-operad in arity $\geq 1$ consisting of complexes of projective $R$-modules. Then the levelwise tensor product of $\ope{P}$ with the (nonunital) Barratt--Eccles operad produces a tamely $\Sigma$-cofibrant replacement $\ope{P}\otimes_{\lev}\ope{E}^\mm{nu}\stackrel{\sim}{\rt}\ope{P}$ of $\ope{P}$. 

On the other hand, the tensor product with the PD surjections operad provides a map $\ope{P}\rt \ope{P}\otimes_{\lev} \Surj^\vee_R$ which is a quasi-isomorphism, but generally not a tame weak equivalence. In fact, the same computation as in Remark \ref{rem:spectral lie operad} shows that for any bounded above complex of projective $R$-modules $V$, one has
$$
\big(\ope{P}\otimes_{\lev} \Surj^\vee_R\big)\circ V\simeq \bigoplus_{r\geq 1} \big(\ope{P}(r)\otimes V^{\otimes r}\big)^{h\Sigma_r}.
$$
The dg-operad $\ope{P}\otimes_{\lev} \Surj^\vee_R$ therefore models a PD $\infty$-operad over $R$ whose algebras can be seen as `$\dgob P$-algebras with divided powers'. 
\end{remark}

\newpage

\section{Simplicial-cosimplicial models for derived PD operads}\label{sec:simp-cosimp point-set models}
Given an ordinary ring $R$, recall that the homotopy theory of simplicial commutative rings over $R$ does not have a good description in terms of chain complexes over $R$, unless $R$ is a $\mathbb{Q}$-algebra. More generally, it is complicated to give chain complex models for derived (PD) $\infty$-operads and their algebras. Instead, we will now introduce simplicial-cosimplicial analogues of the model categories studied in the previous section.

\subsection{Simplicial-cosimplicial models for modules and pro-coherent modules}\label{sec:sc-modules}
Throughout, we fix a discrete ring $R$, a finite group $G$ and $\mc{F}\subseteq\Orb_G$. We will start by introducing simplicial-cosimplicial versions of the model categories of complexes of $G$-representations from Proposition \ref{prop:chain model structures on G-reps}. The main idea will be to build these as some sort of resolution model structures. 
\begin{notation}
For any category $\cat{C}$ with limits, restricting along the Yoneda embedding yields an equivalence of categories $\Fun(\Del^{\op}, \cat{C})\simeq \Fun^{\mm{R}}\big(\cat{sSet}^{\op}, \cat{C})$. We can therefore evaluate a simplicial object $X$ in $\cat{C}$ on a simplicial set $K$, and denote the resulting object in $\cat{C}$ by $X(K)$; it can be computed explicitly as the limit of $X$ over the category of simplices of $K$.
\end{notation}
\begin{definition}\label{def:Kan fibration}
Let $f\colon Y\rt X$ be a map of simplicial chain complexes of $R[G]$-modules. We will say that $f$ is an $\mc{F}$-tame (resp.\ $\mc{F}$-projective) \emph{Kan fibration} if it satisfies the following two conditions:
\begin{enumerate}
\item it is a Reedy fibration of simplicial objects, i.e.\ each map of chain complexes
\begin{equation}\label{eq:matching map}
Y(\Delta[n])\rt Y(\partial\Delta[n])\times_{X(\partial\Delta[n])} X(\Delta[n])
\end{equation}
induces surjections on $H$-fixed points for all admissible $H$.
\item for each horn inclusion $\Lambda^i[n]\rt \Delta^n$, the map 
$$
Y(\Delta[n])\rt Y(\Lambda^i[n])\times_{X(\Lambda^i[n])} X(\Delta[n])
$$ 
is a fibration of complexes of $R[G]$-modules whose fibre is connective with respect to the $t$-structures of Lemma \ref{lem:t-structure on complexes of G-reps}. Note that working in the tame or projective setting results in two different connectivity conditions (see Lemma \ref{lem:t-structure chain}).
\end{enumerate}
Likewise, $f$ is said to be an $\mc{F}$-tame (resp.\ $\mc{F}$-projective) \emph{acyclic Kan fibration} if it is a Reedy fibration and each map \eqref{eq:matching map} has a connective fibre.
\end{definition}
If $X$ is a simplicial chain complex, we will write $\Tot_\oplus(X)$ for the total complex of the corresponding bicomplex, using the direct sum.
\begin{lemma}\label{lem:acyclicity for Kan complexes}
Let $f\colon Y\rt X$ be an $\mc{F}$-tame Kan fibration of simplicial chain complexes of $R[G]$-modules. Then $f$ is an $\mc{F}$-tame acyclic Kan fibration if and only if $\Tot_\oplus(Y)\rt \Tot_{\oplus}(X)$ is an $\mc{F}$-tame weak equivalence. The same statement holds in the $\mc{F}$-projective case.
\end{lemma}
\begin{proof}
Since $\Tot_{\oplus}$ is exact, it suffices to verify that a Kan fibrant simplicial chain complex $X$ is acyclic if and only if $\Tot_{\oplus}(X)\simeq 0$ is weakly equivalent to zero. Note that $\Tot_{\oplus}(X)$ is the total complex of the bicomplex $\big[\dots \rightarrow F(2)\rightarrow F(1)\rightarrow X(0)\big]$, where $F(n)$ is the kernel of $X(n)\rt X(\Lambda^0[n])$.

Assume that $X$ is acyclic, so that $X(0)$ and all $F(n)$ are connective by assumption. This means that $\Tot_{\oplus}(X)$ is connective as well and the spectral sequence associated to $\Tot_{\oplus}(X)$ converges to $\pi_*(X)$ and has $E_1$-page $\big[\dots \rightarrow \pi_*F(2)\rightarrow \pi_*F(1)\rightarrow \pi_*X(0)\big]$. The fact that $X$ is an acyclic Kan fibrant object implies that this is exact, i.e.\ the $E_2$-page vanishes. It follows that $\pi_*\Tot_\oplus(X)=0$, and since $\Tot_\oplus(X)$ was connective it follows that $\Tot_\oplus(X)\simeq 0$.

Conversely, suppose that $\Tot_{\oplus}(X)\simeq 0$ and let $Z(n)$ be the fibre of each $X(n)\rt X(\partial\Delta[n])$. We will prove by induction that each $Z(n)$ is connective. To this end, consider the sub-bicomplexes
$$
C^{(n)}= \Big[\dots \rightarrow F(n+1)\stackrel{\partial_{n+1}}{\rt} Z(n)\rt 0\rt \dots \rightarrow 0\Big].
$$
Since $X$ is Reedy fibrant, the map $\partial_{n+1}$ is surjective (on $H$-fixed points for admissible $H<G$), with fibre given by the kernel of $X(n+1)\rt X(\partial\Delta[n+1])$. It follows that there are fibre sequences $C^{(n+1)}\rt C^{(n)}\rt Z(n)[n, n+1]$, so that an inductive argument shows that $\Tot_{\oplus}(C^{(n)})\simeq \Tot_{\oplus}(X)\simeq 0$ for each $n$. Now note that we have a cofibre sequence of complexes
$$
Z(n)[n]\rt \Tot_{\oplus}(C^{(n)})\rt \Tot_{\oplus}\big[\dots \rightarrow F(n+1)\rightarrow 0\rightarrow \dots \rightarrow 0\big].
$$
Since $X$ is Kan fibrant, each $F(k)$ is connective so that the cofibre is $(n+1)$-connective; since the middle term is contractible, $Z(n)$ is connective.
\end{proof}
We    turn to simplicial-cosimplicial modules, which we will also call sc-modules.
\begin{notation}
Write $\cMod_{R[G]}$ and $\scMod_{R[G]}$ for the (ordinary) categories of cosimplicial and simplicial-cosimplicial $R[G]$-modules. By the classical Dold--Kan correspondence, the normalised chains functor identifies these categories with the categories of nonpositively graded chain complexes and second quadrant bicomplexes. Write $\Tot_\oplus\colon \scMod_{R[G]}\rt \Ch_{R[G]}$ for the functor sending an sc-module to the total complex of the associated bicomplex.
\end{notation}
We will say that a map of sc-modules is an $\mc{F}$-tame ($\mc{F}$-projective) Kan fibration if taking normalised chains in the \emph{cosimplicial} direction yields a Kan fibration between simplicial chain complexes in the sense of Definition \ref{def:Kan fibration}, and similarly for acyclic Kan fibrations. 

\begin{theorem}\label{thm:simp-cosimp model structure}
Let $R$ be a ring, $G$ a finite group and $\mc{F}\subseteq \Orb_G$ a full subcategory. Then the category $\scMod_{R[G]}$ can be endowed with the following two cofibrantly generated, simplicial model structures:
\begin{enumerate}
\item the \emph{$\mc{F}$-projective model structure}, whose (trivial) fibrations are the $\mc{F}$-projective (acyclic) Kan fibrations. Furthermore, a map is a weak equivalence if and only if its image under $\Tot_\oplus$ induces quasi-isomorphisms on $H$-fixed points for all admissible $H<G$.

\item the \emph{$\mc{F}$-tame model structure}, whose (trivial) fibrations are the $\mc{F}$-tame (acyclic) Kan fibrations. Furthermore, a map is a weak equivalence if and only if its image under $\Tot_\oplus$ is an $\mc{F}$-tame weak equivalence.
\end{enumerate}
Furthermore, the total complex functor determines a Quillen equivalence
$$\begin{tikzcd}
\Tot_{\oplus}\colon \scMod_{R[G]}\arrow[r, yshift=1ex] & \Ch_{R[G]}\colon \mm{Res}\arrow[l, yshift=-1ex]
\end{tikzcd}$$
between the $\mc{F}$-tame (resp.\ $\mc{F}$-genuine) model structures.
\end{theorem}
\begin{proof}[Proof of \Cref{thm:simp-cosimp model structure}]
Throughout, we will work with nonpositively graded chain complexes of $R[G]$-modules instead of cosimplicial $R[G]$-modules for simplicity; the two are equivalent by the Dold--Kan correspondence. We will start by constructing the two desired model structures on the category $\mathbf{sCh}_{R[G], \leq 0}$ of simplicial diagrams of nonpositively graded chain complexes of $R[G]$-modules. Given a set $K$ of maps of simplicial sets and a set $L$ of maps of chain complexes,  write $K\boxtimes L$ for the set of maps
$$
T_+\wedge M\cup_{S_+\wedge M}S_+\wedge N\rt T_+\wedge N \qquad\qquad S\stackrel{\in K}{\rt} T, \quad M\stackrel{\in L}{\rt} N,
$$
where $\wedge$ is the evident tensoring of $\mathbf{sCh}_{R[G], \leq 0}$ over pointed simplicial sets. Both model structures have sets of generating (trivial) cofibrations of the form
\begin{equation}\label{eq:generating cofs simcos}
\begin{split}
I&=\big\{\partial\Delta[n]\rightarrow \Delta[n]\big\}\boxtimes \big\{P\rightarrow P[0, 1]\}\\
J&=\big\{\Lambda^i[n]\rightarrow \Delta[n]\big\}\boxtimes \big\{P\rightarrow P[0, 1]\} \cup \big\{\partial\Delta[n]\rightarrow \Delta[n]\big\}\boxtimes \big\{0\rightarrow P[0, 1]\big\}
\end{split}
\end{equation}
for a certain set of dg-$R[G]$-modules $P$. In the $\mc{F}$-projective case, we take the set of shifted representations $P=R[G/H][k]$ with $H<G$ admissible and $k<0$. In the $\mc{F}$-tame case, we use the set of complexes of finite $\mc{F}$-admissible $G$-representations concentrated in degrees $\leq -1$. It follows from Lemma \ref{lem:t-structure chain} that a map is an acyclic Kan fibration if and only if it has the right lifting property against $I$, and a Kan fibration if and only if it has the right lifting property against $J$. 

To see that these generating sets determine model structures, we have to verify that iterated pushouts of maps in $J$ are weak equivalences and that a Kan fibration is a weak equivalence if and only if it is an acyclic Kan fibration.   The form of the sets $I$ and $J$ then implies that both of these model structures are simplicial. 
To see that the maps in $J$ are weak equivalences, note that $\Tot_\oplus$ preserves colimits and that it sends each map in $J$ to a pushout of a map $0\rt P[n, n+1]$. For any $\mc{F}$-quasiprojective complex $P$, the map $0\rt P[n, n+1]$ is a trivial cofibration in both the $\mc{F}$-projective and $\mc{F}$-tame model structure on $\Ch_{R[G]}$. Consequently, iterated pushouts of maps in $J$ are weak equivalences in both the $\mc{F}$-projective and $\mc{F}$-tame case.
Lemma \ref{lem:acyclicity for Kan complexes} shows that an $\mc{F}$-tame (resp.\ $\mc{F}$-projective) Kan fibration is acyclic if and only if it is an $\mc{F}$-tame ($\mc{F}$-projective) weak equivalence.
 
It remains to show that $\Tot_{\oplus}$ is part of a Quillen equivalence. To see this, note that its right adjoint $\mm{Res}$ sends a complex $X$ to the simplicial chain complex given in degree $n$ by the degree $\leq 0$ part of $\Hom(C_*(\Delta[n]), X)$. It follows from this that $\mm{Res}$ preserves cofibrant objects. Furthermore, unraveling the definitions shows that $\Tot_{\oplus}\mm{Res}(X)$ can be identified (up to signs) with the complex
$$\begin{tikzcd}[ampersand replacement=\&, column sep=1pc]
\dots \arrow[r] \& \bigoplus\limits_{n\geq 0}X_1\oplus \bigoplus\limits_{n\geq 1}X_{0}\arrow[rr, "\scalebox{0.8}{$\arraycolsep=3pt \left(\begin{array}{cc} \partial & 1 \\ 0 & \partial\end{array}\right)$}"] \& \&  \bigoplus\limits_{n\geq 0}X_0\oplus \bigoplus\limits_{n\geq 1}X_{-1}\arrow[rr, "\scalebox{0.8}{$\arraycolsep=3pt\left(\begin{array}{cc} \partial & 1 \\ 0 & \partial\end{array}\right)$}"] \& \& \bigoplus\limits_{n\geq 0}X_{-1}\oplus \bigoplus\limits_{n\geq 1}X_{-2}\arrow[r] \& \dots
\end{tikzcd}$$
The counit map $\Tot_{\oplus}\mm{Res}(X)\rt X$ is simply the projection onto the zeroth summand, which is an acyclic fibration. We conclude that the derived counit map is an equivalence. 

To see that $\Tot_\oplus$ is part of a Quillen equivalence, it remains to check that it detects equivalences. In the $\mc{F}$-projective case this is obvious. In the $\mc{F}$-tame case, it suffices to verify that every fibrant-cofibrant object $X$ such that $\Tot_\oplus(X)$ is acyclic is itself acyclic. This follows from Lemma \ref{lem:acyclicity for Kan complexes}.
\end{proof}
\begin{remark}\label{rem:generating cofibrations simplicial-cosimplicial}
Let us spell out the sets of generating (trivial) cofibrations \eqref{eq:generating cofs simcos} of the $\mc{F}$-tame model structure more explicitly in simplicial-cosimplicial terms. To this end, let us write $\wedge$ for the tensoring of $\scMod_{R[G]}$ over pointed simplicial sets. Furthermore, let $\tilde{C}^*(\Delta[1])$ denote the cosimplicial $R$-module of reduced cochains on $\Delta[1]$  and let $\tilde{C}^*(S^1)\rt \tilde{C}^*(\Delta[1])$ be the natural restriction map along $\Delta[1]\rt S^1=\Delta[1]/\partial\Delta[1]$. Now consider the following three types of inclusions of simplicial-cosimplicial $R$-modules (with trivial $G$-action):
$$\begin{tikzcd}[row sep=0.2pc]
i\colon \partial\Delta[n]_+\wedge \tilde{C}^*(\Delta[1]) \amalg_{\partial\Delta[n]_+\wedge  \tilde{C}^*(S^1)} \Delta[n]_+\wedge \tilde{C}^*(S^1)\arrow[r] & \Delta[n]_+\wedge \tilde{C}^*(\Delta[1])\\
j_1\colon \Lambda^i[n]_+\wedge \tilde{C}^*(\Delta[1]) \amalg_{\Lambda^i[n]_+\wedge  \tilde{C}^*(S^1)} \Delta[n]_+\wedge  \tilde{C}^*(S^1)\arrow[r] & \Delta[n]_+\wedge \tilde{C}^*(\Delta[1])\\
\hphantom{\Delta[n]_+\wedge \tilde{C}^*(S^1) \amalg_{\partial\Delta[n]_+\wedge  \tilde{C}^*(S^1)}} j_2\colon \partial\Delta[n]_+\wedge \tilde{C}^*(\Delta[1]) \arrow[r] & \Delta[n]_+\wedge \tilde{C}^*(\Delta[1]).
\end{tikzcd}$$
Then the cofibrations are generated by the set of maps of simplicial-cosimplicial $R[G]$-modules of the form $i\otimes_R P$, where $P$ is a cosimplicial diagram of finite $\mc{F}$-admissible $G$-representations. Likewise, the trivial cofibrations are generated by the set of maps $j_1\otimes_R P$ and $j_2\otimes_R P$.
 
\end{remark}
For later purposes, let us record some further properties of the cofibrations in the $\mc{F}$-projective and $\mc{F}$-tame model structures on simplicial-cosimplicial $G$-representations.
\begin{lemma}\label{lem:cof-triv fib in simcos tame}
A map of simplicial-cosimplicial $R[G]$-modules $f\colon X\rt Y$ is:
\begin{enumerate}
\item an $\mc{F}$-projective cofibration if and only if it is a split monomorphism in each simplicial-cosimplicial bidegree and in each simplicial degree $n$, the cosimplicial $R[G]$-module $(Y/X)_n$ corresponds to an $\mc{F}$-projectively cofibrant chain complex under the Dold--Kan correspondence.

\item an $\mc{F}$-tame cofibration if and only if it is a split monomorphism in each simplicial-cosimplicial bidegree and $(Y/X)$ is given in each simplicial-cosimplicial bidegree by the retract of an $\mc{F}$-admissible $G$-representation.
\end{enumerate}
\end{lemma}
\begin{proof}
Let us again identify simplicial-cosimplicial $R[G]$-modules with of nonpositively graded chain complexes of $R[G]$-modules using the Dold--Kan correspondence to $\scMod_{R[G]}\simeq \mathbf{sCh}_{R[G], \leq 0}$. Under this identification, a map that is a split monomorphism in each simplicial-cosimplicial bidegree corresponds to a map that is a split monomorphism in each simplicial-chain bidegree.

Consider the sets of maps $\{P\to P[0, 1]\}$ in $\Ch_{R[G], \leq 0}$ used in the definition of the set of generating cofibrations \eqref{eq:generating cofs simcos} in the proof of Theorem \ref{thm:simp-cosimp model structure}. In the $\mc{F}$-projective case, we took the objects $P$ to be the shifted representations $R[G/H][k]$ with $H<G$ admissible and $k<0$, and in the tame case we took the objects $P$ to be all complexes of finite $\mc{F}$-admissible $G$-representations in degrees $\leq -1$. One readily verifies that the weakly saturated class generated by these sets consists exactly of the $\mc{F}$-projective and $\mc{F}$-tame cofibrations in $\Ch_{R[G], \leq 0}$; In the $\mc{F}$-tame case, this follows from the proof of Proposition \ref{prop:chain model structures on G-reps}. Note that a map in $\Ch_{R[G], \leq 0}$ is an $\mc{F}$-projective cofibration if and only if it is an $\mc{F}$-injection with $\mc{F}$-cofibrant cokernel, and likewise in the $\mc{F}$-tame case.

The weakly saturated class of the set $I$ from \eqref{eq:generating cofs simcos} then consists of the associated class of Reedy cofibrations, i.e.\ of maps $X\rt Y$ in $\mathbf{sCh}_{R[G], \leq 0}$ such that each relative latching map $X_n\amalg_{L_nX} L_nY\rt Y_n$ is an $\mc{F}$-projective (resp.\ $\mc{F}$-tame) cofibration in $\Ch_{R[G], \leq 0}$. Because each $L_nX\rt X_n$ is the inclusion of a direct summand, by the Dold--Kan correspondence, one readily sees that this is equivalent to $X\rt Y$ being an $\mc{F}$-projective (resp.\ $\mc{F}$-tame) cofibration in each simplicial degree. Assertions (1) and (2) follow directly from this.
\end{proof}

\begin{remark}
Consider the full subcategory $\mathbf{sMod}_{R[G]}\hookrightarrow \scMod_{R[G]}$ of sc-modules that are constant in the cosimplicial direction. This carries an induced model structure, whose weak equivalences and fibrations are the maps inducing weak equivalences and Kan fibrations on $H$-fixed points, for every admissible subgroup $H$. The resulting $\infty$-category can be identified with the $\infty$-category $\Mod_{R[\mc{F}], \geq 0}$ of connective modules over $\mc{F}$. In the case where $\mc{F}=\Orb_G$, this can be thought of as an $R$-linear version of Elmendorf's theorem \cite{elmendorf1983systems}.
\end{remark}
\begin{remark}[Geometric realisations and totalisations]\label{rem:realisations and totalisations in simp-cosimp modules}
The model structure on $\scMod_{R[G]}$ is tensored over the Kan--Quillen model structure in the obvious way. The homotopy colimit of a pointwise cofibrant simplicial diagram $X\colon \Del^{\op}\rt \scMod_{R[G]}$ can therefore be computed by the diagonal of $X$ in the simplicial direction. The image of the diagonal under the functor $\Tot_\oplus$ simply computes the total complex (using the direct sum).

The homotopy limit of a cosimplicial diagram $X\colon \Del^{\op}\rt \scMod_{R[G]}$ cannot be computed by the cosimplicial diagonal in general. However, this does hold when $X$ is a diagram of cosimplicial $R[G]$-modules (constant in the simplicial direction). Indeed, in this case the Quillen equivalence $\Tot_\oplus$ sends the cosimplicial diagonal to the total complex with respect to the sum; since all sums involved are finite this coincides the the total complex using the direct product, which computes the homotopy limit in $\Ch_{R[G]}$.
\end{remark}
\begin{remark}\label{rem:finite totalisation is cof}
Consider a cosimplicial diagram $X\colon \Del\rt R[\mc{F}]$ of finite $\mc{F}$-admissible $R[G]$-modules. Then $X$ defines a cofibrant object of $\scMod_{R[G]}$ with respect to the $\mc{F}$-tame model structure. If $X$ is furthermore finite, i.e.\ $n$-coskeletal for some $n$, then $X$ is also cofibrant in the $\mc{F}$-projective model structure. Indeed, $X$ corresponds under the Dold--Kan correspondence to a complex of finite $\mc{F}$-admissible $R[G]$-modules in degrees $[-n, 0]$. Starting in degree $-n$, such a complex can be obtained by cell attachments from the generating $(\mc{F}$-projective) cofibrations used in the proof of \Cref{thm:simp-cosimp model structure}.
\end{remark}

\subsection{Explicit derived operads, PD operads, and their algebras}\label{explicitderivedoperads}
We will now use the homotopy theory described in \Cref{thm:simp-cosimp model structure} to describe explicit simplicial-cosimplicial models for derived $\infty$-operads, derived PD $\infty$-operads, \mbox{and their algebras.}

\subsubsection*{Simplicial-cosimplicial symmetric sequences}
Write $\mathbf{sSeq}_{R}^{\simcos}$ for the category of symmetric sequences of simplicial-cosimplicial $R$-modules. 
This can be endowed with a Day convolution product $\otimes$, a composition product, and restricted composition product
$$
X\circ Y=\bigoplus_r \big(X(r)\otimes Y^{\otimes r}\big)_{\Sigma_r} \qquad\qquad X\bcirc Y=\bigoplus_r \big(X(r)\otimes Y^{\otimes r}\big){}^{\Sigma_r}.
$$
These operations can be computed in each simplicial-cosimplicial degree. In particular, an algebra with respect to the composition product is simply a \emph{simplicial-cosimplicial operad} (or \textit{sc-operad}) over $R$, while an algebra with respect to the restricted composition product is a \emph{simplicial-cosimplicial restricted operad} over $R$ (or \textit{sc-restricted operad}), see e.g.\ \cite{fresse2000homotopy, ikonicoff, cesarothesis}. 

Recall that the norm map $\mm{Nm}\colon X\circ Y\rt X\bcirc Y$ makes the identity a lax monoidal functor and is an equivalence when $Y$ is in arity $\geq 1$ (cf.\ \Cref{not:restricted composition}). In particular, every sc-restricted operad has an underlying sc-operad. This defines an equivalence between sc-restricted operads and sc-operads in arities $\geq 1$.

The category $\mathbf{sSeq}_{R}^{\simcos}$ acts on $\scMod_R$ by both the composition product $\circ$ and the restricted composition product $\bcirc$. If $\ope{P}$ is an sc-operad, a left module over it in $\scMod_R$ is simply a simplicial-cosimplicial $\ope{P}$-algebra. Likewise, if $\ope{P}$ is an sc-restricted operad, a left $\ope{P}$-module in $\scMod_R$ is simply a \emph{restricted} $\ope{P}$-algebra. A restricted $\ope{P}$-algebra $A$ is in particular an algebra over the operad underlying $\ope{P}$.

\begin{notation}
We will write $\mathbf{Op}^{\simcos}_R$ and $\mathbf{Op}^{\simcos,\res}_R$ for the categories of sc-operads and sc-restricted operads over $R$. Furthermore, we denote by $\mathbf{Alg}^{\simcos}_{\ope{P}}$ and $\mathbf{Alg}^{\simcos,\res}_{\ope{P}}$ the categories of algebras and restricted algebras over an sc-operad, respectively sc-restricted operad $\ope{P}$.
\end{notation}
\begin{remark}[Restricted algebras]\label{rem:restricted algebras}
Suppose that $k$ is a field and $\ope{P}$ is a (restricted) operad in arity $\geq 1$ coming from an operad $\ope{S}$ in sets as $\ope{P} =k[\ope{S}]$. In this case, Ikonicoff \cite[\S 3.1]{ikonicoff} has given an explicit description of restricted $\ope{P}$-algebras as $k$-vector spaces with  operations and relations. Indeed, a restricted $\ope{P}$-algebra $A$ is an ordinary $\ope{P}$-algebra (via the norm), and there are additional operations
$$ \gamma_s(a_1,\ldots a_r) \in A$$
for all $s\in \ope{S}(r)$ and all tuples $\mathbf{a} = (a_1, \ldots, a_r) $ in $A$, which satisfy various properties which generalise the axioms of a divided power algebra (we use a slightly more efficient labelling convention than \cite{ikonicoff}).
To define the element $ \gamma_s(a_1,\ldots a_r)$, we simply apply the structure map $(\ope{P}(r) \otimes A^{\otimes r})^{\Sigma_r} \rightarrow A$ to the element $$\sum_{\rho \in \Sigma_r / \Sigma_{\mathbf{a},s}} [\rho(s)] \otimes a_{\rho(1)} \otimes \ldots  \otimes a_{\rho(r)}$$
where $[\rho(s)]$ is the basis element in $ k[ \ope{S}(r)]$ corresponding to $\rho(s)$ and $\Sigma_{\mathbf{a},s} \subset \Sigma_r$ consists of all permutations which fix both $s$ and the tuple  $\mathbf{a}$.
 These operations are compatible with sums, scalar multiplication and composition in the way one might expect from this equation.
In \Cref{derivedPLAexplicit}, we will use this strategy to  make our point set models for derived partition Lie algebras more explicit.
\end{remark}

\subsubsection*{Model structures}
We now apply the Constructions from Section \ref{sec:sc-modules} in the case of the symmetric groups, taking $\mc{F}=\Orb_{\Sigma_r}$ to be the full orbit category. Along the lines of Definition \ref{def:admissible subgroup}, we will say that a symmetric sequence of discrete $R$-modules is \emph{admissible} if it arises as the $R$-linearisation of a symmetric sequence of sets.
\begin{definition}\label{def:model structure simp-cosimp sym seq}
The \emph{tame model structure} on $\scsSeq_R$ is the simplicial, cofibrantly generated model structure whose:
\begin{itemize}
\item cofibrations are the split monomorphisms in each simplicial-cosimplicial degree, with cokernel given by the retract of an admissible symmetric sequence.
\item (trivial) fibrations are (acyclic) Kan fibrations in each arity, for $\mc{F}=\Orb_{\Sigma_r}$.

\item weak equivalences are maps inducing $\mc{F}$-tame weak equivalences upon applying $\Tot_{\oplus}$.
\end{itemize}
This admits a right Bousfield localisation given by the \emph{projective model structure}, whose weak equivalences are maps $X\rt Y$ such that for every $H<\Sigma_r$, the map $X(r)^H\rt Y(r)^H$ induces a quasi-isomorphism on total complexes (using direct sums).
\end{definition}
Proposition \ref{prop:chain model for genuine G-reps}, Corollary \ref{cor:chain models for pro-coh G-reps} and \Cref{thm:simp-cosimp model structure} then have the following consequence:
\begin{corollary}\label{cor:scsSeq models}
Let $R$ be a commutative ring. Then the projective model structure on $\scsSeq_R$ presents the $\infty$-category $\sSeq^{\gen}_{\ulR}$ of derived symmetric sequences, i.e.\ $\scsSeq_R[W_{\mm{proj}}^{-1}]\simeq \sSeq^{\gen}_{\ulR}$.

If $R$ is furthermore coherent, then the fully faithful left adjoint of $\infty$-categories
$$
\scsSeq_R[W_{\mm{proj}}^{-1}]\hookrightarrow \scsSeq_R[W_{\mm{tame}}^{-1}]
$$
is naturally equivalent to the fully faithful inclusion $\upiota\colon \sSeq^{\gen}_{\ulR}\hookrightarrow \sSeq^{\gen, \vee}_{\ulR}$ into the pro-coherent derived symmetric sequences.
\end{corollary}
\begin{remark}[Connective objects] \label{rem:connective simplicial-cosimplicial}
A map between simplicial symmetric sequences (constant in the cosimplicial direction) is a tame weak equivalence if and only if it induces a weak equivalence on all $H$-fixed points; these simplicial symmetric sequences model the (equivalent) connective parts of $\sSeq^{\gen}_{\ulR}$ and $\sSeq^{\gen, \vee}_{\ulR}$. The functor $\sSeq_R^{\heartsuit}\rt \sSeq^{\gen}_{\ulR, \geq 0}$ from \Cref{ex:algebraic operads} then simply sends a symmetric sequence of discrete $R$-modules to the corresponding constant simplicial symmetric sequence. Note that $\sSeq_R^{\heartsuit}$ is not the heart of the $t$-structure on $\sSeq^{\gen}_{\ulR}$.
\end{remark}
\begin{lemma}\label{lem:simp-cosimp day convolution}
The projective and tame model structure on $\scsSeq_R$ both satisfy the pushout-product axiom with respect to the Day convolution product and levelwise tensor product. The induced symmetric monoidal structures on $\sSeq^{\gen}_{\ulR}$ and (when $R$ is coherent) $\sSeq^{\gen, \vee}_{\ulR}$ coincide with the ones from Construction \ref{con:derived products} and \Cref{prop:derived pro-coh composition}.
\end{lemma}
\begin{proof}
The pushout-product axiom is readily checked using the sets of generating cofibrations and trivial cofibrations from Remark \ref{rem:generating cofibrations simplicial-cosimplicial}. The resulting symmetric monoidal structures on $\sSeq^{\gen}_{\ulR}$ and (when $R$ is coherent) $\sSeq^{\gen, \vee}_{\ulR}$ restrict to the Day convolution and levelwise tensor product on $R[\Orb_\Sigma]$. In the projective case, Remark \ref{rem:realisations and totalisations in simp-cosimp modules} and Remark \ref{rem:finite totalisation is cof} show that they preserve finite totalisations of cosimplicial diagrams in $R[\Orb_\Sigma]$. In the tame case, Remark \ref{rem:realisations and totalisations in simp-cosimp modules} implies that they preserve all totalisations of cosimplicial diagrams in $R[\Orb_\Sigma]$. The result then follows from Remark \ref{rem:left right extended functors}.
\end{proof}
\begin{proposition}\label{prop:restricted composition product compatible with tame}
Let $\bcirc$ denote the restricted composition product on $\scsSeq_R$, computed levelwise. If $X$ is a tamely cofibrant simplicial-cosimplicial symmetric sequence, then the following assertions hold:
\begin{enumerate}
\item The functor $(-)\bcirc X\colon \scsSeq_R\rt \scsSeq_R$ preserves tame cofibrations and trivial cofibrations and the associated functor between $\infty$-categories $(-)\bcirc^{\LL} X\colon \scsSeq_R[W_{\mm{tame}}^{-1}]\rt \scsSeq_R[W_{\mm{tame}}^{-1}]$ preserves all colimits.
\item The functor $X\bcirc(-)\colon \scsSeq_R\rt \scsSeq_R$ preserves tamely cofibrant objects, as well as tame cofibrations and tame weak equivalences between them, and the associated functor between $\infty$-categories $X\bcirc^{\LL} (-)\colon \scsSeq_R[W_{\mm{tame}}^{-1}]\rt \scsSeq_R[W_{\mm{tame}}^{-1}]$ preserves sifted colimits.
\end{enumerate}
The same assertions apply to the usual composition product.
\end{proposition}
Following the same sort of strategy as in the proof of Proposition \ref{prop:composition product compatible with tame}, we will decompose the functors $\bcirc$ and $\circ$ into several functors. Most importantly, we will need a slightly more refined understanding of the functor sending a simplicial-cosimplicial symmetric sequence $X$ to its $r$-fold Day convolution product $X^{\otimes r}$, together with its $\Sigma_r$-equivariant structure. To this end, let us write
$$
\scMod_{R[\Sigma_r\times \Sigma]}=\prod_{q\geq 0} \scMod_{R[\Sigma_r\times \Sigma_q]}
$$
for the category of simplicial-cosimplicial symmetric sequences with an additional action of $\Sigma_r$, equipped with the tame model structure for the family $\mc{F}=\mc{O}_{\Sigma_r\times\Sigma_q}$ of all subgroups.
\begin{lemma}\label{lem:r-fold tensor product}
Let $r\geq 0$ and consider the functor $\scsSeq_R\rt \scMod_{R[\Sigma_r\times \Sigma]}$ sending $X\longmapsto X^{\otimes r}$ to its $r$-fold Day tensor product. This functor preserves tamely cofibrant objects, as well as tame cofibrations and weak equivalences between tamely cofibrant objects. 
\end{lemma}
\begin{proof}
Suppose that $X$ is a tamely cofibrant object. Then $X$ is given in each simplicial-cosimplicial bidegree $(a, b)$ by a retract of the free $R$-linear symmetric sequence $R[S(a, b)]$ generated by a symmetric sequence of sets $S(a, b)$. It follows that $X^{\otimes r}$ is given in simplicial-cosimplicial bidegree $(a, b)$ by a retract of the $R$-linearisation of the $\Sigma_r$-equivariant symmetric sequence $S(a, b)^{\otimes r}$ given by the $r$-fold Day tensor product of symmetric sequences of sets. By Lemma \ref{lem:cof-triv fib in simcos tame}, this implies that $X^{\otimes r}$ is cofibrant in the tame model structure on $\scMod_{R[\Sigma_r\times \Sigma]}$.

To study the behaviour of the functor $X\longmapsto X^{\otimes r}$, we will use a filtration argument: If $X\rt Y$ is a map in $\scsSeq_R$ that is a split monomorphism in each simplicial-cosimplicial bidegree, then there is a natural filtration
$$\begin{tikzcd}
X^{\otimes r}=F_0\arrow[r, rightarrowtail] & F_1\arrow[r, rightarrowtail] & \dots\arrow[r, rightarrowtail] & F_r=Y^{\otimes r}
\end{tikzcd}$$
in $\scMod_{R[\Sigma_r\times \Sigma]}$ obtained by declaring $X$ to be in filtration weight $0$ and $Y$ to be in weight $1$. This filtration has the property that each $F_{m-1}\rt F_{m}$ is a split monomorphism in each simplicial-cosimplicial bidegree and that the associated graded is $\bigoplus_{m=0}^r \big(F_m/F_{m-1}\big)\cong (X\oplus Y/X)^{\otimes r}$. Here $F_m/F_{m-1}$ corresponds to the summands in $(X\oplus Y/X)^{\otimes r}$ containing exactly $m$ tensor factors of $Y/X$.

Recall from Lemma \ref{lem:cof-triv fib in simcos tame} that $X\rt Y$ is a tame cofibration if and only if it is a split monomorphism in each simplicial-cosimplicial bidegree, with tamely cofibrant cokernel. If $X$ is tamely cofibrant, it follows that $(X\oplus Y/X)^{\otimes r}$ is tamely cofibrant, so that each summand $F_m/F_{m-1}$ is cofibrant. Consequently, each $F_{m-1}\rt F_m$ is a tame cofibration in $\scMod_{R[\Sigma_r\times \Sigma]}$ and $X\longmapsto X^{\otimes r}$ indeed preserves cofibrations between cofibrant objects.

Finally, we will show that $X\longmapsto X^{\otimes r}$ preserves trivial cofibrations, and hence weak equivalences, between cofibrant objects. Since it preserves transfinite compositions, the small object argument implies that it will suffice to verify the following: For any pushout square in $\scsSeq_R$
$$\begin{tikzcd}
A\arrow[r]\arrow[d] & B\arrow[d]\\
X\arrow[r] & Y
\end{tikzcd}$$
where $A\rt B$ is a generating trivial cofibration as in Remark \ref{rem:generating cofibrations simplicial-cosimplicial} and $X$ is cofibrant, the map $X^{\otimes r}\rt Y^{\otimes r}$ is a trivial cofibration in $\scMod_{R[\Sigma_r\times \Sigma]}$. Using that $Y/X\cong B/A$, we obtain a filtration on $Y^{\otimes r}$ whose associated graded is $(X\oplus B/A)^{\otimes r}$. Remark \ref{rem:generating cofibrations simplicial-cosimplicial} shows that $B/A$ either admits a simplicial homotopy equivalence to $0$ (for the generating cofibrations arising from $j_1$) or a cosimplicial homotopy equivalence to $0$ (for $j_2$). Since the functor $(X\oplus -)^{\otimes r}$ is ``prolongation'' in both the simplicial and the cosimplicial direction of its restriction to discrete symmetric sequences, a result of Dold \cite[Theorem 5.6]{dold1958homology} implies that $X\rt (X\oplus B/A)^{\otimes r}$ is a simplicial or cosimplicial homotopy equivalence. Consequently, the associated graded of the filtration on $Y^{\otimes r}$ is acyclic and $X^{\otimes r}\rt Y^{\otimes r}$ is indeed a trivial cofibration.
\end{proof}
\begin{proof}[Proof of \Cref{prop:restricted composition product compatible with tame}]
We decompose the composition products $\bcirc$ and $\circ$ into several functors, along the lines of Proposition \ref{prop:composition product compatible with tame}. First, let $r\geq 0$ and consider the functor
$$\begin{tikzcd}
\scMod_{R[\Sigma_r]}\times \scMod_{R[\Sigma_r\times \Sigma]}\arrow[r] & \scMod_{R[\Sigma\times \Sigma_r]}; & (V, W)\arrow[r, mapsto] & V\otimes W
\end{tikzcd}$$
sending a simplicial-cosimplicial $\Sigma_r$-representation and a $\Sigma_r$-equivariant simplicial-cosimplicial symmetric sequence to their tensor product, with diagonal $\Sigma_r$-action. This is readily checked to be a left Quillen bifunctor for the tame model structures (where all subgroups are admissible). 

Next, consider the functor $(-)^{\Sigma_r}\colon \scMod_{R[\Sigma\times \Sigma_r]}\rt \scsSeq_R$ taking $\Sigma_r$-fixed points. This preserves transfinite compositions and pushouts along monomorphisms that are split in each simplicial-cosimplicial degree. Using this and the description of the generating (trivial) cofibrations from Remark \ref{rem:generating cofibrations simplicial-cosimplicial}, one sees that $(-)^{\Sigma_r}$ preserves cofibrations and trivial cofibrations ; This follows from a similar, but much easier, argument as in the proof of Lemma \ref{lem:r-fold tensor product}. The same argument applies to the functor $(-)_{\Sigma_r}\colon \scMod_{R[\Sigma\times \Sigma_r]}\rt \scsSeq_R$.

For (1), notice that $X^{\otimes r}$ defines a tamely cofibrant in $\scMod_{R[\Sigma_r\times \Sigma]}$ by Lemma \ref{lem:r-fold tensor product}. The above two observations then show that the functors sending a simplicial-cosimplicial symmetric sequence $Y$ to $Y\bcirc X = \bigoplus_r (Y(r)\otimes X^{\otimes r})^{\Sigma_r}$ and $Y\circ X = \bigoplus_r (Y(r)\otimes X^{\otimes r})_{\Sigma_r}$ both preserve cofibrations and trivial cofibrations. Since these functors preserve pushouts along cofibrations and direct sums, the induced functors of $\infty$-categories preserve all colimits.

For (2), let $X$ be a tamely cofibrant object. Lemma \ref{lem:r-fold tensor product} and the above two observations then show that for each $r\geq 0$, the functors sending $Y\in \scsSeq_R$ to $(X(r)\otimes Y^{\otimes r})^{\Sigma_r}$ and $(X(r)\otimes Y^{\otimes r})_{\Sigma_r}$ preserve tamely cofibrant objects, as well as tame cofibrations and weak equivalences between tamely cofibrant objects. Taking the sum over all $r$ then shows that $X\bcirc (-)$ and $X\circ(-)$ also preserve tamely cofibrant objects, as well as cofibrations and weak equivalences between them.

It remains to verify that the induced functors of $\infty$-categories preserve sifted colimits. First, note that $X\bcirc^{\LL} (-)$ and $X\circ^{\LL} (-)\colon \scsSeq_R[W_{\mm{tame}}^{-1}]\rt \scsSeq_R[W_{\mm{tame}}^{-1}]$ preserve geometric realisations of simplicial objects, since these can simply be computed by the simplicial diagonal (Remark \ref{rem:realisations and totalisations in simp-cosimp modules}). It remains to verify that $X\bcirc^{\LL} (-)$ and $X\circ^{\LL} (-)\colon \scsSeq_R[W_{\mm{tame}}^{-1}]\rt \scsSeq_R[W_{\mm{tame}}^{-1}]$ preserve colimits of diagrams indexed by filtered posets \cite[Proposition 5.3.1.16, Corollary 5.5.8.17]{lurie2009higher}. 

Let us first verify this for $\omega_1$-filtered posets. For this case, note that the weak equivalences in $\scsSeq_R$ are closed under $\omega_1$-filtered colimits. Indeed, the functor $\Tot_{\oplus}\colon \scsSeq_R\rt \Ch_{R[\Sigma]}$ taking total complexes preserves filtered colimits and detects tame weak equivalences, and the weak equivalences in $\Ch_{R[\Sigma]}$ are closed under $\omega_1$-filtered colimits by Remark \ref{rem:O-tame we detected by finite stuff}. This implies that $\omega_1$-filtered colimits in $\scsSeq_R$ are already homotopy colimits, and the result follows from the fact that $X\bcirc (-), X\circ (-)\colon \scsSeq_R\rt \scsSeq_R$ preserve filtered colimits.

Finally, we need to check that $X\bcirc^{\LL} (-)$ and $X\circ^{\LL} (-)\colon \scsSeq_R[W_{\mm{tame}}^{-1}]\rt \scsSeq_R[W_{\mm{tame}}^{-1}]$ preserve colimits indexed by filtered posets with \emph{countably} many objects. Every such countable poset admits a cofinal functor from $\mathbb{N}$, so we have to show that $X\bcirc (-), X\circ (-)\colon \scsSeq_R\rt \scsSeq_R$ preserve \emph{homotopy} colimits of sequences. This follows immediately from the fact that both functors preserve sequences of cofibrant objects and cofibrations between them, as well as sequential colimits.
\end{proof}
\begin{lemma}\label{lem:composition preserves totalisation}
The composition products $\circ^{\LL}, \bcirc^{\LL}\colon \scsSeq_R[W_{\mm{tame}}]\times \scsSeq_R[W_{\mm{tame}}]\rt \scsSeq_R[W_{\mm{tame}}]$ preserve totalisations of cosimplicial diagrams in $R[\mc{O}_\Sigma]\hookrightarrow \scsSeq_R[W_{\mm{tame}}]$.
\end{lemma}
\begin{proof}
By Remark \ref{rem:realisations and totalisations in simp-cosimp modules}, the totalisation of a cosimplicial object $X^\bullet\colon \Del\rt R[\mc{O}_\Sigma]$ can simply be presented by the cosimplicial symmetric sequence $X^\bullet\in \scsSeq_R$, which is tamely cofibrant. Since the composition products of simplicial-cosimplicial symmetric sequences are computed levelwise, the result follows.
\end{proof}
\begin{proposition}\label{prop:restricted composition product compatible with projective}
If $X$ and $Y$ are projectively cofibrant simplicial-cosimplicial symmetric sequences, then $X\circ Y$ and $X\bcirc Y$ are projectively cofibrant. Consequently, the composition product and the restricted composition product induce monoidal structures $\circ^{\LL}, \bcirc^{\LL}\colon \sSeq^{\gen}_{\ulR}\times \sSeq^{\gen}_{\ulR}\rt \sSeq^{\gen}_{\ulR}$ preserving colimits in the first variable and sifted colimits in the second variable.
\end{proposition}
\begin{proof}
The proof from Proposition \ref{prop:composition product compatible with projective} carries over to this setting. The projective model structure is a right Bousfield localisation of the tame model structure and exhibits 
$$
\sSeq^{\gen}_{\ulR}\simeq \scsSeq_R[W_{\mm{proj}}^{-1}]\subseteq \scsSeq_R[W_{\mm{tame}}^{-1}]
$$ 
as the full subcategory generated under colimits and finite totalisations by the symmetric sequences $R[\Sigma_r/H]$ for all subgroups $H<\Sigma_r$. Verifying that $\circ$ and $\bcirc$ preserve projectively cofibrant objects therefore comes down to verifying that $\sSeq^{\gen}_{\ulR}$ is stable under the (restricted) composition products from Proposition \ref{prop:restricted composition product compatible with tame}. If this is the case, the restriction will automatically preserve colimits in the first variable and sifted colimits in the second variable.

To see that $\sSeq^{\gen}_{\ulR}$ is stable under $\circ^{\LL}$ and $\bcirc^{\LL}$, note that both functors preserve colimits and (hence) finite totalisations in the first variable. It therefore suffices to verify that for all $H<\Sigma_r$, the objects $R[\Sigma_r/H]\circ^{\LL} Y\cong (Y^{\otimes r})_H$ and $R[\Sigma_r/H]\bcirc^{\LL} Y\cong (Y^{\otimes r})^H$ are contained in $\sSeq^{\gen}_{\ulR}$ whenever $Y\in \sSeq^{\gen}_{\ulR}$. Since both functors preserve sifted colimits by Proposition \ref{prop:restricted composition product compatible with tame}, we just need to verify that $(Y^{\otimes r})_H$ and $(Y^{\otimes r})^H$ are contained in $\sSeq^{\gen}_{\ulR}$ whenever $Y=\Tot(Y^\bullet)$ is a the totalisation of a finite cosimplicial diagram $Y^\bullet\colon \Del\rt R[\mc{O}_\Sigma]$. Let us write $D_\bullet\colon \Del^{\op}\rt R[\mc{O}_\Sigma]$ for the $R$-linear dual simplicial object, so that $Y\simeq |D_\bullet|^\vee$ in $\scsSeq_R$. Using Remark \ref{rem:realisations and totalisations in simp-cosimp modules}, one can then identify $$
(Y^{\otimes r})_H\simeq \big|(D_\bullet^{\otimes r})^H\big|^\vee\qquad\qquad \text{and}\qquad\qquad (Y^{\otimes r})^H\simeq \big|(D_\bullet^{\otimes r})_H\big|^\vee.
$$
The simplicial symmetric sequences $\big|(D_\bullet^{\otimes r})^H\big|$ and $\big|(D_\bullet^{\otimes r})_H\big|$ are both perfect, because the functors $X\longmapsto (X^{\otimes r})_H$ and $X\longmapsto (X^{\otimes r})^H$ are $r$-excisive (see the proof of \Cref{cor:polynomial functors between coherent}). In other words, both simplicial symmetric sequences are equivalent to the geometric realisation of a finite simplicial object in $R[\mc{O}_\Sigma]$. This implies that their $R$-linear duals $(Y^{\otimes r})_H$ and $(Y^{\otimes r})^H$ are equivalent to the totalisation of a finite cosimplicial object in $R[\mc{O}_\Sigma]$. Remark \ref{rem:finite totalisation is cof} now implies that $(Y^{\otimes r})_H$ and $(Y^{\otimes r})^H$ are contained in $\sSeq^{\gen}_{\ulR}$.
\end{proof}
\begin{corollary}\label{cor:simp-cosimp model for exotic composition}
Let $R$ be a ring. Then the monoidal structures $\circ^{\LL}$ and $\bcirc^{\LL}$ on $\sSeq^{\gen}_{\ulR}$ from \Cref{prop:restricted composition product compatible with projective} are naturally equivalent to the composition product $\circ$ and the restricted composition product $\bcirc$ from \Cref{con:derived products}.
\end{corollary}
\begin{proof}
The (restricted) composition product from Proposition \ref{prop:restricted composition product compatible with projective} preserves colimits in the first variable and sifted colimits in the second variable. By Remark \ref{rem:finite totalisation is cof} and Remark \ref{rem:realisations and totalisations in simp-cosimp modules}, it furthermore preserves finite totalisations of cosimplicial diagrams in $R[\mc{O}_\Sigma]$. Since it restricts to the usual (restricted) composition product on $R[\mc{O}_{\Sigma}]$, Remark \ref{rem:left right extended functors} then implies that it is naturally equivalent to the (restricted) composition product from \Cref{con:derived products}.
\end{proof}
\begin{theorem}\label{thm:simp-cosimp model for exotic composition}
Let $R$ be a coherent ring. The composition product $\circ$ and the restricted composition product $\bcirc$ on the tame model category $\scsSeq_R$ induce the monoidal structures $\circ$ and $\bcirc$ of  \Cref{prop:derived pro-coh composition} on the underlying $\infty$-category $\sSeq^{\gen,\vee}_{\ulR}$.
\end{theorem}
\begin{proof}
By (1) and (2) of Proposition \ref{prop:composition product compatible with tame}, $\circ$ and $\bcirc$ restrict to monoidal products on the full subcategory of tamely cofibrant symmetric sequences, which preserve weak equivalences in each variable. By part (1) and (3) of Proposition \ref{prop:composition product compatible with tame}, the resulting monoidal structures $\circ^{\sLL}$ and $\bcirc^{\sLL}$ on the $\infty$-category $\sSeq^{\gen, \vee}_{\ulR}$ preserve sifted colimits. Furthermore, the restrictions of these monoidal structures to $R[\Orb_\Sigma]$ coincide with the usual composition product and restricted composition product. Finally, $\circ^{\sLL}$ and $\bcirc^{\sLL}$ preserve totalisations of cosimplicial diagrams in $R[\Orb_\Sigma]$ (Lemma \ref{lem:composition preserves totalisation}), so that both are obtained by left-right extension (Remark \ref{rem:left right extended functors}) and hence coincide with the monoidal structures from Proposition \ref{prop:derived pro-coh composition}
\end{proof}

\subsubsection*{Rectification of derived operads and algebras}
\Cref{cor:simp-cosimp model for exotic composition} and \Cref{thm:simp-cosimp model for exotic composition} show that the composition product $\bcirc$ on $\sSeq^{\gen}_{\ulR}$ and $\sSeq^{\gen, \vee}_{\ulR}$ can be presented by the restricted composition product on $\scsSeq_R$. We will now show how this can be used to produce point-set models for the $\infty$-categories of derived $\infty$-operads and derived PD $\infty$-operads.
 
\begin{proposition}\label{prop:model cats of simp-cosimp operads}
The following categories carry cofibrantly generated semi-model structures whose weak equivalences and fibrations are tame weak equivalences and fibrations on the underlying objects:
\begin{enumerate}
\item the category $\mathbf{Op}^{\simcos}_{R}$ of sc-operads.\vspace{2pt}
\item the category $\mathbf{Op}^{\simcos,\res}_R$ of sc-restricted operads over $R$.\vspace{2pt}
\item the category $\mathbf{Alg}^{\simcos}_{\ope{P}}$ of algebras over an sc-operad $\ope{P}$ whose underlying symmetric sequence is tamely cofibrant.\vspace{2pt}
\item the category $\mathbf{Alg}^{\simcos,\res}_{\ope{P}}$ of algebras over an sc-restricted operad $\ope{P}$ whose underlying symmetric sequence is tamely cofibrant.
\end{enumerate}
One can also endow these categories with a semi-model structure whose weak equivalences and fibrations are the ones from the projective model structure.
\end{proposition}
\begin{proof}
We only treat the case (4) of restricted algebras over an sc-restricted operad $\ope{P}$ whose underlying symmetric sequence is tamely cofibrant; the other cases are similar. By \cite[Theorem 12.1.4]{fresse2009modules}, it suffices to verify the following condition: for a cofibrant $\ope{P}$-algebra $A$, a generating trivial cofibration $V\rt W$ in $\scMod_R$ and a map $f\colon V\rt A$ in $\scMod_R$, the map $A\rt A\amalg_{\ope{P}\bcirc V} \ope{P}\bcirc W$ is a trivial cofibration in $\scMod_R$. To prove this, we use that the pushout carries an exhaustive increasing filtration
\begin{equation}\label{eq:filtration on pushout}\begin{tikzcd}
A=F^0\arrow[r, hook] & F^1\arrow[r, hook] & \dots \arrow[r, hook] & A\amalg_{\ope{P}\bcirc V} \ope{P}\bcirc W
\end{tikzcd}\end{equation}
where each map is a split monomorphism in each simplicial-cosimplicial degree. To see this, consider the category $\Fun(\mathbb{Z}_{\geq 0}, \scMod_R)$ of increasing sequences in $\scMod_R$, with the Day convolution product and the Reedy model structure. We can consider $A$ as a $\ope{P}$-algebra in $\Fun(\mathbb{Z}_{\geq 0}, \scMod_R)$ given by $A$ in each filtration weight.
Likewise, consider $V$ as a constant sequence and let $W'$ denote the sequence given by $V$ in weight $0$ and by $W$ in weight $\geq 1$. Then the pushout $A\amalg_{\ope{P}\bcirc V} \ope{P}\bcirc W'$ of $\ope{P}$-algebras in $\Fun(\mathbb{Z}_{\geq 0}, \scMod_R)$ will produce the desired filtration.

Indeed, note that $A$ being cofibrant implies that it is given in each simplicial-cosimplicial degree $(i, j)$ by the retract of a free $\ope{P}_i^j$-algebra on a projective $R$-module $X_i^j$. Because $V\rt W$ is a split monomorphism in each simplicial-cosimplicial degree, we can identify $A\amalg_{\ope{P}\bcirc V} \ope{P}\bcirc W'\cong \ope{P}_i^j\bcirc \big(X_i^j)\oplus (W/V)_i^j\big)$, where $(W/V)_i^j$ has weight $1$. This shows that the inclusions in \eqref{eq:filtration on pushout} are split injections in each simplicial-cosimplicial degree.

It therefore suffices to show that the associated graded of the filtration \eqref{eq:filtration on pushout} consists of acyclic tamely cofibrant sc-$R$-modules in weight $\geq 0$. The associated graded can be identified with the coproduct $A\amalg \ope{P}\bcirc (W/V)$. It therefore suffices to prove that for any cofibrant $\ope{P}$-algebra $A$ and any contractible cofibrant sc-module $Z$, the map $A\rt A\amalg \ope{P}\bcirc Z$ is an acyclic cofibration of sc-modules. Using the small object argument and a similar filtration argument to the one given above, this can be reduced to the assertion that for any cofibrant sc-module $X$, the map $\ope{P}\bcirc X\rt \ope{P}\bcirc (X\oplus Z)$ is an acyclic cofibration of sc-$R$-modules. This follows from Proposition \ref{prop:restricted composition product compatible with tame}.
\end{proof}
\begin{remark}\label{rem:forget preserves cofibrations}
The proof shows that the forgetful functors $\mathbf{Op}^{\simcos}_R\rt \scsSeq_R$, $\mathbf{Op}^{\mm{res}}_R\rt \scsSeq_R$ and $\mathbf{Alg}^{\simcos}_{\ope{P}}\rt \scMod_{R}$ preserve cofibrations between cofibrant objects.
\end{remark}
\begin{lemma}\label{lem:preservation of diagonal}
Consider the forgetful functors $\mathbf{Op}^{\simcos}_R\rt \scsSeq_R$, $\mathbf{Op}^\mm{\simcos, \res}_R\rt \scsSeq_R$ and $\mathbf{Alg}^{\simcos}_{\ope{P}}\rt \scMod_R$, where $\ope{P}$ is a simplicial-cosimplicial (restricted) operad whose underlying symmetric sequence is tamely cofibrant. Each of these functors is right Quillen for the model structures from Proposition \ref{prop:model cats of simp-cosimp operads}, and the induced functor of $\infty$-categories preserves geometric realisations.
\end{lemma}
\begin{proof}
The forgetful functors are right Quillen functors that preserve weak equivalences by construction. It remains to prove that the induced functor of $\infty$-categories preserves geometric realisations, which we will only do in the case of $\ope{P}$-algebras (the other cases are exactly the same). Recall that $\scMod_R$ is a simplicial model category and note that the cotensoring over simplicial sets preserves $\ope{P}$-algebras. Using this, there is an adjoint pair
$$\begin{tikzcd}
\delta^*\colon \Fun(\Del^{\op}, \mathbf{Alg}^{\simcos}_{\ope{P}})\arrow[r, yshift=1ex] & \mathbf{Alg}^{\simcos}_{\ope{P}}\colon \delta_*\arrow[l, yshift=-1ex]
\end{tikzcd}$$
where $\delta^*$ takes the diagonal in the simplicial direction, and $\delta_*$ sends a simplicial-cosimplicial algebra $A$ to the simplicial diagram $A^{\Delta[\bullet]}$. This adjoint pair is a Quillen pair when $\Fun(\Del^{\op}, \mathbf{Alg}^{\simcos}_{\ope{P}})$ is endowed with the Reedy (semi-)model structure; furthermore, the right adjoint $\delta_*$ sends every every fibrant object $A$ to a simplicial diagram that is homotopically constant on $A$. It follows that the left derived functor of $\delta^*$ computes the homotopy colimit of a simplicial diagram in $\mathbf{Alg}^{\simcos}_{\ope{P}}$. Since the forgetful functor $\mathbf{Alg}^{\simcos}_{\ope{P}}\rt \scMod_R$ commutes with taking the diagonal and preserves Reedy cofibrant simplicial diagrams by Remark \ref{rem:forget preserves cofibrations}, the result follows.
\end{proof}
\begin{theorem}[Rectification of derived $\infty$-operads and derived PD $\infty$-operads]\label{thm:simp-cosimp models for derived PD operads}
Let $R$ be a coherent ring. Then the underlying $\infty$-category of $\mathbf{Op}^{\simcos}_R$ with the projective semi-model structure is equivalent to the $\infty$-category $\Op^{\gen}_{\ulR}$ of derived $\infty$-operads over $R$. Likewise, the underlying $\infty$-category of $\mathbf{Op}^{\simcos,\res}_R$ with the tame semi-model structure is equivalent to the $\infty$-category $\Op^{\gen, \pd}_{\ulR}$ of derived PD $\infty$-operads over $R$. More precisely, there are commuting squares
$$\begin{tikzcd}
\mathbf{Op}^{\simcos}_R[W_{\mm{proj}}^{-1}]\arrow[d]\arrow[r, "\simeq"] & \Op^{\gen}_{\ulR}\arrow[d] & \mathbf{Op}^{\simcos,\res}_R[W_{\mm{tame}}^{-1}]\arrow[d]\arrow[r, "\simeq"] & \Op^{\gen, \pd}_{\ulR}\arrow[d]\\
\scsSeq_R[W_{\mm{proj}}^{-1}]\arrow[r, "\simeq"] & \sSeq^{\gen}_{\ulR} & \scsSeq_R[W_{\mm{tame}}^{-1}]\arrow[r, "\simeq"] & \sSeq^{\gen, \vee}_{\ulR}.
\end{tikzcd}$$
\end{theorem}
\begin{proof}
Following the same argument as in   \Cref{thm:chain models for PD operads}, using \Cref{thm:simp-cosimp model for exotic composition}, Proposition \ref{prop:model cats of simp-cosimp operads} and Lemma \ref{lem:preservation of diagonal}.
\end{proof}
\begin{theorem}[Rectification of algebras: derived setting]\label{thm:simp-cosimp models for derived algebras}
Let $R$ be a coherent ring and $\ope{P}$ a simplicial-cosimplicial (restricted) operad over $R$ whose underlying symmetric sequence is tamely cofibrant. Then the underlying $\infty$-category of the tame semi-model structure on simplicial-cosimplicial algebras over $\ope{P}$ is equivalent to the $\infty$-category $\Alg_{\PP}(\QC^\vee_R)$ of pro-coherent algebras over the associated derived (PD) $\infty$-operad $\PP$. In other words, there are commuting squares
$$\begin{tikzcd}
\mathbf{Alg}^{\simcos}_{\ope{P}}[W_{\mm{tame}}^{-1}]\arrow[d]\arrow[r, "\simeq"] & \Alg_{\PP}^{\gen}(\QC^\vee_R)\arrow[d] & \mathbf{Alg}^{\simcos,\res}_{\ope{P}}[W_{\mm{tame}}^{-1}]\arrow[d]\arrow[r, "\simeq"] & \Alg^{\gen,\pd}_{\PP}(\QC^\vee_R)\arrow[d]\\
\scMod_{R}[W_{\mm{tame}}^{-1}]\arrow[r, "\simeq"] & \QC^\vee_R & 
\scMod_{R}[W_{\mm{tame}}^{-1}]\arrow[r, "\simeq"] & \QC^\vee_R.
\end{tikzcd}$$
\end{theorem}
\begin{proof}
Exactly as in \Cref{thm:chain models for PD algebras}, we combine \Cref{thm:simp-cosimp model for exotic composition}, Proposition \ref{prop:model cats of simp-cosimp operads} and Lemma \ref{lem:preservation of diagonal}.
\end{proof}
\begin{remark}[Rectification of derived $\infty$-operads and algebras over non-coherent rings]\label{rem:rectification algebras non-coherent}
Similarly to \Cref{rem:chain models simplicial}, the parts of \Cref{thm:simp-cosimp models for derived PD operads} and \Cref{thm:simp-cosimp models for derived algebras} pertaining to the projective model structure also apply when $R$ is not a coherent ring. More precisely, for any commutative ring, there is an equivalence of $\infty$-categories $\mathbf{Op}^{\simcos}_R[W_{\mm{proj}}^{-1}]\simeq \Op^{\gen}_{\ulR}$ and for every simplicial-cosimplicial (restricted) operad $\ope{P}$ over $R$ whose underlying symmetric sequence is projectively cofibrant, there are commuting squares
$$\begin{tikzcd}
\mathbf{Alg}^{\simcos}_{\ope{P}}[W_{\mm{proj}}^{-1}]\arrow[d]\arrow[r, "\simeq"] & \Alg_{\PP}^{\gen}(\Mod_R)\arrow[d] & \mathbf{Alg}^{\simcos,\res}_{\ope{P}}[W_{\mm{proj}}^{-1}]\arrow[d]\arrow[r, "\simeq"] & \Alg^{\gen,\pd}_{\PP}(\Mod_R)\arrow[d]\\
\scMod_{R}[W_{\mm{proj}}^{-1}]\arrow[r, "\simeq"] & \Mod_R & 
\scMod_{R}[W_{\mm{proj}}^{-1}]\arrow[r, "\simeq"] & \Mod_R.
\end{tikzcd}$$
Indeed, Lemma \ref{lem:preservation of diagonal} shows that the left vertical functors are monadic right adjoints, and \Cref{cor:simp-cosimp model for exotic composition} identifies the corresponding monads on $\Mod_R$ with the monads $\PP\circ (-)$ and $\PP\bcirc (-)$ whose categories of algebras are precisely $\Alg_{\PP}^{\gen}(\Mod_R)$ and $\Alg^{\gen,\pd}_{\PP}(\Mod_R)$, respectively.
\end{remark}
In particular, we deduce: 
\begin{corollary}[Derived rings as simplicial-cosimplicial commutative rings] \label{derivedringssimplicialcosimplicial}
For any $R$-algebra,   the $\infty$-category $\Alg_{\Com^{ }_R}(\Mod_R)$ 
of derived commutative $R$-algebras 
is obtained from the model category $\mathbf{Alg}^{\simcos}_{\Com_R}$ of simplicial-cosimplicial commutative $R$-algebras by inverting projective weak equivalences.
\end{corollary}

\subsection{Partition Lie algebras}
We will now construct an explicit cosimplicial model for the derived PD $\infty$-operad (cf.\ \Cref{def:derived partition lie operad})  which parametrises derived partition Lie algebras, freely using the techniques developed in \cite{arone2018action}. This cosimplicial model is the linear dual of the subdivided simplicial bar construction of the commutative operad. As the bar construction can already be computed   in pointed simplicial sets, we start by working with symmetric sequences and \mbox{operads in this setting.}

Let $\Com^{\mm{nu}}$ denote the (nonunital) commutative operad, given by $S^0$ with trivial $\Sigma_r$-action in each arity $\geq 1$ and by a point in arity $0$. This is an augmented operad, and we can consider the simplicial bar construction
$$\begin{tikzcd}
{\Barr_\bullet(\mb{1}, \Com^{\mm{nu}}, \mb{1})= \dots }\arrow[r, yshift=1.5ex]\arrow[r, yshift=0.5ex] \arrow[r, yshift=-0.5ex]\arrow[r, yshift=-1.5ex] & \Com^{\mm{nu}}\circ \Com^{\mm{nu}}\arrow[r, yshift=1ex]\arrow[r, yshift=-1ex]\arrow[r] & \Com^{\mm{nu}}\arrow[r, yshift=0.5ex]\arrow[r, yshift=-0.5ex] & \mb{1}.
\end{tikzcd}$$
This simplicial bar construction has a well-known description in terms of  partition complexes, which we will now recall.
\begin{notation}
Write  $\mc{P}_r$ for the poset of partitions of $\ul{r}=\{1, \dots, r\}$, ordered by coarsening; the initial and terminal partitions are given by
$$
\hat{0}=\Boxed{1}\hspace{-0.3 pt} \Boxed{2} \hspace{-0.3 pt} \Boxed{3} \ldots \Boxed{r}\qquad\qquad\qquad\qquad \hat{1}=\Boxed{123\ldots r}.
$$
Write $N_\bullet(\mc{P}_r)^{-[\hat{0}<\hat{1}]}\subseteq N_\bullet(\mc{P}_r)$ for the simplicial subset spanned by the chains of partitions that do not contain $[\hat{0}<\hat{1}]$.
\end{notation}
One can then identify $\Barr_\bullet(\mb{1}, \Com^{\mm{nu}}, \mb{1})$ in arity $r$ with the pointed \mbox{simplicial set}
\begin{equation}\label{eq:simplicial bar of com}
\Barr_\bullet(\mb{1}, \Com^{\mm{nu}}, \mb{1})(r)=\frac{N_\bullet(\mc{P}_r)}{N_\bullet(\mc{P}_r)^{-[\hat{0}<\hat{1}]}}.
\end{equation}
For $r=0$, this is the basepoint, and for $r=1$ this is $S^0$ by convention. The non-basepoint $n$-simplices of $\Barr_\bullet(\mb{1}, \Com^{\mm{nu}}, \mb{1})(r)$ are then given by chains of partitions $[\hat{0}=x_0\leq x_1\leq\dots \leq x_{n-1}\leq x_n=\hat{1}]$; the simplicial structure maps simply insert identities or remove elements from such chains, and give the basepoint if the resulting chain no longer begins at $\hat{0}$ and ends at $\hat{1}$.
\begin{remark}[Levelled trees]\label{rem:levelled tree}
A chain of partitions $[x_0\leq \dots \leq x_t]$ can be viewed as a \emph{levelled forest}, where each leaf is labelled by a subset of $\ul{r}$. Indeed, each leaf is labelled by a subset of $\ul{r}$ corresponding to a class in  $x_0$, and each class of $x_t$ determines a tree. In particular, chains $[\hat{0}=x_0\leq \dots\leq x_t=\hat{1}]$ correspond to \emph{levelled trees} with leaves labelled by the elements of $\ul{r}$.\vspace{5pt}

In these terms, the non-basepoint simplices of $\Barr_\bullet(\mb{1}, \Com^{\mm{nu}}, \mb{1})(r)$ correspond to levelled trees, and  the simplicial face maps are given by contracting edges between two levels or removing the root or leaf vertices; this produces the basepoint if the result is no longer a tree with $r$ leaves.
\end{remark}
\begin{notation}[Barycentric subdivision]\label{not:barycentric com}
We will denote the barycentric subdivision of a simplicial set $X$ by $\sd(X)$. We   abbreviate $\sd(N_\bullet(\mc{P}_r)^{-[\hat{0}<\hat{1}]})$ as $\sd(\mc{P}_r)^{-[\hat{0}<\hat{1}]}$. Finally, we define
$$
\sdBar(\Com^{\mm{nu}})(r):=\frac{\sd(\mc{P}_r)}{\sd(\mc{P}_r)^{-[\hat{0}< \hat{1}]}}
$$
as the barycentric subdivision of the simplicial bar construction \eqref{eq:simplicial bar of com}. Explicitly, this is the quotient of the nerve of the poset of nondegenerate \textit{chains} of partitions $\sigma=[x_0< \dots< x_t]$ by the full subcomplex spanned by chains \mbox{with $x_0\neq \hat{0}$ or $x_t\neq \hat{1}$.}

Even more explicitly,  $d$-simplices in $\sd(\mc{P}_r)\big/\sd(\mc{P}_r)^{-[\hat{0}<\hat{1}]}$   corresponds to  pairs
$$
\big(\sigma, S\big)=\big([\hat{0}=x_0<\dots <x_t=\hat{1}], S_0\subseteq \dots\subseteq S_d\big),
$$
where $\sigma$ is a nondegenerate chain of partitions of $\ul{r}$ and $S_0\subseteq \dots\subseteq S_d=\{0, \dots, t\}$ is an increasing set of subsets. We allow $t=-1$ in this definition, which corresponds to the basepoint. We will refer to such tuples as \emph{nested chains of partitions} of $\underline{r}$.
\end{notation}

Our goal will be to endow the barycentric subdivision $\sdBar(\Com^{\mm{nu}})$ with the structure of a cooperad in pointed simplicial sets. It will be convenient to describe such cooperads as symmetric sequences of pointed simplicial sets $\CC$ together with the following kind of cocomposition maps: for every partition $y$ of\vspace{-3pt} the form $$\ul{r}\cong \ul{r}_1\sqcup\dots\sqcup\ul{r}_b,$$ there is a total \vspace{-3pt}cocomposition map
$$\begin{tikzcd}
\Delta_y\colon \CC(r)\arrow[r] & \CC(b)\wedge \CC(r_1)\wedge\dots\wedge \CC(r_b)
\end{tikzcd}$$ 
which is equivariant with respect to the stabiliser $\Sigma_y<\Sigma_r$ of $y$  and satisfies obvious associativity and unitality constraints.

To define these cocomposition maps $\Delta_y$, we will need some terminology:
\begin{definition}
Let $y$ be a partition of $\ul{r}$ of the form $\ul{r}\cong \ul{r}_1\sqcup \dots \sqcup\ul{r}_b$. 

A nondegenerate chain of partitions $\sigma=[x_0<\dots < x_{t}]$ of $\ul{r}$ is said to be \emph{$y$-branched} if every class $\ul{r}_i$ of the partition $y$ arises as a class in \mbox{some partition $x_\alpha$ in $\sigma$.}

Furthermore, $\sigma$ is said to be \emph{$y$-subbranched} if it is contained in a $y$-branched chain of partitions. Write $\mm{Sub}(y)\subseteq \sd(\mc{P}_r)$ for the simplicial subsets spanned by all 
$$\big(\sigma, S\big)=\big([x_0<\dots <x_t], S_0\subseteq \dots\subseteq S_d\big)$$
for which $[x_0<\dots <x_t]$ is   $y$-subbranched.

Write $\mm{Unbr}(y)^{-{[\hat{0}<\hat{1}]}}\subseteq \sd(\mc{P}_r)$ for the simplicial subsets spanned by all nested chains $\big(\sigma, S\big)$ such that  either $\sigma$ is not $y$-branched or $\sigma$ does not contain $[\hat{0}<\hat{1}]$.
\end{definition}
\begin{notation}
If $\sigma$ is a degenerate chain of partitions,  let $((\sigma))$ be the minimal nondegenerate chain with a map $\sigma\rt \nondeg{\sigma}$; it is obtained by deleting repetitions.
\end{notation}
\begin{cons}[Ungrafting map]\label{cons:cutting trees}
Let $y\colon \ul{r}\cong \ul{r}_1\sqcup \dots\sqcup \ul{r}_b$ be a partition. 

We will define an order-preserving map
$$\begin{tikzcd}
\phi_y\colon \mm{Sub}(y)\arrow[r] & \sd(\mc{P}_b)\times \sd(\mc{P}_{r_1})\times \dots\times \sd(\mc{P}_{r_b})
\end{tikzcd}$$
which we will refer to as the \emph{ungrafting map} (along $y$). 

Start with the   map $\mm{Sub}(y)\rightarrow   \sd(\mc{P}_{r, \geq y})\times \sd(\mc{P}_{r, \leq y})$   induced by the map of posets  
$$ 
 \sigma \mapsto  \big(\nondeg{\sigma\vee y}, \nondeg{\sigma\wedge y}\big).
 $$
Next, note that there are isomorphisms $\mc{P}_{r, \geq y}\cong \mc{P}_{b}$ and $\mc{P}_{r, \leq y}\cong \mc{P}_{r_1}\times \dots \times \mc{P}_{r_b}$. On subdivisions, this induces an isomorphism $\sd(\mc{P}_{r, \geq y})\cong \sd(\mc{P}_{b})$ and a map
$$\begin{tikzcd}
\sd(\mc{P}_{r, \leq y})\arrow[r, "\cong"] & \sd\big(\mc{P}_{r_1}\times \dots \times \mc{P}_{r_b}\big)\arrow[r] & \sd(\mc{P}_{r_1})\times \dots \times \sd(\mc{P}_{r_b})
\end{tikzcd}$$
where the second map sends a nondegenerate tuple $(\sigma_1,\ldots ,\sigma_b)$ of chains to the tuple of nondegenerate chains $\nondeg{\sigma_i}$. 

Combining these two maps, we can assign to each nondegenerate chain $\sigma$ in $\mm{Sub}(y)$ a tuple of chains $(\sigma_{\ul{b}}, \sigma_{\ul{r}_1}, \dots, \sigma_{\ul{r}_b})$. We obtain $\phi_y(\sigma)$ from this tuple by 
removing the maximal partition $\hat{1}$  (if it appears)  from $\sigma_{\ul{r}_i}$ for each class $\ul{r}_i$ of   $y$ that does not appear anywhere  in  $\sigma$.  Note that this map preserves subchain inclusions.
\end{cons}
\begin{remark}[Description via levelled trees]\label{rem:cutting trees}
Let $\sigma$ be a chain of partitions of $\ul{r}$, corresponding to a levelled forest where each leaf is labelled by a subset of $\ul{r}$. Then $\sigma$ is $y$-branched if and only if for every class $\ul{r}_i$ of the partition $y$, there is a branch in the forest whose leaves are precisely labelled by subsets with union $\ul{r}_i$. These various branches may be of different height.

If $\sigma$ is $y$-branched, unraveling the definitions shows that $\phi_y(\sigma)$ is given as follows. The resulting chain $\sigma_{\ul{b}}$ in $\mc{P}_b$ corresponds to the forest obtained by cutting off all $\ul{r}_i$-labelled branches and inserting just enough degeneracies on the top to make the result a levelled forest. The chains of partitions $\sigma_{\ul{r_i}}$ of each $\ul{r}_i$ corresponds to the $\ul{r}_i$-labelled branch, with all of its degeneracies removed.
\begin{figure}[h]
\includegraphics[width=0.9\textwidth]{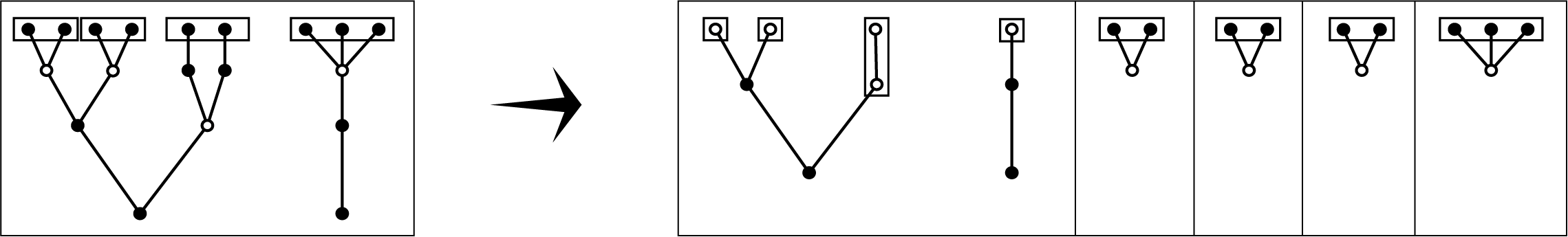}
\end{figure}

For a subchain $\tau\subseteq \sigma$ of such a $y$-branched chain $\sigma$, $\phi_y(\tau)$ simply takes the corresponding subchains in $\mc{P}_b$ and $\mc{P}_{r_i}$ and furthermore removes the endpoint $\hat{1}$ from $\sigma_{\ul{r}_i}$ whenever $\ul{r}_i$ does not appear as a class anywhere in $\tau$.
\end{remark}
\begin{lemma}
Let $y\colon \ul{r}\cong \ul{r}_1\sqcup \dots \sqcup \ul{r}_b$ be a partition. The inclusion $\mm{Sub}(y)\hookrightarrow \sd(\mc{P}_r)$ induces an isomorphism
$$
\frac{\mm{Sub}(y)}{\mm{Sub}(y)\cap \mm{Unbr}(y)^{-[\hat{0}<\hat{1}]}}\stackrel{\ \cong\ }{\rt} \frac{\sd(\mc{P}_r)}{\mm{Unbr}(y)^{- {[\hat{0}<\hat{1}]}}}
$$
and the map $\phi_y$ of  \Cref{cons:cutting trees} descends to a map
$$
\frac{\mm{Sub}(y)}{\mm{Sub}(y)\cap \mm{Unbr}(y)^{-[\hat{0}<\hat{1}]}}\rt \frac{\sd(\mc{P}_b)}{\sd(\mc{P}_b)^{-[\hat{0}< \hat{1}]}}\wedge \frac{\sd(\mc{P}_{r_1})}{\sd(\mc{P}_{r_1})^{-[\hat{0}< \hat{1}]}} \wedge \dots \wedge \frac{\sd(\mc{P}_{r_b})}{\sd(\mc{P}_{r_b})^{-[\hat{0}< \hat{1}]}}.
$$
\end{lemma}
\begin{proof}
For the first assertion, note that the map is injective by construction. To see that it is surjective, note that the non-basepoint simplices of $\sd(\mc{P}_r)\big/\mm{Unbr}(y)^{- {[\hat{0}<\hat{1}]}}$ correspond to nested families of nondegenerate chains $[\sigma_0\leq \dots \leq \sigma_n]$, where $\sigma_n$ is a nondegenerate chain of partitions which is $y$-branched and contains $\hat{0}$ and $\hat{1}$. This implies that all $\sigma_i$ are $y$-subbranched.

For the second assertion, first suppose that $\sigma$ is a $y$-branched chain. If $\sigma$ does not end in $\hat{1}$, then $\phi_y(\sigma)$ defines a chain in $\mc{P}_b$ not ending in $\hat{1}$. On the other hand, if $\sigma$ does not start at $\hat{0}$, then $\phi_y(\sigma)$ defines a chain in at least one $\mc{P}_{r_i}$ that does not start at $\hat{0}$. Both assertions are readily seen by the description of $\phi_y$ in terms of levelled trees (Remark \ref{rem:cutting trees}). Furthermore, if $\tau\subseteq \sigma$ is a subchain of a $y$-branched chain which is not itself $y$-branched, then $\phi_y(\tau)$ contains at least a chain in one $\mc{P}_{r_i}$ which does not end at $\hat{1}$ (by construction). Hence $\phi_y$ descends to the quotient.
\end{proof}
\begin{cons}[Cocomposition]\label{cons:cocomp for partitions}
Let $y: \ul{r}\cong \ul{r}_1\sqcup \dots\sqcup \ul{r}_b$ be a partition. We then define $\Delta_y$ to be the composite map
$$\begin{tikzcd}
\frac{\sd(\mc{P}_r)}{\sd(\mc{P}_r)^{-[\hat{0}<\hat{1}]}}\arrow[r] & \frac{\sd(\mc{P}_r)}{\mm{Unbr}(y)^{-[\hat{0}<\hat{1}]}}\arrow[r, "\cong"] &  \frac{\mm{Sub}(y)}{\mm{Sub}(y)\cap \mm{Unbr}(y)^{-[\hat{0}<\hat{1}]}}\arrow[d, "\phi_y"]\\
& & \frac{\sd(\mc{P}_b)}{\sd(\mc{P}_b)^{-[\hat{0}< \hat{1}]}}\wedge \frac{\sd(\mc{P}_{r_1})}{\sd(\mc{P}_{r_1})^{-[\hat{0}< \hat{1}]}} \wedge \dots \wedge \frac{\sd(\mc{P}_{r_b})}{\sd(\mc{P}_{r_b})^{-[\hat{0}< \hat{1}]}}.
\end{tikzcd}$$

Explicitly, $\Delta_y$ can be described in simplicial degree $d$ as follows. Following \Cref{not:barycentric com}, a $d$-simplex in $\sd(\mc{P}_r)\big/\sd(\mc{P}_r)^{-[\hat{0}<\hat{1}]}$ corresponds to a nested chains of partitions
$$
\big(\sigma, S\big)=\big([\hat{0}=x_0<\dots <x_t=\hat{1}], S_0\subseteq \dots\subseteq S_d\big).
$$
Then $\Delta_y(\sigma, S)$ is the basepoint if $\sigma$ is not $y$-branched. If $\sigma$ is $y$-branched, then it can be ungrafted along $y$, as in Construction \ref{cons:cutting trees}. Write $\sigma_{\ul{b}}$ and $\sigma_{\ul{r}_i}$ for the resulting nondegenerate chains of partitions of $\ul{b}$ and $\ul{r}_i$. These chains are indexed by quotients of $S_d$ (because we have divided out degeneracies), which we will denote by $\pi_{\ul{b}}\colon S_d\twoheadrightarrow S_{\ul{b}, d}$ and $\pi_{\ul{r}_i}\colon S_{d}\twoheadrightarrow S_{\ul{r}_i, d}$. The subsets $S_\alpha\subseteq S_d$ are then sent to subsets $S_{\ul{b}, \alpha}\subseteq S_{\ul{b}, d}$ and $S_{\ul{r}_i, \alpha}\subseteq S_{\ul{r}_i, d}$ as follows: $S_{\ul{b}, \alpha}$ is simply the image $\pi_{\ul{b}}(S_\alpha)\subseteq S_{\ul{b}, d}$. On the other hand, for each $\ul{r}_i$,   $S_{\ul{r}_i, \alpha}\subseteq S_{\ul{r}_i, d}$ is the subset of all $\pi_{\ul{r}_i}(x)$ with $x\in S_\alpha$ such that the partition $x$ has $\ul{r}_i$ as a union of some of its classes.\vspace{4pt}

In terms of trees, $\sigma$ determines a levelled tree and each $S_\alpha$ marks some of its levels. We then mark the same levels in each of the branches and the trunk obtained by ungrafting along $y$ (and removing degeneracies). Furthermore, one marks the top level of the trunk if in the original tree, there was a level marked by $S_\alpha$ entirely above the ungrafting line. 
\end{cons}
\begin{proposition}\label{prop:subdivided bar is a cooperad}
The maps $\Delta_y$ endow $\sdBar(\Com^{\mm{nu}})$ with a cooperad structure.
\end{proposition}
\begin{proof}
We have to check that for any two partitions $z\leq y$, cocomposition along $y$ and $z$ is coassociative. For any $(\sigma, S)$ the image under the two cocompositions $(1\circ \Delta_z)\Delta_y$ and $(\Delta_y\circ 1)\Delta_z$ is the basepoint unless $\sigma$ corresponds to a tree that can both be ungrafted along $y$ and $z$. If $\sigma$ can both be ungrafted along $y$ and $z$, $(1\circ \Delta_z)\Delta_y(\sigma)$ and $(\Delta_y\circ 1)\Delta_z(\sigma)$ have the same underlying chains, by coassociativity of ungrafting. 
Furthermore, the explicit description of the subsets of `marked levels' in each of these trees shows that these marked levels are the same when we first ungraft along $y$ and then along $z$  or vice versa.
\end{proof}
The $R$-linearisation  of a cooperad in pointed simplicial sets is a cooperad in simplicial $R$-modules, and taking $R$-linear duals gives a cosimplicial restricted operad.
\begin{definition}\label{def:simplicial-cosimplicial part lie}
Let $R$ be a coherent ring. Write $\partLieR$ for the cosimplicial restricted operad over $R$ given by
$$
\partLieR=\Map_\ast(\sdBar(\Com^{\mm{nu}}), R).
$$
\end{definition}
In particular, we see that the cosimplicial $R$-module $\partLieR(r)^d$ has a basis given by nested chains of partitions of $\underline{r}$:
$$ 
\big(\sigma, S\big)=\big([\hat{0}=x_0<\dots <x_t=\hat{1}], S_0\subseteq \dots\subseteq S_d\big).
$$
Given $\rho\in \Sigma_r$, we will write $\rho^{-1}(\sigma, S)$ for the nested chain of partitions of $\ul{r}$ obtained by \emph{restricting} each partition $x_i$ along $\rho\colon \ul{r}\to \ul{r}$.
\begin{remark}
Note that $\sdBar(\Com^{\mm{nu}})(r)$ is given in each simplicial degree by a finite $\Sigma_r$-set. Consequently, the cosimplicial restricted operad $\partLieR(r)$ has an underlying symmetric sequence that is tamely cofibrant, so that the category of algebras over $\partLieR$ carries a semi-model structure.
\end{remark}
\begin{theorem}[Simplicial-cosimplicial models for partition Lie algebras]\label{thm:simp-cosimp partition lie}
The cosimplicial restricted operad $\dgob{Lie}^{\pi}_{R,\Delta}$ is a model for the derived partition Lie PD $\infty$-operad $\mm{Lie}^\pi_{R,\Delta}$ of Definition \ref{def:derived partition lie operad}. Consequently, there is an equivalence of $\infty$-categories
$$\begin{tikzcd}
\mathbf{Alg}_{\partLieR}^{\simcos,\res}[W_\mm{tame}^{-1}]\arrow[r, "\simeq"] & \Alg_{\Lie^\pi_{R,\Delta}}^{\gen,\pd}(\QC^\vee_R).
\end{tikzcd}$$
In particular, when $R=k$ is a field, the localisation of the category of simplicial-cosimplicial restricted algebras over $\mathbf{Lie}_{k,\Delta}^\pi$ at the weak equivalences is equivalent to the $\infty$-category of partition Lie algebras from \cite[Definition 5.47]{brantner2019deformation}.
\end{theorem} 

\begin{proof}
The fact that $\partLieR$ is a cosimplicial model for   $\mm{Lie}^\pi_{R,\Delta}$ will follow from our point-set description of Koszul dual PD $\infty$-operads below (\Cref{thm:simp-cosimp KD}). The assertion about algebras then follows from \Cref{thm:simp-cosimp models for derived algebras} and the final conclusion follows from Corollary \ref{cor:derived partition lie monad}.
\end{proof}

Adapting \cite{ikonicoff} to our setting, we can make the  simplicial-cosimplicial  $\partLieR$-algebras appearing in \Cref{thm:simp-cosimp partition lie} more explicit:

\begin{cons}\label{derivedPLAexplicit}
 Let $R$ be a field. A simplicial-cosimplicial restricted $\partLieR$-algebra is a simplicial object in cosimplicial restricted $\partLieR$-algebras.
To equip a cosimplicial $R$-module 
$$\begin{tikzcd}
 \mathfrak{g}^0\arrow[r, yshift=0.5ex]\arrow[r, yshift=-0.5ex] & \mathfrak{g}^1 \arrow[r, yshift=1ex]\arrow[r, yshift=-1ex]\arrow[r] & \mathfrak{g}^2\arrow[r, yshift=1.5ex]\arrow[r, yshift=0.5ex] \arrow[r, yshift=-0.5ex]\arrow[r, yshift=-1.5ex]  & \ldots 
\end{tikzcd}$$
with the structure of a  \textit{cosimplicial restricted $\partLieR$-algebra}, we must first define a $\partLieR$-algebra structure on $\mathfrak{g}$.

This means that for any nested chain $$(\sigma, S)= \big([\hat{0}=x_0<\dots <x_t=\hat{1}], S_0\subseteq \dots\subseteq S_d\big) \in \partLieR(r)^d$$ and any tuple $\mathbf{a} = (a_1,\ldots,a_r)$ in $\mathfrak{g}^d$, we must specify an element
$$ \{a_1, \dots, a_r\}_{(\sigma, S)} \in \mathfrak{g}^d,$$
depending linearly on each entry of $\mathbf{a}$, satisfying the following properties:
\begin{enumerate}
\item The canonical generator $(\sigma, S)\in \partLieR(1)^d\cong R$, given by $\sigma=[\hat{0}=x_0=\hat{1}]$ and subsets $S_0=\dots=S_d=\{0\}$, acts as $\{a\}_{(\sigma, S)} = a$ for any $a \in \mathfrak{g}^d$;
\item Given a permutation $\rho \in \Sigma_r$ and a tuple $\mathbf{a} = (a_1,\ldots,a_r)$ as above, we have 
$$  \{a_{\rho(1)}, \dots, a_{\rho(r)}\}_{(\sigma, S)}  =  \{a_1, \dots, a_r\}_{\rho^{-1}(\sigma, S)} $$
\item Given a partition $y: \underline{r} \cong \underline{r}_1 \sqcup \ldots \sqcup \underline{r}_b$, nested chains $(\tau,T)$ of $\underline{b}$,   \mbox{$(\sigma^i,S^i)$ of $\underline{r}_i$,}
and  tuples $\mathbf{a}_i = (a_1^i,\ldots a_{r_i}^i)$ for all $i=1,\ldots b$ corresponding to a tuple \mbox{$\mathbf{a} = (a_1,\ldots,a_r)$} under $y$, we have
$$ \left\{\{\mathbf{a}_1\}_{(\sigma^1,S^1)}, \ldots , \{\mathbf{a}_b\}_{(\sigma^b,S^b)}\right\}_{(\tau,T)} = \sum_{\Delta_y(\sigma,S) = ((\tau,T), (\sigma^1,S^1), \ldots, (\sigma^b,S^b))} \{ \mathbf{a}\}_{(\sigma,S)}$$
\item Let $\psi_*\colon \mf{g}^d\rt \mf{g}^{d'}$ be the cosimplicial structure map induced by $\psi\colon [d]\to [d']$ in $\Del$. Then
$$
\psi_*\{a_1, \dots, a_r\}_{(\sigma, S)} = \sum_{\psi^*(\tau, T)=(\sigma, S)}\left\{\psi_*a_1, \dots, \psi_*a_r\right\}_{(\tau, T)}.
$$
Here the sum runs over all nested chains $(\tau, T)=(\tau, T_0\subseteq \dots\subseteq T_{d'})$ in $\partLieR(r)^{d'}$ such that $(\sigma, S)=\big(\tau\big|_{T_{\psi(d)}}, T_{\psi(0)}\subseteq \dots\subseteq T_{\psi(d)}\big)$.
\end{enumerate}

Moreover, for any tuple $\mathbf{a}=(a_1, \dots, a_r)$ in $\mf{g}^d$ with stabiliser group $\Sigma_{\mathbf{a}} \leq  \Sigma_r$ and any nested chain $(\sigma,S) \in \partLieR(r)^d$, we must specify an element
$$\gamma_{(\sigma, S)}(a_1, \dots, a_r) \in \mathfrak{g}^d$$
These `divided operations' must satisfy the following properties:
\begin{enumerate} \setcounter{enumi}{4}
\item  Let $\Sigma_{\mathbf{a}, \sigma}=\Sigma_{\mb{a}}\cap \mm{Stab}(\sigma)< \Sigma_r$ be the group of symmetries fixing $\mathbf{a}$ and $\sigma$. Then
$$
\{a_1, \dots, a_r\}_{(\sigma, S)} = \big|\Sigma_{\mathbf{a}, \sigma}\big|\cdot\gamma_{(\sigma, S)}(a_1, \dots, a_r).
$$
\item For any permutation $\rho \in \Sigma_r$, we have
$$\gamma_{(\sigma, S)}(a_{\rho(1)}, \dots, a_{\rho(r)})= \gamma_{\rho^{-1}(\sigma, S)}(a_1, \dots, a_r).$$
\item Suppose that $\mb{a}$ contains (at least) $i$ copies of the same element $a$, indexed by a subset $\ul{i}\subseteq \ul{r}$ and let $\mb{a}_{(\lambda, \ul{i})}$ be the tuple obtained from $\mb{a}$ by scaling each of these $i$ copies by $\lambda\in R$. Writing $\Sigma_{\mb{a}, \mb{a}_{(\lambda, \ul{i})}, \sigma}<\Sigma_r$ for the subgroup of permutations fixing $\mb{a}, \mb{a}_{(\lambda, \ul{i})}$ and $\sigma$, we have
$$
\frac{|\Sigma_{\mb{a}_{(\lambda, \ul{i})}, \sigma}|}{|\Sigma_{\mb{a}, \mb{a}_{(\lambda, \ul{i})}, \sigma}|}\gamma_{(\sigma, S)}(\mathbf{a}_{(\lambda, \ul{i})}) =\lambda^i \cdot \frac{|\Sigma_{\mb{a}, \sigma}|}{|\Sigma_{\mb{a}, \mb{a}_{(\lambda, \ul{i})}, \sigma}|}\gamma_{(\sigma, S)}(\mathbf{a}).
$$
\item Suppose that $\mb{a}$ contains (at least) $i$ copies of the same element $a=b+c$, indexed by a subset $\ul{i} \subseteq \ul{r}$. For each decomposition $i=j+k$,   form a new tuple
$\mathbf{a}_{(j,k)}$ from $\mathbf{a}$ by replacing the first $j$ copies of $a$ by $b$ and the last $k$ copies of $a$ by $c$. Write $\Sigma_{\mb{a},\ul{i}, \sigma}<\Sigma_r$ for the subgroup of permutations fixing $\mb{a}, \sigma$ and the subset $\ul{i}$, and $\Sigma_{\mb{a}, \ul{j}, \ul{k},\sigma}<\Sigma_r$ for the permutations fixing $\mb{a}, \sigma$ and the two disjoint subsets of $i$ containing its first $j$ and last $k$ elements. We then have
$$
\frac{|\Sigma_{\mb{a}, \sigma}|}{|\Sigma_{\mb{a}, \ul{i}, \sigma}|}\gamma_{(\sigma, S)}(\mathbf{a})  = \sum_{i = j+k}  \sum_{\rho \in _{\Sigma_{\mb{a}, \sigma}}\backslash \Sigma_{\mathbf{a}}/_{\Sigma_{\mathbf{a}_{(j,k)}}}}\frac{|\Sigma_{\mb{a}_{(j, k)}, \rho^{-1}(\sigma)}|}{|\Sigma_{\mb{a}, \ul{j}, \ul{k},\rho^{-1}(\sigma)}|}\gamma_{\rho^{-1}(\sigma,S)}(\mathbf{a}_{(j,k)}). $$
 \item Fix a partition $y\colon \underline{r} \cong \underline{r}_1 \sqcup \ldots \sqcup \underline{r}_b$, as well as nested chains $(\sigma^i,S^i) \in  \partLieR(r_i)^d $  and 
$(\tau,T) \in  \partLieR(b)^d $. Given a tuple $\mbox{$\mathbf{a} = (a_1,\ldots,a_r)$}$, let $\mathbf{a}^i = (a_1^i,\ldots a_{r_i}^i)$ denote the tuples corresponding to the classes $i=1, \dots, b$ of the partition $y$. Let us consider the following subgroups of permutations of $\ul{b}$
$$
\Sigma_{(\mathbf{a}^i), (\sigma^i, S^i), \tau} < \Sigma_{(\gamma(\mb{a}^i)), \tau} < \Sigma_b
$$
Here the left subgroup contains permutations that fix $\tau$, the family of tuples $(\mb{a}^i)_{i\in \ul{b}}$ and the family of nested partitions $(\sigma^i, S^i)$. The second subgroups contains the permutations that fix $\tau$ and the tuple $\big(\gamma_{(\sigma^i, S^i)}(\mb{a}^i)\big)$. 

We then have
$$ \frac{|\Sigma_{(\gamma(\mathbf{a}^i)), \tau}|}{|\Sigma_{(\mathbf{a}^i), (\sigma^i, S^i), \tau}|} \cdot  \gamma_{(\tau,T)} \left(\gamma_{{(\sigma^1,S^1)}}(\mathbf{a}_1), \ldots , \gamma_{{(\sigma^b,S^b)}}(\mathbf{a}_b)\right) \ \ \ \  \ \ \ \ \ \ \  \ \ \ \  \ \ \ \ \ \ \  \ \ \ \  \ \ \ \ \ \ \ $$ $$ \ \ \ \  \ \ \ \ \ \ \  \ \ \ \  \ \ \ \ \ \ \ =\sum_{\Delta_y(\sigma,S)=  ((\tau,T), (\sigma^1,S^1), \ldots, (\sigma^b,S^b))} \frac{| \Sigma_{\mb{a},\sigma}|}{| \Sigma_{(\mathbf{a}^i), (\sigma^i, S^i), \tau}\ltimes \prod_{i\in \ul{b}} \Sigma_{\mb{a}_i, \sigma_i}|} \cdot 
 \gamma_{(\sigma,S)}(\mathbf{a}).$$
\item Let $\psi_*\colon \mf{g}^d\rt \mf{g}^{d'}$ be the cosimplicial structure map induced by $\psi\colon [d]\to [d']$ in $\Del$. Then
$$
\frac{|\Sigma_{\psi_*(\mb{a}), \sigma}|}{|\Sigma_{\mb{a}, \sigma}|}\psi_*\gamma_{(\sigma, S)}(a_1, \dots, a_r) = \sum_{\psi^*(\tau, T)=(\sigma, S)}\frac{|\Sigma_{\psi_*(\mb{a}), \sigma}|}{|\Sigma_{\psi_*(\mb{a}), \tau}|} \gamma_{(\tau, T)}(\psi_*a_1, \dots, \psi_*a_r).
$$
Here the sum runs over all nested chains $(\tau, T)=(\tau, T_0\subseteq \dots\subseteq T_{d'})$ in $\partLieR(r)^{d'}$ such that $(\sigma, S)=\big(\tau|_{T_{\psi(d)}}, T_{\psi(0)}\subseteq \dots\subseteq T_{\psi(d)}\big)$. Note that for each such nested pair $(\tau, T)$, there are subgroup inclusions $\Sigma_{\mb{a}, \sigma} < \Sigma_{\psi_*(\mb{a}), \sigma}$ and $\Sigma_{\psi_*(\mb{a}), \tau}<\Sigma_{\psi_*(\mb{a}), \sigma}$
induced by the inclusions $\Sigma_{\mb{a}}<\Sigma_{\psi_*(\mb{a})}$ and $\mm{Stab}(\tau)<\mm{Stab}(\sigma)$.
\end{enumerate}
\end{cons}
\begin{remark}
The description of a cosimplicial restricted $\ope{Lie}^\pi_{R, \Delta}$-algebras in Construction \ref{con:derived products} is an application of the explicit description of restricted algebras over operads given in \cite[Definition 4.1.1, 4.1.4]{ikonicoff}. To translate between the two descriptions, observe that any tuple of elements $(a_1, \dots, a_r)$ is uniquely determined by a partition $y\colon \ul{r}\cong \ul{r}_1\sqcup\dots\sqcup \ul{r}_b$ on $\ul{r}$, together with a $b$-tuple of \emph{mutually distinct} elements $(a_{\ul{r}_1}, \dots, a_{\ul{r}_b})$ of $A$. 

The element denoted $\gamma_{(\sigma, S)}(a_1, \dots, a_r)$ above then corresponds to $\gamma_{[\sigma, S]_{y}, y}(a_{\ul{r}_1}, \dots, a_{\ul{r}_b})$ in the notation of Ikonicoff. By \cite[Definition 4.1.1(3)]{ikonicoff}, these $\gamma_{[\sigma, S]_{y}, y}(a_{\ul{r}_1}, \dots, a_{\ul{r}_b})$ with all $a_{\ul{r}_i}$ mutually distinct determine all other operations appearing in loc.\ cit.
\end{remark}

\subsection{Explicit Koszul duality}
We will now give a simplicial-cosimplicial model for the Koszul dual of more  general augmented derived $\infty$-operads, similarly to the simplicial-cosimplicial model for the partition Lie PD $\infty$-operad in Definition \ref{def:simplicial-cosimplicial part lie}.

\begin{definition}
An augmented (simplicial-cosimplicial) operad $\ope{P}$ over $R$ is called \emph{reduced} if $\ope{P}(0)=0$ and $\ope{P}(1)=R\cdot 1$.
\end{definition}
We start by fixing a reduced operad $\ope{P}$ in discrete $R$-modules and consider the simplicial bar construction
$$\begin{tikzcd}
{\Barr(\mb{1}, \ope{P}, \mb{1})= \dots }\arrow[r, yshift=1.5ex]\arrow[r, yshift=0.5ex] \arrow[r, yshift=-0.5ex]\arrow[r, yshift=-1.5ex] & \ope{P}\circ \ope{P}\arrow[r, yshift=1ex]\arrow[r, yshift=-1ex]\arrow[r] & \ope{P}\arrow[r, yshift=0.5ex]\arrow[r, yshift=-0.5ex] & \mb{1}.
\end{tikzcd}$$
This is a simplicial symmetric sequence of $R$-modules which can be written explicitly as a direct sum
$$
\Barr(\mb{1}, \ope{P}, \mb{1})(r)_d=\bigoplus_{\sigma\in \Barr(\mb{1}, \Com^{\mm{nu}}, \mb{1})(r)_d} \ope{P}(\sigma).
$$
Here the direct sum is indexed by the non-basepoint $d$-simplices of the partition complex, i.e.\ by chains of partitions $\sigma=[\hat{0}=x_0\leq \dots\leq x_d=\hat{1}]$. Each such chain of partitions determines a levelled tree and we denote by $\ope{P}(\sigma)$ the $R$-module of labellings of this tree by elements of $\ope{P}$; in other words, it is a certain tensor product of $\ope{P}(r_\alpha)$ indexed by the vertices of the tree. The simplicial structure is obtained by removing levels and composing operations in $\ope{P}$, and produces zero if the result is no longer a tree with $r$ leaves.

As in the previous section, the barycentric subdivision of the simplicial bar construction can be endowed with the structure of a cooperad.
\begin{notation}[$R$-linear barycentric subdivision]\label{not:linearised barycentric subdivision}
Recall that there is an adjoint pair $\sd\colon \mm{sSet}\leftrightarrows \mm{sSet} \colon \mm{Ex}$ given by the barycentric subdivision and Kan's Ex-functor. The functor $\mm{Ex}$ preserves simplicial symmetric sequences of $R$-modules, and we will write $\sd\colon \mm{ssSeq}_R\rt \mm{ssSeq}_R$ for its left adjoint; in other words, this is the $R$-linear extension of the usual barycentric subdivision. The natural transformation $\mm{id}\rt \mm{Ex}$ is adjoint to an augmentation $\sd\rt \mm{id}$.

For any reduced operad in discrete $R$-modules, we will then write $\sdBar(\ope{P})$ for the $R$-linear barycentric subdivision of $\Barr(\mb{1}, \ope{P}, \mb{1})$. Explicitly, this is given by
$$
\sdBar(\ope{P})(r)_d=\bigoplus_{(\sigma, S)\in \sdBar(\Com^{\mm{nu}})(r)_d} \ope{P}(\sigma)
$$
where the sum runs over all simplices in the (set-valued) barycentric subdivision from Notation \ref{not:barycentric com}. Such simplices correspond to nested nondegenerate chains of partitions $\sigma=[\hat{0}=x_0<\dots<x_t=\hat{1}]$ with $S_0\subseteq \dots\subseteq S_d=\{0, \dots, t\}$. Here we allow $t=-1$, corresponding to the basepoint in $\sdBar(\Com^{\mm{nu}})(r)$; the corresponding summand $\ope{P}(\sigma)$ is zero in this case by definition.
\end{notation}
\begin{remark}\label{rem:linear barycentric subdivision}
For any simplicial set $K$ and any cosimplicial symmetric sequence $P$, one has that $\sd(K_+\wedge P)\cong \sd(K)_+\wedge P$. Because the tame model structure on simplicial-cosimplicial symmetric sequences is simplicial, one sees that for each cofibrant cosimplicial symmetric sequence $P$, the augmentation 
$$
\sd(K_+\wedge P)\cong \sd(K)_+\wedge P\rt K_+\wedge P
$$
is a tame weak equivalence. Furthermore, the description of the generating (trivial) cofibrations for the tame (or projective) model structure on $\mathbf{sSeq}_{R}^{\simcos}$ (see Remark \ref{rem:generating cofibrations simplicial-cosimplicial}) shows that $\sd\colon \mathbf{sSeq}_{R}^{\simcos}\rt \mathbf{sSeq}_{R}^{\simcos}$ is a left Quillen functor. Combining these two observations, one sees that the augmentation $\sd(X)\rt X$ is a weak equivalence for every tamely cofibrant sc-symmetric sequence. By adjunction, this means that $Y\rt \mm{Ex}(Y)$ is a weak equivalence for every fibrant sc-symmetric sequence.
\end{remark}
\begin{cons}[Cocomposition on the subdivided bar construction]\label{cons:cocomp for partitions general}
Let $\ope{P}$ be a reduced operad in discrete $R$-modules and let $y$ be a partition of the form $\ul{r}\cong \ul{r}_1\sqcup \dots\sqcup \ul{r}_b$. We define a map
\begin{equation}\label{eq:ungrafting for P}\begin{tikzcd}
\Delta_y\colon \sdBar(\ope{P})(r)\arrow[r] & \sdBar(\ope{P})(b)\otimes \sdBar(\ope{P})(r_1)\otimes \dots\otimes \sdBar(\ope{P})(r_b)
\end{tikzcd}\end{equation}
as follows. Note that the domain is a direct sum indexed by $\sdBar(\Com^{\mm{nu}})(r)$ while the target in a direct sum indexed by $\sdBar(\Com^{\mm{nu}})(b)\wedge \sdBar(\Com^{\mm{nu}})(r_1)\wedge\dots\wedge \sdBar(\Com^{\mm{nu}})(r_b)$. Then $\Delta_y$ sends the summand indexed by $(\sigma, S)$ to the summand indexed by $\Delta_y(\sigma, S)$, for the comultiplication of Construction \ref{cons:cocomp for partitions}. In particular, $\Delta_y$ sends the summand by $(\sigma, S)$ to zero if $\sigma$ is not $y$-branched. 

If $\sigma$ is $y$-branched, it corresponds to a tree than can be ungrafted along $y$ and $\Delta_y(\sigma, S)$ is the tuple consisting of the branches and trunk of this ungrafted tree. Note that the branches and trunk of the tree associated to $\sigma$ together contain exactly the same vertices labelled by non-identity operations in $\ope{P}$ as the tree $\sigma$ itself. Consequently, the $(\sigma, S)$-summand in $\sdBar(\ope{P})(r)$ is naturally isomorphic to the $\Delta_y(\sigma, S)$-summand in $\sdBar(\ope{P})(b)\otimes \sdBar(\ope{P})(r_1)\otimes \dots\otimes \sdBar(\ope{P})(r_b)$. We then define $\Delta_y$ to be this natural isomorphism.

In other words, the map \eqref{eq:ungrafting for P} sends a levelled tree $\sigma$ with vertices labelled by $\ope{P}$ (and a nested family $S_0\subseteq \dots \subseteq S_d$ of marked levels) to zero if it cannot be ungrafted along $y$, and to the ungrafting if it can.
\end{cons}
As in Proposition \ref{prop:subdivided bar is a cooperad}, the associativity of ungrafting then shows:
\begin{corollary}
Let $\ope{P}$ be a reduced operad in discrete $R$-modules. The operations $\Delta_y$ endow $\sdBar(\ope{P})$ with the structure of a simplicial cooperad in $R$-modules.
\end{corollary}
\begin{remark}
When $\ope{P}=\Com^{\mm{nu}}_R$ is the nonunital $R$-linear commutative operad, $\sdBar(\Com^{\mm{nu}}_R)\cong R\wedge \sdBar(\Com^{\mm{nu}})$ is simply the $R$-linear extension of the cooperad in pointed simplicial sets from Proposition \ref{prop:subdivided bar is a cooperad}.
\end{remark}
\begin{definition}[Subdivided bar construction of reduced simplicial-cosimplicial operads]
Suppose that $\ope{P}$ is a reduced simplicial-cosimplicial $R$-linear operad. We will write $\sdBar(\ope{P})$ for the simplicial-cosimplicial cooperad given by
$$
\sdBar(\ope{P})_{d}^n = \sdBar(\ope{P}_d^n)_d
$$
with cocomposition given in each simplicial-cosimplicial degree as in Construction \ref{cons:cocomp for partitions general}.
The simplicial-cosimplicial \emph{dual restricted operad} is the $R$-linear dual
$$
\mb{D}_{\Delta}(\ope{P})=\sdBar(\ope{P})^\vee.
$$
\end{definition}
\begin{theorem}\label{thm:simp-cosimp KD}
Let $\ope{P}$ be a reduced sc-operad with projectively cofibrant underlying symmetric sequence   (Definition \ref{def:model structure simp-cosimp sym seq}). Let $\PP$ denote the associated augmented derived $\infty$-operad. Then the sc-restricted  operad $\mb{D}_\Delta(\ope{P})$ is a model for the   derived PD $\infty$-operad $\KD^{\pd}(\cat{P})$.
\end{theorem}
The proof requires a preliminary construction:
\begin{cons}[Simplicial-cosimplicial Koszul complex]
Suppose that $\ope{P}$ is a reduced operad in discrete $R$-modules. We will define the \emph{subdivided Koszul complex} of $\ope{P}$ to be the symmetric sequence
$$
\sdK(\ope{P})=\sd\big(\Barr(\mb{1}, \ope{P}, \ope{P})\big).
$$
Explicitly, $\sdK(\ope{P})(r)_d\cong \bigoplus_{(\sigma, S)}\ope{P}(\sigma)$ where the sum runs over chains of partitions $\sigma=[x_{-1}=\hat{0}\leq x_0<\dots <x_t=\hat{1}]$ with a nested family of subsets $S_0\subseteq \dots \subseteq S_d=\{-1, \dots, t\}$ \emph{which all contain $-1$}. The simplicial structure maps act on the nested family of subsets $S_0\subseteq \dots \subseteq S_d$ in the evident way.

In terms of levelled trees, a $d$-simplex of $\sdK(\ope{P})(r)$ consists of a levelled tree with vertices marked by $\ope{P}$, together with a nested family of marked levels, which all contain the top level (i.e.\ the leaves). These levelled trees are almost nondegenerate: one only allows the leaf vertices to all be equal to the identity.

The bar construction $\Barr(\mb{1}, \ope{P}, \ope{P})$ carries a natural right $\ope{P}$-action. This induces a right $\ope{P}$-module structure on $\sdK(\ope{P})$. In terms of levelled trees, this action simply precomposes the leaf vertices labelled by $\ope{P}$ with operations from $\ope{P}$.

On the other hand, $\sdK(\ope{P})$ is a left comodule over $\sdBar(\ope{P})$. Indeed, for every partition $y$ of the form $\ul{r}\cong \ul{r}_1\sqcup \dots \sqcup \ul{r}_b$, there is a comultiplication map
$$\begin{tikzcd}
\Delta_y\colon \sdK(\ope{P})(r)\arrow[r] & \sdBar(\ope{P})(b)\otimes \sdK(\ope{P})(r_1)\otimes\dots\otimes \sdK(\ope{P})(r_b)
\end{tikzcd}$$
defined in exactly the same way as in Construction \ref{cons:cocomp for partitions general}: a levelled tree with vertices marked by $\ope{P}$ is sent to its ungrafting along $y$ if this is possible, and to zero otherwise. Furthermore, the subsets $S_0\subseteq \dots \subseteq S_d$ give rise to subsets of levels for each of the branches and the trunk of the resulting ungrafted tree. Note that the left comodule structure and right $\ope{P}$-module structure commute.

More generally, if $\ope{P}$ is an sc-operad, we define $\sdK(\ope{P})$ by the diagonal
$$
\sdK(\ope{P})_{d}^n=\sdK(\ope{P}_d^n)_d.
$$
This carries a commuting left comodule structure of $\sdBar(\ope{P})$ and a right $\ope{P}$-module structure.
\end{cons}
\begin{proof}[Proof of \Cref{thm:simp-cosimp KD}]
Consider the natural map of right $\ope{P}$-modules $\pi\colon \sdK(\ope{P})\rt \mb{1}$ sending all summands indexed by $(\sigma, S)$ with $\sigma=[\hat{0}=x_{-1}\leq x_0<\dots <x_t=\hat{1}]$ to zero, except the summand in arity $1$ and simplicial-cosimplicial degree zero corresponding to $\hat{0}=x_{-1}\leq x_0=\hat{1}$; this summand is given by $\ope{P}(1)\cong R\cdot 1$. We claim this $\pi$ is a weak equivalence and that the map
\begin{equation}\label{eq:subdivided koszul quotient}\begin{tikzcd}
\sdK(\ope{P})\circ^h_{\ope{P}} \mb{1} \arrow[r] & \sdBar(\ope{P})\circ \sdK(\ope{P})\circ^h_{\ope{P}} \mb{1} \arrow[r] & \sdBar(\ope{P})\circ \mb{1}=\sdBar(\ope{P})
\end{tikzcd}\end{equation}
is a weak equivalence as well. By Proposition \ref{prop:bar construction as coendomorphisms}, this implies that $\sdBar(\ope{P})$ is a model for the $\infty$-categorical bar construction of $\ope{P}$. If $\ope{P}$ is projectively cofibrant as an sc-symmetric sequence, then $\sdBar(\ope{P})$ is easily seen to be tamely cofibrant as an sc-symmetric sequence. It then follows from Lemma \ref{lem:simp-cosimp day convolution} that the $R$-linear dual $\mb{D}_\Delta(\ope{P})$ is a model for the Koszul dual derived PD $\infty$-operad $\KD^{\pd}(\PP)$.

It remains to verify the claim, for which it suffices to treat the case where $\ope{P}$ is a projectively cofibrant cosimplicial operad; the case of a general sc-operad follows by taking geometric realisations (Remark \ref{rem:realisations and totalisations in simp-cosimp modules}). For a reduced cosimplicial operad $\ope{P}$, let $\RMod_{\ope{P}}(\mathbf{sSeq}_{R}^{\simcos})$ be the category of right $\ope{P}$-modules in sc-symmetric sequences. Because $\ope{P}$ is cofibrant as a symmetric sequence, this carries a simplicial model structure whose fibrations and weak equivalences are fibrations and weak equivalences on the underlying symmetric sequences, as in Definition \ref{def:model structure simp-cosimp sym seq}. Now note that the map $\pi$ factors into natural maps of sc-symmetric sequences (equipped with a right $\ope{P}$-action)
$$\begin{tikzcd}
\pi\colon \sdK(\ope{P})\arrow[r] & \Barr(\mb{1}, \ope{P}, \ope{P})\arrow[r] & \mb{1}.
\end{tikzcd}$$
The map $\Barr(\mb{1}, \ope{P}, \ope{P})\rt \mb{1}$ is the usual augmentation of the bar construction, which gives a cofibrant replacement in the simplicial model category $\RMod_{\ope{P}}(\mathbf{sSeq}_{R}^{\simcos})$; in particular, it is a weak equivalence. 

The first map is the canonical augmentation of the linearised barycentric subdivision (Notation \ref{not:linearised barycentric subdivision}). Recall that this map is adjoint to a natural transformation $\theta\colon \mm{id}\rt \mm{Ex}$ of functors $\RMod_{\ope{P}}(\mathbf{sSeq}_{R}^{\simcos})\rt \RMod_{\ope{P}}(\mathbf{sSeq}_{R}^{\simcos})$. The map $\theta$ is a natural transformation between right Quillen functors which is an weak equivalence on fibrant objects, since it is at the level of the underlying symmetric sequences (Remark \ref{rem:linear barycentric subdivision}). This implies that $\sd\rt \mm{id}$ is a natural transformation of left Quillen functors $\RMod_{\ope{P}}(\mathbf{sSeq}_{R}^{\simcos})\rt \RMod_{\ope{P}}(\mathbf{sSeq}_{R}^{\simcos})$, which is a weak equivalence on cofibrant objects. In particular, the map $\sdK(\ope{P})\rt \Barr(\mb{1}, \ope{P}, \ope{P})$ is a weak equivalence between cofibrant left $\ope{P}$-modules.

This shows that $\sdK(\ope{P})\rt \mb{1}$ is a cofibrant resolution of $\mb{1}$ as a right $\ope{P}$-module. The map \eqref{eq:subdivided koszul quotient} can then be identified with the map
$$\begin{tikzcd}
\sdK(\ope{P})\circ_{\ope{P}} \mb{1} \arrow[r, "\Delta"] & \sdBar(\ope{P})\circ \sdK(\ope{P})\circ_{\ope{P}} \mb{1} \arrow[r] & \sdBar(\ope{P})\circ \mb{1}=\sdBar(\ope{P}).
\end{tikzcd}$$
Since $\sdK(\ope{P})\circ_{\ope{P}} \mb{1}\cong \sdBar(\ope{P})$, this map is readily verified to be an isomorphism.
\end{proof}
\begin{corollary}\label{cor:KD partition Lie}
Let $R$ be a coherent ring and let $A$ be a nonunital simplicial commutative $R$-algebra which is cofibrant as a simplicial $R$-module. Then the Koszul dual partition Lie algebra $\KD(A)\in \Alg_{\Lie^\pi_{R,\Delta}}^{\gen,\pd}(\QC^\vee_R)$ can be modelled by the cosimplicial restricted $\partLieR$-algebra
\begin{equation}\label{eq:KD partition Lie explicit}
\big(\sdK(\Com^{\mm{nu}}_R)\circ_{\Com^{\mm{nu}}_R} A\big)^\vee
\end{equation}
whose restricted $\partLieR$-algebra structure arises by $R$-linear duality from the $\sdBar(\Com^{\mm{nu}}_R)$-coalgebra structure on $\sdK(\Com^{\mm{nu}}_R)\circ_{\Com^{\mm{nu}}_R} A$ .
\end{corollary}
\begin{proof}
The proof of \Cref{thm:simp-cosimp KD} shows that $\sdK(\Com^{\mm{nu}}_R)$, equipped with its commuting left $\sdBar(\Com^{\mm{nu}}_R)$-comodule structure and its right $\Com^{\mm{nu}}_R$-module structure, is a model for the Koszul complex of the derived $\infty$-operad $\mm{Com}^\mm{nu}$ (\Cref{cons:koszul complex}). Since $\sdK(\Com^{\mm{nu}}_R)$ is cofibrant as a right $\Com^{\mm{nu}}_R$-module and $A$ is cofibrant as an $R$-module, it then follows that the $\sdBar(\Com^{\mm{nu}}_R)$-coalgebra $\sdK(\Com^{\mm{nu}}_R)\circ_{\Com^{\mm{nu}}_R} A$ is a model for the bar construction of $A$, whose underlying simplicial $R$-module is cofibrant. The $R$-linear dual is then a cosimplicial restricted $\partLieR$-algebra that models $\KD(A)$.
\end{proof}
\begin{remark}
Let $R$ be a coherent ring and let $A$ be a discrete nonunital commutative $R$-algebra that is projective as an $R$-module. Then the Koszul cosimplicial restricted $\partLieR$-algebra $\mf{g}$ from \eqref{eq:KD partition Lie explicit} can be identified explicitly as follows. 
The $R$-module $\mf{g}^d$ consists of families of $R$-linear maps 
$$
\alpha=(\alpha_{(\tau, T)}), \qquad\qquad\qquad \alpha_{(\tau, T)}\colon A^{\otimes j}\rt R
$$
where $(\tau, T)$ runs over the set of nested chains of partitions $\tau=[\hat{0}=x_0<\dots<x_t=\hat{1}]$, $T_0\subseteq \dots \subseteq T_d=\{1, \dots, t\}$ of the finite set $\ul{j}$, for all $j\geq 1$. Furthermore, this family of maps has to be invariant in the following sense: if $(\tau, T)$ and $(\tau', T')$ are two nested chains of partitions of $\ul{j}$ related by the action of some $\gamma\in \Sigma_j$, then $\alpha_{(\tau, T)}(a_1, \dots, a_j)=\alpha_{(\tau', T')}(a_{\gamma(1)}, \dots, a_{\gamma(j)})$.

To describe the cosimplicial structure maps, let us fix a map $f\colon [d]\to [d']$. Suppose that $\tau'=[\hat{0}=x_0<\dots<x_t=\hat{1}]$ and $T_0\subseteq \dots\subseteq T_{d'}=\{0, \dots, t\}$ form a nested chain of partitions of $\ul{j}$, and consider the partition $x_{\min(T_{f(d)})}$ of $\ul{j}$.
This partition induces a quotient map $q\colon \ul{j} \twoheadrightarrow \ul{b}$ and for each $i\in [d]$, the restricted chain $\sigma\big|T_{f(i)}$ induces a chain of partitions of on $\ul{b}$. Let us write $f^*(\tau, T)$ for the resulting nested chain of partitions of $\ul{b}$. For any element $\alpha\in \mf{g}^d$, its image $f_*\alpha\in \mf{g}^{d'}$ is then given by the family of maps
$$
(f_*\alpha)_{(\tau, T)} = \alpha_{f^*(\tau, T)}\big(\prod_{j\in q^{-1}(1)} a_j, \prod_{j\in q^{-1}(2)} a_j, \dots, \prod_{j\in q^{-1}(b)} a_j\big).
$$
Finally, the restricted $\partLieR$-algebra structure on $\mf{g}$ can be described as follows. For $(\sigma, S)\in \partLieR(r)^d$ and $\alpha_1, \dots, \alpha_r\in \mf{g}^d$, the element $\gamma_{(\sigma, S)}(\alpha_1, \dots, \alpha_r)$ is a tuple of maps
$$
\big(\gamma_{(\sigma, S)}(\alpha_1, \dots, \alpha_r)\big)_{(\tau, T)}\colon A^{\otimes j}\rt R.
$$
If $(\tau, T)$ can be cocomposed into $(\sigma, S)$ and branches $(\tau_1, T_1), \dots, (\tau_r, T_r)$, then the above map is a suitable symmetrisation of the product $(\alpha_1)_{(\tau_1, T_1)}\dots (\alpha_r)_{(\tau_r, T_r)}$. If there is no such cocomposition, the map is zero.
\end{remark}

\newpage
\appendix
\section{The PD surjections operad}\label{app:surjections cooperad}
The commutative operad $\Com_R$ admits various well-known explicit resolutions by $\Sigma$-cofibrant dg-operads, like the Barratt--Eccles operad and the surjections operad \cite{mcclure2003multivariable, bergerfresse}. In contrast, the dual problem of finding an explicit (combinatorial) $\Sigma$-cofibrant resolution of the nonunital cocommutative dg-cooperad $\dgob{coCom}_R^{\mm{nu}}$ has not yet been addressed in the literature.

The significance of such a $\Sigma$-cofibrant resolution comes from the Koszul dual problem of trying to find a cofibrant chain model for the Lie operad (this problem seems to be folklore, and is raised for instance in \cite{DehlingVallette15}). Indeed, a $\Sigma$-cofibrant model for the (non-counital) cocommutative cooperad gives rise to a cofibrant model for the Lie operad by the cobar construction.
A partial result in this direction appears in \cite[Proposition 2.3]{dehling2017weak}, 
where the author constructs  a certain $\Sigma$-cofibrant cooperad $\dgob{Lie}_3^\diamond$ and a map  $ \dgob{Lie}_3^{\diamond, s} \to \dgob{coCom}_R^{\mm{nu}}$ 
from its operadic suspension
which is a resolution in low degrees. 

The goal of this section is to present a solution to this problem by giving a construction of $\Surj_R$, the surjection dg-cooperad over a ring $R$, inspired by the surjections cooperad of McClure--Smith.

\begin{theorem}\label{thm:surjections cooperad}
Let $R$ be a commutative ring. There exists an explicit $\Sigma$-cofibrant dg-cooperad $\Surj_R$ in non-negative degrees, together with a quasi-isomorphism of dg-cooperads $\Surj_R \to \dgob{coCom}^{\mm{nu}}_R$ to the cooperad of nonunital cocommutative coalgebras over $R$. We will refer to $\Surj_R$ as the \emph{surjections cooperad}.
\end{theorem}

 We note that $\Surj_R(0)=0$ and there is no natural way of extending the cooperad structure in order to resolve $\dgob{coCom}_R$.
The rest of this section will be devoted to proving \Cref{thm:surjections cooperad}. We will first describe the underlying symmetric sequence of 
$\Surj_R$ (from which $\Sigma$-cofibrancy will be evident), then define a comultiplication on it, and finally prove that the structure described indeed  forms a cooperad. 
The existence of the quasi-isomorphism $\Surj_R \to \dgob{coCom}^{\mm{nu}}_R$ is then  evident.

\subsection*{The underlying complex}
The symmetric sequence underlying $\Surj_R$
agrees with the symmetric sequence underlying the (nonunital) \emph{surjections operad} of McClure--Smith \cite{mcclure2003multivariable}. We briefly recall its definition following the notation and conventions from Berger--Fresse \cite{bergerfresse} (who denote it by $\mathcal X$).
\begin{definition}\label{def:surjections module}
	Let $\ul{r}$ be a set with $r\geq 1$ elements and let $\big<r+d\big>$ be a linear order with $r+d$ elements; up to unique isomorphism, we   identify $\big<r+d\big>$ with   $\{1, \dots, r+d\}$. A map (of sets) $\mb{u}\colon \big<r+d\big>\rt \ul{r}$ can be identified with an (ordered) \emph{sequence} of elements in $\ul{r}$
	$$
	\mb{u}=\big(u_1, \dots, u_{r+d}\big).
	$$
	Such a sequence is said to be \emph{degenerate} if $\mb{u}\colon \big<r+d\big>\rt \ul{r}$ is not surjective or if it sends two consecutive elements in $\big<r+d\big>$ to the same element in $\ul{r}$.
	
	Let $\Surj_R(\ul{r})_d$ be the quotient of the free $R$-module on such sequences $\mb{u}$, by the submodule generated by the degenerate sequences. In other words, $\Surj_R(\ul{r})$ is freely generated by \emph{non-degenerate sequences}. The symmetric group $\mm{Aut}(\ul{r})$ acts in an obvious way on $\Surj_R(\ul{r})_d$. 
\end{definition}
\begin{remark}
	Non-degenerate sequences are often called  non-degenerate \emph{surjections} in \cite{bergerfresse}; we use the term sequences to  highlight the ordering, which becomes important later. 
\end{remark}
\begin{definition}[Caesuras]\label{def:caesura}
	Following \cite{bergerfresse}, we call an element $u_\alpha$ in a sequence $\mb{u}=\big(u_1, \dots, u_{r+d}\big)$ in $\ul{r}$  a \emph{caesura} if it is \emph{not} the last of occurrence of that element in the sequence. There are exactly $d$ such caesuras in the sequence. 
	We  write $\Caes{\mb{u}}$ for the set of caesuras in $\mb{u}$, with its natural linear order.
\end{definition} 
\begin{sgnrule}\label{sgn:differential sur}
	Let $\mb{u}=(u_1, \dots, u_{r+d})$ be a nondegenerate sequence in $\ul{r}$. We associate a sign $\pmu$ to each $u_\alpha$ in this sequence as follows. First consider all $\alpha$ for which the $u_\alpha$ are caesuras; these are given alternating signs $\pm$, \mbox{starting with $+$.} Next consider all $\alpha$ for which $u_\alpha$ occurs for the last time in the sequence; these $u_\alpha$ are given the sign opposite to the sign associated to the previous copy of the element $u_\alpha\in \ul{r}$ in the sequence (if there is no previous copy, we associate $0$ to it, although this will not play a role). For example, we have the following element in $\Surj_R(4)_4$:
	$$
	\big(\overset{+}{1}\overset{-}{3}\overset{+}{2}\overset{0}{4}\overset{-}{1}\overset{-}{2}\overset{+}{3}\overset{+}{1}\big).
	$$
	Note that the sign associated to $u_\alpha$ only depends on (a) whether $u_\alpha$ is a caesura or not and (b) the subsequence $\big(u_1, \dots, u_\alpha\big)$ of elements preceding it.
\end{sgnrule}
The differential $\partial\colon \Surj_R(\ul{r})_d\rt \Surj_R(\ul{r})_{d-1}$ is then given by removing elements from such a sequence, together with the sign from \Cref{sgn:differential sur}
$$
\partial\big(u_1, u_2, \dots, u_{r+d}\big)=\sum_{\alpha=1}^{r+d} \pmu \big(u_1, \dots, \widehat{u_\alpha}, \dots, u_{r+d}\big).
$$
Note that if the element $u_\alpha$ appears only once in the sequence, then removing it gives zero (since the resulting sequence no longer describes a nondegenerate sequence). 
\begin{proposition}[{\cite[Theorem  2.15c]{mcclure2003multivariable}}]\label{prop:qi of future cooperads}
	For every $r\geq 1$, there is an $\mm{Aut}(\ul{r})$-equivariant quasi-isomorphism to the trivial representation
	$$
	\Surj_R(\ul{r})\rt R.
	$$
\end{proposition}
\begin{proof}
	In degree $0$, the map is the $R$-linear extension of the map sending every sequence $(u_1, \dots, u_r)$ to the unit $1$. To see that this is a quasi-isomorphism, one can realise each $\Surj_{R}(\ul{r}- \{1\})$ as a deformation retract of $\Surj_R(\ul{r})$, from which the result follows by induction.
	Indeed, define $i\colon \Surj_R(\ul{r}-\{1\}) \to \Surj_R(\ul{r})$ by $i(\mb{u})=(1, u_1, \dots, u_{r+d})$ and $r\colon \Surj_R(\ul{r})\to \Surj_R(\ul{r}-\{1\})$ by $r(\mb{u})=(u_2,\dots, u_{r+d})$ if $u_1=1$ is the only occurrence of $1$, and $r(\mb{u})=0$ otherwise.
	It is clear that $ri=\id$; the homotopy $h\colon \Surj_R(\ul{r})_d \to \Surj_R(\ul{r})_{d+1}$ between $\id$ and $ir$ is given by $h(\mb{u})=(1, u_1, \dots, u_{r+d})$.
\end{proof}

\subsection*{The cooperad structure}
We will describe the cooperad structure on the symmetric sequence $\Surj_R$ in terms of partial cocomposition maps. To this end, let $\ul{r}$ and $\ul{s}$ be two nonempty finite sets and let $v\in \ul{r}$. 
We denote by $\ul{r}\sqcup_v \ul{s}=\ul{r}\setminus \{v\}\sqcup \ul{s}$ the set obtained by removing $v$ and adding $\ul{s}$. The cocomposition of an $(r+s-1)$-ary operation along $v$ into an $r$-ary and an $s$-ary operation is then a map of the form
$$
\Delta_v\colon \Surj_R\big(\ul{r}\sqcup_v \ul{s}\big)\rt \Surj_R\big(\ul{r}\big)\otimes \Surj_R\big(\ul{s}\big).
$$
This map acts by replacing the first elements of $\ul{s}$ appearing in a sequence $\mb{u}$ by the element $v$ and removing the remaining ones.

More precisely, given a sequence $\mb{u}=\big(u_1, \dots, u_p\big)$ in $\Surj_R\big(\ul{r}\sqcup_v \ul{s}\big)$ (which is of degree ${p-r-s+1}$), its image under $\Delta_v$ is as follows. 
Let $\big(u_{\alpha(1)}, \dots, u_{\alpha(k)}\big)$ be the subsequence consisting of all elements in $\ul{s}$ (in particular, $k\geq s$). Then we define $\Delta_v\big(u_1, \dots, u_p\big)$ to be\vspace{7pt}
\begin{equation} \vspace{10pt}\hspace{30pt}\label{eq:cocomposition for sur}\sum_{i=1}^k \pmc \Big(u_1, \dots, \overset{\displaystyle v}{\cancel{u_{\alpha(1)}}}, \dots, \overset{\displaystyle v}{\cancel{u_{\alpha(i)}}}, \dots, \widehat{u_{\alpha(i+1)}}, \dots, \widehat{u_{\alpha(k)}}, \dots, u_p\Big)\otimes \Big(u_{\alpha(i)},  \dots, u_{\alpha(k)}\Big).
\end{equation}
This gives a sequence of elements in $\ul{r}$, which may be degenerate in case the original sequence has consecutive elements in $\ul{s}$. Furthermore, one takes the sequence $u_{\alpha(i)}, \dots, u_{\alpha(k)}$ of elements in $\ul{s}$; this may either be degenerate or may not exhaust all of $\ul{s}$. When degenerate or non-exhaustive sequences appear, the corresponding term is zero. This typically means that many terms in the above sum are zero: if $u_{\alpha(1)}, \dots, u_{\alpha(i-1)}$ are not all caesuras, then the second factor is not exhaustive and the term vanishes.
Finally, the sign $\pmc$ is dictated by the following Koszul sign rule for caesuras:

\begin{sgnrule}[Koszul sign rule for caesuras]\label{rem:koszul sign for caesuras}
	We will write $\pmc$ for the sign obtained by the following rule: whenever in a formula a caesura passes along another one, one multiplies by $-1$. 
	Explicitly, consider a term in the cocomposition \eqref{eq:cocomposition for sur} of the form $\mb{v}\otimes \mb{w}$ for certain sequences $\mb{v}$ and $\mb{w}$. Then there is a bijection $\Caes{\mb{u}}\cong \Caes{\mb{v}}\star \Caes{\mb{w}}$ between the linear orders of caesuras in $\mb{u}$ and those in $\mb{v}$ and $\mb{w}$ (where $\star$ denotes the addition of ordinals) and the sign $\pmc$ is the sign of this bijection.
	
	Note that the sign rule for caesuras refines the usual Koszul sign rule, in the sense that under the symmetry isomorphism $\Surj_R(r)\otimes \Surj_R(s)\cong \Surj_R(s)\otimes \Surj_R(r)$, $\big(u_1, \dots, u_{r+d}\big)\otimes \big(v_1, \dots, v_{s+e}\big)$ and $\big(v_1, \dots, v_{s+e}\big)\otimes \big(u_1, \dots, u_{r+d}\big)$ agree up to the sign $\pmc$ given by the number of times two caesuras are interchanged.
\end{sgnrule}

\begin{example}
	Consider the sequence $(\textbf{1},\textbf{2},\textbf{3},1,2,3)\in\Surj_R(\ul{3})_3$ (numbers in bold are the caesuras) and  write $\ul{3} = \{1,v\} \circ_v \{2,3\}$. 
	The partial cocomposition $\Delta_v(1,2,3,1,2,3)$ along $v$ is then given by
	\begin{align*}
		&(1,v,1) \otimes (2,3,2,3) + (1,v,v,1) \otimes (3,2,3) +  (1,v,v,1,v) \otimes (2,3) + (1,v,v,1,v,v) \otimes (3)\\
		&=(1,v,1) \otimes (2,3,2,3) + 0\otimes (3,2,3) + 0 \otimes (2,3) + 0\otimes 0\\
		&=(1,v,1) \otimes (2,3,2,3).
	\end{align*}
	There are no caesuras going over other caesuras and therefore all signs are $+$.
	
	If we consider instead the sequence $(\textbf{1},\textbf{2},\textbf{1},\textbf{3},1,2,3)\in\Surj_R(\ul{3})_4$ and we decompose along the same element $v$, we get
\begin{align*}
		  \Delta_v(1,2,1,3,1,2,3)&=-(1,v,1,1) \otimes (2,3,2,3) + (1,v,1,v,1) \otimes (3,2,3) +  (1,v,1,v,1,v) \otimes (2,3) + 0\\
		&=-0 + (1,v,1,v,1) \otimes (3,2,3) + (1,v,1,v,1,v) \otimes (2,3) + 0.
	\end{align*}
The first sign arises since the first $2$, which is a caesura, went over the second $1$ which is also a caesura.
\end{example}

\begin{remark}
	Any caesura $u_\beta$ in $\big(u_1, \dots, u_p\big)$ will appear as a caesura in exactly one of the two factors in the expression for $\Delta_v(u)$: if $u_\beta\not\in \ul{s}$, it will appear as a caesura in the first factor and if it is one of the $u_{\alpha(i)}, \dots, u_{\alpha(k)}$, it will appear as a caesura in the second factor. Finally, all $u_{\alpha(1)}, \dots , u_{\alpha(i-1)}$ will appear as caesuras in the first factor (namely as all but the last copy of $v$). 
\end{remark}

Observe that the maps $\Delta_v$ are well-defined: such maps send a degenerate sequence to a sum of terms, each of which containing a degenerate sequence and likewise for non-surjective sequences. 

\begin{proposition}\label{prop:surjections cooperad}
	For $r, s\geq 1$, the formulas $ \Delta_v\colon \Surj_R\big(\ul{r}\sqcup_v \ul{s}\big)\rt \Surj_R\big(\ul{r}\big)\otimes \Surj_R\big(\ul{s}\big)$ defined above endow $\{\Surj_R(r)\}$ with the structure of a dg-cooperad.
\end{proposition}

The proof of this proposition is a lengthy verification of all the axioms. Postponing this for the moment, we record some simple consequences.

\begin{corollary}
	The maps from Proposition \ref{prop:qi of future cooperads} induce a quasi-isomorphism of dg-cooperads $\Surj_R\rt \mathbf{coCom}_R^{\mm{nu}}= \pi_0(\Surj_R)$.
\end{corollary}
\begin{proof}
	All we need to check is that the map is compatible with the partial cocompositions in degree $0$.
	One readily checks that any partial composition of a permutation in $\Surj_R(\ul{r})_0$ is a tensor product of two permutations with a $+$ sign (there are no caesuras).
\end{proof}

Notice that while the degree $0$ part of the surjections cooperad is a $\Sigma$-free dg-cooperad with underlying symmetric sequence $\Surj_R(\ul{r})_0=k[\Sigma_r] $, this cooperad is \emph{not} the nonunital coassociative cooperad; this should not be expected, since there is no map $\dgob{coAss}_R^{\mm{nu}} \to \dgob{coCom}^{\mm{nu}}_R$.

\begin{remark}[Surjections cooperad in degree $0$]
	One can show that degree $0$ piece of the surjection cooperad $(\Surj_R)_0$ is isomorphic to the linear dual of the operad $\mathbf{Zinb}$ governing Zinbiel algebras. Recall that such Zinbiel algebras are chain complexes equipped with a binary operation $\prec$ satisfying $(x\prec y)\prec z = x\prec (y\prec z + (-1)^{|y||z|} z \prec y)$ \cite[Section 13.5]{LV}.
	
	Indeed, we can define a map $(\Surj_R)_0 \to \mm{Cofree}_{\mm{Coop}}(R\mu \oplus R\mu^{(12)})$ into the cofree cooperad cogenerated by an arity $2$ element $\mu$ with free $\Sigma_2$ action, using that $\Surj_R(2)_0 \cong R[\Sigma_2]$.
	This map restricts to a map $(\Surj_R)_0\to \dgob{Zinb}^\vee\subseteq \mm{Cofree}_{\mm{Coop}}(R\mu \oplus R\mu^{(12)})$. This map is necessarily injective, since it is injective on cogenerators. Since $\dim((\Surj_R)_0(n)) = \dim(\dgob{Zinb}(n)) = n!$, it is an isomorphism. 
\end{remark}

\begin{proof}[Proof of Proposition \ref{prop:surjections cooperad}]
We start by observing that for nonempty finite sets $\ul{r}$ and $\ul{s}$ and $v \in \underline{r}$, the map $\Delta_v\colon \Surj_R\big(\ul{r}\sqcup_v \ul{s}\big)\rt \Surj_R\big(\ul{r}\big)\otimes \Surj_R\big(\ul{s}\big)$ is equivariant with respect to the group $\mm{Aut}(\ul{r}\setminus v)\times \mm{Aut}(\ul{s})\subseteq \mm{Aut}(\ul{r}\sqcup_v \ul{s})$ of permutations of $\ul{r}$ fixing $v$ and of permutations of $\ul{s}$.
To show that the dg-cooperad axioms are satisfied we need to check counitality, coassociativity (both parallel and sequential), and compatibility with the differential.

\subsubsection*{Counitality} 
Note that $\Surj_R(\ul{1})\cong R$ is spanned by the trivial one term sequence $(1)$; this gives the counit. For $v\in \ul{r}$ and a sequence $\mb{u}$ in $\Surj_R(\ul{r}\sqcup_v \ul{1})$, the formula for the cocomposition has only one term in which the second factor is nonzero, for $i=k$, giving
$$
\Delta_v\big(u_1, \dots, u_p\big)=\big(u_1, \dots, \overset{v}{\cancel{1}}, \dots, \overset{v}{\cancel{1}}, \dots, u_p\big)\otimes (1).
$$
In other words, one just replaces all copies of the element $1\in\ul{1}$ by $v$. This shows that the cocomposition is right counital; the verification of left counitality is similar.

\subsubsection*{Parallel coassociativity} 
Let $v_1, v_2\in \ul{r}$ be two distinct elements, and $\ul{s}_1$ and $\ul{s}_2$ two sets. We consider the set $\ul{r}\sqcup_{(v_1, v_2)} (\ul{s}_1, \ul{s}_2)$ obtained by replacing $v_i\in\ul{r}$ by $\ul{s_i}$. Notice that 
$$
\ul{r}\sqcup_{(v_1, v_2)} (\ul{s}_1, \ul{s}_2)= (\ul{r}\sqcup_{v_1} \ul{s}_1)\sqcup_{v_2} \ul{s}_2 = (\ul{r}\sqcup_{v_2} \ul{s}_2)\sqcup_{v_1} \ul{s}_1.
$$
Consider a sequence $\mb{u}=\big(u_1, \dots, u_p\big) \in \Surj_R(\ul{r}\sqcup_{(v_1, v_2)} (\ul{s}_1, \ul{s}_2))_d$. We have to verify that
$$
\Delta_{v_2}\circ \Delta_{v_1}(\mb{u})\in \Surj_R(\ul{r})\otimes\Surj_R(\ul{s}_2)\otimes \Surj_R(\ul{s}_1)
$$
agrees with $\Delta_{v_1}\circ \Delta_{v_2}(\mb{u})$ upon permuting the $\Surj_R(\ul{s}_1)$ and $\Surj_R(\ul{s}_2)$ pieces. We start by doing this verification up to the Koszul sign induced by the caesuras.

Let $u_{\alpha(1)}, \dots, u_{\alpha(k)}$ be the subsequence of elements in $\ul{s}_1$ and $u_{\beta(1)}, \dots, u_{\beta(l)}$ for the subsequence of elements in $\ul{s}_2$. Then $\Delta_{v_2}\circ \Delta_{v_1}(u_1, \dots, u_p)$ is given by
\begin{align*}
	\sum_{j=1}^l\sum_{i=1}^k &\pmc \Big(u_1, \dots, \overset{\displaystyle v_1}{\cancel{u_{\alpha(1)}}}, \overset{\displaystyle v_2}{\cancel{u_{\beta(1)}}} \dots , \overset{\displaystyle v_1}{\cancel{u_{\alpha(i)}}},  \overset{\displaystyle v_2}{\cancel{u_{\beta(j)}}}, \dots, \widehat{u_{\alpha(i+1)}}, \widehat{u_{\beta(j+1)}},\dots,\widehat{u_{\alpha(k)}}, \widehat{u_{\beta(l)}},\dots, u_p\Big)\\
	&\otimes \Big(u_{\beta(j)}, u_{\beta(j+1)}, \dots, u_{\beta(l)}\Big) \otimes \Big(u_{\alpha(i)}, u_{\alpha(i+1)}, \dots, u_{\alpha(k)}\Big).
\end{align*}
In words, one just replaces all $u_{\alpha(1)}, \dots, u_{\alpha(i)}$ by $v_1$ and removes the $u_{\alpha(i+1)}, \dots, u_{\alpha(k)}$, and similarly one replaces $u_{\beta(1)}, \dots, u_{\beta(j)}$ by $v_2$ and removes the $u_{\beta(j+1)}, \dots, u_{\beta(l)}$. In the first factor, the various $u_{\alpha}$ and $u_{\beta}$ need not appear in the order they are depicted: for instance, $u_{\beta(1)}$ may precede $u_{\alpha(1)}$. This is clearly symmetric upon exchanging $v_1\leftrightarrow v_2$ and $\ul{s}_1\leftrightarrow \ul{s}_2$.
Using \Cref{rem:koszul sign for caesuras} it is immediate that the signs $\pmc$ produced in the computation of $\Delta_{v_2}\circ \Delta_{v_1}$ are also produced in $\Delta_{v_1}\circ \Delta_{v_2}$.

\subsubsection*{Sequential coassociativity} 
We now consider sets $\ul{r},\ul{s}$ and $\ul{t}$ and let $v\in \ul{r}$, $w\in \ul{s}$. 
We will address coassociativity on the total set $(\ul{r}\sqcup_v \ul{s})\sqcup_w \ul{t}= \ul{r}\sqcup_v\big(\ul{s}\sqcup_w\ul{t}\big)$.
Concretely, let $\mb{u}$ be a sequence of $\Surj_R(\ul{r}\sqcup_v\ul{s}\sqcup_w\ul{t})_d$ and we will show that 
$$
\Delta_v\circ\Delta_w(\mb{u})=\Delta_w\circ\Delta_v(\mb{u})\in\Surj_R(\ul{r})\otimes\Surj_R(\ul{s})\otimes\Surj_R(\ul{t}).
$$
Let $\big(u_{\alpha(1)}, \dots, u_{\alpha(k)}\big)$ be the subsequence of elements in $\ul{s}\sqcup_w \ul{t}$ and let $\big(u_{\alpha(i_1)}, \dots, u_{\alpha(i_\omega)}\big)$ be the (sub)subsequence of elements in $\ul{t}$.
Then $\Delta_w\circ\Delta_v(\mb{u})$ is given by
\begin{align*}
	\sum_{i=1}^k \sum_{\lambda\colon i_\lambda\geq i}&\pmc \Big(u_1, \dots, \overset{\displaystyle v}{\cancel{u_{\alpha(1)}}}, \dots, \overset{\displaystyle v}{\cancel{u_{\alpha(i)}}}, \dots, \widehat{u_{\alpha(i+1)}}, \dots, \widehat{u_{\alpha(k)}}, \dots, u_p\Big)\\
	&\otimes\Big(u_{\alpha(i)}, u_{\alpha(i+1)}, \dots, \overset{\displaystyle w}{\cancel{u_{\alpha(i_1)}}}, \dots, \overset{\displaystyle w}{\cancel{u_{\alpha(i_\lambda)}}} , \dots, \widehat{u_{\alpha(i_{\lambda+1})}}, \dots, \widehat{u_{\alpha(i_\omega)}}, \dots, u_{\alpha(k)}\Big)\\
	\otimes &\Big(u_{\alpha(i_\lambda)}, u_{\alpha(i_{\lambda+1})}, \dots, u_{\alpha(i_\omega)}\Big).
\end{align*}
In words, from  $\big(u_1, \dots, u_p\big)$ one first removes the part of the subsequence $\big(u_{\alpha(1)}, \dots, u_{\alpha(k)}\big)$ after step $i$ and replaces the part of the subsequence before step $i$ by copies of $v$. Next, from the sequence $\big(u_{\alpha(i)}, u_{\alpha(i+1)}, \dots, u_{i_k}\big)$ one removes the part of the subsubsequence $\big(u_{\alpha(i_1)}, \dots, u_{\alpha(i_\omega)}\big)$ after the step $\lambda$ and replaces the part of the subsubsequence before step $\lambda$ by copies of $w$.

Going `right-to-left' instead, we see that every summand above is obtained by first picking out a subsequence $\big(u_{\alpha(i_\lambda)}, u_{\alpha(i_{\lambda+1})}, \dots, u(i_{a_{\omega}})\big)$ of $\big(u_1, \dots, u_p\big)$ of arbitrary length, then extending it to a larger subsequence (determined by $\alpha(i)$ and $\alpha(k)$) while picking a number $i_1\leq i_\lambda$.

Note that these are precisely the terms obtained when computing $\Delta_v\circ \Delta_w\big(u_1, \dots, u_p\big)$, except that the latter may also produce terms in which $i_1> i_\lambda$. Those additional terms are all zero, since the middle sequence in $\Surj_R(\ul{s})$ is no longer exhaustive (it does not contain any $w$). Because in both computations, the signs arise from the same permutations of caesuras, they agree and we conclude that $\Delta_v\circ \Delta_w\big(u_1, \dots, u_p\big) =\Delta_w\circ \Delta_v\big(u_1, \dots, u_p\big)$.

\subsubsection*{Compatibility with the differential}
It remains to check that the cocomposition is compatible with the differential $\partial$.
Let $\mb{u}\in \Surj_R(\ul{r}\sqcup_v\ul{s})_d$ and let $\big(u_{\alpha(1)}, \dots, u_{\alpha(k)}\big)$ be the subsequence of elements in $\ul{s}$. 
Up to signs, $\partial\circ \Delta_v(\mb{u})$ is given by
\begin{align*}
	\partial &\sum_{i=1}^k\Big(u_1, \dots, \overset{\displaystyle v}{\cancel{u_{\alpha(i)}}}, \dots, \widehat{u_{\alpha(i+1)}}, \dots, \widehat{u_{\alpha(k)}}, \dots, u_p\Big)\otimes \Big(u_{\alpha(i)}, u_{\alpha(i+1)}, \dots, u_{\alpha(k)}\Big)
	\nonumber\\
	=&\sum_{u_{\beta}\not\in\ul{s}}\sum_{i=1}^k\Big(u_1, \dots, \widehat{u_\beta}, \dots, \overset{\displaystyle v}{\cancel{u_{\alpha(i)}}}, \dots, \widehat{u_{\alpha(i+1)}}, \dots, \widehat{u_{\alpha(k)}}, \dots, u_p\Big)\otimes \Big(u_{\alpha(i)}, u_{\alpha(i+1)}, \dots, u_{\alpha(k)}\Big)
	\nonumber\\
	+&\sum_{i=1}^k\sum_{j\leq i}  \Big(u_1, \dots, \overset{\displaystyle \widehat{\hspace{.1cm}v\hspace{.1cm}}}{\cancel{u_{\alpha(j)}}}, \dots, \overset{\displaystyle v}{\cancel{u_{\alpha(i)}}}, \dots, \widehat{u_{\alpha(i+1)}}, \dots, \widehat{u_{\alpha(k)}}, \dots, u_p\Big)\otimes \Big(u_{\alpha(i)}, u_{\alpha(i+1)}, \dots, u_{\alpha(k)}\Big)
	\nonumber\\
	+&\sum_{i=1}^k \sum_{j\geq i}  \Big(u_1, \dots, \overset{\displaystyle v}{\cancel{u_{\alpha(i)}}}, \dots, \widehat{u_{\alpha(i+1)}}, \dots, \widehat{u_{\alpha(j)}}, \dots, \widehat{u_{\alpha(k)}}, \dots, u_p\Big)\otimes \Big(u_{\alpha(i)}, \dots, \widehat{u_{\alpha(j)}}, \dots, u_{\alpha(k)}\Big)
\end{align*}
Here we have split the result into the three types of summands above corresponding to the three kinds of elements which can be removed by the differential: (1) an element $u_\beta\in \ul{r}\setminus \{v\}$, (2) a copy of $v$ put in the place of $u_{\alpha(j)}\in \ul{s}$, or (3) an element $u_{\alpha(j)}\in \ul{s}$.

On the other hand, we have that $\Delta_v\circ \partial(\mb{u})$ is given up to signs by
\begin{align*}
	&\sum_{u_\beta\not\in\ul{s}} \Delta_v\big(u_1, \dots, \widehat{u_\beta}, \dots, u_p\big)+\sum_{j=1}^k \Delta_v\big(u_1, \dots, \widehat{u_{\alpha(j)}}, \dots, u_p\big)\\
	=&\sum_{i=1}^k\sum_{u_\beta\not\in\ul{s}}\Big(u_1, \dots, \widehat{u_{\beta}}, \dots, \overset{\displaystyle v}{\cancel{u_{\alpha(i)}}}, \dots, \widehat{u_{\alpha(i+1)}}, \dots, \widehat{u_{\alpha(k)}}, \dots, u_p\Big)\otimes \Big(u_{\alpha(i)}, u_{\alpha(i+1)}, \dots, u_{\alpha(k)}\Big)
	\nonumber\\
	+&\sum_{j=1}^k \sum_{i>j} \Big(u_1, \dots, \widehat{u_{\alpha(j)}}, \dots, \overset{\displaystyle v}{\cancel{u_{\alpha(i)}}}, \dots, \widehat{u_{\alpha(i+1)}}, \dots, \widehat{u_{\alpha(k)}}, \dots, u_p\Big)\otimes \Big(u_{\alpha(i)}, u_{\alpha(i+1)}, \dots, u_{\alpha(k)}\Big)
	\nonumber\\
	+& \sum_{j=1}^k\sum_{i<j} \Big(u_1, \dots, \overset{\displaystyle v}{\cancel{u_{\alpha(i)}}}, \dots, \widehat{u_{\alpha(i+1)}}, \dots, \widehat{u_{\alpha(j)}}, \dots, \widehat{u_{\alpha(k)}}, u_p\Big)\otimes \Big(u_{\alpha(i)}, \dots, \widehat{u_{\alpha(j)}}, \dots, u_{\alpha(k)}\Big)
\end{align*}
The first type of summand corresponds to the case where the differential removes an element not in $\ul{s}$ whereas the second and third line describe the cocomposition after one has removed the element $u_{\alpha(j)}$ from the sequence.

We first check that up to signs the two computations agree and we will do a careful sign verification afterwards. 
It is clear that the first type of summands agrees in both formulas. 
The other summands are \emph{almost} the same, except that $\partial\circ \Delta_v\big(u_1, \dots, u_p\big)$ also includes the cases where $i=j$ (twice). 
One easily sees that such terms pairwise cancel each other out. 
To be precise, the difference $\partial\circ \Delta_v\big(u_1, \dots, u_p\big)-\Delta_v\circ \partial\big(u_1, \dots, u_p\big)$ is (as usual up to sign) given by\vspace{5pt} 
\begin{align}\label{eq:extra terms}\ 
\sum_{i=1}^k  \Big(u_1, \dots, \overset{\displaystyle v}{\cancel{u_{\alpha(i-1)}}}, \dots, \overset{\displaystyle \widehat{\hspace{.1cm}v\hspace{.1cm}}}{\cancel{u_{\alpha(i)}}}, \dots, \widehat{u_{\alpha(i+1)}}, \dots, \widehat{u_{\alpha(k)}}, \dots, u_p\Big)\otimes \Big(u_{\alpha(i)}, u_{\alpha(i+1)}, \dots, u_{\alpha(k)}\Big) 
\end{align}
$$+\sum_{i=1}^k \Big(u_1, \dots, \overset{\displaystyle v}{\cancel{u_{\alpha(i)}}}, \dots, \widehat{u_{\alpha(i+1)}}, \dots, \widehat{u_{\alpha(j)}}, \dots, \widehat{u_{\alpha(k)}}, \dots, u_p\Big)\otimes \Big(\widehat{u_{\alpha(i)}}, u_{\alpha(i+1)}, \dots, u_{\alpha(k)}\Big).$$
In the first line, the term corresponding to $i=1$ is zero (the first factor does not contain any $v$) and in the second line, the term corresponding to $i=k$ is zero. For $i>1$, the $i$-th term in the first line is precisely cancelled by the $(i-1)$-st term in the second line; we will verify that the signs match in Case \ref{it:case6} 
below.

\subsubsection*{The signs of $\partial \circ \Delta_v$ and $\Delta_v \circ \partial$}

Recall that the differential $\partial$ acts by removing from a sequence $\mb{u}$ the element $u_{\beta}$ for all $j$, with sign $\pm_{\mb{u}, \beta}$ (\Cref{sgn:differential sur}) determined by (a) whether or not $u_\beta$ is a caesura and (b) the caesuras appearing in $u_1, \dots, u_{\gamma}$, where $u_{\gamma}$ denotes the largest caesura with $\gamma\leq \beta$ and $u_\beta=u_\gamma$.

On the other hand, \Cref{rem:koszul sign for caesuras} dictates that $\Delta_v$ produces a sign which is equal to the sign of an unshuffle of caesuras: cutting $\big(u_1,\dots,u_p\big)$ along $v$ in $\Surj_R(\ul{r}\sqcup_v\ul{s})$ carries the sign $\pmc$ given by the exchange of caesuras in $\big(u_{\alpha(i)}, u_{\alpha(i+1)}, \dots, u_{\alpha(k)}\big)\in \Surj_R(\ul{s})$ and the caesuras among the elements $u_{\alpha(i)+1}, \dots, u_p \in \Surj_R(\ul{r})$.

In order to check the difference in signs between $\Delta_v\circ \partial$ and $\partial\circ\Delta_v$ we will go through various cases, denoting by $u_{\beta}$ the element removed by the differential and by $u_{\alpha(i)}$ the element at which one cuts:

\begin{enumerate}[label=(\roman*)]
	\item\label{it:case1} \emph{If $\beta<\alpha(i)$}. Removing $u_\beta$ does not change the caesuras after $u_{\alpha(i)}$ and cutting at $u_{\alpha(i)}$ does not change the amount of caesuras before $u_\beta$, nor whether $u_\beta$ is a caesura. The signs therefore agree.
	
	\item \emph{If $\alpha(i)<\beta$ and $u_\beta$ is a caesura which is not in $\ul{s}$.} In this case, cutting at $u_{\alpha(i)}$ removes the caesuras $u_{\alpha(c)} \in \ul{s}$ before $u_\beta$ which satisfy $\alpha(i) \leq \alpha(c)<\beta$ (all other caesuras remain the same or are changed into caesuras labelled by $v$); the sign therefore changes by the number of such $u_{\alpha(c)}$.
	
	On the other hand, after removing $u_\beta$, the sign from $\Delta_v$ changes by this same number: indeed, the caesuras $u_{\alpha(c)}\in \ul{s}$ with $\alpha(i)\leq \alpha(c)<\beta$ need to be moved passed one less caesura. We conclude that in this \mbox{case the signs agree.}
	
	\item \emph{If $\alpha(i)<\beta$ and $u_\beta$ is a non-caesura which is not in $\ul{s}$}. Let $\gamma<\beta$ be the largest number such that $u_\gamma=u_\beta$. Since $\pm_{\mb{u}, \beta}=- (\pm_{\mb{u}, \gamma})$ (\Cref{sgn:differential sur}), the previous argument shows that after cutting, the sign of the differential is changed by the number of caesuras $u_{\alpha(c)}\in \ul{s}$ such that $\alpha(i)\leq \alpha(c)<\gamma$. Similarly, after removing $u_\beta$ the element $u_\gamma$ is no longer a caesura, so that the sign from $\Delta_v$ changes by the number of caesuras $u_{\alpha(c)}\in \ul{s}$ with $\alpha(i)\leq \alpha(c)<\gamma$ as well. In total the signs coincide.
	
	\item\label{it:case4} \emph{If $\alpha(i)<\beta$ and $u_\beta$ is a caesura belonging $\ul{s}$.} Since $u_\beta\in \ul{s}$, $\beta = \alpha(b)$ for some $b$. 
	Upon cutting at $\alpha(i)$, there are more caesuras before $u_{\alpha(b)}$: indeed, the caesuras $u_\gamma\in \ul{r}\setminus v$ with $\gamma>\alpha(b)$ now precede $u_{\alpha(b)}$ and the sign of the differential changes by their number.
	
	On the other hand, when cutting at $\alpha(i)$ after having removed $u_{\alpha(b)}$, one no longer has to move $u_{\alpha(b)}$ past the caesuras $u_\gamma\in \ul{r}\setminus v$ with $\gamma>\alpha(b)$. Thus, in this case, the signs also agree.
	
	\item \emph{If $\alpha(i)<\beta$ and $u_\beta\in \ul{s}$ is not a caesura.} We again write $\beta = \alpha(b)$ and suppose that $u_\alpha(c)$ is the preceding copy of that same element in $\ul{s}$ (i.e. $c<b$ is the biggest number such that $u_{\alpha(c)}=u_{\alpha(b)}$). There are two subcases:
	\begin{itemize}
		\item $\alpha(i)<\alpha(c)$: Since the sign of the differential at $u_{\alpha(b)}$ is minus the sign of $u_\alpha(c)$, Case \ref{it:case4} shows that after cutting at $u_{\alpha(i)}$, the sign of the differential changes by the number of caesuras $u_\gamma\in\ul{r}\setminus v$ with $\gamma>\alpha(c)$.
		
		On the other hand, after removing $u_{\alpha(b)}$ the element $u_{\alpha(c)}$ is no longer a caesura, therefore moving it past all caesuras in $\ul{r}\setminus v$ after it, no longer contributes to the sign of $\Delta_v$. In total the sign \mbox{therefore remains the same.}
		
		\item $\alpha(c)<\alpha(i)$. Since the sign of the differential at $u_{\alpha(b)}$ is minus the sign of $u_{\alpha(c)}$, Case \ref{it:case1} shows that the sign of the differential is left unchanged. Similarly, removing $u_{\alpha(b)}$ does not change the caesuras appearing after $u_{\alpha(i)}$, so the signs for $\Delta_v$ do not change either.
	\end{itemize} 
	\item\label{it:case6} With all cases considered, it remains to identify the sign of the $i$-th term in the first row of \eqref{eq:extra terms} with \emph{minus} the sign of the $(i-1)$-st term in the second row. Note that both terms are zero if $u_{\alpha(i-1)}$ is not a caesura. 
	
	Now the sign of the $i$-th term in the first row is given by (a) the number of caesuras $u_\gamma$ with $\gamma\leq \alpha(i-1)$ (coming from the differential) and (b) the sign of the unshuffle of the caesuras in $\ul{s}$ and $\ul{r}\setminus v$ appearing in places $\geq  \alpha(i)$ (coming from $\Delta_v$).
	
	On the other hand, the sign of the $(i-1)$-st term in the second row is given by (a1) the number of caesuras $u_\gamma\in \ul{r}\setminus v$, (a2) the number of caesuras $u_{\alpha(b)}\in \ul{s}$ with $\alpha(b)<\alpha(i-1)$ (together these give the sign of the differential) and (b) the sign of the unshuffle of the caesuras in $\ul{s}$ and $\ul{r}\setminus v$ appearing in places $\geq \alpha(i-1)$ (coming from $\Delta_v$).
	
	The difference between the two different signs (a) coming from the differential is given by \emph{minus} the parity of the number of caesuras $u_\gamma\in \ul{r}\setminus v$ with $\gamma>\alpha(i-1)$. The difference between the signs (b) coming from $\Delta_v$ is the parity of the number of times $u_{\alpha(i-1)}$ is moved past a caesura in $\ul{r}\setminus v$ after it. The signs differ therefore by $-1$ and the two terms indeed cancel.
\end{enumerate}

This concludes the proof of Proposition \ref{prop:surjections cooperad} and therefore also Theorem \ref{thm:surjections cooperad}.
\end{proof}

\subsection*{The PD surjections operad}
For most of our purposes we are more interested in the linear dual of the cooperad $\Surj_R$, we conclude by giving an explicit description of the operad $\Surj_R^\vee$, which we dub the PD surjections operad.
\begin{definition}[PD surjections operad]\label{subsec:App PD surj operad}
The \emph{PD surjections operad} $\Surj_R^\vee$ is the  $R$-linear dg-operad defined as follows:
\begin{itemize}
	\item	For each nonempty finite set $\ul{r}=\{1,\dots,r\} $, let $\Surj_R^\vee(\ul{r})$ be the free graded $R$-module spanned in each degree $-d\leq 0$ by by (ordered) sequences $\textbf{u}=(u_1, \dots, u_{r+d})$ that are non-degenerate in the sense that each $1, \dots, r$ appears in the sequence and $u_\alpha\neq u_{\alpha+1}$ for $\alpha=1,\dots, r+d-1$.
The symmetric group $\Sigma_r$ acts on such nondegenerate sequences by permuting each individual $u_\alpha$.
	
	\item Each $\Surj^\vee_R(\ul{r})$ comes equipped with a differential sending a nondegenerate sequence $\mb{u}$ to the (signed) sum of all nondegenerate sequences $\mb{u}_+$ obtained by adding an element to $\mb{u}$. More precisely, 	
	$$
	\partial\big(u_1, u_2, \dots, u_{r+d}\big)=\sum_{\alpha=1}^{r+d+1} \sum_{u_{\alpha-1} \ne v \ne u_{\alpha}} \pm_{\mb{u}_+, \alpha} \big(u_1, \dots, u_{\alpha-1}, v, u_{\alpha} \dots, u_{r+d}\big).
	$$
	Here the sign $\pm_{\mb{u}_+, \alpha}$ is the sign associated to the element $v$ in $\mb{u}_+=(u_1, \dots, v, \dots, u_{r+d})$, as in Sign Rule \ref{sgn:differential sur}.

	\item The operad structure is determined by partial composition maps
	$$\begin{tikzcd}
	\circ_k \colon \Surj_{R}^\vee(\ul{r})\otimes \Surj_{R}^\vee(\ul{s}) \to \Surj_{R}^\vee\big((\ul{r}-\{k\})\sqcup\ul{s}\big)
	\end{tikzcd}$$
	along $k\in \ul{r}$, defined as follows. For any two sequences $\textbf{u}=(u_1,\dots u_{r+d})$ in $\ul{r}$ and $\textbf{v}=(v_1,\dots,v_{s+e})$ in $\ul{s}$, let $(u_{\alpha(1)}, \dots, u_{\alpha(i)})$ denote the subsequence of $\mb{u}$ with values $k$. We then define
	
	\vspace{3pt}

\noindent\scalebox{0.9}{\parbox{\linewidth}{$$
		\mb{u}\circ_k \mb{v} = \sum \pmc \Big(u_1, \dots, \overset{\displaystyle s_{1}}{\cancel{\hspace{.1cm} u_{\alpha(1)}\hspace{.1cm}}}, \dots, \overset{\displaystyle s_{i-1}}{\cancel{\hspace{.1cm} u_{\alpha(i-1)}\hspace{.1cm}}},\dots, \overset{\displaystyle v_{1}}{\cancel{\hspace{.1cm} u_{\alpha(i)}\hspace{.1cm}}}, u_{\alpha(i)+1},\dots, v_{2}, \dots, u_{\beta}, \dots, v_{s+e}, \dots, u_{r+d}\Big).$$	}}
\\[3pt]

\noindent More precisely, we take the sum of all sequences $\mb{w}$ in $(\ul{r}-\{k\})\sqcup\ul{s}$ obtained from $\mb{u}$ by the following procedure:
\begin{itemize}
\item replace the last occurence of $k$ in the sequence $\mb{u}$ by $v_1$.
\item replace the occurences of $k$ that are caesuras by any choice of elements $s_{1}, \dots, s_{i-1}\in \ul{s}$.

\item shuffle the elements $u_{\alpha(i)+1}, u_{\alpha(i)+2}, \dots, u_{r+d}$ appearing after the last occurence of $k$ and the elements $v_{2}, \dots, v_{s+e}$.
\end{itemize}
The sign $\pmc$ is determined by how many caesuras went past each other to reach the final sequence $\mb{w}$, as in Sign Rule \ref{rem:koszul sign for caesuras}. Explicitly, for any sequence $\mb{w}$ as above, there is a (non-ordered) bijection $\Caes{\mb{w}}\cong \Caes{\mb{u}}\star \Caes{\mb{v}}$ between the linearly ordered sets of caesuras (Definition \ref{def:caesura}) of $\mb{w}$ and the concatenation of the linear orders of caesuras in $\mb{u}$ and $\mb{v}$. Then $\pmc$ is the sign of this bijection.

\end{itemize}
\end{definition}
We conclude with the following result about the Koszul dual of the PD surjections operad:
\begin{theorem}\label{thm:KD of pd surj op}
The cobar construction of the surjections cooperad gives a cofibrant replacement $\Omega(\Surj_R)\stackrel{\sim}{\rt} \ope{Lie}^s_R$ of the $R$-linear shifted Lie operad. Equivalently, there is an equivalence of dg-operads $\mb{KD}(\Surj^\vee_R)\stackrel{\sim}{\rt} \dgob{Lie}^s_R$.
\end{theorem}
\begin{proof}
The first assertion follows from the fact that $\Surj_R$ is a $\Sigma$-cofibrant resolution of $\mb{coCom}^\mm{nu}_R$ and the fact $\Omega(\mb{coCom}^\mm{nu}_R)\simeq \mb{Lie}_R^s$ \cite[Theorem 6.8]{fresse346koszul}. Since $\Surj_R(r)$ is a finite rank free $R$-module in each degree, there is an isomorphism $\mb{KD}(\Surj^\vee_R)\cong \Omega(\Surj_R)$.
\end{proof}

\section{Free algebras in monoidal $\infty$-categories}
\label{app:free algebras}
The purpose of this section is to record an existence result for free associative algebras in monoidal $\infty$-categories where the tensor product preserves colimits in the first variable, but not in the second (such as symmetric sequences with the composition product). This is due to Kelly \cite{kelly1980} in the case of ordinary categories and, as we will show, the argument from loc.\ cit.\ carries over to $\infty$-categories.
\begin{cons}
Let $\cat{C}$ be a monoidal $\infty$-category with coproducts and sequential colimits, which are preserved by $-\otimes X$ for each $X\in \cat{C}$. For each $X\in \cat{C}$, we inductively define a sequence of objects in $\cat{C}$ by
$$
T^{(0)}(X)= 1 \qquad \qquad T^{(n)}(X)= 1\amalg \big(X\otimes T^{(n-1)}(X)\big).
$$
We define maps $i_n\colon T^{(n-1)}(X)\rt T^{(n)}(X)$ by setting $i_1\colon 1\rt 1\amalg X$ to be the obvious inclusion and
$$\begin{tikzcd}[column sep=5pc]
i_n\colon 1\amalg\big(X\otimes T^{(n-2)}(X)\big)\arrow[r, "\mm{id}\amalg (X\otimes i_{n-1})"] & 1\amalg\big(X\otimes T^{(n-1)}(X)\big).
\end{tikzcd}$$
\end{cons}
\begin{theorem}\label{thm:free algebra}
Let $\cat{C}$ be a monoidal $\infty$-category with coproducts and sequential colimits, such that each $(-)\otimes X$ preserves finite coproducts and sequential colimits, while each $X\otimes (-)$ preserves sequential colimits. For every object $X\in \cat{C}$, there then exists a $T(X)\in \Alg(\cat{C})$ together with a map $X\rt T(X)$ in $\cat{C}$ which exhibits $T(X)$ as the \emph{free associative algebra} on $X$. In other words, the forgetful functor
$$\begin{tikzcd}
\mm{forget}\colon \Alg(\cat{C})\arrow[r] & \cat{C}
\end{tikzcd}$$
admits a left adjoint $T$. Furthermore, there is a natural equivalence of objects in $\cat{C}$
$$
T(X)\simeq \colim_n T^{(n)}(X).
$$
\end{theorem}
The remainder of this section is devoted to a proof of \Cref{thm:free algebra}; throughout we assume that $\cat{C}$ is a monoidal $\infty$-category with the properties appearing in the theorem. The main idea of the proof will be to deduce \Cref{thm:free algebra} from a statement about left modules.
More precisely, recall that $\cat{C}$ is the free right $\cat{C}$-module $\infty$-category on a single object (the unit $1$), so that there is an equivalence of monoidal $\infty$-categories from $\cat{C}$ to the $\infty$-category of right $\cat{C}$-linear endofunctors of $\cat{C}$ \cite[\S 4.7.1]{lurie2014higher}
$$\begin{tikzcd}
\cat{C}\arrow[r, "\sim"] & \End_{\cat{C}}(\cat{C}); & X\arrow[r, mapsto] & X\otimes (-).
\end{tikzcd}$$
For an object $X\in \cat{C}$, write $F_X\colon \cat{C}\rt \cat{C}$ for the right $\cat{C}$-linear functor $X\otimes (-)$. We will then denote by $\cat{LAct}_{F_X}(\cat{C})$ the lax equaliser of $F_X$ and the identity, i.e.\ the pullback
$$\begin{tikzcd}[column sep=2.5pc]
\cat{LAct}_{F_X}(\cat{C})\arrow[r]\arrow[d, "\mm{Forget}"{left}] & \Fun(\Delta[1], \cat{C})\arrow[d]\\
\cat{C}\arrow[r, "{(F_X, \mm{id})}"{below}] & \Fun\big(\{0\}, \cat{C}\big)\times \Fun\big(\{1\}, \cat{C}\big).
\end{tikzcd}$$ 
This is a pullback diagram of right $\cat{C}$-module $\infty$-categories. One can identify $\cat{LAct}_{F_X}(\cat{C})$ with the $\infty$-category of objects $M\in \cat{C}$ equipped with an action map $X\otimes M\rt M$ (without further structure); the left vertical functor takes the underlying object in $\cat{C}$.
\begin{proposition}\label{prop:free left modules}
The forgetful functor $\cat{LAct}_{F_X}(\cat{C})\rt \cat{C}$ admits a right $\cat{C}$-linear left adjoint $\mm{Free}$, with the following properties:
\begin{enumerate}
\item There is a natural equivalence of right $\cat{C}$-linear endofunctors of $\cat{C}$
$$
\mm{Forget}\circ \mm{Free}(Y)\simeq \colim_n \big(T^{(n)}(X)\otimes Y\big).
$$
\item The free-forgetful adjunction is a monadic adjunction.
\end{enumerate}
\end{proposition}
\begin{proof}
Write $\cat{D}=\Fun(\mathbb{N}, \cat{C})$ for the category of sequences $Y_0\rt Y_1\rt \dots $ in $\cat{C}$ and let $\cat{LAct}_{F_X}(\cat{D})$ be the lax equaliser of the functors $Y_\bullet\mapsto X\otimes Y_\bullet$ and $Y_\bullet\mapsto Y_{\bullet+1}$. In other words, $\cat{LAct}_{F_X}(\cat{D})$ is the $\infty$-category of sequences $M_\bullet$ equipped with a natural map $X\otimes M_\bullet\rt M_{\bullet +1}$. The forgetful functor $\cat{LAct}_{F_X}(\cat{C})\rt \cat{C}$ then factors as the composite of right $\cat{C}$-linear functors
$$\begin{tikzcd}
\cat{LAct}_{F_X}(\cat{C})\arrow[r, "\mm{cst}"] & \cat{LAct}_{F_X}(\cat{D})\arrow[r, "\mm{ev}_0"] & \cat{C}.
\end{tikzcd}$$
The first functor, taking constant sequences, admits a left adjoint sending $M_\bullet$ to $\colim_n M_n$, since $X\otimes (-)$ preserves sequential colimits. We claim that the second functor admits a left adjoint sending $Y$ to the sequence $T^{(\bullet)}(X)\otimes Y$.

To see this, note that $T^{(\bullet)}(X)\otimes Y$ admits a natural left $X$-module structure, given by the obvious inclusion
$$
\lambda\colon X\otimes \big(T^{(n)}\otimes Y\big) \subseteq \Big(1\amalg \big(X\otimes T^{(n)}(X)\big)\Big)\otimes Y=T^{(n+1)}(X)\otimes Y.
$$
Note that $T^{(\bullet)}(X)\otimes Y$ is naturally equivalent to $Y$ in degree $0$. We therefore need to prove that evaluation at $0$ induces a natural equivalence
\begin{equation}\label{diag:free filtered module}\begin{tikzcd}
\Map_{\cat{LAct}_{F_X}(\cat{D})}\big(T^{(\bullet)}(X)\otimes Y, M_\bullet\big)\arrow[r, "\sim"] & \Map_{\cat{C}}(Y, M_0).
\end{tikzcd}\end{equation}
To see this, note that the $X$-linear mapping space from $T^{(\bullet)}(X)\otimes Y$ to $M_\bullet$ can be described inductively: a map $T^{(\bullet)}(X)\otimes Y\rt M_\bullet$ is given by a sequence of maps $f_n\colon T^{(n)}(X)\otimes Y\rt M_n$ together with commuting cubes
$$\begin{tikzcd}[column sep=1pc, row sep=1.5pc]
& X\otimes T^{(n-1)}(X)\otimes Y\arrow[rr, "\lambda"]\arrow[dd, gray, dashed] \arrow[ld, "i_n"{swap}] & & T^{(n)}(X)\otimes Y\arrow[dd, "f_{n}"]\arrow[ld, "i_n"]\\
X\otimes T^{(n)}(X)\otimes Y\arrow[rr, crossing over, "\lambda"{pos=0.6}]\arrow[dd, "X\otimes f_n"{left}] & & T^{(n+1)}(X)\otimes Y\\
& X\otimes M_{n-1}\arrow[rr, gray, dashed]\arrow[ld, gray, dashed] &  & M_{n}\arrow[ld]\\
X\otimes M_n\arrow[rr] & & M_{n+1}\arrow[uu, "f_{n+1}"{pos=0.7, swap}, crossing over, leftarrow].
\end{tikzcd}$$
Unraveling the definitions and using that $(-)\otimes Y$ preserves coproducts, one sees that the top square is coCartesian. Consequently, given $f_0, \dots, f_n$, there is a contractible space of maps $f_{n+1}$ making the above cube commute. Proceeding inductive, one then finds that the map \eqref{diag:free filtered module} is an equivalence.

The description of the left adjoint as $\colim_n T^{(n)}(X)\otimes Y$ gives property (1) and shows that it is right $\cat{C}$-linear (since the tensor product commutes with sequential colimits in the first variable). For (2), note that the free-forgetful adjunction satisfies the conditions of the Barr--Beck--Lurie theorem \cite[Theorem 4.7.3.5]{lurie2014higher}. Indeed, the forgetful functor clearly detects equivaleces. and if $M_\bullet$ is a simplicial diagram of $X$-modules which is split in $\cat{C}$, then it is also split in $\cat{LAct}_{F_X}(\cat{C})$: this follows immediately from the fact that $X\otimes \big(\colim M_\bullet\big)\simeq \colim (X\otimes M_\bullet)$ for any split simplicial diagram $M_\bullet$ in $\cat{C}$.
\end{proof}
\begin{proof}[Proof (of \Cref{thm:free algebra})]
Fix an object $X\in \cat{C}$ and let $F_X\colon \cat{C}\rt \cat{C}$ be its image under the monoidal equivalence $\cat{C}\simeq \mm{End}_{\cat{C}}(\cat{C})$. We will write $T_X\in \Alg(\mm{End}_{\cat{C}}(\cat{C}))$ for the right $\cat{C}$-linear monad associated to the free-forgetful adjunction $\cat{LAct}_X(\cat{C})\leftrightarrows \cat{C}$ from Proposition \ref{prop:free left modules}. 
Note that there is a natural map $\eta\colon F_X\rt T_X$ in $\mm{End}_{\cat{C}}(\cat{C})$, corresponding to the obvious map 
$$
X\rt T^{(1)}(X)=1\amalg X\rt \colim_n T^{(n)}(X)=T(X)
$$ 
under the monoidal equivalence $\cat{C}\simeq \mm{End}_{\cat{C}}(\cat{C})$. It therefore suffices to show that $\eta$ exhibits $T_X$ as the free algebra on $F_X$ in $\mm{End}_{\cat{C}}(\cat{C})$.

To see this, let $T\in \Alg(\mm{End}_{\cat{C}}(\cat{C}))$ be any right $\cat{C}$-linear monad and denote by $G_T\colon \Alg_{T}(\cat{C})\rt \cat{C}$ the right $\cat{C}$-linear forgetful functor from the $\infty$-category of $T$-algebras. Recall that there is a left action of $\mm{End}_{\cat{C}}(\cat{C})$ on the $\infty$-category $\mm{Fun}_{\cat{C}}(\Alg_T(\cat{C}), \cat{C})$ of right $\cat{C}$-linear functors, given by postcomposition. By the right $\cat{C}$-linear version of \cite[Lemma 4.7.3.1]{lurie2014higher}, the monad $T$ then arises as the endomorphism algebra of $G_T\in \mm{Fun}_{\cat{C}}(\LMod_T(\cat{C}), \cat{C})$. We therefore have to show that restriction along $\eta$ defines an equivalence
$$\begin{tikzcd}
\eta^*\colon \Map_{\Alg(\mm{End}_{\cat{C}}(\cat{C}))}\big(T_X, \mm{End}(G_T)\big) \arrow[r, "\sim"] & \Map_{\mm{End}_{\cat{C}}(\cat{C})}\big(F_X, \mm{End}(G_T)\big).
\end{tikzcd}$$
Using the universal property of the endomorphism algebra $\mm{End}(G_T)$, the domain can be identified with the space of $T_X$-module structures $T_X\circ G_T\rt G_T$. Such a $T_X$-module structure on $G_T$ simply endows each $T$-algebra with a natural $T_X$-algebra structure; in other words, the space of such $T_X$-module structures on $G_T$ is equivalent to the space of right $\cat{C}$-linear factorisations of $G_T$ as
$$\begin{tikzcd}
\Alg_{T}(\cat{C})\arrow[rr, dashed]\arrow[rd, "G_T"{swap}] & & \Alg_{T_X}(\cat{C})\arrow[ld, "G_{T_X}"]\\
& \cat{C}.
\end{tikzcd}$$
Likewise, $\Map_{\mm{End}_{\cat{C}}(\cat{C})}\big(F_X, \mm{End}(G_T)\big)$ can be identified with the space of natural maps $F_X\circ G_T\rt G_T$, i.e.\ with factorisations of $G_T$ over $\cat{LAct}_{F_X}(\cat{C})$. The assertion then follows from the fact that restriction along $\eta$ determines an equivalence $\Alg_{T_X}(\cat{C})\stackrel{\sim}{\rt} \cat{LAct}_{F_X}(\cat{C})$, by Proposition \ref{prop:free left modules}.
\end{proof}
\begin{remark}
The proof of \Cref{thm:free algebra} provides an additional property of the free algebra $T(X)$: there is an equivalence between left $T(X)$-modules in $\cat{C}$ and $X$-modules, i.e.\ objects equipped with a map $X\otimes M\rt M$.
\end{remark}

\newpage
\bibliographystyle{amsalpha}
\bibliography{There}

\end{document}